\documentclass[11pt]{amsart}

\usepackage{amssymb,amsmath,amsthm,amsfonts,amscd,latexsym, dsfont, color,tikz}
\usepackage[dvipsnames]{xcolor}

\newcommand{\Pic}{\mathrm{Pic}}

\newcommand{\NAH}{\mathrm{NAH}}
\usepackage{multirow}

\newcommand{\T}{\mathrm{T}}

\newcommand{\Thick}{\mathrm{Thick}}
\newcommand{\SU}{\mathrm{SU}}
\newcommand{\titsangle}{\angle_{\mathrm{Tits}}}
\newcommand{\image}{\mathsf{image}} 
 
\newcommand{\graph}{\mathsf{graph}}
\newcommand{\p}{\mathfrak{p}}
 
\newcommand{\solie}{\mathfrak{so}}
\usepackage{dynkin-diagrams}
\usepackage[title]{appendix}
\usepackage[american]{babel} 

\usepackage{array}
\usepackage[citecolor=ForestGreen,linkcolor=blue,urlcolor=black]{hyperref} 
 \usepackage{graphics}
\usepackage{wrapfig}
\usepackage{pst-node}

\usepackage{tikz-cd} 
\usepackage{tikz-3dplot}

\hypersetup{colorlinks}
\usepackage[normalem]{ulem}

\usepackage{setspace} 

\usepackage{tcolorbox}

\usepackage{footnote}

\usepackage{cancel}

\newcommand{\dev}{\mathrm{dev}}
\newcommand{\hol}{\mathsf{hol}}

\usepackage{MnSymbol}
\usepackage{enumitem}
\newcommand{\Fix}{\mathsf{Fix}}

\usepackage{pict2e}

\makeatletter

\newcommand{\newparallel}{\mathrel{\mathpalette\new@parallel{0.3}}}
\newcommand{\new@parallel}[2]{%
  \begingroup
  \settoheight{\unitlength}{$#1T$}
  \sbox\z@{\new@parallel@slash{#1}{#2}}%
  \mkern0.5mu\copy\z@\mkern-0.5mu\copy\z@\mkern0.5mu
  \endgroup
}
\newcommand{\new@parallel@slash}[2]{%
  \begin{picture}(#2,1)
  \roundcap
  \new@parallel@linethickness{#1}
  \Line(0,0)(#2,1.2)
  \end{picture}%
}
\newcommand{\new@parallel@linethickness}[1]{%
  \linethickness{%
      \ifx#1\displaystyle \fontdimen8\textfont\else
      \ifx#1\textstyle \fontdimen8\textfont\else
      \ifx#1\scriptstyle \fontdimen8\scriptfont\else
      1.1\fontdimen8\scriptscriptfont\fi\fi\fi 3
  }%
}
\makeatother

\newcommand{\spann}{\mathsf{span}}

\newcommand{\V}{\mathcal{E}} 
\definecolor{bluegreen}{RGB}{5, 170, 204}

\definecolor{defnyellow}{RGB}{255,224,102} 
\definecolor{lightblue}{RGB}{170, 220, 255} 
\definecolor{darkblue}{RGB}{0, 125, 230}

\newcommand{\Gtwo}{\mathsf{G}_2}

\newcommand{\Gtwosplit}{\mathsf{G}_2'}

\newcommand{\Sym}{\mathsf{Sym}}

\newtcolorbox{bluebox}[1]{colback=lightblue,colframe=darkblue,fonttitle=\bfseries,title=#1}
\newtcolorbox{redbox}[1]{colback=red!5!white,colframe=red!75!black,fonttitle=\bfseries,title=#1}

\usepackage{pgfplots}

\usepgfplotslibrary{colormaps}
\pgfplotsset{
    compat=newest,
    colormap={mycolormap}{color=(lightgray) color=(white) color=(lightgray) } }

\definecolor{mydarkblue}{RGB}{37, 42, 200}

\newcommand{\sff}{\mathrm{II}}
\newcommand{\tff}{\mathrm{III}}

\definecolor{mygreen}{RGB}{0, 150, 50}

\newcommand{\vis}{\partial_{\mathrm{vis}}}
\newcommand{\cl}{\mathcal{L}}
 
\newcommand{\Ein}{\mathsf{Ein}}

\newcommand{\imoct}{\mathrm{Im}(\Oct')}

\newcommand{\Mat}{\mathsf{Mat}}

\newcommand{\SO}{\mathrm{SO}}
\newcommand{\Gr}{\mathsf{Gr}} 
\newcommand{\Stab}{\mathsf{Stab}} 
\usepackage{mathrsfs}
\usepackage[utf8]{inputenc}
\usepackage[T1]{fontenc}
\newcommand{\R}{\mathbb R}
\newcommand{\RP}{\mathbb R \mathbb P}

\newcommand{\C}{\mathbb C}
\newcommand{\Z}{\mathbb Z}

\newcommand{\Ha}{\mathbb{H}} 
\newcommand{\F}{\mathbb{F}}
\newcommand{\g}{\mathfrak{g}}
\newcommand{\frakk}{\mathfrak{k}} 
\newcommand{\frakp}{\mathfrak{p}} 

\newcommand{\X}{\mathbb{X}}
\newcommand{\K}{\mathcal{K}}
\newcommand{\End}{\mathsf{End}} 

\newcommand{\sphere}{\mathbb{S}} 
 
\newcommand{\quadric}{\hat{ \mathbb{S}}^{2,4}}
\newcommand{\sllie}{\mathfrak{sl}}
\newcommand{\U}{\mathcal{U}} 
\newcommand{\SL}{\mathsf{SL}}
\newcommand{\PSL}{\mathsf{PSL}}
\newcommand{\Sp}{\mathsf{Sp}}
\newcommand{\GL}{\mathsf{GL}}
\newcommand{\gl}{\mathfrak{gl}}

\newcommand{\tr}{\text{tr}}
\newcommand{\id}{\text{id}}
 
 \newcommand{\Lie}{\mathsf{Lie}}
  \newcommand{\Aut}{\mathsf{Aut}}

 \newcommand{\Hom}{\mathsf{Hom}}
  
 \newcommand{\Fr}{\mathsf{Fr}\,}

 \newcommand{\del}{\partial}
\newcommand{\delbar}{\overline{\partial}}
 \newcommand{\Der}{\mathsf{Der}}

\newcommand{\Ann}{\mathsf{Ann}}
\renewcommand{\O}{\mathrm{O}}
 
\newcommand{\Iso}{\mathsf{Iso}} 
\usepackage{mathtools}

\DeclarePairedDelimiter\abs{\lvert}{\rvert}
 
\newcommand{\Pho}{\mathsf{Pho}}

\newcommand{\zbar}{\bar{z}}

\newcommand{\im}{\text{im}}

\newcommand{\sig}{\mathsf{sig}}

\newcommand{\Flag}{\mathsf{Flag}}
\newcommand{\Eintwothree}{\mathsf{Ein}^{2,3}} 
\newcommand{\Oct}{\mathbb{O}} 
 
\newcommand{\diag}{\mathsf{diag}}

\newcommand{\Hit}{\mathsf{Hit}}

\usepackage{thmtools, thm-restate}

\makeatletter
\newcommand\footnoteref[1]{\protected@xdef\@thefnmark{\ref{#1}}\@footnotemark}
\makeatother

\theoremstyle{plain}

\newtheorem{theorem}{Theorem}[section]
\newtheorem{proposition}[theorem]{Proposition}
\newtheorem{lemma}[theorem]{Lemma}
\newtheorem{corollary}[theorem]{Corollary}

\newtheorem{question}[theorem]{Question}
\newtheorem{example}{Example}[section]

\newtheorem{definition}[theorem]{Definition}

\newtheorem{remark}[theorem]{Remark}
\newcommand{\alignL}{\begin{flushleft}}
\newcommand{\alignLend}{\end{flushleft}}

\newtheorem{mainthm}{Theorem}

\newtheorem{mainthmprime}{Theorem}

\usepackage{soul}

\usepackage{bbm}

\title{$\Gtwosplit-$Slodowy Slices I: Geometric Structures}
\author{Colin Davalo and Parker Evans}
\date{\today}

\usepackage[margin= 1in]{geometry}

\numberwithin{equation}{section}

\begin{document}

\begin{abstract} 
Let $S$ be a closed surface of genus $g \geq 2$. 
We construct locally homogeneous geometric structures on closed 5-manifolds fibering over $S$, modeled on the two partial flag manifolds $\Ein^{2,3}$ and $\Pho^\times$ of the split real form $\Gtwosplit$ of the complex exceptional Lie group $\Gtwo^\C$. To this end, we consider two families of representations $\pi_1S\rightarrow \Gtwosplit$ constructed via the non-abelian Hodge correspondence from cyclic Higgs bundles, one associated with each $\Gtwosplit$-partial flag manifold. Each family includes $\Gtwosplit$-Hitchin representations, but is much more general. 
From the Higgs bundles of the first family, called $\beta$-bundles, we construct $(\Gtwosplit, \Ein^{2,3})$-geometric structures on $\Ein^{2,1}$-fiber bundles over $S$, and from Hodge bundles in the second family, called $\alpha$-bundles, we construct $(\Gtwosplit, \Pho^\times)$-geometric structures on $(\RP^2\times \mathbb{S}^1)$-bundles over $S$. In the case of $\Gtwosplit$-Hitchin Hodge bundles, which belong to both families, 
we show the image of the developing map of the respective geometric structures is exactly the domain of discontinuity defined by Guichard-Wienhard and Kapovich-Leeb-Porti. 

Each construction can be interpreted as converting a family of equivariant $J$-holomorphic curves in the pseudosphere $\quadric$ into geometric structures on fiber bundles $M \rightarrow S$. The approach used to build geometric structures, namely \emph{moving bases of pencils}, gives a unified description of prior analytic geometric structures constructions using Higgs bundles and harmonic maps.  
\end{abstract} 
 
\maketitle
\setstretch{1.2}
\tableofcontents
\setstretch{1.0} 

\vspace{-5ex}
\section{Introduction}

Let $S$ be a closed surface of genus $g \geq 2$. 
Beginning with Teichm\"uller space $T(S)$, incarnated as the moduli space of marked hyperbolic structures on $S$, there are now numerous examples of rich moduli spaces of geometric structures associated to surface group representations into semisimple Lie groups $G$. In the classical case of Teichm\"uller space, the holonomies of these hyperbolic structures, \emph{Fuchsian representations}, lie in the real rank one group $G = \PSL(2,\R)$. Hitchin introduced generalizations of Teichm\"uller space in \cite{Hit92}, now called the \emph{$G$-Hitchin components $\Hit(S,G)$}, for $G$ a split, real, simple Lie group of higher rank. The component $\Hit(S,G) \subset \chi(S,G)$ in the character variety $\chi(S,G)$ of reductive representations up to conjugation is a cell and contains a distinguished embedding of $T(S)$ called the \emph{$G$-Fuchsian-Hitchin locus}. These Fuchsian-Hitchin representations factor through the \emph{principal} embedding $\PSL(2,\R)\hookrightarrow G$. 
In fact, one can view Teichm\"uller space $T(S)$ as the $\PSL(2,\R)$-Hitchin component. 
It is now known $\Hit(S,G)$ is realizable as a moduli space, possibly in many ways, of holonomies of \emph{$(G,X)$-structures} on \emph{some (higher-dimensional) manifold} $M$ \cite{GW12}. Here, $M$ is a fiber bundle over $S$ \cite{AMTW25}, and the holonomy of the $(G,X)$-structure, a priori a map $\pi_1M \rightarrow G$, descends to $\pi_1S$ to be conflated with a surface group representation called the \emph{descended holonomy}. Correspondingly, we shall call such a $(G,X)$-manifold $M \rightarrow S$ to be \emph{fibered}. Despite the general (abstract) existence results of \cite{GW12,KLP18}, there are \emph{explicit} geometric structures descriptions of $\Hit(S,G)$ only in a few cases  \cite{CG93, GW08, Bar10, GW12, CTT19, NR25,Eva25, RT25}.

Broadly speaking, the $G$-Hitchin component is a special locus in the character variety $\chi(S,G)$ with inherent geometric meaning that remains to be understood in detail. However, one can consider other loci of special representations, which are not necessarily a union of connected components, or even open in $\chi(S,G)$. A particularly powerful tool to locate special representations is the \emph{non-abelian Hodge (NAH) correspondence}, which for a given Riemann surface $\Sigma$ on $S$, realizes a  homeomorphism between a moduli space $\mathcal{M}_{G}(\Sigma)$ of holomorphic objects, \emph{polystable $G$-Higgs bundles}, and the moduli space of reductive representations, the character variety $\chi(S,G)$. In particular, the NAH correspondence enables us to probe the character variety in search of special representations by looking for distinguished Higgs bundles.  

Via the NAH correspondence, the $G$-Hitchin component $\Hit(S,G)$ can be identified with a \emph{Slodowy slice}, in this case associated to the \emph{principal} $\sllie_2$-triple in $\g^\C$ that is called \emph{magical} in \cite{BCGBO24}. Across all such magical triples, these Slodowy slices conjecturally yield (nearly) all components of character varieties of surface groups with the same miraculous property as the $G$-Hitchin components: \emph{containing only discrete and faithful representations}. Such components are now often referred to as \emph{higher Teichm\"uller spaces}\footnote{\label{noteTubeType}We exclude from the definition of higher Teichm\"uller spaces components with rigidity, i.e. with no Zarisky dense elements, which is the case of maximal representations into non tube-type Hermitian $G$ \cite[Theorem 5]{BIW10}.} due to the analogy with $T(S)$ \cite{Wie18}. Moreover, there is a 1-1 correspondence between magical triples and the \emph{$\Theta$-positive} structures of \cite{GW25}. Currently, $\Theta$-positive representations, or equivalently, representations associated to magical Slodowy slices, comprise all known higher Teichm\"uller spaces\footnoteref{noteTubeType}\cite{BCGBO24, GW25, BGLPW24}.

However, one can also consider non-magical Slodowy slices, such as those for root $\mathfrak{sl}_2$-triples. These slices seem to provide promising candidates for special representations: in \cite{Fil23, Zha25, DB25, Vir25}, for the Lie groups $\SO_0(2,3)$, $\SL(3,\R)$, $\SU(n,1)$, the associated representations have been found to be \emph{Anosov} and to admit fibered geometric structures under favorable circumstances. 

\subsubsection{Our focus} In this paper, we study representations in the split real form $\Gtwosplit$ of the complex exceptional Lie group $\Gtwo^\C$ arising via the NAH correspondence as real points in the simple root $\sllie_2$-Slodowy slices. Now, denote $\beta$ and $\alpha$ the short and long roots of $\g_2$. The two families of Higgs bundles we study, \emph{$\alpha$-bundles} and \emph{$\beta$-bundles}, are particular cyclic $\Gtwosplit$-Higgs bundles. For each family, the associated representations include Hitchin representations, but are more general. Indeed, for $\sigma \in \{\alpha,\beta\}$ and fixed Riemann surface $\Sigma$ on $S =S_g$, the associated slice of Higgs bundles in $\mathcal{M}_{\Gtwosplit}(\Sigma)$ is closed, unbounded, and disconnected, with components $\mathcal{M}_{\sigma}^d(\Sigma)$ indexed by a discrete invariant $d$. By recent work of Collier-Toulisse-Wentworth \cite{CTW25} and, independently Li-Zhang \cite{LZ26}, we know, at least locally, that the associated representations $\rho$ vary as the Riemann surface $\Sigma$ is deformed. In particular, the subspace $\chi_{\sigma}^{d}(S)\subset \chi(S,\Gtwosplit)$ of all the representations
associated to $\bigsqcup_{\Sigma \in T(S)}\mathcal{M}_{\sigma}^d(\Sigma)$, has the following dimension at smooth points (see Proposition \ref{prop:Dimension count}):
\begin{itemize}
    \item $\dim \chi^d_{\beta}(S_g) = 28g-28-2(6g-6-d)$,  where $0\leq d\leq 6g-6$, 
    \item $\dim \chi^d_{\alpha}(S_g) = 28g-28-4(2g-2-d)$, where $-2g+2\leq d\leq 2g-2$. 
\end{itemize}
For comparison, the character variety $\chi(S,\Gtwosplit)$ has dimension $28g-28$ at smooth points. 

We realize the representations in $\chi_{\sigma}^d(S)$ as holonomies of explicit $(\Gtwosplit,X)$-structures on a closed 5-manifold $M^5 = M(\sigma,d)$ fibered over $S$, where $X=\Gtwosplit/P_{\sigma}$ is the corresponding $\Gtwosplit$-flag manifold. 
In the case of $\Gtwosplit$-Hitchin representations, we relate our construction to the $(\Gtwosplit,X)$-manifolds defined abstractly by domains of discontinuity in \cite{GW12, KLP18}. 

One accomplishment of the construction is that beyond the Hitchin case and a few additional reducible exceptions, the associated representations are not known to be Anosov, or even discrete and faithful. Thus, there is no guarantee a priori of such a geometric structure to exist. 

\emph{Parabolic geometries}, a special case of \emph{Cartan geometries}, have become a popular subject in the last forty years (see \cite{CS09} for history and details), born out of familiar examples such as Riemannian, conformal, and projective geometry. Now, locally homogeneous $(G,G/P)$-manifolds are \emph{flat} examples of $(G,P)$-parabolic geometries. The cases of $(\Gtwosplit, P_{\alpha})$ and $(\Gtwosplit,P_{\beta})$ are two of the most exotic examples, and lead to rich geometry as in \cite{BH92, EN20, MNS21, HN22}. 

For $\Gtwosplit/P_{\beta}$, the 5-manifolds carry a \emph{maximally non-integrable 2-plane distribution} $\mathscr{D}$, also called a \emph{$(2,3,5)$-distribution}, along with a conformal structure of signature $(2,3)$. For $\Gtwosplit/P_{\alpha}$, the geometry is more subtle, but the 5-manifolds modeled on this space carry, in particular, a contact structure $\mathscr{C}$ such that the projective bundle $\mathbb{P}(\mathscr{C})$ holds a special family of cubic curves (see \cite{MNS21}). We believe the explicit examples built presently of 5-manifolds $M$ fibered over closed surfaces $S$ of arbitrary genus $g \geq 2$ with fibers closed 3-manifolds may be worthy of additional investigation, particularly regarding how these geometric features interact with the bundle structure of $M$.

\subsubsection{Strategy}
The main strategy employed presently is largely informed by prior work \cite{Bar10, CTT19, DB25} using Higgs bundles and harmonic maps to build fibered geometric structures. Our main techniques to construct such geometric structures make use of \emph{bases of pencils} along equivariant harmonic maps with a certain \emph{parallelism}. These ideas offer a unified framework for previous analytic constructions of geometric structures, as we explain in Appendix \ref{Appendix:Unified}. 

We now summarize the construction. Fix a Riemann surface $\Sigma= (S,J)$ on $S$ and a Higgs bundle $(\V,\Phi)$ on $\Sigma$ that is a $\sigma$-bundle, for $\sigma \in \{ \alpha, \beta\}$. Let $\rho:\pi_1S\rightarrow \Gtwosplit$ be the associated representation. The NAH correspondence passes from $(\V,\Phi)$ to $\rho$ using the unique $\rho$-equivariant harmonic map $f=f_{\rho}: \tilde{\Sigma} \rightarrow \X$ to the $\Gtwosplit$-Riemannian symmetric space $\X=\X_{\Gtwosplit}$. 
The $\sigma$-bundle condition on $(\V,\Phi)$ allows us to define a distinguished \emph{pencil} along $f$, a smoothly varying family of 2-planes $\mathcal{P} \in \Omega^0(\tilde{\Sigma}, \Gr_{2}(f^*\T\X))$. For each pencil $\mathcal{P}_x\subset \T_{f(x)}\X$, one can form an associated codimension two submanifold $\mathcal{B}_{\sigma}(\mathcal{P}_x) \subset \mathcal{F}_{\sigma}$ of the flag manifold $\mathcal{F}_{\sigma} = \Gtwosplit/P_{\sigma}$, called the \emph{$\sigma$-base of pencil}. One can then assemble these bases of pencil into a manifold $\overline{B}_{\mathcal{P}} \subset \tilde{\Sigma} \times \mathcal{F}_{\sigma}$, fibering over $\tilde{\Sigma}$, mapping tautologically into $\mathcal{F}_{\sigma}$ by $(p,F)\mapsto F$. We  work to verify this map is indeed a local diffeomorphism, producing a fibered $(\Gtwosplit,\mathcal{F}_{\sigma})$-manifold structure on the quotient $B_{\mathcal{P}}=(\pi_1S)\backslash \overline{B}_{\mathcal{P}}$ with descended holonomy $\rho$.

This construction is guided by \cite{Dav25}, where the same procedure is studied for maps $f:\tilde{\Sigma}\rightarrow \X$ that are not necessarily harmonic, but instead \emph{nearly geodesic}, and the pencil $\mathcal{P}_x$ is chosen to be the tangent pencil, namely $\mathcal{P}_x=df(\T_x\tilde{\Sigma})$.
In the present case, we \emph{do not} use the tangent pencil, but instead modify it in Lie-theoretic fashion that happens to be very explicit in the present circumstances. 
Indeed, the pencil $\mathcal{P}$ we use is conveniently described in terms of an auxiliary $\rho$-equivariant map, namely a $J$-holomorphic curve $\nu:\tilde{\Sigma} \rightarrow \quadric$ to the almost-complex pseudosphere $\quadric$, a pseudo-Riemannian $\Gtwosplit$-homogeneous space. As in \cite{Nie24, CT24}, a suitable Gauss map construction applied to $\nu$ recovers $f$. The $\sigma$-bundle conditions manifest geometrically in $\nu$, enabling the construction of the pencil $\mathcal{P}$: the first and third fundamental forms $\mathrm{I}, \tff$ of $\nu$ are non-vanishing for $\beta$-bundles and the second fundamental form $\sff$ of $\nu$ is non-vanishing for $\alpha$-bundles. The pencils of interest, in each respective case, are then built out of these non-vanishing objects.

To implement the strategy of `moving bases of pencils', we must overcome three main obstacles: 
\begin{enumerate}[label=(\alph*)]
    \item What is the \emph{topology} of the base of pencil $\mathcal{B}_{\sigma}(\mathcal{P}_x)$?
    \item What is the \emph{geometry} of the base of pencil $\mathcal{B}_{\sigma}(\mathcal{P}_x)\subset \mathcal{F}_{\sigma}$? 
    \item How do the fibers of the base of pencil $\mathcal{B}_{\sigma}(\mathcal{P}_x)$ vary in $x$? 
\end{enumerate}

In general, (a) is already difficult. Presently, we find \emph{the fibers themselves are fiber bundles} whose fibrations are related to principal bundle fibrations of the flag manifolds $\mathcal{F}_{\sigma} =\Gtwosplit/P_{\sigma}$. For (b), we need an explicit coordinate parametrization of our fibers, which is much harder to obtain than just the topology alone. 
For (c), each base of pencil is understood in (b) relative to different points in the symmetric space. Comparing fibers of the moving base of pencil $\overline{B}_{\mathcal{P}}$ then requires translating between different models for the flag manifold $\mathcal{F}_{\sigma}$. Consequently, lengthy calculations are needed to verify the tautological map is indeed a local diffeomorphism, especially in the case of $\Gtwosplit/P_{\alpha}$.\medskip  

The group $\Gtwosplit$ is related to \emph{exceptional} geometry that is unique and interesting in its own right. However, $\Gtwosplit$ is also situated at an interesting middle ground: it is a low rank Lie group, but has sufficiently large dimension and intricacy so as to intimate general features. In short: $\Gtwosplit$ provides an ideal environment to test the efficacy of Lie-theoretic constructions. To this end, the present work highlights some previously obscured difficulties and illuminates new phenomena; see $\S$\ref{Sec:Technique}. 

The study of geometric structures with Higgs bundles does not end with $G=\Gtwosplit$. 
Rather, this work is part of the beginning of a developing trend, largely promoted by Collier, to study Slodowy slices for Lie groups of any rank and whether the associated representations admit fibered geometric structures and are Anosov. In forthcoming work \cite{DE26Anosov}, we prove many of the representations associated to $\sigma$-bundles are $P_{\sigma}$-Anosov and, moreover, the developing maps of the geometric structures constructed here with bases of pencil have image precisely the domain of discontinuity defined in \cite{GW12, KLP18}.

\subsection{Main Results}\label{Sec:MainResults}
We now describe our main results. 

\subsubsection{Using Higgs Bundles}
In $\Gtwosplit$, we consider two families of cyclic Higgs bundles: $\beta$-bundles and $\alpha$-bundles. 
Each of these families of Higgs bundles are sub-families of the cyclic $\Gtwosplit$-Higgs bundles introduced by Collier-Toulisse in \cite{CT24}. 

For $\Sigma$ a Riemann surface on $S$, we consider cyclic bundles with the following general form:
\begin{equation}\label{G2CyclicGeneral_Intro}
\begin{tikzcd}
\mathcal{L}_2\mathcal{L}_1 \arrow[r, "\beta"] & \mathcal{L}_2\arrow[r, "\alpha"]&\mathcal{L}_1\arrow[r, "-i\sqrt{2}\beta"]&\mathcal{O} \arrow[r, "-i\sqrt{2}\beta"]&\mathcal{L}_1^{-1} \arrow[r, "\alpha"]&\mathcal{L}_2^{-1} \arrow[r, "\beta"]\arrow[lllll,bend left, "\delta"]&\mathcal{L}_1^{-1}\mathcal{L}_2^{-1}\arrow[bend left, lllll, "\delta"]
\end{tikzcd}.
\end{equation}
The particular shape of these bundles allows us to place a $\Gtwosplit$-structure on them. Here, $\alpha, \beta, \delta$ are holomorphic endomorphism valued-one forms twisted by the canonical bundle $\K=(\T^{1,0}\Sigma)^*$. For example, the first object $\beta$ is of the form $\beta \in H^0(\Hom(\mathcal{L}_1\mathcal{L}_2, \mathcal{L}_2)\otimes \K)= H^0(\mathcal{L}_1^{-1}\K)$.  
In this diagram, $\alpha, \beta$ play the respective roles of the long and short simple roots for $\g_2$, and $\delta$ is the highest root. In particular, $\delta$ is a long root. Unlike in the corresponding case of $\SO_0(2,3)$, the long roots $\alpha$ and $\delta$ are not symmetric with respect to $\beta$; the following remarks clarify this essential point. 

Among these Higgs bundles, we impose symmetries that allow $\alpha$ or $\beta$ to be non-vanishing, in which case they are holomorphic sections of a trivial bundle $\mathcal{O}$, and can be re-gauged to simply be written `$1$'. Indeed, we consider 
bundles for which  $\beta=1$ (resp. $\alpha=1$), which are called $\beta$-cyclic (resp. $\alpha$-cyclic). 
Hitchin representations correspond exactly to Higgs bundles that are both $\alpha$-cyclic and $\beta$-cyclic \cite{Hit92}. On the other hand, if we consider the case when $\delta , \beta =1$, the corresponding representations (essentially) factor through 
$\SL(3,\R)$, as explained in \cite{CT24}.
Finally, one may consider bundles with $\alpha, \delta =1$. When $\beta=0$, we obtain a representation factoring through a special copy of $\PSL(2,\R)$, but for $\beta \neq 0$, the obtained representations are irreducible. 
\medskip

We now describe the relevant Higgs bundles in greater detail. Let $\K = \K_{\Sigma}$ be the holomorphic cotangent line bundle on $\Sigma$.  
The $\beta$-cyclic bundles on $\Sigma$ obtain the following form, for some holomorphic line bundle $\mathcal{B}$:
\begin{equation}\label{BetaCyclic}
\begin{tikzcd}
\mathcal{B} \arrow[r, "1"] & \mathcal{B}\mathcal{K}^{-1}\arrow[r, "\alpha"]&\mathcal{K}\arrow[r, "-i\sqrt{2}"]&\mathcal{O} \arrow[r, "-i\sqrt{2}"]&\mathcal{K}^{-1} \arrow[r, "\alpha"]&\mathcal{B}^{-1}\mathcal{K} \arrow[r, "1"]\arrow[bend left, lllll, "\delta"]&\mathcal{B}^{-1} \arrow[bend left, lllll, "\delta"]
\end{tikzcd}.
\end{equation}
 
We show that generically such Higgs bundles do give associated geometric structures. 
\begin{mainthm}[Geometric Structures for $\beta$-Bundles]\label{IntroThm:Ein23GeometricStructures}
    Let $(\V,\Phi)$ be a polystable $\beta$-cyclic bundle such that $\mathrm{div}(\delta) \neq 0$ and $\rho:\pi_1S\to \Gtwosplit$ be the associated representation.
    Then there is a canonically associated closed $(\Gtwosplit, \Ein^{2,3})$-manifold $M^5\rightarrow S$ with fiber $\Ein^{2,1}$ and descended holonomy $\rho$. 
\end{mainthm}

We note that the demand that $\delta$ has nonzero divisor merely avoids the case that $\alpha=0$, $\delta=1$ \cite{CT24}. As noted above, such representations (essentially) factor through $\SL(3,\R)$. In particular, they are $\alpha$-Anosov but never $\beta$-Anosov.\footnote{When $\alpha =\delta =0$ and $\mathcal{B}\cong \K^{1/2}$, we obtain representations factoring through the $\SL(2,\R)$-subgroup of the short root $\beta$, and the conclusion of Theorem \ref{IntroThm:Ein23GeometricStructures} does hold. See Proposition \ref{Prop:BetaBundleStability} and Remark \ref{Remk:BetaSpecialCases} for further details.} It is then expected the theorem does not apply in this case. We emphasize, as alluded to earlier, that this manifold $M^5$ is explicitly constructed, though we avoid the precise details in the introduction.

The use of $\Ein^{2,1}$ for the fiber is motivated by the results of \cite{DE26}, where we show the fiber of the geometric structure constructed via the natural domain $\Omega$ in $\Ein^{2,3}$ for $\Gtwosplit$-Hitchin representations, first defined in \cite{GW12}, is $\Ein^{2,1}$. Note, however, that \cite{DE26} determines the fiber and global topology of the $\Ein^{p-1,p}$-manifolds  $\rho(\pi_1S)\backslash \Omega_{\rho}$ for $\rho$ any $p$-Anosov deformations of an $\iota$-Fuchsian in $\SO_0(p,p+1)$ more generally for all $p \geq 3$. 
 
When a Higgs bundle in \eqref{BetaCyclic} satisfies $\mathcal{B}\cong\mathcal{K}^3$, the corresponding representation is $\Gtwosplit$-Hitchin. In this case, we compare our geometric structures with the ones constructed as quotients of domains of discontinuity in \cite{GW12} and  \cite{KLP18}. 
We quickly introduce notation for the theorem. The fibering of the geometric structures built entails the developing map $\dev: \widetilde{M} \rightarrow \Ein^{2,3}$ descends to the $\pi_1S$-cover $\overline{M}$ of $M$ to a map $\dev: \overline{M}\rightarrow \Ein^{2,3}$. In this notation, the result is as follows.

\begin{mainthm}[$\Ein^{2,3}$-structures for $\Gtwosplit$-Hitchin Representations]\label{Thmbeta:DG=GGT}
Let $(\V,\Phi)$ be a $\beta$-bundle on a Riemann surface $\Sigma$ such that $\mathcal{B}\cong\mathcal{K}^3$, and $\rho: \pi_1 S \rightarrow \Gtwosplit$ the associated representation, which is Hitchin.  
Then the developing map $\dev: \overline{M} \rightarrow \Ein^{2,3}$ of the geometric structure from Theorem \ref{IntroThm:Ein23GeometricStructures} is a 
diffeomorphism onto the Tits metric thickening domain $\Omega_{\rho}$. Consequently, $M \cong_{\mathrm{Diff}}\rho(\pi_1S) \backslash \Omega_{\rho}$. 
\end{mainthm}

Hitchin representations have a canonical $\beta$-bundle attached. Indeed, there is a homeomorphism between the $\Gtwosplit$-Hitchin component and the moduli space of $\beta$-bundles with $\mathcal{B}\cong \K^3$ by \cite{Lab17, CT24}. In particular, the construction in Theorem \ref{IntroThm:Ein23GeometricStructures} provides an explicit geometric structures interpretation of $\Hit(S,\Gtwosplit)$, 
as we now spell out. Given $\rho \in \Hit(S,\Gtwosplit)$, by Labourie's identification $\Hit(S,\Gtwosplit)\cong \mathcal{Q}^6$, where $\mathcal{Q}^6 \rightarrow T(S)$ is the bundle of holomorphic sextic differentials, we can take the associated pair $(\Sigma, q_6)$, where $\Sigma\in T(S), q_6 \in H^0(\K^6_{\Sigma})$. We can then construct the $\beta$-bundle in \eqref{BetaCyclic} with $\mathcal{B}=\K^3$, $\alpha=\beta=1$ and $\delta =q_6$. The construction of Theorem \ref{IntroThm:Ein23GeometricStructures} then builds (differential geometrically) a canonical 5-manifold $M \rightarrow S$ endowed with a $(\Gtwosplit, \Ein^{2,3}$)-structure with descended holonomy $\rho$. Theorem \ref{Thmbeta:DG=GGT} asserts that this manifold $M$ is exactly the quotient $\rho(\pi_1S)\backslash \Omega_{\rho}$ of the \cite{GW12, KLP18}-domain of discontinuity. For further details, see $\S$\ref{Subsec:G2Ein23HitchinCase}.\medskip

We now discuss the other partial flag manifold $\Pho^\times \cong \Gtwosplit/P_{\alpha}$. Unlike its cousin, $\Ein^{2,3}\cong\Gtwosplit/P_{\beta}$, which is also an $\SO_0(3,4)$-flag manifold, $\Pho^\times$ requires inescapable $\Gtwosplit$-geometry. All analogous constructions pursued in this case are an order of magnitude more difficult, due to the intricacy of this homogeneous space.

With a similar approach to \cite{DE26}, we can compute the topology of the fibers of the $\frac{\pi}{2}$-Tits metric thickening domain $\Omega_\rho$ from \cite{KLP18} for $\Gtwosplit$-Hitchin representations. 
We find the following. 

\begin{mainthm}
[$\Pho^\times$-fibers for Hitchin representations]\label{IntroThmPhoXFibers}
Let $\rho \in \Hit(S,\Gtwosplit)$ and $\mathcal{M} = \rho(\pi_1S)\backslash \Omega_{\rho} $ the associated $(\Gtwosplit,\Pho^\times)$-manifold. Then there is a smooth fiber bundle $\mathcal{M} \rightarrow S$ with fiber $\RP^2\times \mathbb{S}^1$.
\end{mainthm}

Motivated by this result, we seek to build geometric structures on $(\mathbb{RP}^2\times \mathbb{S}^1)$-fiber bundles over $S$. To this end, we consider \emph{$\alpha$-Hodge bundles}, which are the sub-family of cyclic $\Gtwosplit$-Higgs bundles of the form \eqref{G2CyclicGeneral_Intro} with $\alpha=1$ that are also $\C^*$-fixed points in the moduli space of Higgs bundles. Up to gauge, they have the following form, where either $\beta=0$ or $\delta=0$: 
\begin{equation}
\begin{tikzcd}
\mathcal{T}^2\mathcal{K} \arrow[r, "\beta"] & \mathcal{T}\mathcal{K}\arrow[r, "1"]&\mathcal{T}\arrow[r, "-\sqrt{2}i \beta"]&\mathcal{O} \arrow[r, "-\sqrt{2}i\beta"]&\mathcal{T}^{-1} \arrow[r, "1"]&\mathcal{T}^{-1}\mathcal{K}^{-1} \arrow[r, "\beta"]\arrow[lllll,bend left, "\delta"]&\mathcal{T}^{-2}\mathcal{K}^{-1}\arrow[lllll,bend left, "\delta"]
\end{tikzcd}.
\end{equation}

For $\alpha$-bundles, we can also construct associated geometric structures on an appropriate fiber bundle over the surface $S$. 

\begin{mainthm}[Geometric Structures for $\alpha$-Bundles]\label{IntroThmPhoXStructure}
    Let $(\V,\Phi)$ be a polystable $\alpha$-bundle with $\beta=0$ or $\delta=0$, and $\rho:\pi_1S\to \Gtwosplit$ be the associated representation.
    Then there is a canonically associated closed $(\Gtwosplit, \Pho^\times)$-manifold $M^5 \rightarrow S$ with fiber $\RP^2\times \sphere^1$ and descended holonomy $\rho$. 
\end{mainthm}

When $\mathcal{T}\cong\mathcal{K}$, the associated representation is Fuchsian-Hitchin. In this case we compare our geometric structures with the ones constructed as quotients of domains of discontinuity in \cite{KLP18}:
\vspace{-3ex}
\begin{mainthm}[$\Pho^\times$-structures for Fuchsian-Hitchin Representations]
Let $\rho:\pi_1S\rightarrow\Gtwosplit$ be the representation associated to an $\alpha$-bundle on a Riemann surface $\Sigma$ such that $\mathcal{T}\cong\mathcal{K}$ and $\delta=0$, so that the associated representation $\rho$ is $\Gtwosplit$-Fuchsian-Hitchin.
Then the developing map $\dev: \overline{M} \rightarrow \Pho^\times$ of the geometric structure from Theorem \ref{IntroThmPhoXStructure} is a finite covering map onto the Tits metric thickening domain $\Omega_{\rho}$. Consequently, $M$ is a finite cover of the quotient $\rho(\pi_1S) \backslash \Omega_{\rho}$.
\end{mainthm}

For a description of $M$ topologically, see Corollary \ref{Cor:QuotientTopologyAlpha}.

\subsubsection{Using Equivariant Minimal Surfaces}
 
A key similarity between the present construction of geometric structures and that of \cite{CTT19} is the use of equivariant harmonic maps. 
In \cite{CTT19}, Collier, Tholozan, and Toulisse used maximal spacelike surfaces $\sigma: \tilde{\Sigma} \rightarrow \Ha^{2,n}$ in pseudo-hyperbolic space $\Ha^{2,n}$ to build fibered geometric structures. 
In particular, they consider \emph{immersed} such maps, which they prove are naturally attached to maximal representations $\rho:\pi_1S\rightarrow \SO_0(2,n+1)$. For each such maximal spacelike surface, they build an $(\SO_0(2,n+1), \Pho(\R^{2,n+1}))$-structure on a compact manifold $M \rightarrow S$ with fiber $\Pho(\R^{2,n})$. In the special case $n=3$ that $\SO_0(2,3)$ is split, with the additional demand that $\sigma$ has \emph{non-vanishing second fundamental form}, they convert $\sigma$ into an $(\SO_0(2,3), \Ein^{1,2})$-structure on $\T^1S \rightarrow S$. 
Of emphasis here is the coarse picture suggested: appropriate adjectives on the map $\sigma$ encode instructions for the correct model space.  

Our work is largely motivated by \cite{CTT19} and \cite{CT24}.
In the present case, we trade $\SO_0(2,n+1)$ for $\Gtwosplit$ and maximal spacelike surfaces for \emph{alternating} \emph{$J$-holomorphic} curves in the pseudosphere $\quadric$. The space $\quadric$, the unit sphere in $\R^{3,4}$, carries a canonical $\Gtwosplit$-invariant almost-complex structure, allowing us to consider $J$-holomorphic curves $\nu:\Sigma \rightarrow \quadric$ from a Riemann surface $\Sigma$. Here, the \emph{alternating} condition, named in \cite{Nie24}, entails that $\nu$ has a generalized $(T,N,B)$-\emph{Frenet frame} analogous to space curves in $\mathbb{E}^3$. 
Our two constructions of geometric structures are equivalently described via two equivariant families of $J$-holomorphic curves, which we call \emph{$\alpha$-curves} and \emph{$\beta$-curves}. 

The $\beta$-curves are \emph{immersed} alternating $J$-holomorphic curves $\nu: \tilde{\Sigma} \rightarrow \quadric$. 
Via these curves, we find a reinterpretation of Theorem \ref{IntroThm:Ein23GeometricStructures} from the perspective of harmonic maps. 

\setcounter{mainthm}{1}
\begin{mainthmprime}[$\beta$-curves to Geometric Structures]
Let $\nu: \tilde{S}\rightarrow \quadric$ be a $\rho$-equivariant $\beta$-curve that is linearly full. 
Then there is a canonically associated closed $(\Gtwosplit, \Ein^{2,3})$-manifold $M^5 \rightarrow S$ with fiber $\Ein^{2,1}$ and descended holonomy $\rho$. 
\end{mainthmprime}
We note here that the \emph{linearly full} condition above, namely that the curve is not contained in any codimension one subspace of $\R^{3,4}$, excludes $\mathsf{div}(\delta)=0$ in Theorem \ref{IntroThm:Ein23GeometricStructures} by \cite{CT24}. In this work, Collier and Toulisse explain how the the moduli space of (polystable) $\beta$-bundles with $\alpha \neq0$ is homeomorphic to the moduli space $\mathcal{H}_{\beta}(S)$ of equivariant alternating $\beta$-curves, namely pairs $(\nu, \rho)$, where $\nu:\tilde{\Sigma} \rightarrow \quadric$ is a $\rho$-equivariant $\beta$-curve and $\rho:\pi_1S\rightarrow \Gtwosplit$ is a representation. 

The geometric structures in $\Pho^\times$ are also built from a different flavor of equivariant alternating $J$-holomorphic curves that we call \emph{$\alpha$-curves}, those that are allowed finite branch points, but have \emph{non-vanishing second fundamental form}. 
Indeed, the $\alpha$-bundles (with $\beta\neq0$) defined previously correspond to curves $\tilde{\Sigma} \rightarrow \quadric$, allowing a reinterpretation of Theorem \ref{IntroThmPhoXStructure}.

\setcounter{mainthm}{4}
\begin{mainthmprime}[$\alpha$-curves to Geometric Structures]
Let $\nu: \tilde{S} \rightarrow \quadric$ be a $\rho$-equivariant $\alpha$-curve that is superminimal. 
Then there is a canonically associated closed $(\Gtwosplit, \Pho^\times)$-manifold $M^5 \rightarrow S$ with fiber $\RP^2\times \mathbb{S}^1$ and descended holonomy $\rho$. 
\end{mainthmprime}

Here, the \emph{superminimal} condition on $\nu$ is a demand on its \emph{harmonic map sequence} (cf. \cite[Section 2.4]{Bar10}) and corresponds to the associated Higgs bundle being a $\C^*$-fixed point.  
Besides the Hitchin case, these $\alpha$-curves have not been studied. 
For fixed Riemann surface $\Sigma$, we study more generally the space of polystable $\alpha$-bundles, which has $4g-3$ connected components labeled by $d=\deg(\mathcal{T})$ (see Proposition \ref{Prop:AlphaStability}). In particular, we describe explicitly the degree constraints on $\mathcal{T}$ for the strata with $\beta\neq0$ and $\delta=0$, corresponding to the above theorem. We consider the \emph{generic} locus of stable $\alpha$-bundles, with the Riemann surface $\Sigma \in T(S)$ varying, in Proposition \ref{prop:Dimension count}.

\subsection{Techniques and Related Work}\label{Sec:Technique}
\subsubsection{Geometric structures via Bases of Pencils} 

The first step of our general construction is to understand two special cases in detail. We consider \emph{$\sigma$-Fuchsian} representations $\rho: \pi_1S\rightarrow \Gtwosplit$, for $\sigma \in \{\alpha, \beta\}$, factoring through the $\SL(2,\R)$-subgroup of the root $\sigma$, described in Appendix \ref{Appendix:SL2}.
This subgroup preserves preserves a totally geodesic hyperbolic plane $\Ha^{2}_{\sigma}$ in the $\Gtwosplit$-symmetric space $\X$, which is also the image of the associated equivariant harmonic map $f: \tilde{\Sigma} \rightarrow \X$ under any uniformizing $\sigma$-bundle associated to $\rho$.
Let $\mathcal{F}_{\sigma} = \Gtwosplit/P_{\sigma}$ be the associated flag manifold. 
In this case, we use the tangent pencils $\mathcal{P}_x = \T_{f(x)}\Ha^2_{\sigma}$ to build the $(\Gtwosplit,\mathcal{F}_{\sigma})$-manifold $B_{\mathcal{P}}$ of interest, or rather, its $\pi_1S$-cover. 
Indeed, $\overline{B}_{\mathcal{P}} \subset \tilde{\Sigma} \times \mathcal{F}_{\sigma} $ is defined with fiber $\overline{B}_{\mathcal{P}}|_x$ the \emph{base of the pencil} $\mathcal{B}_{\sigma}(\mathcal{P}_x)\subset \mathcal{F}_\sigma$. This base $\mathcal{B}_{\sigma}(\mathcal{P}_x)$ is concretely the set of flags $F \in \mathcal{F}_{\sigma}$ reached at infinity in the visual boundary $\vis\X$ by traveling along geodesics in $\X$ emanating from $x$, with initial velocity orthogonal to $\mathcal{P}_x$. By \cite{Dav25}, $\mathcal{B}_{\sigma}(\mathcal{P})$ is a smooth closed codimension two submanifold of $\mathcal{F}_{\sigma}$ in this case. 

Now, still in the $\sigma$-Fuchsian case, the most difficult part (c) of the construction, namely understanding how the fibers vary, comes for free from \cite{Dav25}: the tautological map $\overline{B}_{\sigma} \rightarrow \mathcal{F}_{\sigma}$ by $(p,F)\mapsto F$ is a diffeomorphism onto the cocompact Tits metric thickening domain $\Omega_{\rho} \subset \mathcal{F}_{\sigma}$ defined by \cite{GW12, KLP18}. However, we must still work to understand the geometry and topology of this fiber, problems (a) and (b), even in this  case; see further discussion in $\S$\ref{subsec:fibers}.

Let us now take a general $\sigma$-bundle $(\V,\Phi)$ on a Riemann surface $\Sigma$, with associated representation $\rho: \pi_1S \rightarrow \Gtwosplit$. As we have explained, such a Higgs bundle gives an associated $\rho$-equivariant $J$-holomorphic curve $\nu: \tilde{\Sigma} \rightarrow \quadric$ with special properties. 
It turns out we can use $\mathrm{I}+\tff^{0,1}$ and $\sff$ in the respective cases of $\beta$ and $\alpha$-bundles to define moving pencils $\mathcal{P}_{\sigma} \in \Omega^0(\tilde{\Sigma}, \Gr_{2}(f^*\T\X))$ along the associated harmonic map $f:\tilde{\Sigma} \rightarrow \X$. 
These pencils are very special: pointwise $\mathcal{P}_{\sigma}|_x$ is tangent a sub-symmetric space $\mathcal{H}_{\sigma}|_x$, moving in $x$, that is $\Gtwosplit$-equivalent to a tangent space to $\Ha^2_{\sigma}$. This feature is critical. 
Indeed, problems (a) and (b), the geometry and topology of the base of pencils in the general case, are \emph{exactly} the same as the previous case. 
Thus, in the general case, the problem is complementary to the special cases discussed in the previous paragraph: all that remains is (c), to understand how these fibers $\mathcal{B}_{\sigma}(\mathcal{P}_x)$ of bases of pencil vary in $x \in \tilde{\Sigma}$. Just as before,  $\overline{B}_{\mathcal{P}} \subset \tilde{\Sigma} \times \mathcal{F}$ has a tautological map $\dev: \overline{B}_{\mathcal{P}} \rightarrow \mathcal{F}_{\sigma}$, namely $\dev(p, F)=F$. 
In this case, we do not know a priori whether $\dev$ is a local diffeomorphism, and 
we must verify that the developed fibers actually define a local fibration. This verification is the main challenge and requires a nontrivial calculation. 

\subsubsection{Parallel distribution of planes}

The second main idea is to construct the  distribution of planes $(\mathcal{P}_x)_{x\in \widetilde{\Sigma}}$ to be \emph{parallel} in the following sense. We find a conformal Riemannian metric $g$ on $\tilde{\Sigma}$ and a parametrization of $\mathcal{P}$ by a vector bundle map $\Psi_0 \in \Omega^0(\tilde{\Sigma}, \T^*\tilde{\Sigma}\otimes f^*\T\X)$ such that $\Psi_0$ is parallel. 

The equivariant vector bundle map $\Psi_0:\T\widetilde{\Sigma}\to \T\X$ lifts the harmonic map $f:\widetilde{\Sigma}\to \X$, but is not necessarily equal to $\mathrm{d}f$. Now, the differential $\mathrm{d}f$ can be reconstructed from the Higgs field $\Phi$ as it is identified via the Maurer-Cartan form with $\Phi+\Phi^*$. Crucially, the object  $\Phi+\Phi^* \in \Omega^0(\T^*\Sigma \otimes \End(\V))$, for $\V$ the flat bundle associated to $\rho$, can also be seen as an equivariant object of the form $\Phi \in \Omega^0(\tilde{\Sigma},\T^*\tilde{\Sigma}\otimes f^*\T\X)$.  
For Higgs bundles of the form \eqref{G2CyclicGeneral_Intro}, the Higgs field $\Phi$ decomposes into three parts depending on $\alpha$, $\beta$ and $\delta$. We consider another endomorphism-valued one-form $\Phi_0$ by modifying $\Phi$, keeping only the $\alpha$-part, respectively the $\beta$-part, which is non-vanishing by our selection of Higgs bundles. We then define $\Psi_0:=\Phi_0+\Phi_0^*$. 

The parallelism of the pencils and a maximum principle for the Hitchin system allow us to reduce the desired transversality condition on $\dev$ to an explicit, but tedious, computation.

In fact, the previous work of \cite{CTT19} using maximal spacelike surfaces to build geometric structures can be reinterpreted through a parallel distribution of planes constructed analogously. Similarly, for Hitchin representations in $\SL(3,\R)$, geometric structures modeled on $\mathbb{RP}^2$ or on the space of full flags $\Flag(\R^3)$ can be described by the corresponding affine sphere, and equivalently by parallel distributions of planes. We give more details in Appendix \ref{Appendix:Unified}. 

A similar strategy was used in \cite{DB25}, but in that case problems (a) and (b), the geometry and topology of the fiber, are essentially trivial. 
However, in that work, studying $G=\SL(3,\R)$ and $X=\Flag(\R^3)$, problem (c), understanding how the bases of pencil move, is handled alongside an even stronger property, (d) the Anosov property, by relating the moving codimension two bases of pencils with the \emph{nestedness} of certain codimension one submanifolds. In this way, the Anosov property in \cite{DB25} is verified simultaneously alongside the construction of the geometric structure. Thus, the problems faced and results obtained in that work are different from those presently.

\subsubsection{Geometry and topology of the fibers}\label{subsec:fibers}

Let us revisit the initial problems (a) and (b):  understanding the topology and geometry of the base $\mathcal{B}_{\sigma}(\mathcal{P})$ of a tangent pencil $\mathcal{P}$ to the sub-symmetric space $\Ha^2_{\sigma} \hookrightarrow \X_{\Gtwosplit}$ of the $\SL(2,\R)$-subgroup corresponding to $\sigma \in \{\alpha, \beta\}$. 
It is not enough to know that we abstractly used a base of pencil; we must \emph{explicitly} be able to parametrize this locus to prove the transversality property of our geometric structures for more general cyclic $\beta$ and $\alpha$-bundles.

A difference between the $\Gtwosplit$ case and the case of $\SO_0(2,n+1)$ or $\SL(3,\R)$ is that this parametrization is considerably harder to find and write explicitly. Indeed, the fiber is a $3$-manifold in the present case, rather than circle or a union of circles when $M = \mathrm{T}^1S$ or $\mathbb{P}(\T^1S)$.

The fibers of the $(G,X)$-manifolds we build are naturally realized as fiber bundles themselves. 
In fact, their structure as a fiber bundle comes from the ambient space. 
Each of $\Ein^{2,3}$ and $\Pho^\times$ can be realized as $\mathbb{S}^3$-fiber bundles over $\RP^2$. It is no coincidence that the fibers we use, namely $\Ein^{2,1}$ and $\RP^2 \times \mathbb{S}^1$, are each circle bundles over $\RP^2$; the bases of pencils we use for fibers are sub-fibrations of the ambient flag manifolds. This structural result aids in our solution of problem (a): the topology of the fiber. We observe similar behavior for bases of pencils in \cite{DE26} in the flag manifold $X = \Ein^{p-1,p}$ for $G = \SO_0(p,p+1)$ more generally for $p \geq 3$.

The next complication is the geometry of the fiber, for which a comparison is illuminating. In \cite{CTT19}, the auxiliary maximal spacelike surfaces held a simple and elegant relation to the associated geometric structures. 
To be precise, consider once more the $\pi_1S$-cover $\overline{M} \rightarrow \tilde{S}$ where the developing map naturally lives. For $G=\SO_0(2,n+1)$ and $X=\Pho(\R^{2,n+1})$, the developed fibers $\dev(\overline{M}_p)$ in \cite{CTT19} were $\Pho(\R^{2,n})$ \emph{topologically and geometrically}. Indeed, the developed fiber is precisely $\dev(\overline{M}_p)=\Pho(\sigma(p)^\bot)$, for $\sigma:\tilde{\Sigma} \rightarrow \Ha^{2,n}$ the associated maximal surface. In words, the developed fiber is a geometrically `straight' copy of $\Pho(\R^{2,n})$, and the data of $\sigma(p)$ immediately yields the corresponding developed fiber $\dev(\overline{M}_p)$. 
In our case, the fibers are \emph{twisted}. 
For $X=\Ein^{2,3}$, our developed fibers $F_p=\dev(\overline{M}_p)$, for any $p \in \tilde{S}$, are each \emph{topologically} $\Ein^{2,1}$. 
However, the linear span of $F_p$ is not 5-dimensional, as a normal, or `straight' copy of $\Ein^{2,1}\hookrightarrow \Ein^{2,3}$. Instead, the linear span of the developed fiber is the whole ambient vector space: $\spann (F_p)=\R^{3,4}$. 
The same remarks apply for $X=\Pho^\times$, where the fibers are also twisted in similar fashion. 

Without a more general picture in mind, namely using bases of pencils, the present construction was elusive. Indeed, a $(\Gtwosplit,\Ein^{2,3})$-geometric structures interpretation of $\Hit(S,G)$ was provided in \cite{Eva25}. However, those geometric structures have some pathologies, including $\dev$ being non-proper and strictly locally injective in the Fuchsian-Hitchin case. These strange features arose as a consequence of trying to use a `straight' copy of $\Ein^{2,1}$ for the fiber, with a certain small `bad set' excised out, leaving a leftover $\mathbb{S}^1\times \mathbb{S}^1\times \R$ fiber and hence a non-compact manifold housing the geometric structure. The present work has some clear upsides: it is effective for both $\Ein^{2,3}$ and $\Pho^\times$, it works in far greater generality than just the Hitchin case, and it also elucidates the \cite{GW12, KLP18}-geometric structures. 

\subsection{Organization} 
We now discuss the layout of the paper. 

\begin{itemize}[noitemsep]
    \item In Section \ref{Sec:Prelims}, we introduce preliminaries on $\Gtwosplit$ and non-compact symmetric spaces $\mathbb{X}$ in general, including bases of pencils. We then discuss the geometry of the $\Gtwosplit$-symmetric space and the embeddings of the $\Gtwosplit$-partial flag manifolds in the visual boundary of $\X_{\Gtwosplit}$. 

    \item In Section \ref{Sec:CyclicG2'Higgs}, we consider a family of cyclic $\Gtwosplit$-Higgs from \cite{CT24}. We carefully introduce the $\Gtwosplit$-structure on such bundles, and derive Hitchin's equations for such Higgs bundles general. We then justify the dimension counts of the spaces $\chi_{\sigma}^d(S)$ from the introduction.

    \item In Section \ref{Sec:G2Ein23GeomStructures}, we construct $(\Gtwosplit, \Ein^{2,3})$-geometric structures for $\beta$-cyclic Higgs bundles, proving Theorem \ref{IntroThm:Ein23GeometricStructures}. 
    \item In Section \ref{Sec:PhoXGeometricStructures}, we construct $(\Gtwosplit, \Pho^\times)$-structures for $\alpha$-cyclic Hodge bundles, proving Theorem \ref{IntroThmPhoXStructure}. Along the way, we determine the fibers of the $(\Gtwosplit, \Pho^\times)$-manifolds for $\Gtwosplit$-Hitchin representations to be $\mathbb{RP}^2\times \mathbb{S}^1$, proving Theorem \ref{IntroThmPhoXFibers}.
    \item In Appendix \ref{Appendix:Unified}, we provide a unified construction of geometric structures in rank two, for $G \in \{\SL(3,\R), \SO_0(2,n+1), \Gtwosplit\}$ also using bases of pencils with certain parallelism. 

    \item In Appendix \ref{Appendix:SL2}, we explicitly describe the five distinct $\sllie_2\R$-subalgebras in $\g_2'$ up to the adjoint action of $\Gtwosplit$, which were classified in Lie-theoretically in \cite{Dok88}. The uniformizing Higgs bundles for Fuchsian representations in each of the corresponding analytic subgroups of $\Gtwosplit$ appear as $\beta$ or $\alpha$-bundles.

    \item In Appendix \ref{Sec:RegularityPencil}, we discuss regularity of pencils in the $\Gtwosplit$-symmetric space. This technical material is needed for the proof of the fibers in the $\Pho^\times$ case for Hitchin representations. 
\end{itemize}

\subsection*{Acknowledgments}
We would like to thank Alex Nolte for his valuable comments and suggestions on the exposition of the paper. The second named author thanks Jacob Erickson for discussions on $\Gtwosplit$-parabolic geometries. C. Davalo was funded by the European Union via the ERC 101124349 ``GENERATE.'' Views and opinions expressed are however those of the authors only and
do not necessarily reflect those of the European Union or the European Research
Council Executive Agency. Neither the European Union nor the granting authority
can be held responsible for them.

\section{Preliminaries}\label{Sec:Prelims}

In this section, we discuss necessary background on the exceptional Lie group $\Gtwosplit$, symmetric spaces $\X$ of non-compact type and their visual boundaries both in general and in the particular case $\X = \X_{\Gtwosplit}$, and facts about Anosov representations and nearest point projections from \cite{Dav25}.  

\subsection{The Exceptional Lie Group \texorpdfstring{$\Gtwosplit$}{G2'}}

Let us introduce the essential background on the split real Lie group $\Gtwosplit$ needed for the paper. We focus mostly on the irreducible faithful representation in dimension seven related to the imaginary split octonions $\imoct$. Since we need quite a bit of information, we include an outline below. The reader is encouraged to reference this section for background depending on their interests. 

\begin{itemize}[noitemsep]
    \item Subsection \ref{Subsec:G2Defn}: \emph{Definitions}. A definition of $\Gtwosplit$ in terms of the split octonions $\Oct'$, remarks on equivalent definitions.
    \item Subsection \ref{Subsec:CrossProductBases}: \emph{Cross product bases}. The notion of $\mathbb{F}$-cross product bases, which for $\F \in \{\R,\C\}$ describe the eigenbases of a $\g_2^\F$-Cartan subalgebra that is $\R$-split for $\F= \R$. Hence, $\R$-cross product are related to $\Gtwosplit$-flag manifolds. $\C$-cross product bases are needed to understand the cyclic $\Gtwosplit$-Higgs bundles in Section \ref{Sec:CyclicG2'Higgs}. 
    \item Subsection \ref{Subsec:StiefelModels}: \emph{Stiefel triplet models}. Two different Stiefel triplet models for $\Gtwosplit$ that are $\Gtwosplit$-torsors. These models are needed to understand various $\Gtwosplit$-homogeneous spaces. 
    \item Subsection \ref{Subsec:LieTheory}: \emph{Lie theory}. Basic Lie theory for $\g_2^\F$, including root space decomposition, Weyl group, and root vectors (expressed in $\F$-cross product bases). 
    \item Subsection \ref{Subsec:Annihilators}: \emph{Annihilators}. This crucial notion is needed throughout the whole paper, as it is foundational to the $\Gtwosplit$-flag manifolds $\Ein^{2,3}$ and $\Pho^\times$.
    \item Subsection \ref{Subsec:SubgroupHomogeneous}: \emph{Subgroups \& Homogeneous spaces}. Important subgroups of $\Gtwosplit$, including in particular the maximal compact $K$, the parabolic subgroups $P_{\beta}, P_{\alpha}, P_{\Delta}$ and the corresponding $\Gtwosplit$-flag manifolds $\Ein^{2,3}, \Pho^\times, \mathcal{F}_{1,2}^\times$. Additional geometry of $\Ein^{2,3}$ and $\Pho^\times$ is discussed later in Sections \ref{Sec:Ein23ModelSpace} and \ref{Sec:PhoxModelSpace} before the geometric structures are built.
\end{itemize}
As a general reference about $\Gtwo$, we recommend \cite{Fon18}, which discusses each of the complex, split, and compact groups $\Gtwo^\C, \Gtwo', \Gtwo^c$. For the history of $\Gtwo$, see \cite{Agri08}. For further exposition on octonions and their relation to exceptional Lie groups, see \cite{Bae02, Esc18}. 

\subsubsection{\texorpdfstring{$\Gtwo$}{G2} via  Octonions} \label{Subsec:G2Defn}

Among the Lie algebras $\g_2, \mathfrak{f}_4, \mathfrak{e}_6, \mathfrak{e}_7, \mathfrak{e}_8$ of exceptional type from the Cartan-Killing classification of complex simple Lie algebras, $\g_2$ is the smallest in dimension and rank. 
In this section, we briefly introduce the associated Lie groups of type $\Gtwo$: the complex Lie group $\Gtwo^\C$ and its two real forms $\Gtwo'$, the split real form, and $ \Gtwo^c$, the compact real form. Our main focus is on $\Gtwosplit$. Unlike some of its exceptional companions, $\Gtwo'$-geometry has the admirable feature of being reasonably explicit, i.e., computations in $\Gtwosplit$ are \emph{tractable}, involving only linear algebra in dimension 7, opposed to its next eldest sibling $F_4$, whose smallest irreducible faithful representation is in dimension 26 \cite{Bae02}. 
We describe briefly the three fundamental perspectives on the $\Gtwo$'s -- as algebra automorphisms, cross product automorphisms, and 3-form stabilizers -- which are all linked through octonions. 

Over the reals, the octonions come in two flavors: $\Oct$ and $\Oct'$, the latter being the \emph{split octonions}. Each of these objects are 8-dimensional algebras over $\R$ equipped with a plethora of structure. 
The exceptional complex simple Lie group $\Gtwo^\C$ and its two (adjoint) real forms $\Gtwo^c$ and $\Gtwosplit$ of compact and split type are the $\F$-algebra automorphisms of $(\Oct')^\C \cong \Oct^\C:=\mathbb{O}\otimes_\R \C$, $\Oct$, and $\Oct'$, respectively. 

One way to define $\Oct$ and $\Oct'$ is through the \emph{Cayley-Dickson} and \emph{split Cayley-Dickson} processes, which pleasantly explain how $\Oct, \Oct'$ can be built from $\R$ via an inductive process, in analogy to building $\C$ from $\R$. We refer the reader to \cite[Section 4.1]{ER25} or \cite[Section 2.2]{Bae02} for the construction over $\R$ or to \cite{Sch95} for more general discussion of groups of type $\Gtwo$ and octonions over an arbitrary field. We circumvent Cayley-Dickson and offer a direct definition here. 

The algebra $\Oct'$ is an algebra extension $\Oct' =\Ha[\mathbf{l}]$ of the quaternions $\Ha$ by a new element $\mathbf{l}$ such that $\mathbf{l}^2=+1$ and $\mathbf{l}$ anti-commutes with the imaginary quaternions $\mathbf{i},\mathbf{j},\mathbf{k}$. In particular, $\mathcal{M}=(1, \mathbf{i},\mathbf{j},\mathbf{k},\mathbf{l},\mathbf{l}\mathbf{i},\mathbf{l}\mathbf{j},\mathbf{l}\mathbf{k})$ is a canonical vector space basis of $\Oct'$ from this perspective.   
We shall refer to $\mathcal{M}$ as the \emph{multiplication basis} for $\Oct'$.

\begin{remark}
    When we consider $(\Oct')^\C:=\Oct'\otimes_\R\C$, the standard octonion $\mathbf{i}\in \Oct'$ is different from the standard imaginary number $i\in \C$. After $\S$\ref{Sec:Prelims}, we will rarely use the notation $\mathbf{i}$.
\end{remark}

The algebra $\Oct'$ is not associative, but is \emph{alternative}, meaning in particular that the subalgebra $\mathcal{A}_{x,y}$ generated by any two elements $x,y$ is associative. In other words, the \emph{associator} 
\begin{align}\label{Associator}
    [\cdot,\,\cdot,\,\cdot]:\Oct'\times \Oct'\times\Oct'\rightarrow \Oct', \;\mathrm{by} \; [u,v,w]=u(vw)-(uv)w 
\end{align}
is alternating. 
Any triplet of generators $(x,y,z) \in \mathcal{M}^3$ such that $z \notin \mathcal{A}_{x,y}$ \emph{anti-associate}, meaning $(xy)z=-x(yz)$. 
The two given facts on commutators and associators in $\Oct'$ uniquely describe the algebra multiplication in $\Oct'$ among the basis elements $\mathcal{M}$, as displayed in Table \ref{fig:SplitOctonionTable}. 
When a formal symbol is necessary, we denote $\odot$ for the algebra product $\odot:\Oct'\times \Oct'\rightarrow \Oct'$, though we usually simply write $xy=x\odot y \in \Oct'$ to denote the product of $x,y \in \Oct'$ simply by juxtaposition. 

\begin{table}[ht]
    \centering
    \[
\begin{array}{|c|c|c|c|c|c|c|c|c|}
\hline
 col\odot row& 1 & \mathbf{i} & \mathbf{j} & \mathbf{k} & \mathbf{l} & \mathbf{l}\mathbf{i} & \mathbf{l}\mathbf{j} & \mathbf{l}\mathbf{k} \\ \hline
1 & 1 & \mathbf{i} & \mathbf{j} & \mathbf{k} & \mathbf{l} & \mathbf{l}\mathbf{i} & \mathbf{l}\mathbf{j} & \mathbf{l}\mathbf{k} \\[4pt] \hline
\mathbf{i} & \mathbf{i} & -1 & \mathbf{k} & -\mathbf{j} &- \mathbf{l}\mathbf{i} & \mathbf{l} & -\mathbf{l}\mathbf{k} & \mathbf{l}\mathbf{j} \\[4pt] \hline 
\mathbf{j} & \mathbf{j} & -\mathbf{k} & -1 & \mathbf{i} & -\mathbf{l}\mathbf{j} & \mathbf{l}\mathbf{k} & \mathbf{l} & -\mathbf{l}\mathbf{i} \\[4pt] \hline
\mathbf{k} & \mathbf{k} & \mathbf{j} & -\mathbf{i} & -1 & -\mathbf{l}\mathbf{k} & -\mathbf{l}\mathbf{j} & \mathbf{l}\mathbf{i} & \mathbf{l} \\[4pt] \hline
\mathbf{l} & \mathbf{l} & \mathbf{l}\mathbf{i} & \mathbf{l}\mathbf{j} & \mathbf{l}\mathbf{k} & 1 & \mathbf{i} & \mathbf{j} & \mathbf{k} \\[4pt] \hline
\mathbf{l}\mathbf{i} & \mathbf{l}\mathbf{i} & -\mathbf{l} & -\mathbf{l}\mathbf{k} & \mathbf{l}\mathbf{j} & -\mathbf{i} & 1 & \mathbf{k} & -\mathbf{j} \\[4pt] \hline 
\mathbf{l}\mathbf{j} & \mathbf{l}\mathbf{j} & \mathbf{l}\mathbf{k} & -\mathbf{l} & -\mathbf{l}\mathbf{i} & -\mathbf{j} & -\mathbf{k} & 1 & \mathbf{i} \\[4pt] \hline
\mathbf{l}\mathbf{k} & \mathbf{l}\mathbf{k} & -\mathbf{l}\mathbf{j} & \mathbf{l}\mathbf{i} & -\mathbf{l} & -\mathbf{k} & \mathbf{j} & -\mathbf{i} & 1 \\ \hline 
\end{array}
\]
    \caption{\emph{The multiplication table for $\Oct'$ in the multiplication basis $\mathcal{M}$.}}
    \label{fig:SplitOctonionTable}
\end{table}

The split real Lie group $\Gtwosplit$ is then realized by \cite[Corollary 4.4]{Fon18}
\[ \Gtwosplit = \Aut_{\R-\rm{alg}}(\Oct') = \{ \psi \in \GL(\Oct') \mid \psi(xy) =\psi(x)\psi(y)\}.\]
Analogously, for the complex group $\Gtwo^\C$, we have $\Gtwo^\C = \Aut_{\C-\mathrm{alg}}(\Oct'^\C)$.

There are a few other algebraic structures related to $\Oct'$: a conjugation $*$, a symmetric bilinear form $q$ and a cross product $\times$. First, one defines the \emph{imaginary split octonions} $\mathrm{Im}(\Oct')$ as the span of the non-unital generators in $\mathcal{M}$, so that $\Oct' = \R \oplus \imoct$. This direct sum of vector spaces shall upgrade to an orthogonal splitting momentarily. The \emph{split octonion conjugation} $*: \Oct'\rightarrow \Oct'$ is the involution with $+1$-eigenspace $\R = \R\{1_{\Oct'}\}$ and $(-1)$-eigenspace $\imoct$. 
Moreover, $(xy)^*=y^*x^*$ so $*$ is an algebra anti-involution. The quadratic form $q $ on $\Oct'$ is of split signature (4,4) and is defined by $q(x) = xx^*$, which is automatically real. The orthogonal projections $\mathrm{Re}: \Oct' \rightarrow \R$ and $\mathrm{Im}: \Oct' \rightarrow \imoct$ are given by 
\[ \begin{cases}
    \mathrm{Re}(x) = \frac{1}{2}(x+x^*)\\
    \mathrm{Im}(x) = \frac{1}{2}(x-x^*).
\end{cases}\]

The algebra $\Oct'$ is a \emph{composition algebra}, meaning $q$ is both non-degenerate and multiplicative over the algebra product: $q(xy)=q(x)q(y)$. Such algebras are exceedingly rare -- Hurwitz' theorem classifies all composition algebras over $\R$ to be one of the following: $\R,\,\C, \,\Ha, \,\Oct, \, \C', \,\Ha',\,\Oct'$, where $\mathbb{A}'$ denotes the split counterpart of the classical algebra. The relation $q(xy)=q(x)q(y)$ holds precisely due to alternativity of $\Oct'$. By abuse, we may write $u \cdot v$ to denote $q(u,v)$, the bilinear form induced by the quadratic form $q$. We emphasize: the dot product $\cdot$ is scalar-valued and distinct from the vector-valued algebra product $\odot$.

Note that any transformation $\psi \in \Gtwosplit$ must fix $1_{\Oct'}$. Thus, it is standard to consider the action of $\Gtwosplit$ on $\imoct$, upon which one realizes an irreducible faithful representation $\Gtwosplit \rightarrow \GL(\imoct)$, 
which is one of the two \emph{fundamental representations} of $\Gtwosplit$, the other being the adjoint representation \cite[page 32]{Fon18}. 
As $\Gtwosplit$ preserves $*$ and the algebra product, $\Gtwosplit<  \mathrm{O}(\imoct, q) \cong \mathrm{O}(3,4)$.

Observe that the algebra product $\odot$ does not restrict to a map of type $\imoct \times \imoct \rightarrow \imoct$, e.g. $\mathbf{i}^2 =-1$. To remedy this issue, one defines the \emph{cross product} $\times: \imoct \times \imoct \rightarrow \imoct$ by 
\[ u\times v := \mathrm{Im}(u v) = uv-(u\cdot v)1_{\Oct'}. \]
An elementary but frequently useful observation is that $uv= u\times v$ if $u\bot v$. 
Since $\Gtwosplit$ fixes the real axis pointwise, one finds $\Gtwosplit < \Aut(\imoct, \times)$. In fact, the inclusion is an equality, and this leads to the second (possible) definition, namely  $\Gtwosplit=\Aut(\imoct, \times)$.
Additionally, $\Gtwosplit$ preserves the volume form $\mu$ on $\imoct$, so the whole tuple 
$(\odot, q, \mu ,*)$ is $\Gtwosplit$-invariant. Similarly, the group $\Gtwo^\C$ preserves the corresponding tuple $(\odot, q, \mu,*)$ of complexified objects. Moreover, $\Gtwo^\C = \Aut(\imoct^\C, \times)$. We will use this fact to study $\g_2^\C$ in what follows. We shall also use the notation $\Gtwo^\F$ to denote either $\Gtwo^\R=\Gtwo'$ or $\Gtwo^\C$ for uniformity.

Here, we briefly remark on the cross product. We shall denote $\mathcal{C}_{u}:\imoct \rightarrow \imoct$ as the \emph{cross product endomorphism} $\mathcal{C}_u(x) = u\times x$ of $u \in \imoct$. A very useful identity we shall frequently appeal to is the \emph{double cross product identity} (a non-standard yet descriptive title):
\begin{align}\label{DCP}
    u \times (u \times v) = -q(u) v+(u\cdot v)u.
\end{align}

The final definition of the $\Gtwo$'s is related to 3-forms. Define the \emph{scalar triple product} $\Omega \in \Lambda^3(\imoct^*)$ by 
\[\Omega(u,v,w) = (u\times v)\cdot w.\]
Since $\Gtwosplit$ preserves $\times, q$, we have $\Gtwosplit < \Stab_{\GL(7,\R)}(\Omega)$. Once more, the inclusion is an equality; it is a subtle fact that $\Omega$ determines $\times$ and $\mu$ (see \cite[Lemma 4.12, Proposition 4.14]{Fon18}). Thus, the third common definition of $\Gtwosplit$ is by $\Gtwosplit = \Stab_{\GL(7,\R)}(\Omega)$. It is from this perspective that $\dim_{\R}(\Gtwo') =14$ is most clear. Indeed, the 3-form $\Omega$ is \emph{generic} in the sense that its $\GL(7,\R)$-orbit in $\Lambda^3((\R^7)^*)$ is open, and is one of only two 3-forms with this property, the other one being $\Omega_c$, the scalar triple product on $\mathrm{Im}(\Oct)$, whereby $\Gtwo^c= \Stab_{\GL(7,\R)}(\Omega_c)$. This is usually the starting place for studying $\Gtwo^c$-manifolds in Riemannian geometry. We refer the reader to \cite[Section 4.2]{Fon18} for further details on generic 3-forms in dimension seven and to \cite{Kar20} for an introduction to $\Gtwo^c$.

\subsubsection{Real \& Complex cross product Bases}\label{Subsec:CrossProductBases}

We now introduce a special type of vector space basis for $\imoct^\F := \imoct \otimes_{\R} \mathbb{F}$, for $\F \in \{\R,\C\}$.  

\begin{definition}\label{Defn:FCrossProductBasis}
Let $X= (x_k)_{k=3}^{-3}$ be a vector space basis for $\imoct^\F$. We call $X$ an $\F$-cross product basis when $x_k \times x_l=c_{k,l}x_{k+l}$ for some constants $c_{k,l} \in \F$. 
\end{definition}

Note that $x_{k} =0$ for $|k| > 3$ is implied. The name ``$\F$-cross product basis'' was given in \cite{Eva24a} as a unifying notion for appropriately ordered eigenbasis of Cartan subalgebras of $\g_2^\F$. 

We call the $\C$-cross product basis from \cite[(3.78)]{Bar10} the \emph{model} $\C$-cross product basis. Here, the basis $(e_k)_{k=3}^{-3}$ is given by: 
\begin{align}\label{BaragliaBasis}
    \begin{cases}
    e_{\pm 3} &= \frac{1}{\sqrt{2}}(\mathbf{jl} \pm i\, \mathbf{kl}),\\
    e_{\pm 2} &=\frac{1}{\sqrt{2}}(\mathbf{j} \pm i \,\mathbf{k}) ,\\
    e_{\pm 1} &=\frac{1}{\sqrt{2}}(\mathbf{l} \pm i\, \mathbf{il}),\\\
    e_0 &= \mathbf{i}.
    \end{cases}
\end{align}
The complex split octonions $(\Oct')^\C$ have an additional conjugation, namely complex conjugation relative to the real subspace $\Oct'$, denoted by $z\mapsto \overline{z}$, which is different from the split octonion conjugation $*$. Note that $\overline{e_k}=e_{-k}$ for the basis in \eqref{BaragliaBasis}. Since $\overline{uv} = \overline{u}\,\overline{v}$ for any $u,v \in (\Oct')^\C$, this leads to many symmetries in the multiplication table for the basis $(e_k)_{k=3}^{-3}$, as shown in Table \ref{Table:ComplexCrossProductBasis}, which also confirms that $(e_k)_{k=3}^{-3}$ is a $\C$-cross product basis. For later, it will be useful to record the matrix $[q]$ that encodes the complex bilinear form $q=q_{3,4}^\C$ in the basis $(e_k)_{k=3}^{-3}$. Namely, 
\begin{align}\label{QModelBasis}
    [q] = \begin{pmatrix}
        & & & & & & -1\\
        & & & & & +1& \\
        & & & & -1 & &\\
        & & & +1& & & \\
        & & -1& & & & \\
        & +1& & & & & \\
        -1 & & & & & & \\
    \end{pmatrix}
\end{align}

\begin{table}[h!]
\renewcommand{\arraystretch}{1.3}
\centering
\begin{tabular}{|c|c|c|c|c|c|c|c|}
\hline
col $\times$ row & $e_3$ & $e_2$ & $e_1$ & $e_0$ & $e_{-1}$ & $e_{-2}$ & $e_{-3}$ \\
\hline
$e_3$ & $0$ & $0$ & $0$ & $-ie_3$ & $\sqrt{2}e_2$ & $\sqrt{2}e_1$ & $-ie_0$ \\
\hline
$e_2$ & $0$ & $0$ & $\sqrt{2}e_3$ & $ie_2$ & $0$ & $-ie_0$ & $-\sqrt{2}e_{-1}$ \\ \hline
$e_1$ & $0$ & $-\sqrt{2}e_3$ & $0$ & $ie_1$ & $ie_0$ & $0$ & $-\sqrt{2}e_{-2}$ \\ \hline
$e_0$ & $ie_3$ & $-ie_2$ & $-ie_1$ & $0$ & $ie_{-1}$ & $ie_{-2}$ & $-ie_{-3}$ \\ \hline
$e_{-1}$ & $-\sqrt{2}e_2$ & $0$ & $-ie_0$ & $-ie_{-1}$ & $0$ & $-\sqrt{2}e_{-3}$ & $0$ \\ \hline
$e_{-2}$ & $-\sqrt{2}e_1$ & $ie_0$ & $0$ & $-ie_{-2}$ & $\sqrt{2}e_{-3}$ & $0$ & $0$ \\ \hline
$e_{-3}$ & $ie_0$ & $\sqrt{2}e_{-1}$ & $\sqrt{2}e_{-2}$ & $ie_{-3}$ & $0$ & $0$ & $0$ \\ \hline
\end{tabular}
\caption{The cross product in the basis $(e_k)_{k=3}^{-3}$.}
\label{Table:ComplexCrossProductBasis}
\end{table}

We will consider the following basis a \emph{model} $\R$-cross product basis. 
\begin{align}\label{ModelRCrossProductBasis}
    X^\R = \left( \frac{\mathbf{i}+\mathbf{l}\mathbf{i}}{\sqrt{2}},\, \frac{\mathbf{j}-\mathbf{l}\mathbf{j}}{\sqrt{2}}, \,\frac{\mathbf{k}-\mathbf{l}\mathbf{k}}{\sqrt{2}}, \,\mathbf{l}, \,\frac{\mathbf{k}+\mathbf{l}\mathbf{k}}{\sqrt{2}}, \,\frac{\mathbf{j}+\mathbf{l}\mathbf{j}}{\sqrt{2}}, \,\frac{\mathbf{i}-\mathbf{l}\mathbf{i}}{\sqrt{2}} \right)
\end{align}
For the representation $\g_2' \rightarrow \gl_7\R$ in the basis $X^\R$ (slightly re-normalized), see \cite[page 89, (5.5)]{Bar10}.

$\F$-cross product bases are rather constrained in structure. The quadratic form $q$ on $\imoct^\F$ is automatically \emph{anti-diagonal}, so that $x_k \,\bot\, x_{l}$ unless $k = -l$ \cite[Proposition 2.3.6]{Eva24b}. 
We shall see $\R$-cross product bases relate to annihilators and the $\Gtwosplit$-flag manifolds. See Figure \ref{Fig:Apartment}, which shows how an $\R$-cross product basis $(x_k)_{k=3}^{-3}$ naturally describes an apartment in the visual boundary of the $\Gtwosplit$-symmetric space.

\subsubsection{Stiefel Triplet Models}\label{Subsec:StiefelModels}

In this section, we present two contrasting Stiefel triplet models for $\Gtwo^\F$. 
Each model is useful for understanding (and verifying transitivity) of the action of $\Gtwo^\F$ on different homogeneous spaces. 

The first Stiefel model helps to understand the action on spacelike or timelike vectors. Define 
\begin{align*}
    V_{(+,+,-)}(\imoct^\F) = \{(u,v,w) \in (\imoct^\F)^3\mid q(u)&=q(v)=-q(w)=+1, \\
 u\cdot v&=0,u\cdot w=0, v\cdot w =0, (u\times v) \cdot w =0\}.
\end{align*}
An element $(x,y,z) \in V_{(+,+,-)}(\imoct^{\F})$ is precisely one which is the same as $(\mathbf{i},\mathbf{j},\mathbf{l})$ up to the $\Gtwo^\F$-action. This Stiefel model is a well-known $\Gtwo^\F$-torsor \cite[Remark 5.13]{Fon18}. 
\begin{proposition}[First Stiefel Model]\label{Prop:FirstStiefel}
The group $\Gtwo^\F$ acts simply transitively on $V_{(+,+,-)}(\imoct^{\F})$. 
\end{proposition}

The idea of the proof can be summarized as follows for $\F=\R$, but works for $\F=\C$ as well. 
Any given triple $p =(u,v,w) \in V_{(+,+,-)}(\imoct)$ extends to a basis $\Gtwosplit$-equivalent to the multiplication basis $\mathcal{M}$ by:
\[ B_p := (u, \,v, \,u\times v,\, w,\, w\times u, \,w\times v, \,w\times (u\times v)).\]
Indeed, for $p_0=(\mathbf{i},\mathbf{j},\mathbf{l})$, the multiplication basis is given by $\mathcal{M}=B_{p_0}$. 
Now, for arbitrary $p=(u,v,w) \in V_{(+,+,-)}(\imoct)$, there is a unique transformation $\Psi \in \Gtwosplit$ such that $\Psi \cdot (\mathbf{i},\mathbf{j},\mathbf{l}) = (u,v,w)$ and moreover this transformation $\Psi$ is the unique linear transformation such that $\Psi \cdot \mathcal{M} = B_p$. 

\begin{remark}
The idea of using Stiefel models to understand the action of a Lie group $G$ is not special to $\Gtwosplit$, though it is especially useful in this case. For example, $G =\SO(n)$ acts simply transitively on $V_{n-1}(\R^n)$, the Stiefel manifold of orthonormal $(n-1)$-tuples in $\R^{n,0}$. Each such tuple $\boldsymbol{v}= (v_k)_{k=1}^{n-1}$ extends uniquely to an oriented orthonormal basis $B_{\boldsymbol{v}}=(v_k)_{k=1}^n$ of $\R^{n,0}$ and simple transitivity holds as above: $\psi \in \SO(n)$ satisfies $\psi \cdot \boldsymbol{v}=\boldsymbol{v}'$ if and only if $\psi \cdot B_{\boldsymbol{v}}= B_{\boldsymbol{v}'}$. 
\end{remark}

To understand the action of $\Gtwo^\F$ on isotropic vectors in $\imoct^\F$, one needs another model. The next Stiefel model from \cite{BH14} serves this purpose. Define the following set $\mathcal{N}$ of pairwise orthogonal isotropic triples: 
\[ \mathcal{N}_{\F} = \{ (u,v,w) \in Q_0(\imoct)^{\F}\mid u\cdot v=u\cdot w = v\cdot w=0, \,\Omega(u,v,w)= \sqrt{2}\}.\]
One can replace the constant $\sqrt{2}$ with any fixed non-zero constant and the result remains true.
\begin{lemma}[Null Triplet Model {\cite{BH14}}]\label{Lem:NullStiefel}
The group $\Gtwo^\F$ acts simply transitively on $\mathcal{N_\F}$. 
\end{lemma}
The idea of the proof is very similar to that of Proposition \ref{Prop:FirstStiefel}. Indeed, any triple $n=(u,v,w) \in \mathcal{N}_{\mathbb{F}}$ extends naturally to an $\mathbb{F}$-cross product basis $B_n$, and the bases $B_n$ and $B_{n'}$ are $\Gtwo^\F$-equivalent for any given pair $n,n' \in \mathcal{N}$. Here, we define 
\begin{align}\label{NullStiefelExtension}
    B_{n} = (u\times v,\,u,\,v,\,(u\times v) \times w,\,u \times w, \,v \times w, \,w).
\end{align}
Observe that since $B_n$ is determined from $n$ by the cross product, $\Psi \in \Gtwo^\F$ satisfies $\Psi \cdot n= n'$ if and only if $\Psi \cdot B_n = B_n'$. 

Now, to prove that for any $n, n'\in \mathcal{N}_{\F}$, there is a(n obviously unique) map $\Psi \in \Gtwo^\F$ such that $\Psi\cdot B_n= B_{n'}$, we can make a further observation. Indeed, $B_n$ in \eqref{NullStiefelExtension} is in fact, an $\F$-cross product basis (though not stated in such language in \cite{BH14}). The key to prove such $\Psi \in \Gtwo^{\F}$ exists is to show that the structure constants $(c_{k,l})$ and $(c_{k,l}')$ determined by $x_k \times x_l=c_{k,l}x_{k+l}$ and similarly for $c_{k,l}'$, do in fact satisfy $c_{k,l}= c_{k,l}'$. This is precisely what Baez \& Huerta prove, verifying Lemma \ref{Lem:NullStiefel}. 
See \cite[Section 2.3.2]{Eva24b} for further details on $\F$-cross product bases. 

\subsubsection{Basic Lie Theory of \texorpdfstring{$\g_2'$}{g2'}}\label{Subsec:LieTheory}

Recall from Subsection \ref{Subsec:G2Defn} that the group $\Gtwo^\F$ is described by $\Gtwo^\F = \Aut(\imoct^\F, \times)$. As a consequence, one finds a description of the Lie algebra as infinitesimal symmetries (\emph{derivations}) of the cross product:
\[ \g_2^\F = \Der(\imoct^\F, \times) \ = \{\psi \in \gl(\imoct^\F) \mid \psi(x\times y) = \psi(x) \times y+ x\times \psi(y) \}.\]
Let $X$ be an $\F$-cross product basis for $\imoct^\F$. One then immediately notes that the following is a Cartan subalgebra of the rank two Lie algebra $\g_2^\F$, which is $\R$-split for $\F = \R$:
\begin{align}\label{g2Cartan}
\mathfrak{a}_\F = \{ X = \diag(r+s,r,s,0,-s,-r,-r-s) \mid r,s \in \F\} < \g_2^\F .
\end{align}
 
\begin{figure}[ht]
    \centering
    \includegraphics[width=0.50\linewidth]{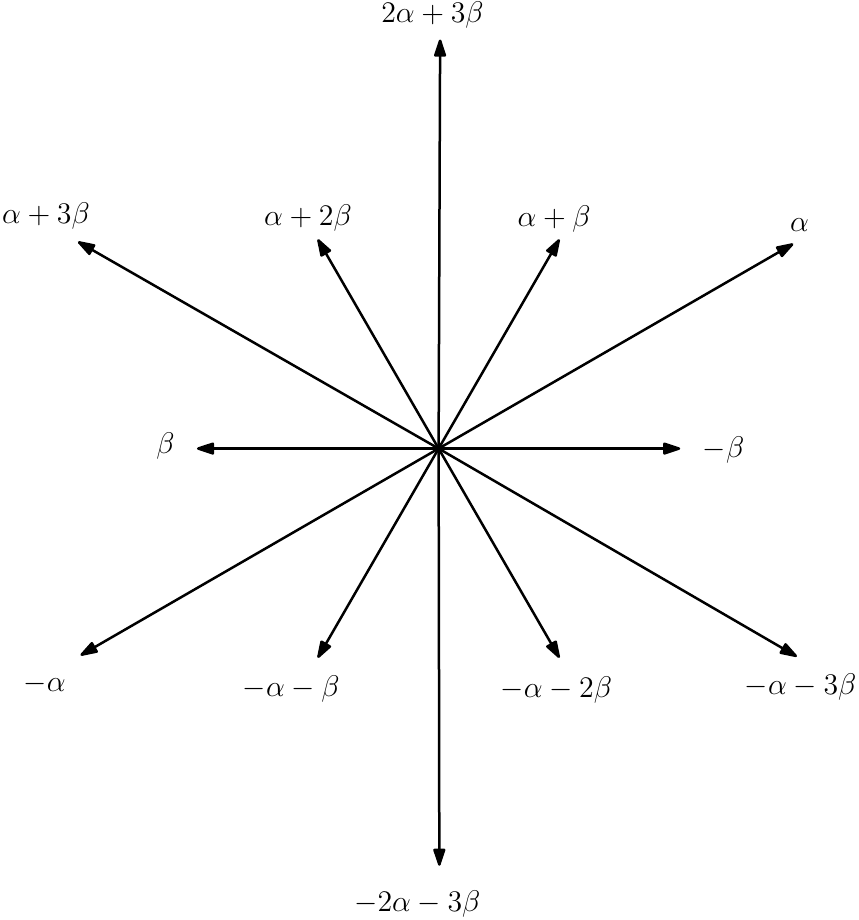}
    \caption{\emph{\small{The $G_2$ root system.}}}
    \label{fig:G2rootsystem}
\end{figure}

The associated $\F$-root system $\Sigma(\g^\F,\mathfrak{a}_\F) \subset \Hom_{\F}(\mathfrak{a}_\F,\F)$ are each of type $G_2$, as displayed in Figure \ref{fig:G2rootsystem}. In particular, we have the root space decompositions:
\[ \g_2^\F = \mathfrak{a}_\F \oplus \bigoplus _{\sigma \in \Sigma}\g_{\sigma}^\F.\]
Here the one-dimensional root spaces $\g_{\sigma}^\F \subset \g_2^\F$ are the usual  simultaneous eigenspaces of $\mathfrak{a}_\F$: 
\[ \g_{\sigma}^\F = \{ X \in \g_2^\F \mid [t,X]=\sigma(t)X, \forall t \in \mathfrak{a}_\F\}.\] 

Parametrizing $\mathfrak{a}_\F$ as in \eqref{g2Cartan}, we will frequently use the convention that $\beta:= s^*$ and $\alpha:= r^*-s^*$ are choices of simple short and long roots of $\mathfrak{a}_\F$, respectively. We will denote $\Delta:= \{\alpha, \beta\}$ as well as $\Sigma= \Sigma^+ \sqcup \Sigma^-$ the induced partition of the root system into positive and negative roots relative to $\Delta$. See Figure \ref{fig:G2regularity} for the corresponding (closed) Weyl chamber relative to this choice. Note that the Weyl group $ W= N(\mathfrak{a})/Z(\mathfrak{a})$ is the dihedral group of order twelve. In particular, $W$ acts transitively on the short roots and on the long roots. 

Next, we shall describe the form of a short root vector $e_{-\beta} \in \g_{-\beta}$ and a long root vector $e_{-\alpha} \in \g_{-\alpha}$ in $\g_2^\C$. Up the action of the Weyl group, this allows one to obtain root vectors in each root space $\g_{\sigma}$. We shall need these two root vectors for the discussion of $\Gtwosplit$-Higgs bundles in Section \ref{Sec:CyclicG2'Higgs}. 

Let $X=(e_i)_{i=3}^{-3}$ be the model $\C$-cross product basis from \eqref{BaragliaBasis}. The following linear transformations $E_{\gamma} \in \g_{\gamma}$, for $\gamma \in \{-\beta, -\alpha \}$, are root vectors. To avoid writing $(7\times 7)$-matrices, we suggestively express these linear transformations as diagrams, encoding their action on the basis $X$.
\begin{align} \label{longroot}
   E_{-\alpha} &= \left[  e_{3} \stackrel{0}{\longrightarrow} e_2 \stackrel{1}{\longrightarrow}  e_1\stackrel{0}{\longrightarrow} e_0\stackrel{0}{\longrightarrow} e_{-1}\stackrel{1}{\longrightarrow} e_{-2}\stackrel{0}{\longrightarrow} e_{-3} \right]\\
   \label{shortroot}
   E_{-\beta} &= \left[  e_{3} \stackrel{1}{\longrightarrow} e_2 \stackrel{0}{\longrightarrow}  e_1\stackrel{-\sqrt{2}i}{\longrightarrow} e_0\stackrel{-\sqrt{2}i}{\longrightarrow} e_{-1}\stackrel{0}{\longrightarrow} e_{-2}\stackrel{1}{\longrightarrow} e_{-3} \right]
\end{align}
Note that the above expressions are meant to imply $e_{-3}$ lies in the mutual kernel of $E_{-\alpha},E_{-\beta}$. 
We now explain a way to verify these linear maps are derivations that serves the dual purpose of describing $\g_2^\F$ via its action on null Stiefel triplets. 

Suppose that $(\varphi_t)_{{t \in (-\varepsilon, \varepsilon)}} \subset \Gtwo^\F$ is a smooth curve. We can describe each map $\varphi_t$ by a triple $(u_t, v_t,w_t) \in \mathcal{N}_{\F}$ relative to a background choice of $n = (u,v,w) \in \mathcal{N}_{\F}$. That is, $\varphi_t$ is the unique $\Gtwo^\F$ transformation such that $\varphi_t\cdot(u,v,w) =(u_t,v_t,w_t)$. We can then differentiate the seven relations of the triple $(u_t,v_t,w_t)$ to obtain infinitesimal analogues: if $(\dot{u},\dot{v},\dot{w})=\frac{d}{dt}\big|_{t=0}(u_t,v_t,w_t)$, then 
\begin{align}\label{InfinitesimalDerivation}
\begin{cases}
    \langle u, \dot{u}\rangle =0\\
    \langle v, \dot{v}\rangle =0\\
    \langle w, \dot{w}\rangle =0\\
    \langle \dot{u}, v\rangle + \langle u, \dot{v}\rangle =0\\
    \langle \dot{u}, w\rangle + \langle u, \dot{w}\rangle =0\\
    \langle \dot{v}, w\rangle + \langle v, \dot{w}\rangle =0\\
    \Omega(\dot{u},v,w)+\Omega(u,\dot{v}, w)+\Omega(u,v,\dot{w})=0.
    \end{cases}
\end{align}
Here, we write $\langle u,v\rangle$ to avoid unsightly expressions such as $u\cdot \dot{u}$. 
Lemma \ref{Lem:DerivationExtension} shows these seven constraints are precisely the requirements for a linear map $\phi:\spann_{\F}(n) \rightarrow \imoct^\F$, described by 
\[ \phi = \begin{cases} u &\longmapsto \dot{u} \\
                        v &\longmapsto \dot{v} \\
                        w &\longmapsto \dot{w} \\
                        \end{cases} ,\]
to extend to a derivation of $\imoct^\F$. Note the derivation condition entails that if $\phi$ extends to $\tilde{\phi} \in \Der(\imoct^\F)$, then the extension $\tilde{\phi}$ is unique and it must be given in the $\F$-cross product basis $B_{n}$ from \eqref{NullStiefelExtension} by:
\begin{align}\label{Derivation_Extension}
\tilde{\phi}=
\begin{cases} u\times v &\longmapsto \dot{u}\times v+ u \times \dot{v} \\
u &\longmapsto \dot{u} \\
v &\longmapsto \dot{v} \\
(u\times v)\times w&\longmapsto (\dot{u}\times v+ u \times \dot{v}) \times w+(u\times v) \times \dot{w}\\
u\times w &\longmapsto \dot{u}\times w + u\times \dot{w} \\
v \times w &\longmapsto \dot{v}\times w+ v\times \dot{w}\\
w &\longmapsto \dot{w} 
\end{cases}
\end{align}
In principle, one must check $\binom{7}{2} =21$ equations to verify the derivation condition. This lemma drastically simplifies the computations necessary. 
\begin{lemma}[$\g_2^\F$ via Null Triples]\label{Lem:DerivationExtension}
Let $n=(u,v,w) \in \mathcal{N}_\F$ be a null triple. Define 
\[\g_n \subset \Hom_{\F}(\spann_{\F}(n), \imoct^\F)\] as the subspace of linear maps satisfying \eqref{InfinitesimalDerivation}. The restriction map $\pi:\g_2^\F \rightarrow \g_n$ by $\Phi \mapsto \Phi|_{\spann(n)}$ is a linear isomorphism. The inverse $\tilde{\phi} = \pi^{-1}(\phi)$ of $\phi \in \g_n$ is given by \eqref{Derivation_Extension}.
\end{lemma}

\begin{proof}
The map $\pi$ is well-defined by the argument preceding the lemma. Any derivation $\psi \in \g_2^\F$ must obtain the form \eqref{Derivation_Extension}, which immediately implies $\pi$ is injective. 

The map $\pi$ is a linear isomorphism if $\dim_{\F}(\g_n) = 14$. Expressing $(\dot{u},\dot{v},\dot{w})$ in the background basis $B_{n}$, one easily verifies this is the case. 
\end{proof}
While the expression \eqref{Derivation_Extension} for the extension may appear nightmarish at first sight, in some  situations it is tractable to compute. Indeed, considering the transformations $E_{\gamma}$, for $\gamma \in \{-\beta, -\alpha\}$, from \eqref{longroot}, \eqref{shortroot}, respectively, one finds that for $n = (e_2,e_1,e_{-3}) \in \mathcal{N}_{\F}$, the map $E_{-\beta}$ corresponds to $(\dot{u},\dot{v},\dot{w})= (0, -\sqrt{2}i e_0, 0)$ and $E_{-\alpha}$ corresponds to $(\dot{u},\dot{v},\dot{w})= (e_1, 0, 0)$. Hence, $E_{\gamma}|_{\spann(n)} \in \g_n$. One also verifies each of $E_{-\alpha}, E_{-\beta}$ is indeed given by the extension formula \eqref{Derivation_Extension}, meaning $E_{-\alpha},E_{-\beta} \in \g_2^\C$. 

\subsubsection{Annihilators}\label{Subsec:Annihilators}
We now introduce \emph{annihilators}, a notion fundamental to $\Gtwosplit$-flag manifolds. 

\begin{definition}
Fix $u \in \imoct$. Then define the \textbf{annihilator} of $u$ as 
\[ \Ann(u) := \ker(\mathcal{C}_u)=\{ v \in \imoct \mid u \times v=0\}.\]
\end{definition}
If $q(u)\neq 0$, then $\Ann(u) = \R\{u\}$, by the double cross product identity \eqref{DCP}. However, if $u$ is isotropic, then its annihilator is three-dimensional. 
We include a proof of the following indispensable result, also shown in \cite{BH14}. Here, we need the model space $\Gr_{(3,0)}^\times(\imoct)$:
\begin{align*}
    \Gr_{(3,0)}^\times (\imoct )= \{ U \in \Gr_{(3,0)}(\imoct) \mid U\times_{\imoct}U=U\}.
\end{align*}
 In fact, we can describe the annihilator as the graph of a unique linear map relative to a point in the above model space. We shall frequently appeal to this result. 
 \begin{proposition}[Annihilator 3-planes]\label{Prop:AnnihilatorGraph}
 Let $x \in Q_0(\imoct)$ and $U \in \Gr_{(3,0)}^\times(\imoct)$.  
 \begin{enumerate}[noitemsep]
     \item $\Ann(x)$ is a three-dimensional, isotropic subspace of $\imoct$. 
     \item The orthogonal projection map $\pi_{U}: \Ann(x) \rightarrow U$ is a linear isomorphism. Moreover, $\Ann(x)$ is the graph of a unique linear map $\phi:U \rightarrow U^\bot$ that is an anti-isometry onto its image. 
     \item Write $cx=u+z$ for some unit spacelike vector $u \in Q_+(U) $, unit timelike vector $z \in Q_-(U^\bot)$, and nonzero scalar $c \in \R_+$. Then for any spacelike vector $v \in U$, 
     \[ \phi(v) = -z(uv).\]
 \end{enumerate}
 
 \end{proposition}
 
 \begin{proof}
 Write $cx= u+z$ as in (3). First, note by \eqref{DCP} that no nonzero element in $U$ or $U^\bot$ can lie in $\Ann(x)$. Hence, any element $y \in \Ann(x)$ has nonzero orthogonal projection to $U$ and to $U^\bot$. 
 
 Now, take any $v \in U$. We wish to find conditions on $w \in U^\bot$ such $(v+w) \in \Ann(x)$. Using the $\Z_2$-cross product grading $\imoct = U\oplus U^\bot$, one finds 
 \begin{align}\label{AnnihilatorEquivalence}
     (u+z) \times (v+w) = 0 \Leftrightarrow \begin{cases} u\times v+ z \times w =0  & (U)\\ 
 z \times v+u \times w= 0 & (U^\bot) \end{cases}. 
 \end{align}
 If $u=v$, then $z= w$ is forced by the equations \eqref{AnnihilatorEquivalence}. 
 
 Otherwise, we may suppose $v \bot u$ up to moving $v$ inside the span of $\lbrace u,v\rbrace$. The double cross product identity \eqref{DCP} 
 allows us to re-write the right-hand side of \eqref{AnnihilatorEquivalence}.
 Apply $\mathcal{C}_u$ to equation $(U^{\bot})$, then use  \eqref{DCP} to simplify: 
 \[ w =_{\eqref{DCP}} -u\times(u\times w)= u\times (z\times v).\] 
 Hence, we observe that 
\[ w \cdot z = \Omega(u, z\times  v, z)=\Omega(z\times v, z, u) = 0.\]
Thus, $w$ is orthogonal to $z$. With similar reasoning on the equation $(U)$, one then concludes $(v+w) \in \Ann(x)$ if and only if $w$ satisfies the seemingly overdetermined system 
\begin{align} \label{AnnihilatorEquations}
    \begin{cases} 
        w = -z\times (u\times v)\\
        w=u\times (z \times v).
    \end{cases} 
\end{align}
However, the two are equal. To see the equality, recall the \emph{associator} \eqref{Associator} is alternating. Hence,
 \begin{align}\label{AssociatorApplication}
     z(uv) + u(zv) = z(uv) - (zu)v + (zu)v + u(zv) = [z,u,v] + [u, z, v] = 0. 
 \end{align}
 Consequently, using that $z(uv) = z\times (u\times v)$ since these elements satisfy $(u,v,z) \in V_{(+,+,-)}(\imoct)$, we obtain
 \[ z\times(u\times v) = z(uv) =_{\eqref{AssociatorApplication}}-u(zv) = -u \times (z\times v).\]
Thus, for any $v \in U$, the equation \eqref{AnnihilatorEquivalence} has a unique solution $w \in U^\bot$. 
 
We have just seen that
 \[ w = -z\times (u\times v)= -z(uv) \] 
when $u \bot v$. However, $w = -z(uv)$ holds as well $v \in \R\{u\}$, so the formula holds 
for all $v\in U$. Thus, $\Ann(x)$ is the graph of the map $\phi: U \rightarrow U^\bot$ given by $\phi(v) = -z(uv)$. Note that $\phi$ is an anti-isometry due to 
 $q(-z(uv)) = q(z)q(u)q(v) = -q(v)$, using multiplicativity of $q$. This verifies (2) and (3). Since $\phi$ is an anti-isometry, $\Ann(x)$ is isotropic, proving (1). 
 \end{proof} 
 Note that any isotropic 3-plane in $\R^{3,4}$ obtains the same form as a graph; the $\Gtwosplit$-geometry facts are (1) and (3) in Proposition \ref{Prop:AnnihilatorGraph}. 
 
Later, we shall study \emph{annihilator photons} $\omega \in \Gr_2(\imoct)$, which are two-planes for which the cross product vanishes identically.
Proposition \ref{Prop:AnnihilatorGraph} immediately implies the following. 
\begin{corollary}[Annihilator Photons as Graphs]\label{Cor:AnnihilatorPhotonGraph}
Let $U \in \Gr_{(3,0)}^\times(\imoct)$. 
\begin{enumerate}
    \item For annihilator photon $\omega$, there is a unique rank two linear map $\varphi_{U,\omega}:U \rightarrow U^\bot$ whose graph over $W:= \pi_{U}(\omega)$ is $\omega$, where $\pi_U:\imoct \rightarrow U$ denotes orthogonal projection.  
    \item Moreover, $\omega \in \Gr_2(\imoct)$ 
is an annihilator photon if and only if it obtains the form $\omega = \spann \{u+z, v+(uv)z\}$ for 
an orthonormal pair $u,v \in Q_+(U)$ and $z \in Q_-(U^\bot)$. 
\end{enumerate}
\end{corollary}

There is a useful relation between $\R$-cross product bases and annihilators. 
\begin{remark}
A convenient feature of an $\R$-cross product basis $(x_k)_{k=3}^{-3}$ is $q(x_3)=0$ and $\Ann(x_3) = \spann_\R \langle x_3, x_2,x_1\rangle$. In fact, $x_k$ is isotropic for $k\neq0$ and its annihilator can be read off from Figure \ref{Fig:Apartment}. 
\end{remark}

\subsubsection{Subgroups of \texorpdfstring{$\Gtwosplit$ \& Homogeneous Spaces}{G2'}}\label{Subsec:SubgroupHomogeneous}

In this section, we describe some important $\Gtwosplit$-subgroups, namely the maximal compact $K$, and the parabolic subgroups $P_{\beta}, P_{\alpha},P_{\Delta}$. We remark on some other subgroups as well. 

We first discuss the maximal compact subgroup.  To this end, it is useful to introduce the following model space of spacelike 3-planes closed under cross product: 
\begin{align}\label{ModelSpaceForX}
    \Gr_{(3,0)}^\times (\imoct )= \{ U \in \Gr_{(3,0)}(\imoct) \mid U\times_{\imoct}U=U\}.
\end{align}
For example, $\spann_{\R} \{ \mathbf{i},\mathbf{j},\mathbf{k}\} \in \Gr_{(3,0)}^\times(\imoct)$.
\begin{proposition}[Maximal Compact $K < \Gtwosplit$]
Fix $U \in \Gr_{(3,0)}^\times(\imoct)$. The stabilizer subgroup $K = \Stab_{\Gtwosplit}(U)$ is a 
maximal compact subgroup of $\Gtwosplit$ and moreover $K \cong \SO(4)$. 
\end{proposition}

In summary, $K < \Gtwosplit$ stabilizes certain space+time splittings $\R^{3,4} = U\oplus U^\bot$ compatible with the cross product. 
For a proof, see \cite[Lemma 3.2]{ER25}. As a consequence, the Riemannian symmetric space $\X = \X_{\Gtwosplit}$ is naturally identified with the model space $\Gr_{(3,0)}^\times(\imoct)$. A more comprehensive description can be found in \cite[Lemma 3.5]{Eva24a}. 
\begin{remark}
The $\Gtwosplit$-symmetric space $\X$ is also naturally identified with the Grassmannian of quaternionic subalgebras $\Gr_{\Ha}(\Oct')$ via the equivariant map $\Gr_{(3,0)}^\times(\imoct) \rightarrow \Gr_{\Ha}(\Oct')$ by $U \mapsto \R\oplus U$. 
\end{remark}

Now, we introduce a certain key splitting of $\imoct$ that we shall repeatedly use. 

\begin{definition}[Frenet frame]
\label{defn:FrenetFrame}
We call a tuple $(\mathscr{L},T,N,B)$ a \textbf{Frenet frame} when 
\begin{itemize}
    \item $\imoct =\mathscr{L} \oplus T \oplus N \oplus B,$
    \item $\mathscr{L},T,N,B$ are subspaces of alternating signatures $(1,0)$, $(0,2)$, $(2,0)$, $(0,2)$,
    \item each of $T,N,B$ are closed under cross product with $\mathscr{L}$.
\end{itemize}
\end{definition}

In fact, a Frenet frame naturally arises
from a choice of multiplication basis for $\imoct$. For convenience, we take the standard multiplication basis $\mathcal{M} = (1,\mathbf{i},\mathbf{j},\mathbf{k},\mathbf{l},\mathbf{li},\mathbf{lj},\mathbf{lk})$, as in $\S$\ref{Subsec:G2Defn}. This yields the model Frenet frame splitting: 
\begin{align}\label{ModelFrenetFrame}
    \begin{cases}
    \mathscr{L} = \R \{ \mathbf{i} \},\\
    T = \spann \{ \mathbf{l}, \mathbf{li} \},\\
    N = \spann \{ \mathbf{j}, \mathbf{k} \},\\
    B = \spann \{ \mathbf{lj}, \mathbf{lk} \}.
\end{cases}
\end{align} 
More generally, if $p=(u,v,z) \in V_{(+,+,-)}(\imoct)$, then we obtain the Frenet frame 
\[ \mathscr{L} = \R \{ u\}, \;T = \spann \{ z, uz \}, \; N= \spann \{ v, uv\}, \; B = \spann \{ zv, z(uv)\}.  \] 
At some points, we may use a slight refinement of a Frenet frame as follows: we consider tuples $(x,T,N,B)$, where $x \in Q_+(\imoct)$ and 
$(\R\{x\}, T,N,B)$ is a Frenet frame in the sense of Definition \ref{defn:FrenetFrame}. We briefly remark that the set of refined Frenet frames is a $\Gtwosplit$-homogeneous space of interest. 
\begin{proposition}[\cite{Eva24a}]
The space of refined Frenet frames $(x,T,N,B)$ is isomorphic as a $\Gtwosplit$-space to $\Gtwosplit/\mathscr{T}$, where $\mathscr{T}< K$ is a maximal torus in the maximal compact subgroup. 
\end{proposition}
Next, we discuss the action of $\Gtwosplit$ on each of its orbits in projective space $\mathbb{P}\imoct$. To this end, partition projective space into three parts, depending on the $q$-signature of the lines:
\[ \mathbb{P}(\R^{3,4}) = \mathbb{S}^{2,4}\sqcup \Ein^{2,3} \sqcup \mathbb{H}^{3,3},\]
where $\mathbb{S}^{2,4}= \mathbb{P}Q_+(\imoct)$, $\Ein^{2,3} = \mathbb{P}Q_0(\imoct)$, $\Ha^{3,3} = \mathbb{P}Q_-(\imoct)$ are the sets of positive, isotropic, and negative lines in $\imoct$, respectively. 

\begin{proposition}[Orbits on Quadrics]
$\Gtwosplit$ acts transitively on $Q_{\epsilon}(\imoct)$ for $\epsilon \in \{+,0,-\}$. 
\end{proposition}

\begin{proof}
The result follows immediately from the Stiefel models in Proposition \ref{Prop:FirstStiefel} and Lemma \ref{Lem:NullStiefel}. 
\end{proof}

The corresponding stabilizer subgroups are as follows: if $x \in Q_\epsilon(\imoct)$, then 
\begin{align}\label{eq:QuadricStabilizers}
    \begin{cases}
    \Stab_{\Gtwosplit}(x) \cong \SU(1,2) & \epsilon =+\\
    \Stab_{\Gtwosplit}(x) \cong \SL(3,\R) & \epsilon =-
\end{cases} 
\end{align}
The former is explained in \cite[Proposition 3.9]{ER25}, where the sub-symmetric space of $\SU(1,2)$ is also studied. The latter is explained in \cite[Corollary 3.11]{CT24}. \medskip 

The quadric hypersurface $\quadric= Q_+(\imoct)$ admits a canonical $\Gtwosplit$-invariant almost-complex structure as follows. For $x \in \quadric$, \eqref{DCP} says $\mathcal{C}_x|_{x^\bot}^{2}=-\id|_{x^\bot}$, so that $J: \T\quadric \rightarrow \T\quadric$ by $J_x = \mathcal{C}_x$ is an almost-complex structure on $\T_x\quadric=[x^\bot \subset \imoct]$. See Figure \ref{Fig:ACStructure}. The almost-complex structure $J$ is not integrable \cite{CT24}. This is exactly analogous to the situation for $\mathbb{S}^6$, which admits a canonical $\Gtwo^c$-invariant almost-complex structure from $(\mathrm{Im}(\Oct),\times)$ that is non-integrable \cite{Leb87}. 

\begin{figure}[ht]
\begin{center}

\tdplotsetmaincoords{70}{110} 
\begin{tikzpicture}[tdplot_main_coords,scale=3]
\def\r{1}

\coordinate (O) at (0,0,0);
\coordinate (N) at (0,0,\r);

\draw[tdplot_screen_coords] (0,0) circle (\r); 

\begin{scope}[tdplot_screen_coords]
  \shade[ball color=blue!30, opacity=0.6] (0,0) circle (\r);
\end{scope}

\tdplotdrawarc[]{(O)}{\r}{-68}{112}{anchor=south}{}
\tdplotdrawarc[dotted]{(O)}{\r}{112}{112+180}{anchor=south}{}

\filldraw[fill=blue!10,opacity=0.6] 
   (-0.8,-0.8,\r) -- (0.8,-0.8,\r) -- (0.8,0.8,\r) -- (-0.8,0.8,\r) -- cycle;
   \node[left] () at (0.8,-0.8,\r) {$\T_{\mathbf{k}}\mathbb{S}^2$};
\fill (N) circle (0.5pt) node[below right]{};
\draw[->,thick,red] (N) -- ++(-0.6,0,0) node[above]{$\mathbf{j}$};
\draw[->,thick,red] (N) -- ++(0,0.6,0) node[above]{$\mathbf{i}$};
\draw[->,thick,red] (N) -- ++(0,0,0.6) node[right]{$\mathbf{k}$};
\draw[->,thick] (-.05,.3,\r) arc (90:170:.3);
\node[above right] at (-.3*0.707, .3*0.707,\r)  ()    {$J|_{\mathbf{k}}$};
    \end{tikzpicture}
\end{center}
    \caption{A cartoon representing the almost-complex structure $J:\T\quadric \rightarrow \T\quadric$ instead on $\mathbb{S}^2=Q_+(\R^3)$. In both settings, the cross product endomorphism $\mathcal{C}_x$ of the position vector $x$ defines a distinguished rotation $J|_x:=\mathcal{C}_x$ by $\frac{\pi}{2}$ in the tangent space.}
\label{Fig:ACStructure}
\end{figure}
\begin{remark}
Note that if $(x,T,N,B)$ is a (refined) Frenet frame, then $\T_x\quadric= T\oplus N\oplus B$ is an orthogonal decomposition into complex lines.  
\end{remark}

Finally, we discuss the $\Gtwosplit$-flag manifolds, the central interest of this work.  The model spaces are as follows: 
\begin{align}
    \Ein^{2,3} &= \{ [x] \in \mathbb{P}(\R^{3,4}) \; | \; q(x) = 0 \},\\ 
    \Pho^\times  &= \{ \omega \in \Pho(\R^{3,4}) \; |\; \omega \times_{3,4} \omega = 0\}, \\ 
    \mathcal{F}_{1,2}^\times &= \{  \,(\ell, \omega) \in \Ein^{2,3} \times \Pho^\times | \; \ell \subset \omega \}.
\end{align}
Recall that a \emph{photon} $\omega \in \Pho(\R^{3,4})$ is an isotropic 2-plane. 
We call an element $\omega \in \Pho^\times$ an \emph{annihilator photon} and a pair $(\ell, \omega) \in \mathcal{F}_{1,2}^\times$ a \emph{pointed annihilator photon}. 

For additional clarity, we realize the identifications with these model spaces. 
For the statement, let $P_{\Theta}$ denote the parabolic subgroup associated to $\Theta \subset \Delta = \{\alpha, \beta\}$, using the notation from $\S$\ref{Sec:SymmSpace} with the same choices $(\mathfrak{a}, \Sigma, \Delta)$. Recall that $P_{\Theta}$ is the normalizer in $\Gtwosplit$ of the corresponding Lie algebra $\mathfrak{p}_{\Theta}$, given by the sum of the root spaces of non-negative $\Theta$-height. 

\begin{proposition}[$\Gtwosplit$-Flag Manifolds]
\label{Prop:G2Flags}
Let $X=(x_k)_{k=3}^{-3}$ be an $\R$-cross product basis for $\R^{3,4}$.  There are $\Gtwosplit$-equivariant diffeomorphisms
\[ \begin{cases}
    \Gtwosplit/P_{\beta} \cong \Ein^{2,3} \\
    \Gtwosplit/P_{\alpha} \cong \Pho^\times \\
    \Gtwosplit/P_{\Delta} \cong \mathcal{F}_{1,2}^\times,
\end{cases} \] 
given by $g P_{\Theta} \mapsto g\cdot p_{\Theta}$ for $\Theta \in \{ \alpha, \beta, \Delta \}$,
where $p_{\beta}= \langle x_3\rangle $, $p_{\alpha} = \langle x_3, x_2\rangle$, and $p_{\Delta} = (p_{\beta}, p_{\alpha})$. 
\end{proposition}

The transitivity of the action of $\Gtwosplit$ on $\mathcal{F}_{1,2}^\times$, and hence on $\Ein^{2,3}, \Pho^\times$ as well, follows immediately from Lemma \ref{Lem:NullStiefel}. The equalities $\mathfrak{p}_{\sigma} = \Stab_{\g_2'}(p_{\sigma})$, for $\sigma \in \{\alpha, \beta, \Delta\}$, on the infinitesimal level, can be verified easily using \cite[page 89, (5.5)]{Bar10}, the representation of $\g_2'$ in an $\R$-cross product basis. The equality on the group level is found with slightly more work. 
These geometric models are well-known  \cite{Bry98, LM03,MNS21, Man25}.

\subsection{Symmetric Spaces} \label{Sec:SymmSpace}
In this section, we introduce some basic terminology for symmetric spaces of non-compact type and their visual boundaries. We then introduce the relevant details on the $\Gtwosplit$-symmetric space. Standard references for symmetric spaces are Helgason \cite{Helg78} and Eberlein \cite{Ebe96}. 

The symmetric space $\X_G$ of a non-compact simple Lie group $G$ is a $\mathrm{CAT}(0)$-space and more specifically a \emph{Hadamard manifold}: a complete, non-positively curved, simply-connected Riemannian manifold. Consequently, $\X_{G}$ is topologically a cell, with $\exp_{x}:\T_x\X\rightarrow \X$ a diffeomorphism for any $x \in \X$. One natural way to compactify $\X_G$ is by attaching the visual boundary as follows. 

\begin{definition}
The \textbf{visual boundary} $\vis\X_G$ of the symmetric space $\X_G$ is the set of (unit speed) geodesic rays $\gamma:\R_{\geq 0}\to \X_G$, where two rays are equivalent if their images are at bounded Hausdorff distance.
\end{definition}

We now highlight some notation to be used frequently.
\begin{definition}
Let $\gamma_{x,v}$ denote the unique geodesic of $\X$ satisfying $\gamma_{x,v}(0) = x$ and $\gamma_{x,v}'(0)=v$. Similarly, for $\xi \in \vis\X$ and $x \in \X$, let $v_{x, \xi} \in \T^1_x\X$ be the unique tangent vector such that $[\gamma_{v_{x,\xi}}]=\xi$. 
\end{definition}Fix a point $x \in \X$. The map $\mathrm{T}^1_x\X \rightarrow \vis\X$ by $v\mapsto [\gamma_{x,v}]$ is a bijection. There is a unique topology on $\vis\X$ for which one (in fact, every) such map is a homeomorphism. The topology on $\overline{\X} = \X \sqcup \vis\X$ is such that a sequence $(x_n) \subset\X$ satisfies $x_n \rightarrow [\gamma]$ if and only for one (in fact, every) basepoint $x \in \X$, the geodesics rays $\gamma_{n}:[0,\infty)$ emanating from $x$ and passing through $x_n$ converge uniformly on compacta to the representative of $[\gamma]$ based at $x$. Going forward, we may write $\gamma(\infty) \in \vis\X$ to denote the equivalence class $[\gamma]$, for $\gamma:[0,\infty) \rightarrow\X$ a geodesic ray.  

One can define flag manifolds through their connection to the visual boundary $\vis\X$.

\begin{definition}
A subgroup $P < G$ is called a \textbf{parabolic} subgroup when $P = \Stab_{G}( [\gamma])$ for some point $[\gamma] \in \vis\X$. The associated homogeneous space $\mathcal{F}= G/P$ to a parabolic subgroup is called a \textbf{flag manifold}. 
\end{definition}

We now introduce terminology to combinatorially classify the points in the visual boundary. 

\begin{definition}
Let $v \in \mathrm{T}_x\X$. We say $v$ \textbf{points towards} $\mathcal{F}=G/P$ when the subgroup $\Stab_{G}([\gamma_{x,v}])$ is conjugate in $G$ to the subgroup $P$. 
\end{definition}

Let $G$ be a non-compact simple Lie group. 
Every parabolic subgroup $P < G$ obtains a standard form $P =P_{\Theta}$ up to conjugation. 
Fix a Cartan subalgebra $\mathfrak{a} < \g$, restricted real root system $\Sigma(\g, \mathfrak{a})$ with simple roots $\Delta$. 
Given any subset $\Theta \subset \Delta$ of simple roots, there is an associated standard parabolic subalgebra $\p_{\Theta}$. In particular, $\p_{\Theta}$ consists of the sum of all root spaces of non-negative $\Theta$-height. 

The associated Lie group $P_{\Theta}$ is the normalizer of $\p_{\Theta}$ in $G$. The parabolic subalgebras $\mathfrak{p}_{\{\alpha_i\}}$, for $\alpha_i \in \Delta$, are precisely the \emph{maximal} (proper) parabolic subalgebras of $\g$ up to conjugacy. 

When $P <G$ is a (proper) maximal parabolic subgroup, $P$ uniquely stabilizes a point $[\gamma] \in \vis\X$ and also $P$ uniquely stabilizes a point $f \in \mathcal{F} =G/P$. Thus, when $\Stab_{G}([\gamma]) = P$ is a maximal parabolic, for $\gamma = \gamma_{x,v}$, we shall say that $v$ \emph{points towards} the flag $f \in G/P$. In other words, the $G$-equivariant embedding $G/P\hookrightarrow \vis\X$ is canonical when $P$ is a maximal parabolic. Such embeddings exist when $P$ is non-maximal, but require a choice. 

The \emph{Tits angle} $\titsangle$ between two geodesic rays $f_1=[\gamma_1],f_2=[\gamma_2] \in \vis\X$ is defined by 
\[ \titsangle( [\gamma_1], [\gamma_2]):= \sup_{x \in \X } \measuredangle_x(v_{x,f_1}, v_{x,f_2}) \]
That is, we take the largest possible Riemannian angle, across all choice of basepoint, subtended by pairs of tangent vectors pointing towards the given points in the visual boundary. 

Fix a basepoint $o \in \X$. Then $K =\Stab_{G}(o)$ is a maximal compact subgroup, leading to a Cartan decomposition $\g = \frakk \oplus \frakp$. We can regard a Cartan subalgebra $\mathfrak{a} \subset \mathfrak{p}$ also as a subspace $\mathfrak{a} \subset \T_o\X$. Fix an (open) model Weyl chamber $\mathfrak{a}^+ \subset \mathfrak{a}$. The \emph{Cartan projection} is the map $\mu: \T\X \rightarrow \overline{\mathfrak{a}}^+$ that sends each tangent vector $X \in \mathrm{T}_x\X$ to the unique element $g\cdot X$ in its $G$-orbit that lies in $\overline{\mathfrak{a}}^+$. 

\begin{definition}
We call a tangent vector $X \in \mathrm{T}_x\X$ to be $\boldsymbol{\Theta}$-\textbf{regular} when $\mu(X)$ has the property that $\theta(\mu(X))\neq 0$ for all $\theta \in \Theta$.
\end{definition}
 In fact, $X$ points towards $\mathcal{F}_{\Theta} = G/P_{\Theta}$ exactly when $\mu(X)$ is $\Theta$-regular and $\alpha(\mu(X))=0$ for all $\alpha \in \Delta\setminus \Theta$ (cf. \cite[Proposition 10.64]{BH99}).

\subsubsection{The \texorpdfstring{$\Gtwosplit$}{G2'}-Symmetric Space}  \label{Subsec:G2SymmetricSpace}

We now discuss the $\Gtwosplit$-symmetric space $\X_{\Gtwosplit}$. We will concretely embed the associated flag manifolds $\Ein^{2,3}$ and $\Pho^\times$ in the visual boundary $\vis \X_{\Gtwosplit}$ and we describe \emph{pointing toward} $\Ein^{2,3}$ and $\Pho^\times$. We also discuss the related symmetric space $\X_{\SO(3,4)}$ of $\SO(3,4)$ and the discrepancy between pointing towards $\Ein^{2,3}$ in $\X_{\SO(3,4)}$ and in $\X_{\Gtwosplit}$. \medskip

The symmetric space $\X_{\SO(3,4)}$ associated with $\SO(3,4)\cong \SO(\imoct)$ can be interpreted as the collection of $3$-dimensional spacelike subspaces in $\R^{3,4}\cong \imoct$:
\[\X_{\SO(3,4)}\cong\Gr_{(3,0)}(\imoct). \]

The symmetric space $\X_{\GL(7,\R)} \cong\GL(7,\R)/\mathrm{O}(7,\R)$ associated with $\GL(7,\R)$ is the space of Euclidean metrics on $\R^7$. The inclusion $\SO(3,4)\subset \GL(7,\R)$ induces a totally geodesic embedding $\X_{\SO(3,4)}\to \X_{\GL(7,\R)}$. Let $q$ denote a fixed signature $(3,4)$-quadratic form. Concretely,  this embedding is the map that associates to $U\in \Gr_{(3,0)}(\R^{3,4})$ the Euclidean metric $h= q|_U \oplus (-q)|_{U^\bot}$.

The exceptional simple Lie group $\Gtwosplit$ is a subgroup of $\SO(3,4)$, hence its symmetric space embeds in a totally geodesic way inside the symmetric space of $\SO(3,4)$. Recall that in Subsection \ref{Subsec:SubgroupHomogeneous}, we interpreted $\X_{\Gtwosplit}$ as $\Gr_{(3,0)}^\times(\imoct) \subset \Gr_{(3,0)}(\imoct)$, which realizes this embedding concretely. 

We now recall a basic equivalence on the Riemannian metrics on $\mathbb{X}_{\Gtwosplit}$ and $\mathbb{X}_{\SO(p,q)}$ under their respective identifications with $\Gr_{(3,0)}^\times(\imoct)$ and $\Gr_{(p,0)}(\R^{p,q})$. 
Below, we use the standard and canonical identification $\T_U\Gr_{(p,0)}(\R^{p,q}) \cong \Hom(U,U^\bot)$. 
 
 \begin{proposition}[Metrics in Model Spaces]
 \label{Prop:ModelSpaceMetric}
Let $Q$ denote a fixed signature $(p,q)$-quadratic form. 
Up to a (universal) constant $c > 0$, the Riemannian metric $g$ on $\Gr_{(p,0)}(\R^{p,q})$ is given by 
\[ c\,g_{P}(\phi, \psi) = -\tr(\phi^{*Q} \circ \psi), \]
where $A^{*Q}$ denotes the $Q$-adjoint. Consequently, the same holds for $\Gr_{(3,0)}^\times(\imoct)$.
 \end{proposition}
 
 \begin{proof}
Write $G: =\SO(p,q)$ and $Q = Q_U \oplus Q_V$ in some orthogonal splitting $\R^{p,q} = U \oplus V$. Let $A^{*Q}$ denote the adjoint of $A: U \rightarrow V$ in the sense of $\langle Ax, y \rangle_Q = \langle x, A^{*Q}y \rangle_Q$. In particular, in coordinates, $A^{*Q} = q_U^{-1}A^Tq_V$. In this splitting, if $X = \begin{pmatrix} 0 & B \\ A & 0\end{pmatrix} \in \mathfrak{so}(p,q)=\g$, then $B =-A^{*Q}$. 

Now, a Cartan decomposition $\g = \frakk \oplus \frakp $ is obtained via $\frakk = \left \{ \begin{pmatrix} * & 0 \\ 0 & * \end{pmatrix} \in \g  \right \}$ and 
$\frakp = \left \{ \begin{pmatrix} 0 & * \\ * & 0 \end{pmatrix} \in \g  \right \}$. By the previous paragraph, $X \in \p$ obtains the form $X = A-A^{*Q}$ for some linear map $A: U \rightarrow U^\bot=V$. Hence, the Killing form $B$ on $\p$ obtains the form $B( X, X ) = -\tr( A^* \circ A)$, up to a positive constant. By the $G$-equivariance of the identifications $(G/K, B) \cong \X_{\SO(p,q)} \cong (\Gr_{(p,0)}(\R^{p,q}), g)$ the claim follows.  
 \end{proof}
 
 \begin{corollary}[Orthogonality in $\mathbb{X}_{\SO(p,q)}$] \label{Cor:Orthogonality}
 Let $\phi, \psi \in \T_{U} \Gr_{(p,0)}(\R^{p,q})$. Then $\langle \phi, \psi \rangle_\X= 0$ if and only if for 
 any orthonormal basis $(u_k)_{k=1}^p$ of $P$, we have $\sum_{k=1}^p \langle \phi(u_k), \, \psi(u_k) \,\rangle_{\R^{p,q}} = 0$.  
 \end{corollary} 
 
 \begin{proof}
Take any orthonormal basis $\{u_k\}_{k=1}^p$ for $P$. Then the claim follows from
 \[\sum_{k=1}^p \langle \psi (u_k), \phi u_k \, \rangle_{\R^{p,q}}= \sum_{k=1}^p\langle u_k, \psi^{*Q} \circ \phi(u_k) \rangle_{\R^{p,q}} =  \tr( \psi^{*Q} \circ \phi) .\]
 \end{proof} 
 
We now turn our attention to flag manifolds of $\Gtwosplit$. As in Proposition \ref{Prop:G2Flags}, we can identify
the flag manifold $\mathcal{F}_{\beta} = \Gtwosplit/P_{\beta}$ with the space $\Ein^{2,3}$ of isotropic lines in $\imoct$ and $\mathcal{F}_{\alpha} = \Gtwosplit/P_{\alpha}$ with the space of annihilator photons $\Pho^\times$ in $\imoct$. We describe concretely the natural $\Gtwosplit$-equivariant inclusions:
\[ \mathcal{F}_\beta,\mathcal{F}_\alpha\hookrightarrow \vis\X_{\Gtwosplit}. \]
Let us use an $\R$-split maximal torus $T$ of diagonal transformations in an $\R$-cross product basis $B = (x_k)_{k=3}^{-3}$. 
Now, form one-parameter subgroups $T_{\beta}, T_{\alpha} < T$ associated to the (Cartan projection of the) co-roots $\tau_{\beta},\tau_\alpha \in \mathfrak{a}$ of the roots $\beta, \alpha$. In particular,
\begin{align}\label{ModelEinGeodesic}
    T_{\beta} &= \{ \diag(e^{2s}, e^{s}, e^{s}, 1, e^{-s}, e^{-s}, e^{-2s}) \}_{s \in \R} \\ \label{ModelPhoGeodesic} 
   T_{\alpha} &= \{ \diag(e^r, e^r, 1, 1, 1, e^{-r}, e^{-r} ) \}_{r \in \R}.
\end{align}
We can choose $B$ in such a way that 
$U_0 = \spann \{ x_3 + x_{-3}, \,x_2 + x_{-2}, \,x_{1} + x_{-1} \} \in \Gr_{(3,0)}^\times(\imoct)$. For example, set $B= X^\R$ from \eqref{ModelRCrossProductBasis}, in which case $p_0 = \spann \{ \mathbf{i},\mathbf{j},\mathbf{k}\}.$ 

\begin{figure}[ht]
    \centering
\includegraphics[width=.7\textwidth]{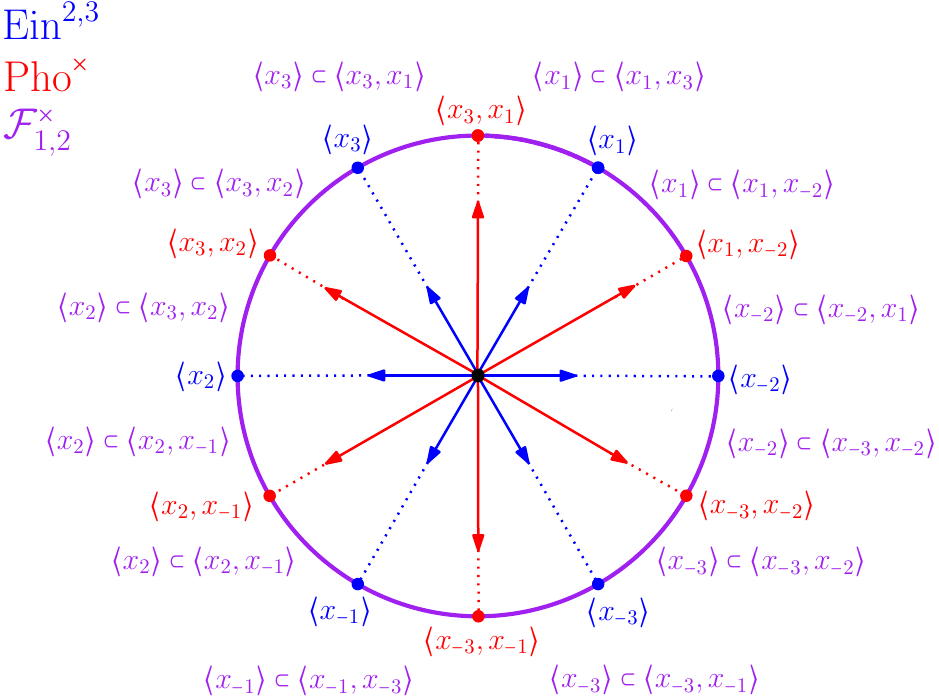}
\caption{ A flat $A\subset \X_{\Gtwosplit}$ and corresponding apartment $\vis A \subset \vis \X$ associated to the $\R$-cross product basis $(x_k)_{k=3}^{-3}$. Each vertex is a partial flag in $\Ein^{2,3}$ or $\Pho^\times$. The open segments of $\titsangle$-length $\frac{\pi}{6}$ correspond to full flags in $\mathcal{F}_{1,2}^\times$. Tangent vectors at $o \in \X$ are drawn non-unitally to visualize the $G_2$ root system. }
\label{Fig:Apartment}
\end{figure}

The geodesic rays $\gamma_{\beta} := T_{\beta} \cdot U_0$ and $\gamma_{\alpha}:= T_{\alpha} \cdot U_0 $, respectively, converge to points in the visual boundary whose stabilizers are respectively the stabilizer of the isotropic line $\langle x_3 \rangle\in \Ein^{2,3}$, which is conjugated to $P_{\beta}$, and the stabilizer of the annihilator photon $\spann \langle x_3,x_2 \rangle\in \Pho^\times$, which is conjugated to $P_{\alpha}$. Therefore, $\Gtwosplit$-equivariant embeddings $\Ein^{2,3} \hookrightarrow \vis\X$ and $\Pho^\times \hookrightarrow \vis\X $ are induced by identifying $\langle x_3\rangle$ with the endpoint of $\gamma_{\beta}$ and $\spann \langle x_3,x_2 \rangle$ with the endpoint of $\gamma_{\alpha}$, respectively. See Figure \ref{Fig:Apartment}. \medskip 

Let us fix some point $U \in \X_{\Gtwosplit}$. For every $\ell\in \Ein^{2,3}$, (resp. $\omega\in \Pho^\times$), there exist a unique
tangent vector $v \in \mathrm{T}^1_U\X$ pointing towards $\ell$ (resp. $\omega$). We give a more concrete characterization of these vectors in Proposition \ref{Prop:PointingTowardsEinsteinPrelim} and \ref{Prop:PointingTowardsPhotonPrelim}, respectively.

To understand pointing towards $\Ein^{2,3}$ in $\X_{\Gtwosplit}$, it is useful to take a detour through $\X_{\SO(3,4)}$. Recall the identification $\mathrm{T}_U\Gr_{(3,0)}(\R^{3,4}) \cong \Hom(U,U^\bot)$, under which $\mathrm{T}_U\X_{\Gtwosplit} \hookrightarrow \mathrm{T}_U\X_{\SO(3,4)}$ includes as the subspace $\mathrm{T}_U\Gr_{(3,0)}^\times(\imoct) = \Hom^\times(U,U^\bot)$ of derivations of $\times_{\imoct}$. We note the relation between these tangent spaces. 
\begin{proposition}[Complementary Subspace to Derivations]\label{Prop:AntiDerivation}
Let $U \in \X_{\Gtwosplit}$. The splitting \[\mathrm{T}_U\X_{\SO(3,4)} = \mathrm{T}_{U}\X_{\Gtwosplit} \oplus \mathcal{C}_{U^\bot}\] 
is orthogonal, where 
\[ \mathcal{C}_{U^\bot} = \{ \mathcal{C}_{z}|_{U} \in \Hom(U,U^\bot) \mid z \in U^\bot \}.\] 
\end{proposition}

Recall that for $z\in \imoct$, $\mathcal{C}_z:\imoct\to\imoct$ is its left cross product endomorphism.

\begin{proof}
First, we show $\mathcal{C}_{U^\bot} \subset \T_U\X_{\SO(3,4)}$. To see this holds, since  $\mathcal{C}_{U^\bot} \subset \Hom(U, U^\bot)$ by definition, it remains only to show $\mathcal{C}_{U^\bot} \subset \solie(3,4)$, which we do now.

By linearity, it suffices to show $\mathcal{C}_z \in \solie(3,4)$ for any $z \in Q_-(U^\bot)$. Thus, let any such $z$ be fixed. Here, it is convenient to work with the Euclidean metric $h =q|_U  \oplus (-q|_{U^\bot})$. 
Now, set $\psi: = \mathcal{C}_z|_{z^\bot}$. 
Note that $\psi \in \O(z^\bot, q)$ and also $\psi \in \O(z^\bot, h)$. Hence, $\psi^{-1} = \psi^{*h}$. However, $\psi =\psi^{-1}$ by \eqref{DCP}. This means $\psi = \psi^{*h}$. Now, 
define $\psi^{*q}$ as the adjoint satisfying the following for all $x, y \in z^\bot$:
\[ \langle \psi(x), y\rangle_q =\langle x, \psi^{*q}(y) \rangle_q \] 
Note that $\psi^{*q} = -\psi^{*h}$ since $\psi$ exchanges $U$ and $U^\bot$. Thus, $\psi=-  \psi^{*q}$, which means $\mathcal{C}_z|_{z^\bot} \in \solie(z^\bot, q)$ and hence $\mathcal{C}_z \in \solie(3,4)$ since $\mathcal{C}_z(z)=0$. We shall use the condition $(\mathcal{C}_z)^{*q} = -\mathcal{C}_z$ going forwards. 

By dimension count, the proposition follows if  $\mathcal{C}_{U^\bot}$ is orthogonal to  $\X_{\Gtwosplit}$. Take $\phi \in \T_U\X_{\Gtwosplit}$ and write $\phi(x)=x', \phi(y)=y', \phi(xy) =x'y+xy'$ for some orthonormal basis $(x,y,xy)$ of $U$. Let $\psi\in \mathcal{C}_{U^\perp}$ be an element of the form ${\mathcal{C}_z}|_{U}$ for some $z\in U^\perp$. By Proposition \ref{Prop:ModelSpaceMetric}, for some $c >0$, we find:
\begin{align*}
    c\langle \psi, \phi\rangle_{\X} &=\tr(\psi^{*h}\circ \phi)= \tr(\mathcal{C}_z|_{U^\bot}\circ \phi) =\langle x,\mathcal{C}_z\circ \psi(x)\rangle_q +\langle y, \mathcal{C}_z\circ \psi(y)\rangle_q+\langle xy, \mathcal{C}_z\circ \psi(xy)\rangle_q \\
    &= \langle x,zx'\rangle_q + \langle y,zy'\rangle_q +\langle xy, z(x'y)\rangle_q +\langle xy, z(xy')\rangle_q =0.
\end{align*}
The final equality holds due to cancellation of terms in pairs. For example, using $\mathcal{C}_y \in \mathrm{O}(y^\bot, q)$,
\[ \langle x, zx'\rangle_q = \langle xy, (zx')y\rangle_q =_{(\star)}\langle xy, -z(x'y)\rangle)_q.\]
Let us explain the equality $(\star)$. Observe that only the orthogonal projection of $x'$ onto $\R \{xz\}$ can contribute to the pairing $\langle xy, (zx')y\rangle_q$. Thus, without loss of generality, we can assume that $x' \in \R\{xz\}$, in which case $y$ is orthogonal to the subalgebra generated by $x',z$. In this case, the triple $(x,y',z)$ anti-associates, 
which finally yields the desired equality. 
\end{proof}

The result of Proposition 
\ref{Prop:AntiDerivation} can be considered a  consequence of the orthogonal decomposition $\solie(3,4)=\g_2'\oplus \mathcal{C}_{\imoct}$ and its interaction with the Cartan decomposition of $\solie(3,4)$. 

The following Proposition \ref{Prop:PointingTowardsEinsteinPrelim} is a special case of the more general result in  \cite[Section 2.2]{DE26} about pointing towards $\Ein^{p-1,p}$ in $\X_{\SO_0(p,p+1)}$. 

\begin{proposition}[Pointing Toward $\Ein^{2,3}$]
    \label{Prop:PointingTowardsEinsteinPrelim}
Let $U\in \X_{\Gtwosplit}$ and let $\ell\in \Ein^{2,3}$. Write $\ell = [u+z]$ for 
$u \in U, \, z \in U^\bot$. Let $\widetilde{\phi}\in \T_{U}\X_{\SO(3,4)}$ be the unique rank one tangent vector $\tilde{\phi}$ such that $\tilde{\phi}(u) =z$. Then $\tilde{\phi}$ points towards $\Ein^{2,3}$ and moreover $\tilde{\phi}$ points towards $\graph^*(\tilde{\phi}):=\graph(\tilde{\phi}_{\mid L})=\ell$ where $L$ is the orthogonal to the kernel of $\tilde{\phi}$.  
\end{proposition}

It will be essential in this work that we consider tangent vectors in $\X_{\SO(3,4)}$ pointing towards $\Ein^{2,3}$ rather than those in $\X_{\Gtwosplit}$. Indeed, the former are rank one maps, while the latter are rank three. Note that $\Ein^{2,3}$ is a flag manifold of both $\SO(3,4)$ and $\Gtwosplit$, but is embedded differently in the visual compactification of the corresponding symmetric spaces.

\begin{proposition}[Projections of Tangent Vectors]
\label{Prop:PointingEinSame}
Suppose $U \in \X_{\Gtwosplit}$ and $\tilde{\phi} \in \mathrm{T}_U\X_{\SO(3,4)}$ points toward $\ell \in \Ein^{2,3}$ in $\X_{\SO(3,4)}$. Let $\pi: \mathrm{T}_U\X_{\SO(3,4)} \rightarrow \mathrm{T}_U \X_{\Gtwosplit}$ be the orthogonal projection. Then $\phi := \pi (\tilde{\phi})$ points towards $\ell \in \Ein^{2,3}$ in $\X_{\Gtwosplit}$. 
\end{proposition}

\begin{proof}
We use Proposition \ref{Prop:AntiDerivation} to compute $\pi(\tilde{\phi})$. Write $\ell = [u+z]$ for unit elements $u \in Q_+(U), \, z\in Q_-(U^\bot)$. Define $Z:=-zu$, so $z=Zu$. Form an orthonormal basis $(z_k)_{k=1}^4$ of $U^\bot$ with $z_1 = Z$. Then $(\mathcal{C}_{z_k})_{k=1}^4$ is an orthonormal basis of $\mathcal{C}_{U^\bot}$. It is obvious that $\langle \tilde\phi, \mathcal{C}_{z_k}\rangle=0$ for $k \in \{2,3,4\}$. Hence, \[\pi(\tilde{\phi})=\tilde{\phi}+a\mathcal{C}_Z\] for the unique constant $a \in \R$ such that the right-hand side is orthogonal to $\mathcal{C}_Z$. One computes $a=-\frac{1}{3}$ directly, and consequently, up to positive scalars, 
$\pi(\tilde{\phi})$ obtains the form 
\begin{align}\label{EinModelExplicit}
   \begin{cases}u &\longmapsto 2Zu \\
    v &\longmapsto -Zv \\
    w &\longmapsto -Zw ,\end{cases}
\end{align}
where $(u,v,w)=(u,v,u\times v)$ is an orthonormal basis for $P$. 

On the other hand, the model geodesic $\gamma_{\beta} = T_{\beta}\cdot U_0$ from \eqref{ModelEinGeodesic} obtains the same form as \eqref{EinModelExplicit}, so that $\pi({\tilde{\phi})}$ points towards $\Ein^{2,3}$. Examining the model geodesic, one finds moreover that \eqref{EinModelExplicit} points towards $\ell = [u+Zu] = [u+z]$.  
\end{proof}

\begin{corollary}
Let $\phi \in \mathrm{T}_U\X_{\Gtwosplit}$. Then $\phi$ points towards $\ell =[u+z]\in \Ein^{2,3}$ if and only if $\phi$ obtains the form \eqref{EinModelExplicit} up to positive scalars. 
\end{corollary}

To underscore the point, we leave the following remark on the discrepancy between convergence in $\vis \mathbb{X}_{\Gtwosplit}$ and $\vis\mathbb{X}_{\SO_0(3,4)}$.

\begin{remark}
Let $\iota: \X_{\Gtwosplit} \hookrightarrow \X_{\SO(3,4)}$ be the inclusion map. If $\phi \in \mathrm{T}_U\X_{\Gtwosplit}$ points towards $\Ein^{2,3}$, then $\iota(\phi)$ points towards $\Iso_{\{1,3\}}(\R^{3,4})$, the flag manifold of pointed maximal isotropic subspaces.  
\end{remark}

We now describe geometrically when tangent vectors point towards annihilator photons. 

\begin{proposition}[Pointing Toward $\Pho^\times$]
    \label{Prop:PointingTowardsPhotonPrelim}
Let $P\in \X_{\Gtwosplit}$. Then $\phi \in \mathrm{T}_U\X_{\Gtwosplit}$ points towards $\Pho^\times$ if and only if $\phi$ obtains the form
\begin{align}\label{ModelAlphaGeodesic}
	 \phi = \begin{cases}
	u &\longmapsto z \\
	v &\longmapsto (uv)z \\
	  uv &\longmapsto 0. \end{cases} 
\end{align}
for some orthonormal triple $(u,v,z)$ such that $u,v \in U$, $z \in U^\bot$.
\end{proposition}

A proof of Proposition \ref{Prop:PointingTowardsPhotonPrelim} is contained in the proof of \cite[Proposition 3.9]{ER25}. We note here a brief alternate argument. First, note that by Proposition \ref{Prop:FirstStiefel}, any two such tangent vectors $\phi$ are equivalent up to the $\Gtwosplit$-action, as they are described precisely by $(u,v,z) \in V_{(+,+,-)}(\imoct)$. Consequently, the claim follows if it holds for the model geodesic \eqref{ModelPhoGeodesic}. One verifies immediately this is the case, with $u,v,z$ given by 
$u = \frac{1}{\sqrt{2}}(x_3+x_{-3})$, $v = \frac{1}{\sqrt{2}}(x_2+x_{-2})$, $z= \frac{1}{\sqrt{2}}(x_3-x_{-3})$. 

\begin{corollary}
The tangent vector $\phi$ in \eqref{ModelAlphaGeodesic} points towards $\omega \in \Pho^\times$ given by 
\[\omega = \graph^*(\phi) = \spann \{ u+z,v+(uv)z\}.\]
\end{corollary}
Here we write again $\graph^*(\phi)$ as the graph of $\phi_{\mid W}$ where $W= U \cap \ker(\phi)^\bot$. Note that $\omega \in \Pho^\times$ by Corollary \ref{Cor:AnnihilatorPhotonGraph}. 

\begin{remark}
A fact we shall employ later on is that a tangent vector $\phi \in \mathrm{T}_U\X_{\Gtwosplit}$ that is rank two must obtain the form \eqref{ModelAlphaGeodesic}. This follows immediately from the derivation condition.
\end{remark}

\subsection{Anosov Representations \& Bases of Pencils}

In this section, we recall the following: the definition of a \emph{base of pencil}, a heuristic definition of Anosov representations, the notion of a domain of discontinuity $\Omega \subset \mathcal{F}$ in a flag manifold $\mathcal{F}$ defined by \emph{Tits metric thickening} as in \cite{KLP18}, and finally how such domains can, under favorable circumstances, be fibered by bases of pencils as in \cite{Dav25}.

\subsubsection{Bases of Pencils}

Let $G$ be a non-compact simple real Lie group. Fix a non-zero unit vector $\tau \in \overline{\mathfrak{a}}^+$ in the model closed Weyl chamber $\overline{\mathfrak{a}}^+ \subset \T_o\X$ at the basepoint $o \in \X$. In practice, we will focus on $\tau$ whose associated $G$-orbit $\mathcal{F}_\tau$ in the visual boundary of the symmetric space is one of the partial flag manifolds $\Ein^{2,3}$ or $\Pho^\times$ for $G = \Gtwosplit$.\medskip

For any point $x\in \X$, we call a plane $\mathcal{P}\subset \mathrm{T}_x\X$ a \emph{pencil of tangent vectors} or a \emph{pencil} for short. Such a plane defines naturally a subset of the flag manifold $\mathcal{F}_\tau$ of expected codimension two, that we call the \emph{$\tau$-base}. We will be especially interested in the case $\mathcal{F}_{\tau} =G/P_{\Theta}$ when $P_{\Theta}$ is maximal parabolic, so the ($G$-equivariant) embedding $\mathcal{F}_{\tau} \hookrightarrow \X_{G}$ is uniquely defined, as in $\S$\ref{Sec:SymmSpace}. 

\begin{definition}[Base of Pencils]\label{Defn:TauBasePencil}
Let $\mathcal{P} \subset \T_{x} \mathbb{X}$ be a pencil. Then the $\boldsymbol{\tau}$\textbf{-base of} $\boldsymbol{\mathcal{P}}$, denoted $\mathcal{B}_{\tau}(\mathcal{P})$, is given by 
\[ \mathcal{B}_{\tau}(\mathcal{P}) = \{ \gamma_{x,v}(\infty) \in \mathcal{F}_{\tau} \; | \; v \in \T_{x}\mathbb{X}, \, v \,\bot \,\mathcal{P}\}\]
\end{definition}
\label{defn:Pencil}
In other words $\mathcal{B}_{\tau}(\mathcal{P})$ contains the $\tau$-flags that can be reached in $\vis\mathbb{X}$ by traveling from $x$ via directions orthogonal to $\mathcal{P}$ in $\mathrm{T}_x\X$. 
\begin{remark}
There are unique $\Gtwosplit$-equivariant embeddings $\Ein^{2,3},\, \Pho^\times \hookrightarrow \vis\X_{\Gtwosplit}$. We will consider only $\mathcal{B}_{\mu(\tau_\alpha)}(\mathcal{P})\subset \Pho^\times$ and $\mathcal{B}_{\mu(\tau_\beta)}(\mathcal{P})\subset \Ein^{2,3}$ where $\tau_\alpha$, $\tau_\beta$ are the respective coroots. Consequently, we may simply refer to these bases respectively as the $\alpha$-base $\mathcal{B}_\alpha(\mathcal{P})$ and the $\beta$-base $\mathcal{B}_\beta(\mathcal{P})$.
\end{remark}

The following Proposition serves as both an example and a clarification. Recall that $\Ein^{2,3}$ is both an $\SO(3,4)$ and a $\Gtwosplit$-flag manifold, associated to a maximal parabolic in each group. 
The corresponding bases of pencils are identified, regardless of which ambient geometry is used. 

\begin{proposition}[Equivalent Bases in $\Ein^{2,3}$]\label{Prop:EquivalentBases}
Suppose that $\mathcal{P} \subset \mathrm{T}_U\X_{\Gtwosplit}$ is a pencil. Denote $\pi:\mathrm{T}_U\X_{\SO(3,4)}\rightarrow \mathrm{T}_U\X_{\Gtwosplit}$ as the orthogonal projection. Let $\tilde{\tau}$ be the Cartan projection $\mu(v)$ of any unit vector $v \in \mathrm{T}_U\X_{\SO(3,4)}$ pointing towards $\Ein^{2,3}$. Set $\tau = \pi(\tilde{\tau})$. Then the bases of pencils ${B}_{\tilde{\tau}}(\mathcal{P})$ and $B_{\tau}(\mathcal{P})$ are naturally identified. 
\end{proposition}

\begin{proof}
The proof follows immediately from Proposition \ref{Prop:PointingEinSame} and the fact that if $X \in \mathrm{T}_U\X_{\Gtwosplit}$ and $Y \in \mathrm{T}_U\X_{\SO(3,4)}$, then $\langle X,Y\rangle_{\X_{\SO(3,4)}}=0$ if and only if $\langle X,\pi(Y)\rangle_{\X_{\Gtwosplit}}=0$. 
\end{proof}

We will take advantage of Proposition \ref{Prop:EquivalentBases} in Section \ref{Sec:G2Ein23GeomStructures} when building geometric structures. 

\subsubsection{Domains of Discontinuity via Tits Metric Thickening}

Recall that for a surface group $\pi_1S$, the \emph{Gromov boundary} $\partial_{\infty}\pi_1S$ is a topological circle that has a canonical $\pi_1S$-action. For any semisimple real Lie group, a representation $\rho: \pi_1S \rightarrow G$ is called $P_{\Theta}$-Anosov, for a self-opposite parabolic subgroup $P_{\Theta}$, when there is a continuous, $\rho$-equivariant \emph{transverse} map $\xi: \partial_{\infty}\pi_1S \rightarrow \mathcal{F}_{\Theta}$ with some additional contraction properties \cite{GW12}. Here, we write $\mathcal{F}_{\Theta} = G/P_{\Theta}$ and call $\xi$ \emph{transverse} when $(\xi(x),\xi(x'))$ lie in the unique open $G$-orbit in $G/P_{\Theta} \times G/P_{\Theta}$ for every pair $x \neq x' \in \partial_{\infty}\pi_1S$. The map $\xi$ is unique if it is exists. 

Given a $P$-Anosov representation $\rho:\pi_1S\rightarrow G$, the idea developed in \cite{GW12} and later extended in \cite{KLP18}, is to build a domain of discontinuity $\Omega_{\rho} \subset \mathcal{F}'=G/Q$ by removing a ``thickening'' of the limit set $\Lambda = \xi(\pi_1S) \subset G/P$ in the following sense. Often in the setup, the flag manifolds $\mathcal{F}=G/P$ and $\mathcal{F}'=G/Q$ are different.  
First, one defines for a single flag $f\in \mathcal{F}$ the notion of a \emph{thickening} $K_f \subset \mathcal{F}'$ that is a finite union of Schubert subvarieties of $\mathcal{F}'$. One then defines a domain of discontinuity for the image of the representation by removing the thickening of the limit set: $\Omega = \mathcal{F}'\backslash \bigcup_{f \in \Lambda} K_{f}$. 
For example, in \cite{GW12}, for $G =\SO(p,q)$, with $1<p < q$, this strategy is applied for $P_1$-Anosov representations with $\mathcal{F}=\Ein^{p-1,q-1}=G/P_1$, the set of isotropic lines in $\R^{p,q}$, and $\mathcal{F}'=\Iso_{p}(\R^{p,q}) = G/P_p$, the set of maximal isotropic planes in $\R^{p,q}$. They also do the same with the roles of $P_1$ and $P_p$ reversed.  

Unlike in most cases, we will set $\mathcal{F} =G/P_{\Theta}=\mathcal{F}'$ to be the same flag manifold, where $\Theta \in \{\alpha, \beta\}$ and $G = \Gtwosplit$. It is a rare feature exhibited by $G=\Gtwosplit$ and its partial flag manifolds that this construction works. The case of interest for us presently is when the thickening is defined via \emph{Tits (angle) metric thickening}. Namely, we shall only consider thickenings of $f \in \mathcal{F}_{\Theta}$ the following form: 
\[ K_{f} = \left\{ f' \in \mathcal{F}_{\Theta} \mid \titsangle(f,f')\leq \frac{\pi}{2} \right \}.\]

Note that to make sense of the Tits angle, one has to embed $\mathcal{F}_{\Theta}$ in the visual boundary $\vis \X$. In the present cases of $\Theta = \{\alpha\}$ or $\Theta =\{\beta\}$, we have already described these (unique) embeddings explicitly in $\S$\ref{Subsec:G2SymmetricSpace}, so there is no ambiguity.

Now, let $\rho:\pi_1S \rightarrow \Gtwosplit$ be a $P_{\Theta}$-Anosov representation for $\Theta \in \{\alpha, \beta\}$ with associated limit map $\xi:\partial_{\infty}\pi_1S \rightarrow \mathcal{F}_{\Theta}$. One can then define an open domain by the same strategy as in the general case: 
\begin{align}\label{Omega_Thick}
    \Omega^{\Thick}_\Theta = \mathcal{F}_{\Theta} \backslash \bigcup_{x \in \partial_{\infty}\pi_1S} K_{\xi(x)}. 
\end{align}
Here, $\Omega = \Omega^{\Thick}_{\Theta}$ is just the complement of the $\frac{\pi}{2}$-neighborhood of $\Lambda = \image(\xi)$ with respect to the Tits angle metric. 

The following result is a consequence of \cite[Theorem 8.6]{GW12} for $\Theta = \{\beta\}$ and \cite[Theorem 1.8]{KLP18} for $\Theta \in \{\alpha,\beta\}$. 

\begin{theorem}[Domains of Discontinuity for Anosov Representations in $\Gtwosplit$]\label{Thm:G2DOD}
Let $\rho:\pi_1S \rightarrow \Gtwosplit$ be $P_{\Theta}$-Anosov. The action of $\rho(\pi_1S)$ on the domain \eqref{Omega_Thick} is properly discontinuous and cocompact. 
\end{theorem}

Although this paper primarily concerns the construction of geometric structures with differential geometry, we shall compare our construction to the quotients $\rho(\pi_1S)\backslash \Omega_{\rho}^{\Thick}$ in Theorems \ref{Thm:BetaTits=Pencil} and \ref{Thm:AlphaHodgeTits=Pencil} when $\rho:\pi_1S\rightarrow \Gtwosplit$ is Hitchin.

\subsubsection{Fibration of Domains of Discontinuity via Bases of Pencils}

In this section, we discuss how to construct fibrations of some domains of discontinuity constructed by metric thickening, for representations preserving a totally geodesic copy of $\mathbb{H}^2$ in the symmetric space, following \cite{Dav25}. \medskip

The first relevant notion is that of regularity, stated presently for $\X=\X_{\Gtwosplit}$. 
\begin{definition}
We call $X \in \T_x\X_{\Gtwosplit}$ to be $\boldsymbol{\tau_\alpha}$\textbf{-regular} or just $\boldsymbol{\alpha}$\textbf{-regular}, respectively  $\boldsymbol{\tau_\beta}$\textbf{-regular} or just $\boldsymbol{\beta}$\textbf{-regular} when $\alpha(\mu(X))\neq 0$, respectively $\beta(\mu(X))\neq 0$. 

Call an immersion $u: M \rightarrow \X_{\Gtwosplit}$ to be $\boldsymbol{\alpha}$\textbf{-regular}, respectively $\boldsymbol{\beta}$\textbf{-regular} when $du(X)$ is respectively $\alpha, \beta$-regular for all non-zero $X \in \T M$. 
\end{definition}

Next, we recall how the notion of $\tau$-bases for $\tau\in \overline{\mathfrak{a}}^+$ relates to fibrations of cocompact domains of discontinuity for Anosov representations. Note that we will only be interested in the case $\tau=\tau_\alpha$ or $\tau=\tau_ \beta$. Let $f: \tilde{S} \rightarrow \mathbb{X}$ be a totally geodesic embedding that is $\tau$-regular. 
Fixing an arbitrary basepoint $o \in \X$, we can define a domain $\Omega_{f}^{\tau} $ in the flag manifold $\mathcal{F}_{\tau}$ using Busemann functions by 
\begin{align}\label{Busemann_Domain}
 \Omega_{f}^{\tau} := \{a \in \mathcal{F}_{\tau} \mid b_{a,o} \circ f\; \text{is proper, bounded below} \}.
\end{align}
Here, the \emph{Busemann} function $b_{a,o}$ measures the relative distance of points $x \in \X$ to $a \in \vis\X$ from the point of view of $o$ by 
\[ b_{a,o}(x) := \lim_{t \rightarrow \infty} d_{\X}(\gamma_{o,a}(t),x)-t.\]
By the triangle inequality, the definition of $b_{a,o}$ is well-defined. The map $b_{a,o}$ is smooth in the case of $\X$ a symmetric space. 

There is now a natural projection $\Omega^{\tau}_{f} \rightarrow \tilde{S} $ as follows. 

\begin{lemma}[Nearest Point Projection]\label{Lem:NearestPointProjection}
Let $f: \tilde{S} \rightarrow \mathbb{X}$ be totally geodesic and $\tau$-regular. There is a natural projection $\pi: \Omega^{\tau}_f \rightarrow \tilde{S}$, where $\pi(a) =x$ for the unique point $x \in \tilde{S}$ such that $b_{a,o} \circ f$ has a critical point at $f(x)$. Moreover:
\begin{enumerate}[noitemsep]
    \item $\Omega^{\tau}_{f}$ is open, 
    \item $\pi$ is a fiber bundle projection,
    \item $(\Omega^{\tau}_f)|_x = \mathcal{B}_{\tau}(\mathcal{P}_x)$, where $\mathcal{P}_x \subset \T_{f(x)}\mathbb{X}$ is the pencil $ df(\T_{x}\tilde{S})$. 
\end{enumerate}
\end{lemma} 

\begin{proof}
The definition of $\pi$ is well-defined and $\Omega^{\tau}_f$ is open by \cite[Lemma 7.2]{Dav25}. Then \cite[Theorem 7.3]{Dav25} settles points (2) and (3).   
\end{proof}

To make the final link, the domain \eqref{Busemann_Domain} defined via Busemann functions is the same as the domain \eqref{Omega_Thick} defined via Tits metric thickening in the cases of interest. 

\begin{proposition}[{\cite[Theorem 7.11]{Dav25}}]
\label{prop:NearestPointProjection}
Let $\rho: \pi_1S\rightarrow (\mathrm{P})\SL(2,\R) \hookrightarrow \Gtwosplit$ be $P_{\alpha}$-Anosov, respectively $P_\beta$-Anosov, and assume that the corresponding $\rho$-equivariant totally geodesic surface $f: \tilde{S} \hookrightarrow \X$ is $\alpha$-regular, respectively $\beta$-regular. Then the corresponding domains in $\mathcal{F}_\alpha$, respectively $\mathcal{F}_\beta$ from \eqref{Omega_Thick} and \eqref{Busemann_Domain} coincide. 
\end{proposition}

\section{Cyclic \texorpdfstring{$\Gtwosplit$}{G2'}-Higgs Bundles}\label{Sec:CyclicG2'Higgs}
In this section, we recall some basic facts on the non-abelian Hodge (NAH) correspondence, and describe the general cyclic $\Gtwosplit$-Higgs bundles of \cite{CT24} in terms of $\C$-cross product bases for $\imoct^\C$. We then derive Hitchin's equations for the Higgs bundles of interest.

\subsection{Non-Abelian Hodge Correspondence}

We now provide a brief review of the non-abelian Hodge correspondence, as developed by Hitchin, Simpson, Corlette, and Donaldson \cite{Hit87, Sim88, Cor88,Don87}. We introduce the correspondence centered around the case $G = \SL(n,\C)$, but also discuss $G$-Higgs bundles in some greater generality in preparation for the discussion on $\Gtwosplit$-Higgs bundles in the next subsection. We refer the reader to the following surveys on Higgs bundles and the non-abelian Hodge correspondence for more details: \cite{Gothen14, Gui18, Li19}.

Let $S$ be a closed surface and $\Sigma =(S,J)$ a Riemann surface on $S$. The non-abelian Hodge correspondence $\NAH_{\Sigma}$ provides a dictionary to translate between surface group representations $\rho: \pi_1S \rightarrow G$ and holomorphic objects called $(G-$)\emph{Higgs bundles} on $\Sigma$. For a fixed real semisimple Lie group $G$, the map $\NAH_{\Sigma}$ is a homeomorphism between the moduli space $\mathcal{M}_{G}(\Sigma)$ of \emph{polystable} $G$-Higgs bundles on $\Sigma$, up to gauge equivalence, and the moduli space $\chi(S,G) = \Hom^{\mathrm{red}}(\pi_1S, G)/G$ \emph{reductive} representations $\pi_1S \rightarrow G$, up to conjugation. We now recall the details. 

We shall opt to describe $G$-Higgs bundles, for many  examples including $\Gtwosplit$, in terms of vector bundles. 
In the simplest setting, for $G=\GL(n,\C)$, a Higgs bundle on $\Sigma$ is a pair $(\V,\Phi)$ consisting of a rank $n$ holomorphic vector bundle $\V$ on $\Sigma$ and a holomorphic $\End(\V)$-valued (1,0)-form on $\Sigma$, denoted $\Phi \in H^0(\End(\V)\otimes \K)$, called a \emph{Higgs field}. Here, $\K=\K_{\Sigma}$ is the holomorphic cotangent line bundle. Now, we present both the general definition and some relevant examples expressed in relation to this vector bundle formulation. 

Let $G$ be a non-compact semisimple real Lie group and $K < G$ be a maximal compact subgroup. The Lie subalgebra $\frakk$ of $K$ has an orthogonal complement $\mathfrak{p}$ under the Killing form, which provides a $\Z_2$-Lie algebra grading $\g=\frakk\oplus \frakp$ called a Cartan decomposition. We may complexify to obtain $\g^{\C}=\frakk^\C \oplus \mathfrak{p}^\C$. A $G$-Higgs bundle $(\mathcal{P}, \Phi)$ on a Riemann surface $\Sigma$ then consists of a holomorphic principle $K^\C$-bundle $\mathcal{P}\rightarrow \Sigma$ and a Higgs field $\Phi \in H^0(\mathcal{P}\times_{Ad} \mathfrak{p}^{\C} \otimes \K)$. An isomorphism $g:(\mathcal{P},\Phi)\rightarrow (\mathcal{P}',\Phi')$ of Higgs bundles is a holomorphic principal bundle isomorphism that pulls back $\Phi'$ to $\Phi$. The moduli space $\mathcal{M}_G(\Sigma)$ is the space of polystable $G$-Higgs bundle up to isomorphism. We refer to \cite{GPGMiR12} for the general notion of stability of $G$-Higgs bundles. 

In the present context, we shall only be interested in $G$-Higgs bundles for simple Lie groups $G \in \{ \SL(n,\C), \SL(n,\R), \SO_0(3,4), \Gtwosplit\}$. We will use the common strategy to factor the principal bundle $\mathcal{P}$ as a sub-bundle of the frame bundle $\mathrm{Fr}(\V)$ of some holomorphic vector bundle $\V$ equipped with additional structure that $\mathcal{P}$ preserves. In this way, $G$-Higgs bundles can be thought of as `decorated' $\GL(n,\C)$-Higgs bundles. To make this idea precise, consider the following examples in ascending order of structure. See \cite{Col19} for further details on these examples and more. 

\begin{example}
    For $G=\GL(n,\C)$, first observe that the principal and vector bundle formulations of Higgs bundles are equivalent. Indeed, 
    regarding $G$ as a real Lie group, then $\frakk^\C = (\mathfrak{u}_n)^\C = \gl_n\C=\frakp^\C$. Hence, if $\V$ is a holomorphic rank $n$ vector bundle, then $\mathcal{P}=\Fr(\V)$ is the $G=K^\C$-principal bundle and $\mathcal{P} \times_{Ad} \p^{\C} \cong \End(\V)$, so that the two notions of Higgs fields coincide. 
\end{example}

Incorporating a volume form, we can refine the previous example. 
    \begin{example} For $G=\SL(n,\C)$, a $G$-Higgs bundle is equivalently described by a tuple $(\V,\Phi, \omega)$ such that $\det(\V) \cong \mathcal{O}$ is holomorphically trivial, $\omega :\det(\V) \rightarrow  \mathcal{O}$ is a fixed holomorphic trivialization, and $\tr(\Phi) =0$. Note that $\g=\sllie_{n}\C=(\mathfrak{su}_n)^\C = \mathfrak{p}^\C$. The identification $\End_0(\V)\cong \mathcal{P}\times_{Ad}\p^\C$, where $\End_{0}(\V)$ denotes trace-free endomorphisms, allows us to identify the two Higgs field perspectives. 
    \end{example}

Moving from complex to real groups requires further decorations on the Higgs bundle, as the following example illustrates.
\begin{example} 
    For $G = \SL(n,\R)$, recall the Cartan decomposition $\g = \frakk\oplus \frakp$, where $\frakk =\mathfrak{so}_n\R$ and $\frakp=\mathrm{Sym}_0(\R^n)$, the trace-free symmetric endomorphisms. In particular, $\frakp^\C = \Sym_0(\C^n)$. Hence, an $\SL(n,\R)$-Higgs bundle can be regarded as a tuple $(\V, \Phi, \omega, Q)$, where $Q: \V \times \V \rightarrow \mathbb{C}$ is a holomorphic non-degenerate bilinear form on $\V$ such that $\Phi$ is trace-free and $Q$-symmetric.
\end{example} 

Similarly, considering Higgs bundles of Lie subgroups requires further decorations still.

\begin{example} For $G= \SO_0(3,4)$, 
    we refine the previous example. In this case, $K = \SO(3)\times \SO(4)$ is the subgroup preserving a splitting $\R^{3,4} = \R^{3,0} \oplus \R^{0,4}$ into space + time. Relative to this splitting, the Cartan decomposition is 
    \[ \g = \frakk \oplus \mathfrak{p} = \begin{pmatrix} * &\\ & *\end{pmatrix} \oplus \begin{pmatrix}  &*\\ *& \end{pmatrix} .\] 
    Transformations $X \in \frakp$ obtain the form $X=\eta - \eta^{*q}$ for $\eta: \R^{3,0} \rightarrow \R^{0,4}$, and $q= q_{3,4}$, thus identifying $\mathfrak{p} \cong \Hom(\R^{3,0}, \R^{0,4})$. Then $\g^{\C} = \frakk^\C \oplus \frakp^\C$ admits a similar description, now relative to the non-degenerate complex bilinear form $Q = q_{3,4}^\C$. As a consequence, a $G$-Higgs bundle may be described by a tuple $(\mathcal{U}, \mathcal{V}, Q_{\mathcal{U}}, Q_\mathcal{V}, \omega, \eta)$, where 
    \begin{itemize}
        \item $\mathcal{U}, \mathcal{V}$ are rank 3 and 4 holomorphic vector bundles, respectively,  
        \item $Q_{\mathcal{U}}$ and $ Q_{\mathcal{V}}$ are non-degenerate holomorphic bilinear forms on $\mathcal{U}, \mathcal{V}$, respectively, 
        \item The object $\eta$ is of the form $\eta \in H^0( \Hom(\mathcal{U}, \mathcal{V})\otimes \K)$, and determines a Higgs field $\Phi = \eta-\eta^{*Q}$, 
        \item $\omega = \omega_{\mathcal{U}} \oplus \omega_{\mathcal{V}}$ holomorphically trivializes\footnote{For $G=\SO(3,4)$, one only assumes that we have a trivialization of the tensor product $\det(\mathcal{U})\otimes\det(\mathcal{V})$.}  $\det(\mathcal{U})$ and $\det(\mathcal{V})$.
    \end{itemize} 
    In this case, the Higgs field $\Phi = \eta - \eta^{*Q}$ satisfies $\Phi \in H^0(\mathcal{P}\times_{Ad}\mathfrak{p}^\C \otimes \K)$, where $\mathcal{P}$ is the holomorphic $K^{\C} = \SO(3,\C) \times \SO(4,\C)$-frame bundle of $(\V, Q_{\mathcal{U}}\oplus (-Q_{\mathcal{V}}), \omega)$. 
\end{example}

Recall from $\S$\ref{Subsec:G2Defn} that $\Gtwosplit < \SO_0(3,4)$ is exactly the subgroup of $\SO_0(3,4)$ additionally preserving the cross product $\times_{3,4}: \R^{3,4} \rightarrow \R^{3,4} \rightarrow \R^{3,4}$. Thus, a $\Gtwosplit$-Higgs bundle is a certain refinement of an $\SO_0(3,4)$-Higgs bundle, with a holomorphic cross product $\times_{\V}: \V \times \V \rightarrow \V$ now incorporated and $\Phi \in H^0(\End(\V)\otimes \K)$ a derivation of $\times_{\V}$. We explain the details in the following subsection.  

Having given some relevant examples, we now step back to describe the non-abelian Hodge correspondence. We will focus on the case of $G = \SL(n,\C)$, which highlights all the essential ideas. 

\subsubsection{Higgs Bundles to Representations}  
The idea is to search for a hermitian metric $h$ on $\V$ that is \emph{harmonic}, as clarified by the following discussion. 
Fix an $\SL(n,\C)$-Higgs bundle $(\V, \Phi, \omega)$. Write $\V= (E,\delbar_{E})$ for $E$ a smooth vector bundle and $\delbar_{E}: \Omega^{0}(E) \rightarrow \Omega^{0,1}(E)$ the Dolbeault operator of $\V$. For a given hermitian metric $h$ on $\V$, form the associated connection
\[\nabla:=\nabla(h)= \nabla_{\delbar_{E}, h}+ \Phi+\Phi^{*h}, \]where $\nabla_{\delbar_{E}, h}$ is the Chern connection of the Hermitian, holomorphic vector bundle $(\V,h)$ and $\Phi^{*h} \in \Omega^{0,1}(\Sigma, \End(\V))$ is the $h$-adjoint of $\Phi$. We shall call $h$ \emph{harmonic} when $\nabla(h)$ is flat. In the case $\nabla(h)$ is flat, let $\rho: \pi_1S\rightarrow \SL(n,\C)$ denote the holonomy, which is contained in $\SL(n,\C)$ since $\omega$ is $\nabla$-parallel. The flat bundle $(E,\nabla)$ is isomorphic to $\tilde{\Sigma}\times_{\rho} \C^n$ by the Riemann-Hilbert correspondence. Now, the hermitian metric $h$ defines a reduction of structure of the frame bundle from $\SL(n,\C)$ to $\SU(n)$. Correspondingly, $h$ is a section $h \in \Gamma(\Sigma, \Fr(\V)\times_{\rho} \SL(n,\C)/\SU(n))$ and thus also a $\rho$-equivariant map $h: \tilde{\Sigma} \rightarrow \SL(n,\C)/\SU(n) $ by standard bundle identifications. With these remarks in place, we note the following. 

\begin{lemma}\label{Lem:HarmonicMap}
Let $\Sigma= (S,J)$ be a Riemann surface. 
If a hermitian metric $h$ on an $\SL(n,\C)$-Higgs bundle $(\V,\Phi,\omega)$ on $\Sigma$ is harmonic, in the sense that $\nabla(h)$ is flat, then $h: \tilde{\Sigma} \rightarrow \mathrm{SL}(n,\C)/\mathrm{SU}(n)$ is a harmonic map of Riemannian manifolds. 
\end{lemma}

The idea of the proof is that the flatness of $\nabla$ is determined by the following system of PDE, called \emph{Hitchin's equations}, which decomposes the curvature of $\nabla$ into $\frakk$ and $\frakp$ parts as follows:
\begin{align}
    F_{\nabla_{\delbar_{\V},h}}+[\Phi,\Phi^*]&=0 \label{HitchinPDE},\\
    \nabla^{0,1}\Phi&=0 \label{HiggsHolomorphic}.
\end{align}
In this setting, the Higgs field $\Phi$ reinterprets as $\Phi = h^{-1}\del h$, up to a constant. Here, we can view $h: \tilde{\Sigma} \rightarrow \Mat_n(\C)$ as a matrix in coordinates, whereby $h^{-1}dh \in \Omega^1(\tilde{\Sigma}, \End(\underline{\C}^n))$ makes sense. By equivariance, the object $h^{-1}dh$ descends to $\Omega^1(\Sigma, \End(\V))$. See \cite{Li19} for further details on the identifications. 
Set $\X = \SL(n,\C)/\SU(n)$. Now, if $D$ is the connection induced in $\T^*\tilde{\Sigma} \otimes h^*T\X$, then the harmonicity equation $D^{0,1}\del h=0$ for the map $h:\tilde{\Sigma} \rightarrow\X$ corresponds exactly to the equation \eqref{HiggsHolomorphic} for the holomorphicity of the Higgs field $\Phi$, explaining the main idea of Lemma \ref{Lem:HarmonicMap}. 

The question of interest now is when a harmonic metric exists. The Hitchin-Simpson theorem states that under certain \emph{stability} conditions, this is always the case. 

\begin{definition}[Stability]\label{Defn:StabilitySL}
Let $\mathcal{H} = (\V,\Phi,\omega)$ be an $\SL(n,\C)$-Higgs bundle, so that $\deg(\V)=0$. 
\begin{itemize}
    \item $\mathcal{H}$ is \textbf{stable} when $\deg(\mathcal{E'}) <0$ for any proper holomorphic $\Phi$-invariant subbundle $\V'$ of $\V$.
    \item $\mathcal{H}$ is \textbf{polystable} when $(\V,\Phi) = \bigoplus_{i=1}^k (\V_i,\Phi_i)$ for degree zero $\Phi$-invariant sub-bundles $\V_i$ such that $\V = \bigoplus_{i=1}^k \V_i$ and Higgs sub-fields $\Phi_i \in H^0(\End(\V_i)\otimes \K)$ such that $\Phi = \sum_{i=1}^k \Phi_i$. 
\end{itemize}
\end{definition}
The $\Phi$-invariance condition above means that $\Phi(\mathcal{E}')\subset \mathcal{E}'\otimes \K$. 

Hitchin and Simpson proved the following fundamental result on existence of harmonic metrics. 
\begin{theorem}[\cite{Hit87,Sim88}]
Let $\mathcal{H}=(\V,\Phi, \omega)$ be an $\SL(n,\C)$-Higgs bundle. Then $\mathcal{H}$ admits a harmonic metric of $\omega$-unit volume, meaning $h =1$ on $\mathcal{O} \cong \det(\V)$, if and only if $\mathcal{H}$ is polystable. Moreover, when $\mathcal{H}$ is stable, then the metric is unique. 
\end{theorem}

The non-abelian Hodge correspondence $\NAH_{\Sigma}: \mathcal{M}_{\Sigma}(G) \rightarrow \chi(S,G)$ is induced by the following associations: $[(\V,\Phi)] \mapsto [ (E,\nabla)]\mapsto [\hol(\nabla)]$ from Higgs bundles to flat bundles to representations, up to equivalence. We refer the reader to the aforementioned references for further details.

\subsubsection{Representations to Higgs Bundles}
The inverse construction is due to Donaldson for $G =\SL(2,\C)$ and Corlette, who generalized his results to $G$ a semisimple complex Lie group. We focus on $G=\SL(n,\C)$ here. In particular, for $\rho: \pi_1S \rightarrow \SL(n,\C)$ a reductive representation, one can first form the associated flat bundle $(E, D)$, with $E_{\rho} = \tilde{\Sigma} \times_{\rho} \C^n$ and $D$ the flat connection on $E_{\rho}$ induced by the trivial connection on $\underline{\C}^n$. One then searches for a harmonic metric $h$ on the flat bundle, which again encodes a $\rho$-equivariant harmonic map $h:\tilde{\Sigma} \rightarrow \X$. To this end, split $D=D_h+\Psi_h$ uniquely into an $h$-connection $D_h$ and a 1-form $\Psi_h \in \Omega^1(\Sigma, \End(E))$ that is $h$-self-adjoint. Next, define an \emph{energy functional} on the space of hermitian metrics on $(E,D)$ as follows: $\mathscr{E}_{\rho}(h) = \int_{\Sigma} \tr(\Psi_h \wedge \Psi_h)$. We call $h$ a \emph{harmonic metric} on $(E,D)$ when $h$ is a critical point of the energy functional $\mathscr{E}_{\rho}$. Writing out the Euler-Lagrange equations for $h$ to be critical for $\mathscr{E}_{\rho}$, one re-obtains Hitchin's equations. The following result then completes the converse association of $\NAH_{\Sigma}$.\footnote{Here, we suppress the volume form $\omega$, but it is induced by $E_{\rho}$ due to the representation taking values in $\SL(n,\C)$.}

\begin{theorem}[{\cite{Don87, Cor88}}]
Let $\rho: \pi_1S \rightarrow \SL(n,\C)$ be a  representation. The flat bundle $(E_{\rho},D)$ admits a harmonic metric if and only if $\rho$ is reductive, and moreover $h$ is unique if $\rho$ is irreducible. When $h$ is harmonic, the map $h: \tilde{\Sigma} \rightarrow \X$ is a harmonic map of Riemannian manifolds. Moreover, in this case, $(E, D^{0,1}, \Psi_h^{1,0})$ is an $\SL(n,\C)$-Higgs bundle. 
\end{theorem}

The whole discussion is perhaps clarified by introducing the notion of a \emph{harmonic bundle} $(E,D,h)$, which is a flat bundle equipped with a harmonic metric. In the above, we quickly sketched the bijective correspondence between Higgs bundles, harmonic bundles, and representations, up to their respective equivalences, where the harmonic bundle is equivalently encoded by the associated equivariant harmonic map. Altogether, this comprises the non-abelian Hodge correspondence. 

In the case that we have $G$-Higgs bundles, rather than just $\GL(n,\C)$ or $\SL(n,\C)$-Higgs bundles, one can further demand harmonic metrics that are compatible with the additional structures (e.g. $h$ having $\omega$-unit volume). In this case, the flat connection $\nabla$ has holonomy in $G$ rather than just $\GL(n,\C)$ and moreover the harmonic map $h: \tilde{\Sigma} \rightarrow \GL(n,\C)/\mathsf{U}(n)$ takes values in the totally geodesic submanifold $\X = G/K$, the sub-symmetric space of $G$. The more general notion of stability and polystability for $G$-Higgs bundles leads to the generalization of Hitchin-Simpson, and to the more general non-abelian Hodge correspondence $\NAH_{\Sigma,G}: \mathcal{M}_{G}(\Sigma) \rightarrow \chi(S,G)$. We refer the reader to \cite{GPGMiR12} for the full details.

\subsection{\texorpdfstring{$\Gtwosplit$-Higgs Bundles}{G2'}}\label{Sec:G2'Higgs}
In this section, we discuss the technical details surrounding $\Gtwosplit$-Higgs bundles, recall some information from \cite{CT24} regarding harmonic metrics on such Higgs bundles.  

Recall that $\Gtwo^\C = \Aut(\imoct^\C, \times)$. Thus, reducing the structure group of a rank seven holomorphic vector bundle $\V$ from $\GL(7,\C)$ to $\Gtwo^\C$ amounts to placing a holomorphic cross product on $\V$ that is fiberwise isomorphic to that of $\imoct^\C$. This leads to the following definition. 

\begin{definition}
A $\Gtwosplit$-Higgs bundle structure on a $\SO_0(3,4)$-Higgs bundle $(\mathcal{U},\mathcal{V}, Q_\mathcal{U},Q_\mathcal{V}, \eta)$ is a holomorphic bundle map $\times_\V : \bigwedge^2 \V \to \V$ that defines a $Q$-cross product fiberwise, for $Q= Q_{\mathcal{U}}\oplus (-Q_{\mathcal{V}})$, such that $\eta \in \Der(\times_{\V})$.
\end{definition}
Going forward, the volume forms $\omega_{\mathcal{U}}$ and $\omega_{\mathcal{V}}$ in the $\SO_0(3,4)$-Higgs bundles will be tautological and hence omitted. 

The Cartan decomposition $\g_2'= \frakk \oplus \frakp $ is induced by the Cartan decomposition of $\solie(3,4)$. Indeed, this follows from general theory, or the explicit description of the maximal compact $\SO(4)\cong K<\Gtwosplit$ in $\S$\ref{Subsec:SubgroupHomogeneous}. Now,  associated to $(\mathcal{U},\mathcal{V}, Q_{\mathcal{U}}, Q_{\mathcal{V}}, \times_{\V})$ is the holomorphic $K^{\C}$-frame bundle for $K < \Gtwosplit$, given by intersecting the $\Gtwo^\C$-frame bundle $\mathrm{Fr}^\times(\V)$ of cross product preserving frames with the $\SO(3,\C)\times \SO(4,\C)$-frame bundle of $\V$. In particular, note that $\Phi \in H^0(\mathcal{P}\times_{Ad}\p^\C \otimes \K)$ holds automatically for $\Phi \in \Der(\times_{\V})$ such that $\Phi$ is an $\SO_0(3,4)$-Higgs field.

\begin{remark}
Let $Q$ be any non-degenerate complex bilinear form on $\C^7$. Then there is only one holomorphic $Q$-cross product on $\C^7$ up to isomorphism \cite{Fon18}. Thus, it is automatically true that $(\V|_p, \times_{\V})$ is fiberwise isomorphic to $(\imoct^\C, \times)$.  
\end{remark}

Next, we introduce the $\Gtwosplit$-Higgs bundles of interest. The notion of $\C$-cross product bases from Definition \ref{Defn:FCrossProductBasis} will play a key role. Indeed, our ambition is to leverage the particular shape of a holomorphic bundle to build form a vector bundle version of a $\C$-cross product basis. We shall again use the notation $(e_k)_{k=3}^{-3}$ for the model such basis from \eqref{BaragliaBasis} and $(c_{i,j})$ for the structure constants, i.e., $x_i\times x_j = c_{i,j} x_{i+j}$.  

Now, let us consider a vector bundle $\V = \bigoplus_{k=3}^{-3} \mathcal{L}_k$ that is a direct sum of holomorphic line subbundles with once-and-for-all fixed isomorphisms 
\begin{align}\label{FixedProductIsomorphism}
    \begin{cases}
        \mathcal{L}_0= \mathcal{O}\\
        \mathcal{L}_{-k}\cong \mathcal{L}_{k}^{-1}, \;\; k \in \{1,2,3\} \\
    \mathcal{L}_1 \mathcal{L}_2 \cong \mathcal{L}_{3}.
    \end{cases}
\end{align}
We can then leverage these isomorphisms to define a global holomorphic cross product $\times_{\V}: \V\times \V \rightarrow \V$ as follows. We decompose $\times_{\V} \in H^0(\Sigma, \Lambda^2(\V^*) \otimes \V)$ as a global map into sub-tensors of the form $\times_{\V} := \sum_{I}( \times_{\V})_{i,j,k} $, where $I : = \{ (i,j,k) \in \Z^3 \mid-3 \leq i , j, k \leq 3, \; i+j=k\}$, and $(\times_{\V})_{i,j,k}$ is a map
\[ (\times_{\V})_{i,j,k}: \mathcal{L}_i\otimes\mathcal{L}_j \rightarrow \mathcal{L}_k.\]
For indices $(i,j,k) \in I$, we define 
\[ (\times_\V)_{i,j,k} := c_{i,j}.\] 
Note that $c_{i,j} =0$ if the indices $(i,j,k)$ are not pairwise distinct. 
Concretely, the bundle cross product $\times_{\V}$ is described as follows: for any elements $x_i \in \mathcal{L}_i|_p, \;y_j \in \mathcal{L}_{j} |_p$, we have 
\begin{align}\label{BundleCrossProduct}
 x_i \times_{\V} y_j = c_{i,j} (x_i \otimes y_i) \in \mathcal{L}_{i+j}|_p. 
\end{align}
In other words, $\times_{\V}$ is nothing more than a carefully chosen linear combination of ordinary tensor products of line bundle factors. 
We now explain how to ensure that the map $\times_\V$ from \eqref{BundleCrossProduct} is compatible with an appropriate $\SO_0(3,4)$-structure and hence a legitimate bundle cross product.

Take $\V =\bigoplus_{k=3}^{-3}\mathcal{L}_i$ a holomorphic vector bundle satisfying \eqref{FixedProductIsomorphism}. We build the underlying $\SO_0(3,4)$-bundle data $(\U, \mathcal{V}, Q)$ as follows. Set :
\begin{align}
\begin{cases}\label{SO34Structure}
    \mathcal{U} &= \cl_{2} \oplus \cl_{0} \oplus \cl_{-2}, \\
    \mathcal{V} &= \cl_{3} \oplus \cl_{1} \oplus \cl_{-1} \oplus \cl_{-3}, \end{cases}
\end{align}
and then define the bilinear form $Q$ on $\mathcal{U}\oplus \mathcal{V}$ by 
\begin{align}\label{QonG2Higgs}
    Q=\sum_{i=3}^{3}(-1)^{i}Q_{i,-i}=\begin{pmatrix}
        & & & & & & -1\\
        & & & & & +1& \\
        & & & & -1 & &\\
        & & & +1& & & \\
        & & -1& & & & \\
        & +1& & & & & \\
        -1 & & & & & & \\
    \end{pmatrix},
\end{align}
where $Q_{i,-i}: \mathcal{L}_i\times \mathcal{L}_{-i}\rightarrow \C$ is the natural dual pairing. Observe that $Q = Q_{\mathcal{U}}\oplus -Q_{\mathcal{V}}$ obtains the desired direct sum form across $\mathcal{U} \oplus \mathcal{V}$.  
We then achieve a cross product $\times_{\V}$ on $\V$ by the above procedure, as clarified by the following proposition. 

\begin{proposition}[$\Gtwosplit$-Structure on Vector Bundle]\label{Prop:HiggsBundleCrossProduct}
Let $(\mathcal{E} =\mathcal{U}\oplus\mathcal{V}, Q)$ be a tuple satisfying \eqref{FixedProductIsomorphism},  \eqref{SO34Structure}, \eqref{QonG2Higgs}. Then the map $\times_{\V}$ defined by \eqref{BundleCrossProduct} is a $Q$-cross product. 

Moreover, we have the following coordinate description of the product $\times_{\V}$. 
Fix any point $p \in \Sigma$. Choose nonzero elements $x_2 \in \cl_{2}|_p$, $x_1 \in \cl_{1}|_p$. Define 
\begin{itemize}[noitemsep]
    \item $x_{3} := (x_2\otimes x_1) \in (\cl_{3})|_p$, 
    \item $x_{-k} :=x_k^{*} \in (\cl_{i}^{-1})|_p = \mathcal{L}_{-i}|_p$, for $i \in \{1,2,3\}$,\footnote{Here, we insist on taking duals using our fixed background isomorphism $\mathcal{L}_{-i}\cong \mathcal{L}_i^{-1}$ and \emph{not} with the non-degenerate pairing $Q$. These two identifications do not agree for all $i$, and this distinction appears in the proof below.} 
    \item $x_0 :=1|_p \in \mathcal{O}|_p$. 
\end{itemize}
The linear map $\Xi:(\V_{\mid p}, \times_\V) \rightarrow (\imoct^\C,\times)$ satisfying $x_k \mapsto e_k$, with $(e_k)_{k=3}^{-3}$ the model basis \eqref{BaragliaBasis}, is cross product preserving. That is, the structure constants of $(x_k)_{k=3}^{-3}$ and that of $(e_{k})_{k=3}^{-3}$ agree. 
\end{proposition}

\begin{proof}
Let $(x_k)_{k=3}^{-3}$ be the basis of $\V_{\mid p}$ extending $(x_2, x_1)$ according to the hypotheses. By definition of $\times_{\V}$, note that $(x_k)_{k=3}^{-3}$ is a $\C$-cross product basis, with some structure constants $(C_{i,j})$ satisfying $x_{i}\times_{\V}x_j=C_{i,j}x_{i+j}$. It remains, however, to show that $\times_{\V}$ is a $Q$-cross product basis. 

Note that we have $Q(x_k, x_{-j})={\mathrm{sign}(k)}\delta_{k,j}$ by definition of $Q$ in \eqref{QonG2Higgs}. Hence, the map $\Xi$ is orthogonal by \eqref{QModelBasis}, the expression of $q$ in the basis $(e_k)_{k=3}^{-3}$. Now, one verifies by bare hands means that by definition of $\times_{\V}$ and of $(x_k)_{k=3}^{-3}$ that the respectively structure constants $(C_{i,j})$ of $(x_k)_{k=3}^{-3}$ and $(c_{i,j})$ of $(e_k)_{k=3}^{-3}$ agree. 
The following calculation illustrates the general idea: to check whether $x_3 \times_{\V}x_{-2}\stackrel{?}{=}c_{3,-2}x_1$, we find by definition unraveling
\[ x_3\times_{\V} x_{-2} = (x_2\otimes x_1)\times_{\V}(x_2^*)= c_{3,-2}(x_2\otimes x_2^*\otimes x_1)=c_{3,-2}x_1.\]

One can easily check the remaining cases by similar arguments, eventually concluding $(C_{i,j}) =(c_{i,j})$. Thus, the linear map $\Xi$ is cross product preserving. 
Since $\Xi: (\V|_p, Q) \rightarrow (\imoct^\C, q)$ is orthogonal, we conclude $\times_{\V}$ is a $Q$-cross product since $\times_{\imoct^\C}$ is a $q$-cross product. 
\end{proof}

\begin{remark}
Note that the holomorphic vector bundle 
\[ \V = \mathcal{L}_2\mathcal{L}_1 \oplus \mathcal{L}_2 \oplus \mathcal{L}_1 \oplus \mathcal{O} \oplus \mathcal{L}_1^{-1}\oplus \mathcal{L}_2^{-1}\oplus (\mathcal{L}_1\mathcal{L}_2)^{-1}\] has transitions
pointwise of the form $g = \diag(wz, w,z,1, z^{-1},w^{-1},z^{-1}w^{-1})$, for $w,z \in \C^*$, which lie in a maximal torus $T^\C < \Gtwo^\C$ (in the basis \eqref{BaragliaBasis}). Hence, one can also introduce $\times_{\V}$ in local coordinates as in Proposition \ref{Prop:HiggsBundleCrossProduct} and argue $\times_{\V}$ is induced globally since $T^\C \subset \Gtwo^\C$ respects these identifications. 
\end{remark}

We now introduce the cyclic $\Gtwosplit$-Higgs bundles of interest, inspired by \cite{CT24}, as well as the two sub-families for which we will build geometric structures. 

\begin{definition}\label{Defn:G2Cyclic}
We call $(\V,\Phi, Q,\times_{\V})$ to be a \textbf{cyclic} $\Gtwosplit$-Higgs bundle when $\V = \bigoplus_{k=3}^{-3} \mathcal{L}_k$ satisfies \eqref{FixedProductIsomorphism}, \eqref{SO34Structure}, \eqref{QonG2Higgs}, is equipped with the holomorphic cross product $\times_{\V}$ in Proposition \ref{Prop:HiggsBundleCrossProduct}, and has Higgs field $\Phi$ of the form
\begin{equation}\label{CyclicG2'Higgs}
\begin{tikzcd}
 \mathcal{L}_3 \arrow[r, "\beta"] & \mathcal{L}_2\arrow[r, "\alpha"]& \mathcal{L}_1 \arrow[r, "-i\sqrt{2}\beta"]& \mathcal{L}_0 \arrow[r, "-i\sqrt{2}\beta"]&\mathcal{L}_{-1}\arrow[r, "\alpha"]&\mathcal{L}_{-2} \arrow[r, "\beta"]\arrow[lllll,bend left, "\delta"]&\mathcal{L}_{-3} \arrow[bend left, lllll, "\delta"]
\end{tikzcd}.
\end{equation}
Furthermore: 
\begin{itemize}[noitemsep]
    \item When $\mathcal{L}_{-1} \cong \K^{-1}$ and $\beta = 1$, we call the bundle $\boldsymbol{\beta}$-\textbf{cyclic}.
    \item When $\mathcal{L}_{2}\cong \mathcal{L}_1\K$ and $\alpha = 1$, we call the bundle $\boldsymbol{\alpha}$-\textbf{cyclic}.
\end{itemize} 
\end{definition}

By our study of root vectors in $\S$\ref{Subsec:LieTheory}, the Higgs field $\Phi$ in \eqref{CyclicG2'Higgs} satisfies $\Phi \in \Der(\times_{\V})$. \medskip 

\begin{remark}
Going forward, for the Higgs bundles in Definition \ref{Defn:G2Cyclic}, we will always discuss stability for the underlying $\SL(7,\C)$-Higgs bundles, in the sense of  Definition \ref{Defn:StabilitySL}. Stability for $\beta$-bundles and $\alpha$-bundles is considered in Sections \ref{Subsec:Ein23Structures} and Section \ref{Sec:StabilityAlplhaBundles}, respectively. 
\end{remark}

Since we have made a uniform convention for the structures $Q, \times_{\V}$ on a cyclic $\Gtwosplit$-Higgs bundle, we shall notationally suppress this data going forward. 

We note that only Higgs bundles that are both $\alpha$-cyclic and $\beta$-cyclic correspond to $\Gtwosplit$-Hitchin representations, namely when $\mathcal{L}_1 \cong \mathcal{K}$ and $\mathcal{L}_3 \cong \mathcal{K}^3$. We construct fibered geometric structures from $\beta$-bundles and $\alpha$-bundles in Sections \ref{Sec:G2Ein23GeomStructures} and \ref{Sec:PhoXGeometricStructures}, respectively.

We now consider the shape of the harmonic metric for cyclic $\Gtwosplit$-Higgs bundles. 
Here are two properties of the harmonic metric from \cite{CT24}. 

\begin{proposition}[$\Gtwosplit$-Harmonic Metric] \label{Prop:HarmonicMetricDiagonal}
    The harmonic metric $h = \diag(h_i)_{i=3}^{-3}$ on a polystable cyclic $\Gtwosplit$-Higgs bundle obtains the form 
    \[ h = \diag \left( h_1^{-1}h_2^{-1}, \; h_2^{-1}, \; h_1^{-1}, 1, \; h_1, \; h_2, \; h_1h_2 \right).\]
\end{proposition}

\begin{remark}\label{Remk:UnitaryCrossProd}
We will use frequently that for any point $p \in \Sigma$, or any open set $U \subset \Sigma$, we can find a basis for $\V|_p$ or a local frame for $\V|_U$, respectively, of the form $(x_k)_{k=3}^{-3}$ that is $h$-unitary and satisfies the multiplication Table \ref{Table:ComplexCrossProductBasis}. This holds for the basis $(x_k)_{k=3}^{-3}$ produced in Proposition \ref{Prop:HiggsBundleCrossProduct} when taking $x_1$ and $x_2$ to have unit $h$-norm by the symmetries of the harmonic metric $h$ in Proposition \ref{Prop:HarmonicMetricDiagonal}. 
\end{remark}

\begin{remark}Note that the $\Gtwosplit$-structure we place on our Higgs bundle, and the one from \cite{CT24} are different, but they only differ by a constant diagonal unitary gauge transformation, so the proposition still holds in our case.
\end{remark}

\subsubsection{Equivariant Frenet Frames}\label{Sec:AssociatedCurves}
In the present case of cyclic $\Gtwosplit$-Higgs bundles as in Definition \ref{Defn:G2Cyclic}, the harmonic metric $h$ yields a harmonic map we now denote $f: \tilde{\Sigma} \rightarrow \X_{\Gtwosplit}$. 
The harmonic metric $h$ also yields a $\nabla$-parallel real locus $\V^\R \subset \mathcal{E}$, where $\nabla = \nabla_{h} + \varphi + \varphi^{*h}$ is the associated flat connection. Here, $\V^\R$ is the fixed point set of an $\R$-linear involution $\lambda: \V \rightarrow \V $; see \cite[Proposition 4.16]{CT24}. The real structure $\lambda$ preserves $\mathcal{L}_0$ and each of $\mathcal{L}_i\oplus \mathcal{L}_{-i}$ for $i=1,2,3$. We write $\mathscr{L},T,N,B$ for the following sub-bundles of $\V^\R$: 
\begin{align}\label{LTNB}
    \begin{cases}
    \mathscr{L} &= \Fix(\lambda|_{\mathcal{L}_0})\\
    T &= \Fix(\lambda|_{\mathcal{L}_1\oplus \mathcal{L}_{-1}}) \\
    N&= \Fix(\lambda|_{\mathcal{L}_2\oplus \mathcal{L}_{-2}}) \\
    B &= \Fix(\lambda|_{\mathcal{L}_3\oplus \mathcal{L}_{-3}}).\end{cases}
\end{align} 

When $\beta\neq 0$, the names $\mathscr{L},T,N,B$ come from the relationship to an associated (equivariant) \emph{alternating almost-complex curves} $\nu: \tilde{\Sigma} \rightarrow \mathbb{S}^{2,4}$. These curves have been especially studied in the case of $\beta$-cyclic bundles, as in \cite{Bar10, CT24, Eva24a, Eva25}, for which they are immersed. However, the $\alpha$-bundles also yield alternating almost-complex curves $\tilde{\Sigma} \rightarrow \mathbb{S}^{2,4}$, as long as $\beta \neq 0$. In this case, the curve $\nu$ is not necessarily immersed and instead has the property that the second fundamental form $\sff$ is non-vanishing, evidenced by $\alpha$. Later, we will remark on the relation between the geometric structures built and this associated $J$-holomorphic curve. 

We now offer a brief summary of how the curve $\nu:\tilde{\Sigma} \rightarrow \quadric$ arises from a cyclic $\Gtwosplit$-Higgs bundle with $\beta\neq0$. This idea originates in Baraglia's thesis \cite[Section 3.6]{Bar10}, though he considered only the $\Gtwosplit$-Hitchin Higgs bundles. The key is to consider the tautological section $s\in \Omega^0(\Sigma, \V)$ given by $s(p) = 1 \in \mathcal{O}|_p$. Now, recall that the tuple $(Q, \times_{\V}, \lambda)$ is $\nabla=\nabla(h)$-parallel. The section $s$ is $\lambda$-real and satisfies $Q(s) = +1$ and hence $s \in \Omega^0(\Sigma, Q_+(\V^\R))$ corresponds to a $\rho$-equivariant map $\nu :\tilde{\Sigma}\rightarrow \quadric$. One can verify $s \times \nabla_{z}s= i \,s $, which holds by the definition of $\times_{\V}$. 
This condition corresponds to $\nu$ being $J$-holomorphic. Finally, the tuple $(\mathscr{L},T,N,B)$ corresponds to the so-called \emph{Frenet frame} of the curve $\nu$ (again, assuming $\alpha \neq 0$, otherwise $\nu$ is totally geodesic), with $T,N,B$ the tangent, normal, and binormal subspaces to $\nu$. These subspaces are real 2-planes, and also complex lines in $\nu^*\T\quadric$, that yield an orthogonal decomposition $\nu^*\T\quadric =T \oplus N \oplus B$ satisfying the hypothesis from Definition \ref{defn:FrenetFrame}. The $\rho$-equivariant harmonic map $f:\tilde{\Sigma} \rightarrow \X$ associated to an $\alpha$ or $\beta$-bundle, is described by $f(p) = \mathscr{L}_p \oplus N_p \in \Gr_{(3,0)}^\times(\imoct)$. 
In particular, $f$ is a generalized Gauss map of $\nu$, created by the spacelike components of the Frenet frame. We refer the reader to \cite{CT24} for an extensive background on the differential geometry of $J$-holomorphic curves in $\quadric$, focused on the case of $\beta$-bundles. 
\begin{remark}
We will frequently exchange the sub-bundles $(\mathscr{L},T,N,B)$ of $\V^\R$ as equivariant objects on $\tilde{\Sigma}$. For example, $T$ can be viewed as of type $T \in \Omega^0(\Sigma, \Gr_{(0,2)}(\V^\R))$ in the Higgs bundle, $T \in \Omega^0(\tilde{\Sigma}, \Gr_{(0,2)}(\underline{\R^{3,4}}))$ in the trivial bundle $\underline{\R^{3,4}} = \tilde{\Sigma} \times \R^{3,4}$, or as a map $T: \tilde{\Sigma} \rightarrow \Gr_{(0,2)}(\R^{3,4})$. 
\end{remark}

\subsubsection{Hitchin's Equations} 
We now derive Hitchin's equations for cyclic $\Gtwosplit$-Higgs bundles in general. The specific cases of $\beta$-cyclic and $\alpha$-cyclic Higgs bundles will be treated later on. 

Before the proof, we set some conventions. Let $g$ be a conformal metric on $\Sigma$ locally written 
$g_0|dz|^2$ and $\omega_g$ the volume form locally written $\omega_g = i \frac{g_0}{2}dz \wedge d\zbar$. We write $\Lambda$ for the contraction by $\omega_g$ and $\Delta_g$ for the Laplacian $\Delta_g = i \Lambda \del \delbar$, which is locally $\frac{2}{g_0} \partial_{z} \partial_{\zbar}$. We denote $F_{\mathcal{L}}$ for the curvature of a hermitian holomorphic line bundle $\mathcal{L}$ as well as $\kappa_g= i\Lambda F_{\K^{-1}}$ for the curvature function.

Note that by Proposition \ref{Prop:HarmonicMetricDiagonal}, there is no loss of generality in the imposed hypotheses on $h$ in the lemma below.

\begin{lemma}[Hitchin's Equations for Cyclic $\Gtwosplit$-Higgs bundles]\label{Lem:HitchinsEquationsGeneral}
Let $g$ be a conformal metric on $\Sigma$ and $(\V,\Phi)$ a cyclic $\Gtwosplit$-Higgs bundle on $\Sigma$. A hermitian metric $h = \diag(h_i)_{i=3}^{-3}$ on $\V$ with $h_0 = 1$, $h_{-i} = h_i^*$, and $h_3 = h_1h_2$ is harmonic if and only if the following equations hold:
\begin{align}
\begin{cases}\label{HitEquations_Primitive}
F_{{h}_3}&= \beta \wedge \beta^*-\delta\wedge \delta^*,\\
F_{{h}_2}&=\alpha\wedge \alpha^*-\beta \wedge \beta^*-\delta\wedge \delta^*, \\
F_{{h}_1}&=-\alpha \wedge \alpha^* +2\beta \wedge \beta^*. \\
\end{cases}
\end{align}
Additionally, $\alpha, \beta, \delta$ satisfy the following equations on $S$, away from their respective zeros: \begin{align}\label{HitchinEquationsBeta_FinalForm}
    \begin{cases}
   \Delta_g \log\lVert \alpha \rVert^2&= 2\lVert \alpha \rVert ^2-3\lVert \beta \rVert^2-\lVert \delta \rVert ^2 +\kappa_g,\\
 \Delta_g \log\lVert \beta \rVert^2&= 2\lVert \beta\rVert ^2-\lVert \alpha\rVert^2 + \kappa_g,\\
 \Delta_g \log\lVert \delta \rVert^2&=2\lVert \delta \rVert ^2-\lVert  \alpha\rVert ^2 + \kappa_g.
\end{cases}
\end{align}
\end{lemma}

\begin{proof}
Hitchin's equations for the harmonic metric can be written as:
\begin{align}\label{HitchinEquationGeneral}
     F_{\nabla^h}+ [\Phi, \Phi^*]=0. 
\end{align}
Let us write $\beta ' =-\sqrt{2} i\beta$. Each of the line bundles $\mathcal{L}_i$ inherits a metric $h_i$ with curvature 
$F_{\mathcal{L}_i}$. One finds that $[\Phi,\Phi^*]=\Phi\wedge \Phi^*+\Phi^*\wedge \Phi$ is diagonal, with entries indexed from $3$ to $-3$:
\begin{align*}
    [\Phi, \Phi^*]_{33}&= \delta\wedge\delta^*+\beta^*\wedge \beta, \\
    [\Phi, \Phi^*]_{22}&= \alpha^*\wedge \alpha+\beta \wedge \beta^* +\delta\wedge \delta^*\\
    [\Phi, \Phi^*]_{11}&= \alpha \wedge \alpha^*+(\beta')^*\wedge \beta'. 
\end{align*}
Note that $F_{-\mathcal{L}_i}=-F_{\mathcal{L}_i}$ for all $0\leq i\leq 3$, by the hypothesis $h_{-i}= h_i^*$. Then Hitchin's equation \eqref{HitchinEquationGeneral} combined with the above immediately yields \eqref{HitEquations_Primitive}.

Now, for a holomorphic hermitian line bundle $(\mathcal{L},h)\rightarrow \Sigma$ with a local non-vanishing holomorphic section $\eta \in H^0(\mathcal{L}|_O)$ over an open set $O$, the curvature $F_h$ of the Chern connection of $(\mathcal{L},h)$ is locally written 
\begin{align}\label{CurvatureEquation}
    F_h = i(\Delta_g \log||\eta||_h^2) \omega_g.
\end{align}
Observe that $\mathcal{L}_i\mathcal{K}$ and $\mathcal{L}_i\mathcal{L}_j\mathcal{K}$ are now hermitian holomorphic line bundles, where the conformal metric $g$ on $\T\Sigma \cong \K^{-1}$ induces a hermitian metric on $\K$. 
Now, we apply \eqref{CurvatureEquation} to the following holomorphic sections: $\beta \in H^0(\mathcal{L}_{-1}\mathcal{K})$, $\alpha \in H^0(\mathcal{L}_1\mathcal{L}_{-2}\mathcal{K})$,  $\delta\in H^0(\mathcal{L}_2\mathcal{L}_3\mathcal{K})$, and contract by the volume form $\omega_g$ on both sides. 
Using the relation $\eta \wedge \eta^*=-i ||\eta||^2\omega_g$, for $\eta \in H^0(\mathcal{L}\otimes \K)$, one obtains the equations \eqref{HitchinEquationsBeta_FinalForm}. 
\end{proof}

Let us briefly discuss the dimensions of the moduli space of `generic' cyclic $\sigma$-bundles, which comprise an open subset of the moduli of all polystable $\sigma$-bundles. 
For that purpose, we consider a parametrization $\mathcal{H}$ of a dense set of cyclic bundles, and show that this set is immersed in the character variety through the holonomy map.

\begin{proposition}
\label{prop:Dimension count}
    Let $S$ be a closed surface of genus $g \geq 2$ and let $0\leq d\leq 6g-6$. Consider the moduli space $\mathcal{H}_{\beta,d}^{gen}(S)$ of quadruples $(\Sigma, \mathcal{B},\alpha, \delta)$, where $\Sigma = (S,J)$ is a Riemann surface structure on $S$, $\mathcal{B}\rightarrow \Sigma$ is a holomorphic line bundle degree $d$ and $\alpha, \delta$ are respectively non-zero sections of $H^0(\mathcal{K}^3\mathcal{B}^{-1}), H^0(\mathcal{B}^2)$, that we quotient by the $\C^*$-action $\lambda\cdot (\alpha, \delta)=(\lambda^{-1}\alpha, \lambda^{2}\delta)$.

Similarly, let  $-2g+2\leq d\leq 2g-2$ and consider the moduli space $\mathcal{H}_{\alpha,d}^{gen}(S)$ of quadruples $(\Sigma, \mathcal{B},\beta, \delta)$, where $\Sigma = (S,J)$ is a Riemann surface structure on $S$, $\mathcal{T}\rightarrow \Sigma$ is a holomorphic line bundle degree $d$ and $\beta, \delta$ are respectively non-zero sections of $H^0(\mathcal{K}\mathcal{T}^{-1}), H^0(\mathcal{K}^3\mathcal{T}^3)$, that we quotient by the $\C^*$-action $\lambda\cdot (\beta, \delta)=(\lambda^{-1}\beta, \lambda^{3}\delta)$.

    Then :
    \begin{itemize}
        \item $\mathcal{H}_{\beta,d}^{gen}(S)$ is a smooth manifold of (real) dimension $28g-28-2(6g-6-d)$. 
        \item $\mathcal{H}_{\alpha,d}^{gen}(S)$ is a smooth manifold of (real) dimension $28g-28-4(2g-2-d)$. 
    \end{itemize}
    Moreover, the map $\hol: \mathcal{H}_{\sigma,d}^{gen}(S)\rightarrow \chi(S,\Gtwosplit)$ given by the holonomy of the corresponding $\Gtwosplit$-Higgs bundles is an immersion into the smooth locus of the character variety, which has dimension $28g-28$.
\end{proposition}

In other words, there are canonical maps from $\mathcal{H}_{\sigma,d}^{gen}(S)$ to the joint moduli space $\mathcal{M}_{\Gtwo^{\C}}(S)$ of $\Gtwo^{\C}$-Higgs bundles built in \cite{CTW25}. For example, there is a smooth map
\[\Xi_{\alpha}:\mathcal{H}_{\alpha,d}^{gen}(S) \rightarrow \mathcal{M}_{\Gtwo^{\C}}(S), \qquad  \big[(\Sigma,\mathcal{T},\beta,\delta)\big] \longmapsto \big[(\Sigma,\V(\mathcal{T}),\Phi(\beta,\delta))\big] , \]
where $(\V(\mathcal{T}),\Phi(\beta,\delta))$ is given by \eqref{AlphaCyclic}. Then $\Xi_{\beta}$ can be defined analogously via \eqref{BetaCyclicG2'Higgs}.

\begin{proof}
    First note by Propositions \ref{Prop:BetaBundleStability} and  \ref{Prop:AlphaStability} that the non vanishing of the sections $\alpha,\delta,\beta$ imply that the corresponding Higgs bundles are stable, and hence their holonomies belong to the smooth part of the character variety. The dimension $\dim_{\R}(\chi^{irr}(S,\Gtwosplit)) = \chi(S_g)\dim_{\R}(\Gtwosplit)$ is well known.
    
    For $\sigma=\beta$ the dimension count as well as the proof of the fact that this space is a smooth manifold is done in \cite[{Theorem 5.9}]{CT24}. 
    For $\sigma=\alpha$ we explain how to proceed similarly.  \medskip 
    
    \textbf{Case 1}: $-g+1<d\leq 2g-2$. Let us fix a Riemann surface structure $\Sigma$ on $S$. Choosing a line bundle $\mathcal{T} \in \Pic^d(\Sigma)$ together with a nonzero projective class $[\beta]\in \mathbb{P}(H^0(\mathcal{T}^{-1}\mathcal{K}))$ is equivalent to choosing a divisor $D$ on $\Sigma$ of degree $2g-2-d$. The space of such pairs $(\mathcal{T}, [\beta])$ is biholomorphic to the symmetric product $\mathrm{Sym}^{2g-2-d}(\Sigma)$ of $\Sigma$, which has real dimension $4g-4-2d$.
    The space of triples $(\Sigma, \mathcal{T}, [\beta])$ is a complex manifold, as the universal symmetric power $\Sym^{2g-2-d}(\mathcal{C}/T)$. Here, $\mathcal{C} \rightarrow T$ is the universal curve over Teichm\"uller space $T=T(S)$, with fiber $\mathcal{C}|_{[\Sigma]}$ biholomorphic to $\Sigma $ (see \cite[Section 6.8]{Hub06}), and we denote $\Sym^k(Y/X)$ for the relative $k$-fold symmetric product of a proper holomorphic submersion $Y \rightarrow X$ of complex manifolds. That is, $\Sym^k(Y/X) \rightarrow X$ has fiber $\Sym^k(Y/X)_{\mid x}\cong \Sym^k(Y_{\mid x})$.  
    In our case, $\Sym^k(\mathcal{C}/T)$ is defined as the quotient of the fiber product\footnote{The \emph{fiber product} is defined in the holomorphic category just as in the smooth category, see \cite[Corollary 0.32]{Fis}.} $\mathcal{C}^k = \mathcal{C}\times_{T} \times \cdots \times_{T} \mathcal{C}$ by the fiberwise action of the symmetric group $\mathfrak{S}_k$. Hence, $\Sym^k(\mathcal{C}/T) \rightarrow T$ has fiber $\Sym^k(\mathcal{C}/T)_{\mid {\Sigma}}\cong \Sym^k(\Sigma)$. 
    For completeness, we briefly explain how to see the manifold structure on $\Sym^k(\mathcal{C}/T) $.
    
    Now, let $(\Sigma, D) \in \Sym^k(\mathcal{C}/T)$ be arbitrary.
    Let us write $D= \sum_{i=1}^n k_i(p_i)$, where $p_i \in \Sigma$ are pairwise distinct and $\sum_{i=1}^nk_i=k$. 
    Then there is a smooth open neighborhood of $x=(\Sigma, D)\in \Sym^k(\mathcal{C}/T)$ of the following form:
    \[ U\times \prod_{i=1}^n\Sym^{k_i}(\mathbb{D}), \]
    where $U$ is an open set in $T$, $\mathbb{D} \subset\mathbb{C}$ denotes the open unit disk, and we choose trivializations of the form $U\times \mathbb{D}$ of \emph{local} neighborhoods of $(\Sigma,p_i)$ in $\mathcal{C}|_{U}$, for $1\leq i \leq n$.
    Each of the spaces $\Sym^{k_i}(\mathbb{D})$ admits a holomorphic structure via elementary symmetric polynomials. This structure is independent of the chosen trivialization of the neighborhood of $x$, as the holomorphic structure on $\Sym^k(\mathbb{D})$ is preserved by biholomorphisms.

    We have fixed $d > -g+1$ in this case. For such $d$, the space $\mathcal{H}_{\alpha, d}^{gen}(\Sigma)$ of $\alpha$-bundles on the fixed Riemann surface $\Sigma$ with $\beta,\delta\neq 0$ fibers over $\Sym^{2g-2-d}(\Sigma)$, with fiber $H^0(\mathcal{T}^3\mathcal{K}^3)\setminus \lbrace 0\rbrace$. By Riemann-Roch, $\dim_{\R}(H^0(\mathcal{T}^3\mathcal{K}^3)\setminus \lbrace 0\rbrace ) = 10g-10+6d$, independent of $\mathcal{T}$. 

    Now we check that the projection $\mathcal{H}_{\alpha, d}^{gen}(\Sigma)\rightarrow \Sym^{2g-2-d}(\Sigma)$  defines a smooth fiber bundle. Indeed, given $(\mathcal{T},[\beta])$, once we fix a representative $\beta$ of $[\beta]$, choosing $\delta\in H^0(\mathcal{T}^3\mathcal{K}^3)\setminus \lbrace 0\rbrace$ amounts to selecting a holomorphic sextic differential $q_6 \in H^0(\K_{\Sigma}^6)$ whose divisor $D'$ satisfies $D'-3D\geq 0$. More precisely, a triple $(\mathcal{T},[\beta],q_6)$, where $q_6 \in H^0(\K_{\Sigma}^6(-3D))$, determines $(\mathcal{T},\beta,\delta)$ up to the prescribed $\C^*$-action. 
    Note that each pair $(\mathcal{T},[\beta])$, or equivalently, each divisor $D \in \Sym^{2g-2-d}(\Sigma)$, defines a finite dimensional subspace $H^0(\K_{\Sigma}^6(-3D))$ of the finite dimensional space $H^0(\K^6_{\Sigma})$, with constant dimension, that varies smoothly with respect to $D$. This hence defines the fiber bundle structure $\mathcal{H}_{\alpha, d}^{gen}(\Sigma) \rightarrow \Sym^{2g-2-d}(\Sigma)$. 
    
    We now vary the Riemann surface structure. Recall that there is a holomorphic vector bundle $\mathcal{Q}^6 \rightarrow T$ over Teichm\"uller space with fiber $\mathcal{Q}^6|_{\Sigma} \cong H^0(\K^6_{\Sigma})$ \cite{Ber61}. 
    We can check in the above local coordinates for $\Sym^{2g-2-d}(\mathcal{C}/T)$ that the projection $\mathcal{H}_{\alpha, d}^{gen}(S)\rightarrow \Sym^{2g-2-d}(\mathcal{C}/T)$ defines a smooth fiber bundle, whose fiber is a fiber sub-bundle of the pullback of $\mathcal{Q}^6$ over $\Sym^{2g-2-d}(\mathcal{C}/T)$. In total, we have then verified the manifold structure on $\mathcal{H}_{\alpha, d}^{gen}(S)$ for $-g+1<d\leq 2g-2$. 
   
    Adding the dimensions counts of $(\mathcal{T},[\beta])$ and $\delta$ together with the $6g-6$ dimensions of Teichm\"uller space for the choice of $\Sigma$, we obtain a smooth manifold of total (real) dimension $28g-28-4(2g-2-d)$.\medskip 

    \textbf{Case 2}: $-2g+2\leq d\leq -g+1$. We proceed similarly. The space of pairs $(\mathcal{T}, [\delta])$ of $\mathcal{T} \in \Pic^d(\Sigma)$ equipped with a projective class $[\delta] \in \mathbb{P}(H^0(\mathcal{T}^3\mathcal{K}^3))$ is biholomorphic to the space $\Sym^{3d+6g-6}(\Sigma)$ of divisors of degree $3d+6g-6$ on $\Sigma$. For such $d$, the space $\mathcal{H}_{\alpha,d}^{gen}(S)$ is the disjoint union over $(\Sigma, \mathcal{T}, [\delta])\in \Sym^{3d+6g-6}(\mathcal{C}/T)$, of the space $H^0(\mathcal{T}^{-1}\mathcal{K})\setminus \lbrace 0\rbrace$ parametrizing the choice of $\beta$, whose dimension is $2g-2-2d$ by Riemann-Roch. With the same argument as before we check that $\mathcal{H}_{\alpha,d}^{gen}(S)$ is a smooth manifold, and a fiber bundle over $\Sym^{3d+6g-6}(\mathcal{C}/T)$. We obtain the same formula for the total dimension as in Case 1. \medskip

    Next, we check that the holonomy map is an immersion when the Riemann surface structure is fixed. We argue for Case 1, with a similar argument for Case 2.
    
    Now, let $\Sigma \in T(S)$ be fixed. We claim the projection $\pi: \mathcal{H}_{\alpha,d}^{gen}(\Sigma)\rightarrow \mathcal{M}_{\Gtwosplit}^{\text{s},6}(\Sigma)$ to the moduli space of 6-cyclic stable $\Gtwosplit$-Higgs bundles over $\Sigma$ is an immersion.
    Given a stable $6$-cyclic $\Gtwosplit$-Higgs bundle $(\V,\Phi)$, one can associate its $6$-cyclic splitting $\V=\bigoplus_{i \in \Z_6} \V_i$, which is uniquely defined. Indeed, if a given stable bundle admits two such splittings, there must be two $\Gtwosplit$ gauge transformations $g_1,g_2$ such that $g_1\cdot(\mathcal{E},\Phi)=g_2\cdot (\mathcal{E},\Phi)=(\mathcal{E}, e^{2i\pi/6}\Phi)$. This implies that $g_1g_2^{-1}$ is an automorphism of the stable bundle $(\mathcal{E},\Phi)$, but $\SL(7,\C)$-stable bundles have no non-trivial automorphisms apart from multiples of the identity, and $\Gtwosplit$ has no center, so $g_1=g_2$ and the splittings agree. 
    This allows to extract the divisor of $\beta$. Note also that the Hitchin fibration defines a smooth map $\mathcal{H}_{\alpha,d}^{gen}(\Sigma) \rightarrow H^0(\K^{6}_{\Sigma})$ given by $[(\mathcal{T},\beta, \delta)]\mapsto  \beta^3\delta$. The map $[(\mathcal{T},\beta, \delta)]\mapsto (\text{div}(\beta), \beta^3\delta)$ factors through $\pi$, and is an immersion. This proves that for a fixed $\Sigma$ the projection to the moduli space of Higgs bundles is an immersion.
     
    We now vary $\Sigma$. The differential of the holonomy map when moving the Riemann surface structure does not vanish since it factors smoothly through the holonomy map on the 6-cyclic part of the joint moduli space of Higgs bundles constructed in \cite{CTW25} (denoted there by $\mathbf{Y}_6 \cap \mathbf{W}_0$), which is an immersion by \cite[Theorem G]{CTW25}. 
\end{proof}

\section{\texorpdfstring{$(\Gtwosplit, \Ein^{2,3})$}{(G2',Ein23)}-Geometric Structures}\label{Sec:G2Ein23GeomStructures}

In this section, we consider the representations $\rho:\pi_1S\rightarrow \Gtwosplit$ associated to the $\beta$-cyclic Higgs bundles defined in Definition \ref{Defn:G2Cyclic} via the non-abelian Hodge correspondence, including, but not limited to, the $\Gtwosplit$-Hitchin case; see Remark  \ref{Remk:BetaSpecialCases}.  When $\rho$ is Hitchin, it is $\beta$-Anosov and the fiber of the Tits metric thickening domain $\Omega \subset \Ein^{2,3}$ is $\Ein^{2,1}$ by \cite{DE26}. 
Motivated by this result, for each representation $\rho$ associated to a $\beta$-cyclic bundle, we build a 5-manifold $M$ that is an $\Ein^{2,1}$-fiber bundle $M \rightarrow S$ with a fibered $(\Gtwosplit, \Ein^{2,3})$-structure whose holonomy descends to $\pi_1S$ as $\rho$. 
When the representation $\rho$ is Hitchin, we then show the resulting geometric structures agree with those defined by domains of discontinuity in \cite{GW12, KLP18}. 

\subsection{The Flag Manifold \texorpdfstring{$\Ein^{2,3}$}{Ein23} }\label{Sec:Ein23ModelSpace}

We now recall relevant features of $\Ein^{2,3}$ as a $\Gtwosplit$-homogeneous space. The most essential property will be its realization as a principal $\sphere^3$-bundle over $\RP^2$, and the flexibility of such fibrations.\medskip

The Einstein (2,3)-universe $\Ein^{2,3}$ is the projective null quadric in $ \imoct \cong \R^{3,4}$, namely,
\[ \Ein^{2,3} = \{ [x] \in \mathbb{P}\imoct \mid q(x)=0\}.\]

The first geometric description of $\Ein^{2,3}$ that we shall frequently appeal to is as follows. 

\begin{proposition}[Spacelike 3-plane Models]\label{Prop:Ein23EasyModel}
Fix $U \in \Gr_{(3,0)}(\imoct)$. There is a 2-1 smooth covering map $F_U : Q_+(U)\times Q_-(U^\bot) \rightarrow \Ein^{2,3}$ by $(u,v) \mapsto [u+v]$. Hence, $\Ein^{2,3}\cong (\mathbb{S}^2\times \mathbb{S}^3)/ (-\id, -\id)$. 
\end{proposition}

\begin{proof}
If $\ell \in \Ein^{2,3}$, then $\pi_{U}(\ell) \neq 0$ and $\pi_{U^\bot}(\ell) \neq0$. Since $q(\ell)=0$, we have $\ell = [u+v]$ for $u \in Q_+(U)$, $v \in Q_-(U^\bot)$. As $u \in \pi_U(\ell), \, v\in \pi_{U^\bot}(\ell)$, the (evidently smooth) map $F_U$ is clearly 2-1. 
\end{proof}

The map $F_U$ in Proposition \ref{Prop:Ein23EasyModel} is far from $\Gtwosplit$-equivariant, but is $K_U=\Stab_{\Gtwosplit}(U)$-equivariant. The most important feature of Proposition \ref{Prop:Ein23EasyModel} is the flexibility in the choice of $U$. In the construction of geometric structures in $\S$\ref{Subsec:Ein23Structures}, for points $p \in \widetilde{M}$, we shall understand the developed fibers $\dev(\widetilde{M}_{p}) \subset \Ein^{2,3}$ of the $(\Gtwosplit, \Ein^{2,3})$-manifold $M \rightarrow S$ through the models $F_{U(p)}$ of different 3-planes $U(p)$. 

Now, let us take a fixed identification $F_U: \Ein^{2,3} \rightarrow (\mathbb{S}^2 \times \mathbb{S}^3)/\sim$ from Proposition \ref{Prop:Ein23EasyModel}, with respect to $U \in \X_{\Gtwosplit}$. 
Our next description of $\Ein^{2,3}$ uses the cross product to kill off the $\Z_2$-quotient from the model $F_U$. For the result, set $ K_U = \Stab_{\Gtwosplit}(U)$. Also, note the extra demand on $U$. 

\begin{proposition}[\cite{BH14}]\label{Prop:Ein23ModelExplicit}
For $U \in \Gr_{(3,0)}^\times(\imoct)$, there is a $K_U$-equivariant diffeomorphism from $\RP^2 \times \mathbb{S}^3$ onto $\Ein^{2,3}$ given by $G_U: \mathbb{P}(U) \times Q_-(U^\bot) \rightarrow \Ein^{2,3}$ via $( [u], v) \mapsto [ u +u\times v]$. 
\end{proposition} 

\begin{proof}
First, note that the pre-image of a point $\ell =[u+v]\in \Ein^{2,3}$ obtains the form $F_U^{-1}([u+v] = \{(u,v),(-u,-v)\}$. 
The map $G_U$ removes the $\Z_2$-ambiguity.

Fix $u \in Q_+(U)$. Since the map $\mathcal{C}_u: U^\bot \rightarrow U^\bot $ via $\mathcal{C}_u(w) = u\times w$ satisfies $\mathcal{C}_u^{\circ 2} = -\id_{U^\bot} $ by the double cross product identity \eqref{DCP}, the inverse of $G_U$ is the well-defined map 
\[ [x + y] \mapsto \pm (x,y) \mapsto ( [x],\,-x\times y). \] 
The map $G_U$ is a diffeomorphism because it is $K:= \Stab_{\Gtwosplit}(U)$-equivariant and $K$ acts transitively on $\Eintwothree$ by Proposition \ref{Prop:KEinsteinModel}. 
\end{proof} 

Next, we discuss the structure of $\Eintwothree$ as a homogeneous $K$-space. 

\begin{proposition}\label{Prop:KEinsteinModel}
Let $U \in \mathbb{X}_{\Gtwosplit}$ and set $K:= K_U =\Stab_{\Gtwosplit}(U)$. There $K$-equivariant diffeomorphism $\Eintwothree \cong K/( \SO(2) \times \Z_2)$.
\end{proposition} 

\begin{proof}
The maximal compact subgroup $K:= \Stab_{\Gtwosplit}(U)$ is identified in the Stiefel triplet model with 
\begin{align}\label{KStiefelModel}
V_{(+,+,-)}(U) := \{ (u,v,z) \in \imoct^3 \; | \; u, v \in Q_+(U), \; u \cdot v = 0, \; z \in Q_-(U^\bot) \}.
\end{align}
That is, all transformations in $K$ are uniquely prescribed by their action on any triple $p_0 \in V_{(+,+,-)}(U)$ by Proposition \ref{Prop:FirstStiefel}, and $\varphi \in K$ necessarily has $\varphi \cdot p_0 \in V_{(+,+,-)}(U)$. 
Hence, the map $K \rightarrow V_{(+,+,-)}(U)$ by $\varphi \mapsto \varphi \cdot p_0$ is a $K$-equivariant diffeomorphism. 

With this out of the way, the rest is simple. Fix $\ell \in \Eintwothree$ and write $F_U^{-1}(\ell)=\pm (u,v)$. 
Then define $H:= \Stab_{K}(\ell)$. 
Fix a basepoint $p_0$ of the form $p_0=(u, y, v) \in V_{(+,+,-)}(U)$. Then for $\varphi \in H$, there is $w \in Q_+(U \cap u^\bot) \cong \mathbb{S}^1$
such that either $\varphi( u , y, v) = (u, w, v)$ or $\varphi( u , y, v) = (-u, w, -v)$. Hence, define $\eta \in \Gtwosplit $ by $\eta \cdot (u, y,  v) = (-u, y, -v)$. 
Then $H \cong \SO(2) \times \Z_2$ is a direct product, with $\Z_2$-factor generated by $\eta$ and $\SO(2)$-factor given by $ \Stab_{\Gtwosplit}(u) \cap \Stab_{\Gtwosplit}(v) \cap K$. 
\end{proof} 

We now describe realizations of $\Ein^{2,3}$ as a principal bundle. To this end, define the pointwise stabilizer  $H_U:=\Stab^{pt}_{\Gtwosplit}(U)$ that fixes $U$ identically. The subgroup $H_U$ is isomorphic to $\Sp(1)$ by Proposition \ref{Prop:FirstStiefel}; any orbit map $H_U \rightarrow Q_-(U^\bot)$ by $\psi \mapsto \psi\cdot z$, for $z \in Q_-(U^\bot)$, is a diffeomorphism. 

\begin{proposition}[$\Ein^{2,3}$ as a Principal Bundle]\label{Prop:Ein23Principal}
Fix a point $U \in \mathbb{X}_{\Gtwosplit}$. Then the projection map $\pi_U: \Ein^{2,3} \rightarrow \Gr_1(U)$ realizes $\Ein^{2,3}$ as a principal $H_U:=\Stab_{\Gtwosplit}^{pt}(U)$-bundle.
\end{proposition}

\begin{proof}
Given any point $[u] \in \Gr_1(U)$, choose $\hat{u} \in Q_+([u])$. The fiber $\pi^{-1}([u])$ consists of all isotropic lines of the form $[\hat{u}+\hat{v} ] \in \Ein^{2,3}$ with $ \hat{v} \in Q_-(U^\bot)$. Hence, the fiber is topologically $\mathbb{S}^3 \cong Q_-(U^\bot)$. It is evident that $H_U$ preserves the fibers. 
In the Stiefel model $V_{(+,+,-)}(U)$ for $K_U:=\Stab_{\Gtwosplit}(U)$, as in \eqref{KStiefelModel}, one sees the subgroup $H_U$ acts simply transitively on the fibers of $\pi_U$. The claim follows. 
\end{proof}

\begin{remark}
The space $\Ein^{2,3}$ is also an $\SO(3,4)$-homogeneous space. However, there is additional structure that is only $\Gtwosplit$-invariant, namely a (2,3,5)-distribution $\mathscr{D} \subset \T\Ein^{2,3}$. See \cite[Section 8.2]{Eva24a} for relations between $\mathscr{D}$ and annihilators. 
This distribution $\mathscr{D}$ is related to one of the first descriptions of $\g_2'$ given by Cartan in \cite{Car10}; see \cite{Agri08, BH92}. 
\end{remark}

\subsection{\texorpdfstring{$\Eintwothree$}{Ein23}-structures for \texorpdfstring{$\beta$}{beta}-cyclic Bundles}\label{Subsec:Ein23Structures}  

We now construct fibered geometric structures for the representations corresponding to $\beta$-bundles. 
More specifically, we construct $(\Gtwosplit, \Ein^{2,3})$-structures on $\Ein^{2,1}$-fiber bundles over closed surfaces $S$. 

A $\beta$-bundle on a Riemann surface $\Sigma$ has the form \eqref{CyclicG2'Higgs} with $\mathcal{T} \cong \mathcal{K}$ and $\beta=1$. We are then reduced to considering the following $\Gtwosplit$-Higgs bundles in the present section: 
\begin{equation}\label{BetaCyclicG2'Higgs}
\begin{tikzcd}
\mathcal{B} \arrow[r, "1"] & \mathcal{B}\mathcal{K}^{-1}\arrow[r, "\alpha"]&\mathcal{K}\arrow[r, "-i\sqrt{2}"]&\mathcal{O} \arrow[r, "-i\sqrt{2}"]&\mathcal{K}^{-1} \arrow[r, "\alpha"]&\mathcal{B}^{-1}\mathcal{K} \arrow[r, "1"]\arrow[bend left, lllll, "\delta"]&\mathcal{B}^{-1} \arrow[bend left, lllll, "\delta"]
\end{tikzcd}.
\end{equation}
Here, $\mathcal{B} $ is a holomorphic line bundle on $\Sigma$, $\alpha \in H^0(\Sigma, \mathcal{B}^{-1} \mathcal{K}^3)$, and $\delta \in H^0(\Sigma, \mathcal{B}^2)$. 

The relevant (poly)stability considerations for $\beta$-bundles was given in \cite{CT24}, up to one small but noteworthy change we highlight here. 
\begin{remark}[Totally Geodesic $\beta$-Curves]
\label{Remk:TotallyGeodesic}
Suppose $(\V,\Phi)$ is a $\beta$-cyclic bundle on $\Sigma$, as in \eqref{BetaCyclicG2'Higgs}, that is polystable. The associated $J$-holomorphic curve $\nu:\tilde{\Sigma} \rightarrow \quadric$ is totally geodesic if and only if $\alpha = 0$  \cite[Proposition 3.23, Theorem 3.24]{CT24}. 
\end{remark}
In \cite{CT24}, they exclude the polystable bundles with $\alpha=0$ due to the focus on non-trivial, i.e., non-totally geodesic, $J$-holomorphic curves. In our case, the strictly polystable case $\alpha =\delta=0$ is one where our method to build geometric structures \emph{does apply}, it is even the simplest case in which our techniques apply, so we wish to include this case as a possibility. 

\begin{proposition}[Stability for $\beta$-bundles {\cite{CT24}}]\label{Prop:BetaBundleStability}
Let $\mathcal{H} = (\V,\Phi)$ be a $\beta$-cyclic Higgs bundle on $\Sigma$ as in Definition \ref{Defn:G2Cyclic}. If $\mathcal{H}$ is polystable, then $0 \leq \deg(\mathcal{B}) \leq 6g-6$. Moreover, if $\alpha\neq 0$, we have:
\begin{enumerate}[noitemsep]
    \item If $\mathcal{H}$ is stable, then $\deg(\mathcal{B})> 0$. 
    \item If $g-1 <\deg(\mathcal{B}) \leq 6g-6$, then $\mathcal{H}$ is stable. 
    \item If $0 < \deg(\mathcal{B}) \leq g-1 $, then $\mathcal{H}$ is stable if and only if it is polystable, which occurs exactly when  $\delta \neq 0$. 
    \item If $\deg(\mathcal{B})=0$, then $\mathcal{H}$ is polystable if and only if $0\neq\delta \in H^0(\mathcal{O})$, which entails $\delta$ is non-vanishing.
\end{enumerate}
Now if $\alpha=0$,
\begin{enumerate}[noitemsep]
    \item[(5)] If $\delta=0$, the bundle is polystable if and only if $\deg(\mathcal{B})=g-1$. 
    \item[(6)] If $\delta\neq 0$, the bundle is polystable if and only if $0 \leq \deg(\mathcal{B}) <2g-2$ and additionally stable when $0<\deg(\mathcal{B})$.
\end{enumerate}
\end{proposition}

\begin{proof}
Note that \cite[Proposition 5.6]{CT24} addresses (1)--(3), while \cite[Theorem 5.14]{CT24} addresses (4). Points (5)--(6) are similar and can be checked directly.
\end{proof}

\begin{remark}[Special $\beta$-bundles] \label{Remk:BetaSpecialCases}
Here are some noteworthy cases of polystable $\beta$-bundles:
\begin{itemize}
    \item Then $\deg(\mathcal{B})=6g-6$ if and only if $\mathcal{B} \cong \K^3$ if and only if the corresponding representation is $\Gtwosplit$-Hitchin \cite{CT24}. Here, the `if' comes from Hitchin's parametrization of $\Hit(S,\Gtwosplit)$ and the `only if' comes from polystability and the fact that $0\neq \alpha \in H^0(\K^3\mathcal{B}^{-1})$ is a holomorphic section of a degree zero bundle and hence must be non-vanishing.
    \item The Higgs bundle associated to a $\beta$-Fuchsian representation $\pi_1S\rightarrow \SL(2,\R)\hookrightarrow \Gtwosplit$ in the $\SL(2,\R)$-subgroup of the short root $\beta$, described in Appendix \ref{Subsec:SL2ShortRoot}, 
    obtains the following form: 
    \[ \mathcal{K}^{1/2} \stackrel{1}{\longrightarrow}\mathcal{K}^{-1/2} \; \; \;\;\;\mathcal{K}\stackrel{-i\sqrt{2}}{\longrightarrow}\mathcal{O} \stackrel{-i\sqrt{2}}{\longrightarrow}\mathcal{K}^{-1}\;\;\;\;\;\mathcal{K}^{1/2} \stackrel{1}{\longrightarrow}\mathcal{K}^{-1/2}. \]
    \item If $\alpha = 0$ and $\delta =0$, then the representation factors through the $\SL(2,\R)$ subgroup of the short root $\beta$, twisted by a copy of $\SO(2)$ in the centralizer. These Higgs bundles can be found by twisting $\mathcal{B} \cong \K^{1/2}$ by a line bundle $\mathcal{L}$ of degree $0$.
    \item If $\alpha =0$, then the corresponding representation factors through $\SO(2,2)$. 
    \item If $\mathcal{B}^2 \cong \mathcal{O}$ and $\alpha , \delta\neq 0$, then $\delta $ is pointwise non-vanishing and the corresponding representation factors through (a $\Z_2$-cover of) $\SL(3,\R)$ and is (essentially) a non-Fuchsian $\SL(3,\R)$-Hitchin representation. See \cite[Theorem 5.14]{CT24} and the surrounding discussion for a precise statement. 
    \item When $\mathcal{B}\cong\mathcal{O}$ and $\delta =1$, $\alpha=0$, we obtain the uniformizing representation of a Fuchsian-Hitchin $\SL(3,\R)$-representation included in $\Gtwosplit$ through the $\PSL(2,\R)$-subgroup in Appendix \ref{Appendix:SL2Number4}. A uniformizing Higgs bundle is given as follows: 
    \[\begin{tikzcd}
\mathcal{O} \arrow[r, "1"] & \mathcal{K}^{-1} &\mathcal{K}\arrow[r, "-i\sqrt{2}"]&\mathcal{O} \arrow[r, "-i\sqrt{2}"]&\mathcal{K}^{-1} &\mathcal{K} \arrow[r, "1"]\arrow[bend left, lllll, "1"]&\mathcal{O} \arrow[bend left, lllll, "1"]
\end{tikzcd}.\]
\end{itemize}
\end{remark}
We emphasize that $\beta$-bundles give representations more general than $\Gtwosplit$-Hitchin representations: by Proposition \ref{Prop:BetaBundleStability} and Remark \ref{Remk:BetaSpecialCases}, if at least one of $\alpha, \delta$ is not zero and $0 < \deg(\mathcal{B}) < 6g-6$, then $\rho$ is irreducible and non-Hitchin. \medskip 

Despite the lengthy remark above, there remains one important case to further clarify.
\begin{remark}\label{Remk:SL3Beta}
In the case $\delta$ is non-vanishing, then $\mathcal{B}^2\cong\mathcal{O}$, which entails the corresponding representation $\rho$ factors through a certain $\Z_2$-cover of $\SL(3,\R)$ by Remark \ref{Remk:BetaSpecialCases}. This $\Z_2$-cover is explicitly $\Stab_{\Gtwosplit}([y])$ for a timelike line $[y] \in \mathbb{P}Q_-(\imoct)$. In particular, this implies $\rho$ is \textbf{not} $P_{\beta}$-Anosov. (However, it is possible $\rho$ is $P_{\alpha}$-Anosov.)
\end{remark}
The invariant $\deg(\mathcal{B})$ detects whether the associated curve $\nu$ is \emph{linearly full}.
\begin{remark}
Let $(\V,\Phi)$ be a polystable $\beta$-bundle with $\alpha \neq 0$. The associated $J$-holomorphic curve $\nu:\tilde{\Sigma}\rightarrow \quadric$ is \textbf{not} linearly full if and only if $\deg(\mathcal{B})=0$. Moreover, in this case, $\delta\neq 0$ by Proposition \ref{Prop:BetaBundleStability} and $\nu$ is linearly full in a copy of $\R^{3,3}$. See \cite[Theorem 5.14]{CT24}. 
\end{remark}

We now determine the relevant maximum principles for $\beta$-bundles for the holomorphic differentials $\alpha$ and $\delta$. These bounds will be essential for verifying the developing map in the next subsection.  
\begin{lemma}\label{Lem:MaximumPrincipleBeta}
Let $(\V, \Phi)$ be a cyclic $\beta$-bundle as in
\eqref{BetaCyclicG2'Higgs}. Then we have the following:
\begin{enumerate}[noitemsep]
    \item $\lVert \alpha\rVert_h \leq  \sqrt{2} \lVert 1 \rVert_h$. Moreover, $\lVert\alpha(p_0) \rVert=\lVert 1(p_0) \rVert$ at some point $p_0$ if and only if $||\alpha||=||1||$ globally. In particular, equality at one point implies $\alpha$ is non-vanishing. 
    \item  $\lVert \delta\rVert_h \leq  \lVert 1 \rVert_h$. Moreover, $||\delta(p_0)||= ||1(p_0)||$ at some point $p_0$ if and only if $||\delta ||=||1||$ globally. In particular, equality at one point implies $\delta$ is non-vanishing.  
    \item When $\delta=0$, then $\lVert \alpha\rVert_h <  \frac{\sqrt{5}}{\sqrt{3}} \lVert 1 \rVert_h$.
\end{enumerate} 
\end{lemma}

\begin{proof}
First, we apply Lemma \ref{Lem:HitchinsEquationsGeneral},
Hitchin's equations for cyclic $\Gtwosplit$-Higgs bundles in general to the case of $\beta$-cyclic bundles. We shall write $\beta =1$. Now, away from the zeros of $\alpha$ and $\delta$, we obtain: 
\begin{align}\label{BetaCyclicHitchinEquations}
    \Delta_g\log \left( \frac{||\alpha||^2}{\lVert 1\rVert^2} \right) &= 3\lVert\alpha\rVert^2-5\lVert1\rVert^2- \lVert \delta \rVert ^2,\\
     \Delta_g\log \left( \frac{||\delta||^2}{\lVert 1\rVert^2} \right) &= 2\lVert\delta\rVert^2-2\lVert1\rVert^2.
\end{align}
We use these equations to prove (1)--(3). 

(1) Consider a point $x_0 \in S$ where the ratio $\frac{\lVert \delta \rVert}{\lVert \beta \rVert}$ achieves its maximum. Such a point exists since the tautological section $\beta$ does not vanish. At a local maximum, the Laplacian is non-positive. 
Thus, 
\[0 \geq \Delta_g \log\left(\frac{\lVert \delta \rVert}{\lVert \beta \lVert}\right)(x_0) =2\lVert \delta \rVert^2(x_0)-2\lVert \beta \rVert^2(x_0). \]
In particular $\left(\frac{\lVert \delta \rVert}{\lVert \beta \rVert}\right) (x_0)\leq 1$. As $x_0$ is a maximum, this inequality holds globally. 

(2) Consider a point $x_1 \in S$ where the ratio $\frac{\lVert \alpha \rVert}{\lVert \beta \rvert}$ achieves its maximum. 
Using the previous bound on $\lVert \delta\rVert$ from (1), along with similar reasoning as in (1), we obtain 
\[0 \geq \Delta_g \log\left(\frac{\lVert \alpha \rVert}{\lVert \beta \lVert}\right)(x_1) =3\lVert \alpha \rVert^2(x_1)-5\lVert \beta \rVert^2(x_1)-\lVert \delta \rVert^2(x_1)\geq 3\lVert \alpha \rVert^2(x_1)-6\lVert \beta \rVert^2(x_1). \] 
Hence, $\left(\frac{\lVert \alpha \rVert}{\lVert \beta \rVert}\right) (x_1)\leq\sqrt{2}$. As $x_1$ is a maximum, this inequality holds globally. Setting $\delta=0$ yields (3) immediately. 

The strong maximum principle (cf. \cite{Jos13}) implies that if $\left(\frac{\lVert \delta \rVert}{\lVert \beta \rVert}\right)(x_0)=1$, then $\left(\frac{\lVert \delta \rVert}{\lVert \beta \rvert}\right)$ is constant, and in particular $\delta$ has no zeroes. Similarly, if $\left(\frac{\lVert \delta \rVert}{\lVert \beta \rVert}\right)(x_1)=\sqrt{2}$, then $\left(\frac{\lVert \alpha \rVert}{\lVert \beta \rvert}\right)$ is constant, and in particular $\alpha$ has no zeroes. 
\end{proof}

\begin{remark}
    The case of equality for both $\alpha, \beta$ is attained for cyclic Higgs bundles over the complex plane $\C$ for which $\alpha,\beta,\delta=1$. Such Higgs bundles are considered in \cite{Eva24a} in relation to polynomial almost-complex curves $\C \rightarrow \mathbb{S}^{2,4}$. The associated minimal surface $\C \rightarrow \X_{\Gtwosplit}$ is a flat. 
\end{remark}

\subsubsection{Bases of pencils in \texorpdfstring{$\Ein^{2,3}$}{Ein23}}  \label{Sec:BasesPencilEin23}

Fix a polystable $\beta$-bundle $(\V,\Phi)$ on a Riemann surface $\Sigma$ with corresponding representation $\rho:\pi_1S\rightarrow \Gtwosplit$ via $\NAH_{\Sigma}$. We will construct the desired $(\Gtwosplit, \Ein^{2,3})$-manifolds with bases of pencils.
Here is the broad idea. Associated to the $\beta$-bundle is a pair of $\rho$-equivariant-objects: a conformal harmonic map $f: \tilde{\Sigma} \rightarrow \X$ and an immersed alternating $J$-holomorphic curve $\nu: \tilde{\Sigma} \rightarrow \quadric$, whose `spacelike Gauss map' is $f$. The fact that $\nu$ is immersed entails $\mathrm{I}$ and $\mathrm{III}$ are both pointwise non-vanishing, which allows us to define a distinguished pencil $\mathcal{P}$ of tangent vectors along $f$. That is, $\mathcal{P} \in \Omega^0(\tilde{\Sigma}, \Gr_{2}(f^*\T\X))$. The $\pi_1S$-cover of our 5-manifold is a fiber bundle over $\tilde{\Sigma}$, with fiber at $p$ given by the $\beta$-base $\mathcal{B}_{\beta}(\mathcal{P}|_p)$. \medskip 

We now describe the construction in more detail. To start, the distribution of planes $\mathcal{P}$ can be parametrized by a $\rho$-equivariant 1-form $\Psi_0 \in \Omega^1(\tilde{\Sigma},f^*\T\X)$ whose image is $\mathcal{P}$. 
Recall $\nu$ has Frenet frame $\mathscr{F}_{\nu} =(\mathscr{L},T,N,B)$. Using this splitting, the flat bundle $(\underline{\R^{3,4}}, D)$, where $D$ is the trivial connection on $\underline{\R}^{3,4} = \tilde{\Sigma} \times \R^{3,4}$, decomposes under the Frenet frame splitting as follows:  
\[ D = \begin{pmatrix} \nabla_{\mathscr{L}} & -\mathrm{I}^* & & \\
\mathrm{I} & \nabla_T & -\sff^* & \\
 &  \sff & \nabla_{N}& -\tff^* \\
 &  &  \tff&\nabla_{B} \end{pmatrix}. \]
Here, the adjoints are with respect to $q = q_{3,4}$. The objects $\sff \in \Omega^1(\tilde{\Sigma}, \Hom(T,N))$ and $\tff \in \Omega^1(\tilde{\Sigma}, \Hom(N,B))$ are the second and third fundamental forms of $\nu$ and $\mathrm{I}$ corresponds to $d\nu$. 

We shall want a variant of $\tff$, defined as follows. Using the cross product with $\nu$, each of $N$ and $B$ are complex lines. In particular, for we can decompose the $\R$-linear map $\tff_p(X): N_p\rightarrow B_p$ into two parts: the $\C$-linear part, denoted $\tff^{1,0}(X)$ and the $\C$-anti-linear part $ \tff(X)^{0,1}$.
Identifying $f(p) = \mathscr{L}_p \oplus N_p$ and $\mathrm{T}_U\X \cong \Hom^\times(U,U^\bot)$, it turns out that $(\mathrm{I}+\tff^{0,1}) \in \Omega^1(\tilde{\Sigma}, \End(\underline{\R^{3,4}}))$ becomes an object of the form $(\mathrm{I}+\tff^{0,1} )\in \Omega(\tilde{\Sigma}, f^*\T\X)$. We shall call this object by 
\[ \Psi_0 := \mathrm{I}+\tff^{0,1} \]
going forward. Now, the 5-manifold manifold $B_{\Psi_0} \rightarrow \Sigma$ of interest has $\pi_1S$-cover $\overline{B}_{\Psi_0} \rightarrow \tilde{\Sigma}$, where
\begin{align}\label{BetaBaseAlong_f}
    \overline{B}_{\Psi_0} = \{ (p, \ell) \in \tilde{\Sigma} \times \Ein^{2,3} \mid \ell \in \mathcal{B}_{\beta}(\Psi_0|_p) \}. 
\end{align}
Here, treating $\Psi_0|_p$ as a pencil in $\T_{f(p)}\X$ we may define the base $\mathcal{B}_{\beta}(\Psi_0|_p)$ as in Definition \ref{Defn:TauBasePencil}.

In the remainder of this subsection, we describe the fibers of $\overline{B}_{\Psi_0}$ geometrically in three different ways. In $\S$\ref{Subsec:DevEin23}, we study the tautological map $\dev: \overline{B}_{\Psi_0} \rightarrow \Ein^{2,3}$ by $(p,\ell)\mapsto \ell$, and verify it is indeed a local diffeomorphism.  \medskip 

Next, we translate the harmonic map picture described above back into the Higgs bundle language. The bundle $(\underline{\R^{3,4}}, D) \rightarrow \tilde{\Sigma}$ is isomorphic to the flat bundle $(\pi^*\V^\R, \pi^*\nabla)$, where $\pi: \tilde{\Sigma} \rightarrow \Sigma$ is the universal covering map. 
In particular, the object $\Psi_0 \in \Omega^1(\tilde{\Sigma}, \End(\underline{\R^{3,4}}))$ defined above descends to obtain the form $\Psi_0 \in \Omega^1(\Sigma, \End(\V^\R))$. In fact, the object $\Psi_0$ obtains a simpler description in the Higgs bundle: we have $\Psi_0 = \Phi_0+\Phi_0^*$, where $\Phi_0$ is the \emph{$\beta$-Fuchsian} part of the Higgs field, given by the following diagram:
\begin{equation}\label{Phi0}
\begin{tikzcd}
\mathcal{B} \arrow[r, "1"] & \mathcal{B}\mathcal{K}^{-1} &\mathcal{K}\arrow[r, "-i\sqrt{2}"]&\mathcal{O} \arrow[r, "-i\sqrt{2}"]&\mathcal{K}^{-1}  &\mathcal{B}^{-1}\mathcal{K} \arrow[r, "1"]&\mathcal{B}^{-1} 
\end{tikzcd}.
\end{equation}
To see $\Psi_0=\Phi_0+\Phi_0^*$, we observe that $\tff^{0,1}= 1+1^*$ and $\tff^{1,0}  =\delta + \delta^*$ due to the fact that $\mathcal{C}_{\nu}= \diag(i,-i,-i,0,i,i,-i)$ in the line bundle splitting $\bigoplus_{i=3}^{-3}\mathcal{L}_i$ by Proposition \ref{Prop:HiggsBundleCrossProduct}.
In other words, writing $\Phi = \Phi_{-\beta} + \Phi_{-\alpha} +\Phi_{\delta}$ as the decompositions into root spaces, we have $\Phi_0 = \Phi_{-\beta}$. 
If we write $\Psi_0$ as a matrix in the Higgs bundle Frenet splitting $\V^\R = \mathscr{L}\oplus T\oplus N \oplus B$, rather than the line bundle splitting $\bigoplus_{i=3}^{-3} \mathcal{L}_i$, then $\Psi_0$ obtains the form:
\[ \Psi_0 = \begin{pmatrix}  & -\mathrm{I}^* & & \\
\mathrm{I} &  &  & \\
 & & & -(\tff^{0,1})^* \\
 &  &  \tff^{0,1}& \end{pmatrix}.\]
Recall that for $U \in \X_{\Gtwosplit}$ fixed, there is an associated Cartan decomposition $\g_2' = \frakk(U)\oplus \p(U)$. Hence, using the identification $\T_U\X_{\Gtwosplit}=\Hom^\times(U,U^\bot) \cong \mathfrak{p}(U)$ by $A \mapsto (A-A^{*q})$, we see that 
$\Psi_0$ identifies with $\mathrm{I}+\tff^{0,1}$ as defined earlier.

For our first description of the base $\mathcal{B}_{\beta}(\Psi_0) \subset \Ein^{2,3}$, we use the model from Proposition \ref{Prop:Ein23EasyModel}. 

\begin{proposition}[$\Psi_0$-base via $\mathcal{R}$] 
\label{Prop:Ein23Psi0BaseViaR}
Fix $p \in \tilde{\Sigma}$ and set $U = f(p)$.
Define the trivial rank two vector bundle $\mathcal{R}_{\Psi_0} \rightarrow Q_+(U)$, a sub-bundle of the trivial bundle $\underline{U^\bot}\rightarrow Q_+(U)$, by 
\[ \mathcal{R}_{\Psi_0}|_u=\{\Psi_0(X)(u) \mid X \in \mathrm{T}_p\tilde{\Sigma} \}\]
Then writing $\underline{U^\bot} = \mathcal{R}_{\Psi_0}\oplus \mathcal{R}^\bot_{\Psi_0}$, the base $\mathcal{B}_{\beta}(\Psi_0|_p) \subset \Ein^{2,3}$ of the pencil $\Psi_0$ is given by 
\[ \mathcal{B}_{\beta}(\Psi_0|_p) =\{ [u+z] \in \Ein^{2,3}\mid u \in Q_+(U),  z \in Q_-(\mathcal{R}_{\Psi_0}^\bot|_u)\}, \]
In particular, $\mathcal{B}_{\beta}(\Psi_0|_p)\rightarrow \mathbb{P}(U)$ defines a circle bundle. 
\end{proposition}

\begin{proof}
While \cite{DE26} proves the result more generally, we recall the idea for completeness. 

Let $(\psi_1,\psi_2)$ be a basis for the plane $\Psi_0|_p$. We regard tangent vectors in $\T_U\X_{\Gtwosplit}$ as linear maps $U\to U^\perp$. Define $s_i:Q_+(U)\rightarrow \mathcal{R}_{\Psi_0}$ by $s_i(u) = \psi_i(u)$, for $i \in \{1,2\}$. These sections produce a non-vanishing frame $(s_1,s_2)$ of $\mathcal{R}_{\Psi_0}$. 
Indeed $\Psi_0$ decomposes into two parts, denoted previously $\mathrm{I}$ and $\tff^{0,1}$, which satisfy  $\image(\mathrm{I}_p)=T_p$ and $\image(\tff^{0,1}_p)=B_p$. Since $u$ has nonzero projection on at least one of $\mathscr{L}_p$ or $N_p$, it follows that $(s_1,s_2)$ is a non-vanishing frame. 
This shows $\mathcal{R}_{\Psi_0} \cong \underline{\R}^2$. This property is also a more general consequence of the fact that the pencil ${\Psi_0}_{|p}$ is ${\beta}$-regular, see \cite[Section 3.1]{DE26}. 

Next, take $\ell \in \Ein^{2,3}$ and write $\ell=\pm(u,z)$ as in Proposition \ref{Prop:Ein23EasyModel}, using the model $F_U$ built from $U=f(p)$. By Proposition \ref{Prop:PointingTowardsEinsteinPrelim}, the unique tangent vector $v_{U,\ell} \in \T^1_U\X$ pointing towards $\Ein^{2,3}$, up to positive scalars, is the unique rank one map $v_{U,\ell}:U\rightarrow U^\bot$ satisfying $v_{U,\ell}(u) =z$. Then $\ell \in \mathcal{B}_{\beta}(\Ein^{2,3})$ by definition exactly when $v_{U,\ell} \, \bot_{\T_U\X} \,\Psi_0|_p$. Take $\psi \in \Psi_0|_p$ and up to some nonzero scalar $c \in \R^*$, 
\[ c\langle v_{P,\ell}, \psi \rangle_{\X}= \langle z, \psi(u) \rangle_q\]
by Corollary \ref{Cor:Orthogonality}. It follows immediately that $[u+z] \in \mathcal{B}_{\beta}(\Ein^{2,3})$ if and only if 
$z \in \mathcal{R}_{\Psi_0}^\bot|_u$. Moreover, 
$\mathcal{B}_{\beta}(\Psi_0)$ is the total space of the circle bundle $\mathcal{B}_{\beta}(\Ein^{2,3})\rightarrow \mathbb{P}(U)$ with $\mathbb{S}^1$-fiber $Q_-(\mathcal{R}_{\Psi_0}^\bot|_u)$. 
\end{proof}
Recall from Proposition \ref{Prop:Ein23Principal}
that the projection $\Ein^{2,3} \rightarrow \mathbb{P}(U)$ defines a principal $\mathbb{S}^3$-bundle. Proposition \ref{Prop:Ein23Psi0BaseViaR} shows that the $\Psi_0$-base defines a slice of this fiber bundle, which turns out to be an $\mathbb{S}^1$-sub-fibration.

Next, we describe a Higgs bundle analogue of the previous description of the $\Psi_0$-base.  
The construction of $\mathcal{R} \rightarrow Q_+(U)$, in Proposition \ref{Prop:Ein23Psi0BaseViaR}, writing $\mathcal{R}=\mathcal{R}_{\Psi_0}$, occurred relative to a single point $p \in \tilde{\Sigma}$. 
The Higgs bundle analogue $\mathscr{R}$ of the construction yields a rank two vector bundle $\mathscr{R} \rightarrow Q_+(\mathscr{U})$, where $\mathscr{U} \subset \V^\R$ is the rank three real sub-bundle corresponding to $\mathscr{L} \oplus N$. At $(p,u) \in Q_+(\mathscr{U})$, the bundle $\mathscr{R}$ has fiber  
\[ \mathscr{R}|_{(p,u)}= \{\Psi_0(X)(u) \in \mathscr{U}|_p^\bot \mid X \in \mathrm{T}_p\Sigma\}.\]
Local triviality of $\mathscr{R}$ is checked in a local coordinate. We then have the option to view the base $\mathcal{B}_{\beta}(\Psi_0)$ in the associated bundle $\tilde{\Sigma} \times_{\rho} \Ein^{2,3} \cong \V[\Ein^{2,3}] \subset \mathbb{P}(\V^\R)$. Following this through leads to a global identification of $\mathcal{B}_{\beta}(\Psi_0)\rightarrow Q_+(\mathscr{U})$ of the base of pencil as a circle bundle over $Q_+(\mathscr{U})$, as a bundle version of Proposition \ref{Prop:Ein23Psi0BaseViaR}. We will use both perspectives on $\mathcal{B}_{\beta}(\Psi_0)$ in Subsection \ref{Subsec:DevEin23}. In particular, Proposition \ref{Prop:Ein23BaseViaH} below gives an equivalent description of the base of pencil in terms of the Higgs bundle. 

For later, it will be useful to consider pencils other than $\Psi_0$. This causes no changes to the proof below. Here, as usual, $f: \tilde{\Sigma} \rightarrow \X$ is the associated equivariant harmonic map to $\rho$ via $\NAH_{\Sigma}$. 

\begin{proposition}[Base of Pencil via Harmonic Metric]\label{Prop:HiggsBundleBetaBase}
Let $(\V,\Phi)$ be a $\beta$-bundle associated to a representation $\rho$. For any pencil $\mathcal{P}$ along $f^*\T\X$, the base of pencil $\mathcal{B}_{\beta}(\mathcal{{P}}) \subset \mathbb{P}(\V^\R)$ is given by 
 \begin{align}\label{Prop:Ein23BaseViaH}
     \mathcal{B}_{\beta}(\mathcal{P})|_p=\{ [Z] \in \mathbb{P}(\V^\R)|_p \mid h(Z, \psi(X)(Z))=0, \forall \psi \in \mathcal{P}_p, \, \forall X \in \T_p\Sigma\}.
\end{align}
\end{proposition}

As above, the pencil $\mathcal{P} \in \Omega^1_{\rho}(\tilde{\Sigma},f^*\T\X)$ can be viewed in the form $\mathcal{P} \in \Omega^1(\Sigma, \End(\V^\R))$.
\begin{proof}

To see \eqref{Prop:Ein23BaseViaH} holds, write $Z = u+z$ for $u \in \mathscr{U}$, $z  \in \mathscr{U}^\bot$. Recall that $h|_{\V^\R} = q|_{\mathscr{U}}\oplus (-q|_{\mathscr{U}^\bot})$ and that any endomorphism $\psi \in \mathcal{P}|_p$ is $h$-self-adjoint and exchanges $\mathscr{U}$ and $\mathscr{U}^\bot$. Hence, 
\[ h(Z, \psi(Z)) = h(u,\psi z) +h(z, \psi u)=2h(z,\psi u)= -2q(z,\psi u).\]
The claim follows by the proof of Proposition \ref{Prop:Ein23Psi0BaseViaR}. 
\end{proof}

We now parametrize the $\Psi_0$-base using the Frenet frame splitting $\V^\R =\mathscr{L} \oplus T\oplus N \oplus B$. This parametrization leads to a local coordinate description of the base $\mathcal{B}_{\beta}(\Psi_0)$ in the Higgs bundle, which will be necessary in the proof of the developing map.

\begin{lemma}[$\Psi_0$-base via Frenet Frame]\label{Lem:Psi0BaseFrenetFrame}
Any element $Z \in \mathcal{B}_{\beta}(\Psi_0)$ obtains the form 
\begin{align}\label{FrenetBaseOfPencil}
     \left[ u + \frac{2 \, u_{\mathscr{L}} \times v - u_N \times v}{\sqrt{4 \lVert u_{\mathscr{L}}\rVert_h^2+\lVert u_N\rVert_h^2}} \right]. 
\end{align}
for some $u \in Q_+(\mathscr{U})$ and $ v \in Q_-(B)$. 

If we consider an $h$-unitary complex cross product basis $(e_k)_{k=3}^{-3}$ for $\V$ with $e_k\in \mathcal{L}_k$ and with multiplication table as in Table \ref{Table:ComplexCrossProductBasis}, then every such element $Z$ can be written in the following form for some $z,z_2\in \C$, $x_0\in \R$ such that $x_0^2+2\abs{z_2}^2=1$ and $\abs{z}=\frac{1}{2}$:

\[ Z=\begin{bmatrix}
    \sqrt{2}ix_0z \\
    \lambda z_2\\
    \overline{z_2}z\\
    \lambda x_0\\
    z_2\overline{z}\\
     \lambda \overline{z_2}\\
 -\sqrt{2}ix_0\overline{z}
\end{bmatrix}. \]
Here, $\lambda := \sqrt{{2x_0^2+\lvert z_2\rvert^2}}$. 
 \end{lemma}

\begin{proof}
By Remark \ref{Remk:UnitaryCrossProd}, let us take a local $h$-unitary complex cross product frame $X=(e_k)_{k=3}^{-3}$ for $\V$ such that the cross product obtain the form in Table  \ref{Table:ComplexCrossProductBasis}. In such a basis, elements of the $\Psi_0$-pencil obtains the form $\Psi_0(z)$, for some $z \in \C$, where 
\[ \frac{1}{||\beta||}\Psi_0(z)= \begin{pmatrix}
    0& \overline{z} &0&0&0&0&0\\
    z& 0 &0&0&0&0&0\\
    0& 0 &0&i\sqrt{2}\,\zbar&0&0&0\\
    0& 0 &-i\sqrt{2}\,z&0&i\sqrt{2}\,\zbar&0&0\\
    0& 0 &0&-i\sqrt{2} \,z&0&0&0\\
    0& 0 &0&0&0&0&\overline{z}\\
    0& 0 &0&0&0&z&0
       \end{pmatrix}. \] 
Fix $p \in \Sigma$. In the frame $X$, real elements $x \in \V^\R|_p$ have coordinatization $x = \sum_{i=3}^{-3} x_ie_{i}$ with $x_{-i} = \overline{x_i}$. 
Now, any element $u \in \mathscr{U}|_p=\mathscr{L}_p\oplus N_p$ obtains the form $u=(z_2,x_0,\overline{z}_2)$ in the basis $(e_2, e_0,e_{-2})$. Thus, 
\[ \Psi_0(z)(u) = (\overline{z} z_2, i\sqrt{2}\,\overline{z}\,x_0, -i\sqrt{2} \,zx_0, \,z\overline{z}_2 )\in U^\perp=\left(\mathcal{L}_3\oplus \mathcal{L}_1\oplus \mathcal{L}_{-1}\oplus \mathcal{L}_{-3}\right)_{|\mathbb{R}}.\]
Now, take an element $v \in B_p$ and write $v= (w,0,0,\overline{w})$ in the basis $(e_3, e_1,e_{-1}, e_{-3})$. On the other hand, for $c \in \R$, we have 
\[ c \,u_{\mathscr{L}}\times v -u_N \times v = (c\,ix_0 w, \, \sqrt{2} \overline{z}_2w, \sqrt{2}z_2\overline{w}, \,-c ix_0\overline{w})\]
We then compute that 
\[ h( \Psi_0(z)(u), cu_{\mathscr{L}}\times v-u_N\times v) = (2c-4)x_0\mathrm{Re}(iwz \overline{z}_2).\] 
Hence, for $c =2$, we find $c u_{\mathscr{L}} \times v-u_N \times v \in \mathcal{R}_{\Psi_0}^\bot|_u$. By dimension count, we conclude that every element in $\mathcal{R}_{\Psi_0}^\bot|_u$ obtains this form.

Now, writing $Z \in \mathcal{B}_{\beta}(\Psi_0)$ as $Z = [u+y]$ for $y \in \mathcal{R}_{\Psi_0}^\bot|_u$ and normalizing such that $||u||_h = 1 = ||y||_h$, the equation \eqref{FrenetBaseOfPencil} follows. Similarly, the coordinate expression for $Z$ is obtained by choosing $\lambda$ such that $||\lambda u||_h =||y||_h$, where 
$y = (\sqrt{2}ix_0z,\, \overline{z}_2z,\,z_2\overline{z},\,-\sqrt{2}ix_0\overline{z} ) \in \mathcal{R}_{\Psi_0}^\bot|_u$, and we normalize $x_0, z_2, z$ as in the hypotheses. 
\end{proof}

\subsubsection{Construction of the geometric structure}\label{Subsec:DevEin23}

In this section, we prove the $\Psi_0$-base studied in the previous subsection yields a $(\Gtwosplit, \Ein^{2,3})$-structure on a 5-manifold $M \rightarrow S$ with $\Ein^{2,1}$-fibers. 
Here are the details. We have built a manifold 
\[ \overline{B}_{\Psi_0} \subset \tilde{\Sigma} \times \Ein^{2,3} \] 
with fiber over $p$ given by $\mathcal{B}_{\beta}(\Psi_0|_p)$. 
There is a tautological map $\dev: \overline{B}_{\Psi_0} \rightarrow \Ein^{2,3} $ by $(p, \ell)\mapsto \ell$. Our goal is to verify that $\dev$ is a local diffeomorphism. By definition, $\dev$ is an injective immersion on fibers. The essence of the proof is to understand how $\Psi_0$ varies, and hence how these fibers vary, to confirm we have a local diffeomorphism. A certain parallelism property of $\Psi_0$ plays a crucial role to this end. 
Going forward, we will always use $\overline{\, \cdot \,}$ to denote the corresponding object over $\tilde{S}$. For example,  $\overline{B}_{\Psi_0}$ fibers over $\tilde{S}$ and $B_{\Psi_0}$ is the corresponding $\pi_1S$-quotient fibering over the closed surface $S$. 

There are two possible geometric descriptions of the closed 5-manifold that carries the $(\Gtwosplit,\Ein^{2,3})$-structure of interest. 
We have defined already $B_{\Psi_0}$ from the perspective described above that explicitly involves bases of pencils. Later, we shall introduce $M_{\Psi_0}$, another diffeomorphic model space that is simpler to describe. We first pursue $\dev$ from the point of view of $B_{\Psi_0}$ and discuss $M_{\Psi_0}$ afterwards.

For emphasis, we highlight a certain conformal metric $g$ related to the parallelism of $\Psi_0$, that we will refer to as the \emph{projected metric}. Here, we write $h=\diag(h_i)_{i=3}^{-3}$, where $h_i$ is the harmonic metric $h$ restricted to $\mathcal{L}_i$. 

\begin{definition}[$\beta$-Projected metric]\label{Defn:BetaConformalMetric}
Note that $h_2h_3^{-1}$ defines a metric on $\mathcal{B}^{-1} \otimes \mathcal{B}\mathcal{K}^{-1}\simeq \mathcal{K}^{-1}$. Denote $g_\beta$ the corresponding Riemannian metric on $S$.
\end{definition} 
The metric $g_{\beta}$ is a conformal metric on the Riemann surface $\Sigma$. 
Up to a constant multiplicative factor, the metric $g_\beta$ is the (negative) induced metric of the associated $\beta$-curve $\nu:\tilde{S}\rightarrow \quadric$, namely $-\nu^*g_{\quadric}$. 
Let $f: \tilde{\Sigma} \rightarrow \X$ be the associated harmonic map. 
\begin{remark}
\label{rem:Projected metric}
    Another interpretation of the metric $g_{\beta}$ is that it is, up to a multiplicative constant $C$, the metric is induced by the pullback metric $f^*g_{\X}$ once tangent vectors are projected to the planes defined by $\Psi_0$. That is, $g_\beta(v,v)=C\lVert \Psi_0(v)\rVert^2$. In particular, up to the multiplicative constant $C$, it is always smaller than the metric induced by the harmonic map, namely $g(v,v)=\lVert df(v)\rVert^2=\frac{1}{4}\lVert\Psi(v)\rVert^2$. This interpretation justifies the terminology ``projected metric''.
\end{remark}

The following proposition describes a crucial property of $\Psi_0$. Indeed, it shows that the way it varies under the flat connection is determined by $\Psi_0$, and not by the derivative. Moreover, this property will allow us to reduce verification of our developing map in Theorem \ref{Thm:DevEin23} to a $C^0$-condition on the pair $(\Psi,\Psi_0)$. 

\begin{proposition}[$\Psi_0$-parallelism for $\beta$-bundles]\label{Prop:ParallelismPsi0Alpha1}
Let $\nabla^h$ and $\nabla^g$ be the Chern connections on $\End(\mathcal{E})$ and $\T\Sigma\simeq \mathcal{K}^{-1}$, respectively, where $g=g_\beta$. 
For any vector field $X \in \Omega^0(\Sigma, \T \Sigma)$, 
\begin{align}\label{Parallelism_BetaBundles}
    (\nabla^h \circ \Phi_0)(X)=& \,(\Phi_0\circ \nabla^g)(X),\\
    (\nabla^h \circ \Psi_0)(X)=& \,(\Psi_0\circ \nabla^g)(X).
\end{align}
In other words, $\Phi_0 \in \Omega^0(\K\otimes \End(\V))$ and $\Psi_0 \in \Omega^0(\T^*S\otimes \End(\V))$ are parallel. 
\end{proposition}

\begin{proof}
In this proof, for clarity we denote by $\beta$ the tautological section that was denoted by just $1$ before.
The Chern connection $\nabla^h$ preserves the line decomposition $\V =\bigoplus_{k=3}^{-3} \mathcal{L}_k$ by Proposition \ref{Prop:HarmonicMetricDiagonal}. 
 We denote by $h_k$ the associated metric on $\mathcal{L}_k$. The induced connection of $\nabla^h$ on $\End(\mathcal{E})$ satisfies:
 \[ \nabla^h \Phi_0(X)=\begin{pmatrix}
    0& 0 &0&0&0&0&0\\
    \nabla^{h_{2}h_3^{-1}}\beta(X)& 0 &0&0&0&0&0\\
    0& 0 &0&0&0&0&0\\
    0& 0 &-i\sqrt{2}\nabla^{h_{0}h_1^{-1}}\beta(X)&0&0&0&0\\
    0& 0 &0&-i\sqrt{2}\nabla^{h_{-1}h_0^{-1}}\beta(X)&0&0&0\\
    0& 0 &0&0&0&0&0\\
    0& 0 &0&0&0&\nabla^{h_{-3}h_{-2}^{-1}}\beta(X)&0
       \end{pmatrix}. \]
Each of the four tautological sections $\beta$ exhibit the desired parallelism, as we now explain. For example, the first object entitled $\beta$ in the matrix of $\Phi_0$ is of the form $\beta_{3,2} \in \Omega^0(\Hom(\mathcal{L}_3,\mathcal{L}_2)\otimes \K)$. 
Now, consider the vector bundle $V=\Hom(\K^{-1}, \End(\V))$, with its natural induced connection from $\nabla^g$ on $\K^{-1}$ and $\nabla^h$ on $\End(\V)$.
We then compute 
\[ \nabla^h(\beta_{3,2}(X))= (\nabla^V\beta_{3,2})(X)+\beta_{3,2}(\nabla^gX).\]
Here is the key: the induced connection $\nabla^V$ restricts to $\Hom(\K^{-1}, \Hom(\mathcal{L}_3,\mathcal{L}_2))$ as the trivial connection since $(\Hom(\mathcal{L}_3,\mathcal{L}_2), h)\cong (\mathcal{L}_{3}^{-1}\mathcal{L}_2, h_3^{-1}h_2)\cong (\K^{-1},g)$. In particular, $\nabla^V\beta_{3,2}=0$. Nearly identical reasoning on the other tautological sections leads to $\nabla^h\Phi_0(X)=\Phi_0(\nabla^gX)$, using the symmetries of $h$ from Proposition \ref{Prop:HarmonicMetricDiagonal}. 

The same reasoning holds for $\Phi_0^*$ and therefore also for $\Psi_0$.
\end{proof}

We now prove the main result of this section: the tautological map $\dev$ indeed deserves its title. 

\begin{theorem}[$\Ein^{2,3}$-structures for $\beta$-bundles]\label{Thm:DevEin23}
Let $(\V, \Phi)$ be a polystable $\beta$-cyclic Higgs bundle such that $\mathrm{div}(\delta) \neq 0$. Then the tautological map $\dev: \overline{B}_{\Psi_0} \rightarrow \Ein^{2,3}$ is a local diffeomorphism.
\end{theorem}

In particular, by Proposition \ref{Prop:BetaBundleStability}, $\delta$ has a nonzero divisor if $(\V,\Phi)$ is stable. The hypotheses exclude only the case $\deg(\mathcal{B})=0$ when $\delta$ is non-vanishing.

\begin{proof}

Fix an arbitrary basepoint $b_0 \in \tilde{\Sigma}$. We can then diffeomorphically identify the pullback bundle $\pi^*B_{\beta}(\Psi_0)$, where $\pi: \tilde{\Sigma} \rightarrow \Sigma$ is the universal covering map, defined as in Proposition \ref{Prop:HiggsBundleBetaBase}, with
$\overline{B}_{\Psi_0} \subset \tilde{\Sigma} \times \Ein^{2,3}$ as in \eqref{BetaBaseAlong_f}
by the following map: $(p,\ell)\longmapsto (p, \mathscr{P}_{p,b_0}^{\nabla}(\ell))$. Here, $\mathscr{P}_{p,b_0}: \mathcal{E}^\R_p\rightarrow \mathcal{E}^\R_{b_0}$ denotes $\nabla$-parallel transport, which is path-independent since $\nabla$ is flat, and we make a fixed identification  $\Ein(\mathcal{E}^\R|_{b_0})\cong \Ein^{2,3}$. 
With this perspective,  
$\dev: \pi^*B_{\beta}(\Psi_0)\rightarrow \Ein^{2,3}$ is given by 
\[ \dev(p,\ell) = \mathscr{P}_{p, b_0}^{\nabla}(\ell). \]
We now show this map $\dev$ is a local diffeomorphism. Let us also write $\overline{B} = \pi^*B_{\beta}(\Psi_0)$ for simplicity.\medskip 

\textbf{Step 0: Setup.} Let $p_0\in \widetilde{\Sigma}$ be arbitrary. Fix any point $\ell_0 \in \overline{B}_{\mid p_0}$ in the fiber above $p_0$. To prove $\dev$ is an immersion, it suffices to show the differential $d\dev_{\ell_0}$ surjects. By construction, $\dev$ is an immersion when restricted to the fiber $\overline{B}|_{p_0}=:\mathcal{F}_{p_0}$.  Thus, we need only prove $d\dev_{\ell_0}$ surjects onto $\T\Ein^{2,3}/\T(\dev (\mathcal{F}_{p_0}))$. The following procedure defined next will allow us to prove this. 

Recall from Proposition \ref{Prop:HiggsBundleBetaBase} that $[Z] \in \Ein^{2,3} $ satisfies $[Z] \in \mathcal{F}_{p_0}$ if and only if for all $v\in \T_{p_0}\widetilde{\Sigma}$:
\[h(\Psi_0(v)Z,Z)=0.\] 
 
Fix any $v_0 \in \T^1_{p_0}\tilde{\Sigma}$ be any unit tangent vector. Let $\left(\gamma_t \right)_{t\in (-\varepsilon,\varepsilon)}$ be the geodesic for the metric $g$ from Definition \ref{Defn:BetaConformalMetric} on $\widetilde{\Sigma}$ starting at $x_0$ with $\dot{\gamma}(0) = v_0$.
Let $\left(Z_t\right)_{t\in (-\varepsilon,\varepsilon)}$ be the section of $\pi^*\V^\R$ above $\gamma$ such that $[Z_0]=\ell_0$ and $Z_t$ is parallel for the Chern connection $\gamma^*\nabla^h$. We claim that $[Z_t]$ is a section of $\gamma^*\overline{B}$. That is, we claim that 
\begin{align}\label{BaseOfPencilClaim}
    h(\Psi_0(Y_t)Z_t,Z_t) =0,\;\; \text{for all} \; Y_t \in \T_{\gamma(t)}\tilde{\Sigma}. 
\end{align} To prove the claim, suppose that $\mathscr{Y}$ is any $\nabla^g$-parallel vector field along $\gamma$ and define 
$f_{\mathscr{Y}}:(-\varepsilon, \varepsilon) \rightarrow \R$ by $f(t) = h(\Psi_0(\mathscr{Y}_t)Z_t, Z_t)$. The claim \eqref{BaseOfPencilClaim} holds if $f_{\mathscr{Y}}\equiv 0$. 
Now, by hypothesis, $f_{\mathscr{Y}}(0) =0$. We then compute that $f_{\mathscr{Y}} '\equiv 0$ using that $\mathscr{Y}$ is $\nabla^g$-parallel, $\nabla^hh=0$, along with the parallelism property of $\Psi_0$ from Proposition \ref{Prop:ParallelismPsi0Alpha1}. Hence, $f_{\mathscr{Y}}\equiv0$, which verifies the claim that $[Z_t]$ is a section of $\gamma^*\overline{B}$. 

Observe that for any point $(p_0,\ell_0)=:y \in \overline{B}_{\mid p_0}$, the above process defines a splitting $\T_y \overline{B}_{\Psi_0}=H_y \oplus V_y$, where $V_y$ is the vertical subspace, i.e. the tangent space of the fiber, and $H_y$ is the space of tangent vectors to horizontal lifts just defined.

Let $(Z_t)_{t \in (-\varepsilon, \varepsilon)} \subset \overline{B}$ be the same curve as above. 
The essential task of the proof will be show the following:
\begin{align}\label{KeyDerivativeEin}
    A=\frac{d}{dt}\bigg|_{t=0} \langle \Psi_0(v)\mathscr{P}^{\nabla}_{t,0}(Z_t),\mathscr{P}^{\nabla}_{t,0}(Z_t) \rangle_{h} > 0.
\end{align} 
Here, $\mathscr{P}^{\nabla}_{t,0}: \V^\R _{\gamma(t)}\rightarrow \mathcal{E}^{\R}_{\gamma(0)}$ denotes $\nabla$-parallel transport. We now explain why this inequality completes the proof. 

By parallel transporting $Z_t$ from $\V^\R|_{\gamma(t)}$ to $\V^\R|_{p_0}$, we can compare the developed fibers $\dev(\overline{B}_{\mid p_0})$ and $\dev(\overline{B}_{\mid \gamma(t)})$. 
The inequality \eqref{KeyDerivativeEin} implies that $\dev([Z_t])$ is moving away from $\dev(\overline{B}_{\mid p_0})$, evidenced by the fact that $\dev([Z_t])$ does not satisfy the equation $\langle \Psi_0(v)(\cdot), (\cdot) \rangle =0$ even infinitesimally. Now, the aforementioned process defines for each $v_0 \in \T_{p_0}\tilde{\Sigma}$ and 
$\ell_0 \in \overline{B}_{|{p_0}}$
in the fiber above $p_0$ a tangent vector $V_0 \in H_{\ell_0}$ that projects to $v_0$. The key inequality \eqref{KeyDerivativeEin} implies $d\dev_{\ell_0}(V_0)$ is a non-zero element of $\T\Ein^{2,3}/\T(\dev(\mathcal{F}_{p_0}))$ and hence $d \dev_{\ell_0}:H_{\ell_0} \rightarrow\T\Ein^{2,3}/\T(\dev(\mathcal{F}_{p_0}))$ is a linear isomorphism, meaning $d\dev_{\ell_0}$ surjects. 
The remainder of the proof is to prove  \eqref{KeyDerivativeEin}. \medskip

\textbf{Step 1: Differentiate.}
Let us rewrite \eqref{KeyDerivativeEin}, noting that $\Psi_0(v)$ and $h=h_{p_0}$ are fixed, and $\Psi_0$ is $h$-self-adjoint. We write $D^C_t$ for the covariant differentiation operator along $\gamma$ of a connection $C$.
\[A=2\langle\Psi_0(v_0) Z_0 ,\frac{d}{dt}\big|_{t=0}\mathscr{P}^{\nabla}_{t,0} (Z_t)\rangle_h=2\langle\Psi_0(v_0) Z_0 , D^{\nabla}_{t}(Z_t)(0)\rangle_h= 2\langle\Psi_0(v_0) Z_0 , \Psi(v_0)Z_0\rangle_h,\]
where in the final equality we use that $D^{\nabla}_t=D^{\nabla^h}_t+\Psi(\gamma'(t))$ and that $Z_t$ is $\nabla^h$-parallel. Since $\Psi(v_0)$ is also $h$-self-adjoint, we can write $A$ once more as follows
\[A=2\langle\Psi_0(v_0)\Psi(v_0) Z_0 ,Z_0\rangle_h=\langle\left(\Psi(v_0)\Psi_0(v_0)+\Psi_0(v_0)\Psi(v_0)\right) Z_0 ,Z_0\rangle_h.\]

\vspace{2ex}
\textbf{Step 2: Express $A$ in coordinates.}

Recall we have $p_0\in \tilde{\Sigma}$ fixed. By Remark \ref{Remk:UnitaryCrossProd}, we now take an $h$-unitary complex cross product basis $(e_k)_{k=3}^{-3}$ of $\pi^*\mathcal{E}|_{p}$ satisfying multiplication table \ref{Table:ComplexCrossProductBasis} and such that $e_i\in \mathcal{L}_i$. 
By Lemma \ref{Lem:Psi0BaseFrenetFrame}, an arbitrary element $Z_0$ of the fiber $\overline{B}_{\Psi_0}|_{p_0}$ can be written in the basis $(e_k)_{k=3}^{-3}$ as follows: 

\begin{align}\label{GeneralFiberElementEin23}
  Z_0 =\begin{bmatrix}
    \sqrt{2}ix_0z \\
    \lambda z_2\\
    \overline{z_2}z\\
    \lambda x_0\\
    z_2\overline{z}\\
     \lambda \overline{z_2}\\
 -\sqrt{2}ix_0\overline{z}
\end{bmatrix} .
\end{align}
Here, $z, z_2 \in \C$, $x_0 \in \R$, and we have normalized such that $\lambda := \sqrt{{2x_0^2+\lvert z_2\rvert^2}}$, $2\abs{z_2}^2+x_0^2=1$, and $\abs{z}^2= \frac{1}{2}$.

Let us denote $\psi=\Psi(v_0), \, \psi_0=\Psi_0(v_0)$, from which we build the endomorphism $M=\psi\psi_0+\psi_0\psi$. 
Therefore, the quantity that we want to show does not vanish is 
\begin{align*}
A&= h\left( M  Z_0,  Z_0 \right) .
\end{align*}

Up to a diagonal $h$-unitary gauge transformation of the basis $(e_k)_{k=3}^{-3}$, the transformations $\psi$, $\psi_0$, and $M$ become the following matrices in the given basis:

\[ \frac{1}{||\beta||}\psi_0=\begin{pmatrix}
    0& 1 &0&0&0&0&0\\
    1& 0 &0&0&0&0&0\\
    0& 0 &0&i\sqrt{2}&0&0&0\\
    0& 0 &-i\sqrt{2}&0&i\sqrt{2}&0&0\\
    0& 0 &0&-i\sqrt{2}&0&0&0\\
    0& 0 &0&0&0&0&1\\
    0& 0 &0&0&0&1&0
       \end{pmatrix},\;\; \psi= \frac{1}{||\beta||}\begin{pmatrix}
     0& 1 &0&0&0&\color{olive}{\delta_0}&0\\
    1& 0 &\color{blue}{\overline{\alpha_0}}&0&0&0&\color{olive}{\delta_0}\\
    0& \color{blue}{\alpha_0} &0&i\sqrt{2}&0&0&0\\
    0& 0 &-i\sqrt{2}&0&i\sqrt{2}&0&0\\
    0& 0 &0&-i\sqrt{2}&0&\color{blue}{\overline{\alpha_0}}&0\\
    \color{olive}{\overline{\delta_0}}& 0 &0&0&\color{blue}{\alpha_0} &0&1\\
    0& \color{olive}{\overline{\delta_0}} &0&0&0&1&0
       \end{pmatrix} \] 
    
\[ \frac{1}{||\beta||^2}M=\begin{pmatrix}
    2& 0&\color{blue}{\overline{\alpha_0}}&0&0&0&\color{olive}{2\delta_0}\\
    0& 2  &0&\color{blue}{i\sqrt{2}\overline{\alpha_0}}&0&\color{olive}{2\delta_0}&0\\
    \color{blue}{\alpha_0}& 0 &4&0&4&0&0\\
    0&\color{blue}{-i\sqrt{2}\alpha_0} &0&8 &0&\color{blue}{i\sqrt{2}\overline{\alpha_0}}&0\\
    0& 0 &4&0&4&0&\color{blue}{\overline{\alpha_0}}\\
    0& \color{olive}{2\overline{\delta_0}} &0&\color{blue}{-i\sqrt{2}\alpha_0}&0&2 &0\\
    \color{olive}{2\overline{\delta_0}}& 0 &0&0&\color{blue}{\alpha_0}&0&2
       \end{pmatrix} \] 
       
We can now compute $A'=\frac{1}{||\beta||^2}A$ using the formula \eqref{GeneralFiberElementEin23} for $Z_0$ in local coordinates. We obtain:

\begin{align*}
A'=&4x_0^2+4\lambda^2\abs{z_2}^2+4\abs{z_2}^2+8\lambda^2x_0^2-8\text{Re}(z_2^2\overline{z}^2)\color{blue}{-2\sqrt{2}x_0(2\lambda^2-1)\text{Re}(i\alpha_0z_2)}\color{olive}{+\text{Re}\left(4\delta_0(-2x_0^2\overline{z}^2+\lambda^2\overline{z_2}^2)\right)}. \\
 \geq &4x_0^2+4\lambda^2\abs{z_2}^2+8\lambda^2x_0^2\color{blue}{-2\sqrt{2}x_0(2\lambda^2-1)\text{Re}(i\alpha_0z_2)}\color{olive}{+\text{Re}\left(4\delta_0(-2x_0^2\overline{z}^2+\lambda^2\overline{z_2}^2)\right)}.
 \end{align*}

\begin{remark} \label{rem:explainationComputation}
    In this computation we color the terms that contain $\alpha_0$ or $\delta_0$. We know a priori that the sum of all uncolored terms in the expression of $A$ is always strictly positive, as it corresponds to the case where the $\beta$-regular pencils considered are all tangent to a totally geodesic copy of $\mathbb{H}^2$. In this case, the fibers are disjoint as they are fibers of the extension of the nearest point projection from \cite{Dav25}. The key part of the computation is to control the colored terms to verify that the maximum principle on $\lVert \delta\rVert$ and $\lVert \alpha\lVert$ is strong enough.
\end{remark}

\textbf{Step 3: Prove $A>0$.}

Let us write $x=x_0^2$ so $0\leq x\leq 1$. One has $\abs{z_0}^2=\frac{1-x}{2}$, $\lambda^2=\frac{3x+1}{2}$, and $|z_2|^2= \frac{1-x}{2}$. 
By Lemma \ref{Lem:MaximumPrincipleBeta}, $\frac{||\alpha||}{||\beta||}=\abs{\alpha_0}\leq\sqrt{2}$ and $\frac{||\delta||}{||\beta||}=\abs{\delta_0}\leq 1$. Since $\delta$ must have a zero by hypothesis, we have $\abs{\delta_0}<1$ by Lemma \ref{Lem:MaximumPrincipleBeta}. 
Hence, we obtain the following strict inequality from the triangular inequality:
\begin{align*}
    A'&> 4x+(3x+1)(1-x)+4x(3x+1)\color{blue}{-6\sqrt{2}x\sqrt{x(1-x)}}\color{olive}{-4x-(3x+1)(1-x)}.\\
    A'&>x\left(12x+4-6\sqrt{2}\sqrt{x(1-x)}\right).
\end{align*}

Using the arithmetic-geometric mean inequality, for $0\leq x \leq 1$, we have
\[12x+4=12x+\frac{4}{3}+\frac{4}{3}+\frac{4}{3}\geq 4\sqrt[4]{12x\times \left(\frac{4}{3}\right)^3}=\frac{16}{\sqrt{3}}\sqrt{x}\geq 6\sqrt{2}\sqrt{x(1-x)}.\] 
For the last inequality we used that $\frac{16}{\sqrt{3}}>6\sqrt{2}$. This concludes the proof that $A>0$.
\end{proof}

\begin{remark}
It is noteworthy in this argument that the bound  $\lVert\alpha\rVert\leq \sqrt{2}$ from the maximum principle in Lemma \ref{Lem:MaximumPrincipleBeta} is not strictly needed, meaning that a weaker bound would also suffice, whereas the exact bound $\lVert\delta\rVert< 1$ is strictly needed. 
This seems to be related to the fact that $\beta$-cyclic bundles with $||\alpha|| \equiv \sqrt{2}$ (bundles with $\alpha$ is non-vanishing) correspond to Hitchin representations, which are $P_{\beta}$-Anosov. On the other hand, if $||\delta(p)|| =1$ occurs at a single point, then $\delta$ must be non-vanishing, which entails the corresponding representation is \textbf{not} $P_{\beta}$-Anosov as explained in Remark \ref{Remk:SL3Beta}. This explains why one expects $\dev$ to fail to be an immersion in this case. 
\end{remark}

\subsubsection{Geometric Structures via $J$-holomorphic Curves}
We now describe the construction of the developing map from the point of view of $J$-holomorphic curves. This leads to the previously advertised alternate model $M_{\Psi_0}$ for the manifold carrying the $\Ein^{2,3}$-structures. To start, we define the $(\mathbb{S}^2 \times \mathbb{S}^1)$-bundle $M_{\Psi_0} \rightarrow \Sigma$ via the Higgs bundle as a fiber bundle direct sum:
\begin{align}\label{Mpsi0Ein23}
M_{\Psi_0}:=Q_+(\mathscr{L} \oplus N)\oplus Q_-(B). 
\end{align}
Treating the Frenet frame subspaces as equivariant objects over $\tilde{\Sigma}$ gives the $\pi_1S$-cover $\overline{M}_{\Psi_0} \rightarrow \tilde{S}$ as follows. 
Define $\overline{M}_{\Psi_0} \subset \tilde{\Sigma} \times (\underline{\R}^{3,4})^2$ with fiber 
\begin{align}\label{MPsi0Ein23}
    \overline{M}_{\Psi_0}|_p = \{ (u,v) \in \underline{\R^{3,4}})^2 \mid u \in Q_+(\mathscr{L}_p\oplus N_p) , v \in Q_-(B_p)\}.
\end{align}
There is a natural 2-1 covering map $\sigma: \overline{M}_{\Psi_0} \rightarrow \overline{B}_{\Psi_0}$, through which the natural developing map
$\mathcal{D}: \overline{M}_{\Psi_0} \rightarrow \Ein^{2,3}$ factors. Here, $\mathcal{D}$ is given by 
\begin{align}\label{GlobalEin23DevMap}
 \mathcal{D}(p,u, v) &= \left[ u + \frac{2 \, u_\mathscr{L}\times v - u_N \times v}{\sqrt{4 \lVert u_{\mathscr{L}}\rVert_h^2+\lVert u_N\rVert_h^2}} \right]. 
\end{align}
In \eqref{GlobalEin23DevMap}, we write $ u \in Q_+(\mathscr{L}\oplus N)$ as a sum $u = u_{\mathscr{L}} + u_N$ with $u_\mathscr{L}\in \mathscr{L}$ and $u_N \in N$. That the map $\mathcal{D}$ factors through $\overline{B}_{\Psi_0}$ via $\sigma(p,\ell)=(p,\mathcal{D}(\ell))$, is an immediate consequence of Lemma \ref{Lem:Psi0BaseFrenetFrame}. 

By the work \cite{CT24} on the moduli space of  equivariant $J$-holomorphic $\beta$-curves, we note the following reinterpretation of Theorem \ref{Thm:DevEin23}. In particular, we can convert special equivariant harmonic maps into fibered $(G,X)$-structures.   
\begin{corollary}[$\beta$-curves to Geometric Structures]
Let $\nu: \tilde{\Sigma} \rightarrow \quadric$ be an alternating $\rho$-equivariant $J$-holomorphic curve for some representation $\rho:\pi_1S\to \Gtwosplit$ that is immersed and linearly full. Define $\overline{M}_{\Psi_0} \rightarrow \tilde{S}$ as in \eqref{MPsi0Ein23}. Then the map $\mathcal{D}: \overline{M}_{\Psi_0} \rightarrow \Ein^{2,3}$ is a local diffeomorphism that defines a fibered $(\Gtwosplit, \Ein^{2,3})$-structure on $M_{\Psi_0}$ with descended holonomy $\rho$. 
\end{corollary}

\begin{proof}
By Proposition \ref{Prop:BetaBundleStability} and \cite[Theorem A]{CT24}, there is a bijection between pairs $(\nu, \rho)$ of linearly full $\rho$-equivariant $J$-holomorphic curves $\nu: \tilde{\Sigma} \rightarrow \quadric$ and $\beta$-cyclic Higgs bundles that are stable, each up to appropriate equivalence.\footnote{Here, the `linearly full' condition automatically excludes the case of $\nu$ being totally geodesic, which forces $\alpha \neq 0$ by Remark \ref{Remk:TotallyGeodesic}.} The result then follows from Theorem \ref{Thm:DevEin23}. 
\end{proof}

\subsection{The \texorpdfstring{$\Gtwosplit$}{G2'}-Hitchin Case: Comparison with Tits Metric Thickening}\label{Subsec:G2Ein23HitchinCase} 

We now consider the construction of $(\Gtwosplit, \Ein^{2,3})$-geometric structures from the previous subsection in the case of Hitchin representations. We show $\dev$ is a diffeomorphism onto the Tits metric thickening domain \eqref{Omega_Thick}.\medskip 

We begin with a remark. Let $\rho:\pi_1S\to \Gtwosplit$ be $\beta$-Anosov and $\iota:\Gtwosplit\hookrightarrow \SO_0(3,4)$ be the standard inclusion. 
Then $\iota \circ \rho: \pi_1S \rightarrow \SO_0(3,4)$ is $P_3$-Anosov, where $P_3 <\SO_0(3,4)$ is the stabilizer of an isotropic 3-plane. Write the Anosov limit maps as $\xi_\rho:\partial \pi_1S\to \Ein^{2,3}$ 
and $\xi^3_\rho:\partial \pi_1S\to \mathsf{Iso}_3(\R^{3,4})$. 
We have the equality $\xi^3_\rho = \Ann(\xi_\rho)$ as a direct consequence of \cite[Proposition 4.4]{GW12} applied to $\iota$. 

We now define another domain $\Omega_{\rho}$ that is the same as $\Omega^{\mathrm{Thick}}_{\rho}$ in \eqref{Omega_Thick}, but more explicitly involves the Anosov limit maps of $\rho$ in a way that is convenient presently: 
\begin{align}\label{OmegaEin23_GWPerspective}
\Omega_\rho :=\Ein^{2,3}\setminus \bigcup_{x\in \partial \pi_1S} \mathbb{P}(\xi^3_\rho(x)).
\end{align} 
In fact, this definition of $\Omega_{\rho}$ is exactly the one given by Guichard-Wienhard  \cite{GW12}, but is also equivalent to the \cite{KLP18} perspective, namely $\Omega^{\Thick}_{\rho}$. 

We need one small observation here about the complement of $\Omega$, namely the thickening of the limit set. 
\begin{proposition}[Complement of $\Omega_{\rho}$ in $\Ein^{2,3}$]\label{Prop:OmegaThickComplementEin23}
Let $\rho: \pi_1S \rightarrow \Gtwosplit$ be $P_{\beta}$-Anosov and define $K = \Ein^{2,3}\backslash \Omega_{\rho}$ for $\Omega_{\rho}$ from \eqref{OmegaEin23_GWPerspective}. Then $K$ is homeomorphic to $\sphere^1 \times \RP^2$.
\end{proposition}

\begin{proof}
The group $\Gtwosplit$ has exactly four orbits in $\Ein^{2,3}\times \Ein^{2,3}$, described in \cite{BH14}.\footnote{Figure \ref{Fig:Apartment} suggestively encodes these orbits in terms of an $\R$-cross product basis $(x_i)_{i=3}^{-3}$. Representatives for the four orbits are: $([x_3], [x_3]), \,([x_3],[x_2]), ([x_3],[x_{-1}]), ([x_3],[x_{-3}])$, which have $\titsangle$ given by: $0, \frac{\pi}{3},\frac{2\pi}{3},\pi$, respectively.} 
Fix a background $\R$-cross product basis $(x_i)_{i=3}^{-3}$ for $\imoct$. 
In particular, if $\ell, \ell' \in \Ein^{2,3}$ are transverse, then the pair $(\ell, \ell')$ is $\Gtwosplit$-equivalent to $([x_3], [x_{-3}])$. Hence, $\mathbb{P}\Ann(\ell) \cap \mathbb{P}(\Ann(\ell') = \emptyset$ if and only if $\ell, \ell' \in \Ein^{2,3}$ are transverse. 

Now, let $C  \subset \Ein^{2,3}$ be a transverse circle, meaning $\ell \pitchfork \ell'$ for pair $\ell \neq \ell' \in C$. We consider the vector bundle $A \rightarrow C$ given by $A_{\ell} = \Ann(\ell)$. The vector bundle $A$ is trivial, as we now show. Fix any point $U \in \X_{\Gtwosplit}$. Then by Proposition \ref{Prop:AnnihilatorGraph}, the 3-plane $A|_{\ell} = \Ann(\ell)$ is the graph of a linear map $\phi_{x, U}: U \rightarrow U^\bot$. Hence, there is a vector bundle isomorphism from the trivial bundle $\underline{U} = C \times U$ to $A$ by $(\ell, y) \mapsto (\ell, \phi_{\ell,U}(y))$. Thus, $A \rightarrow C$ is trivial. Now, we conclude that the fiber bundle $\RP^2 \rightarrow \mathbb{P}(A) \rightarrow C$ is also trivial. However, since $C$ is a transverse circle, the inclusion map $\iota: \mathbb{P}(A) \rightarrow \Ein^{2,3}$ by $(\ell,[y])\mapsto [y]$ is injective.

In particular, for $C = \im(\xi^{\beta})$ as the image of the limit set in $\Ein^{2,3}$, we conclude that the thickening $K_{\im(\xi^{\beta})}$ is homeomorphic to the total space of $\mathbb{P}(A) \rightarrow C$, which concludes the proof.
\end{proof}

We now compare our construction of $(\Gtwosplit,\Ein^{2,3})$-manifolds with those defined in \cite{GW12, KLP18} via domains of discontinuity. 

\begin{theorem}[$\Gtwosplit$-Fuchsian-Hitchin Case]\label{Thm:BetaTits=Pencil}
Let $\rho: \pi_1S \rightarrow \Gtwosplit$ be a representation associated to a $\Gtwosplit$-Hitchin Hodge bundle $(\V,\Phi)$. That is, $\mathcal{B}\cong \K^3$, $\alpha=\beta =1$, and $\delta=0$ in \eqref{BetaCyclicG2'Higgs} for some conformal structure $\Sigma$ on $S$. Then the developing map $\dev: \overline{B}_{\Psi_0} \rightarrow \Ein^{2,3}$ from Theorem \ref{Thm:DevEin23} is a diffeomorphism onto $\Omega_\rho$ in \eqref{OmegaEin23_GWPerspective}.
\end{theorem}

\begin{proof}
We begin by proving that the image of our developing map lies inside the domain $\Omega_{\rho}$. In fact, we can prove much more by interpolating the base of pencil used between $\Psi_0$, the tangent pencil to $\Ha^2_{\beta}$, the sub-symmetric space of the $\SL(2,\R)$-subgroup of $\beta$, and $\Psi_{\Delta}$, the tangent pencil to $\Ha^2_{\Delta}$, the sub-symmetric space of the principal $\PSL(2,\R)$-subgroup. We emphasize that the representation $\rho$, and its associated Hodge bundle $(\V, \Phi)$ are fixed here; we vary only the pencils $\Psi_t \in \Omega^1_{\rho}(\tilde{\Sigma}, f^*\T\X)$ along the fixed $\rho$-equivariant harmonic map $f:\tilde{\Sigma} \rightarrow \X$ associated to $(\V, \Phi)$, which realizes a totally geodesic embedding of $\Ha^2_{\Delta}$. \medskip 

Let us now define $\Psi_t \in \Omega^1(\Sigma, \End(\V))$ for $0 \leq t \leq 1$ by 
\[ \Psi_t: = \Psi_0 +t \Psi_{\alpha},\]
where $\Psi_0=\Phi_{-\beta} + \Phi_{-\beta}^*$ as in $\S$\ref{Sec:BasesPencilEin23} and $\Psi_{\alpha} = \Phi_{-\alpha} + \Phi_{-\alpha}^*,$ with $\Phi_{-\alpha}$ explicitly given by
\begin{equation}\label{Phialpha}
\begin{tikzcd}
\mathcal{B} & \mathcal{B}\mathcal{K}^{-1} \arrow[r, "1"] &\mathcal{K}&\mathcal{O} &\mathcal{K}^{-1} \arrow[r, "1"] &\mathcal{B}^{-1}\mathcal{K} &\mathcal{B}^{-1} 
\end{tikzcd}.
\end{equation}
Associated to each base of pencil $\Psi_t$, we can form the same recipe to build a manifold 
$\overline{B}_{\Psi_t} \subset \tilde{\Sigma} \times \Ein^{2,3}$, with fiber 
\[ \overline{B}_{\Psi_t} |_p= \mathcal{B}_{\Ein^{2,3}}(\Psi_t|_p).\]
By $\beta$-regularity, $\overline{B}_{\Psi_t}$ is a smooth 5-manifold for every $0\leq t\leq 1$.

We can read off the Anosov boundary maps of $\rho$ via the embedding $\Ha^{2}_{\Delta} \hookrightarrow \X$. 
Fix a point $x\in \widetilde{\Sigma}$. The isotropic 3-plane $\xi^3_\rho(\zeta)$ for $\zeta\in \partial \pi_1S$ is exactly the sum of the positive eigenspaces for $\Psi(v)$ where $v\in \mathrm{T}_x \widetilde{\Sigma}$ is a tangent vector pointing towards $\zeta$. 
On the other hand, the fiber at $x$ of $\overline{B}_{\Psi_t}$ is mapped by its developing map $\dev_t: \overline{B}_{\Psi_t} \rightarrow \Ein^{2,3}$ to a subset of the isotropic lines of the quadratic form associated to $\Psi_t(v)$ by Proposition \ref{Prop:Ein23BaseViaH}. We will show that $\dev(\overline{B}_{\Psi_t})\subset \Omega_{\rho}$ holds for all $0\leq t \leq 1$. 

Let us fix a background $h$-unitary complex cross product basis for the pullback bundle $\V \rightarrow \tilde{\Sigma}$ at $x$. 
Up to simultaneously conjugating the matrices representing $\Psi_t(v)$ and $\Psi(v)$ by a diagonal unitary gauge, we can represent these endomorphisms by matrices $A_t$ and $A$, respectively, as follows: 
\[ \frac{1}{||\beta||}A_t=\begin{pmatrix}
    0& 1 &0&0&0&0&0\\
    1& 0 &a_t&0&0&0&0\\
    0& a_t&0&\sqrt{2} &0&0&0\\
    0& 0 &\sqrt{2}&0&\sqrt{2}&0&0\\
    0& 0 &0&\sqrt{2} &0&a_t&0\\
    0& 0 &0&0&a_t&0& 1\\
    0& 0 &0&0&0& 1&0
        \end{pmatrix}, \; \;\frac{1}{||\beta||}A=\begin{pmatrix}
    0& 1 &0&0&0&0&0\\
    1& 0 &\sqrt{\frac{5}{3}}&0&0&0&0\\
    0& \sqrt{\frac{5}{3}}&0&\sqrt{2}&0&0&0\\
    0& 0 &\sqrt{2} &0&\sqrt{2}&0&0\\
    0& 0 &0&\sqrt{2} &0&\sqrt{\frac{5}{3}}&0\\
    0& 0 &0&0&\sqrt{\frac{5}{3}}&0& 1\\
    0& 0 &0&0&0&1&0
       \end{pmatrix},  \]
for some real number $0 \leq a_t \leq \sqrt{\frac{5}{3}}$. 
Indeed, the norm of $\alpha$ for the harmonic metric is constant equal to $\sqrt{\frac{5}{3}}$ by Lemma \ref{Lem:MaximumPrincipleBeta}. We now prove the following lemma that is purely elementary algebra.

\begin{lemma}\label{Devt_Einstein}
    The Hermitian matrix $A_t$ is positive on the span of the $3$ eigenvectors with positive eigenvalues of $A$. Thus, $\im(\dev_t) \cap K_{\rho} = \emptyset$.
\end{lemma}

\begin{proof}
The top three eigenvectors $(v_1, v_2, v_3)$ of $\frac{1}{||\beta||}A$, with corresponding eigenvalues 
$\lambda_1=  \sqrt{6}, \, \lambda_2=2\sqrt{\frac{2}{3}},\,\lambda_3 = \sqrt{\frac{2}{3}}$, can be computed explicitly to be the following:
\begin{align*}
    v_1 =& (1,\sqrt{6}, \,\sqrt{15},\, 2 \sqrt{5}, \,\sqrt{15}, \,\sqrt{6}, \,1), \\
    v_2 =& \bigg(-1, \,-2 \sqrt{\frac{2}{3}},\, -\sqrt{\frac{5}{3}},\,0, \, \sqrt{\frac{5}{3}},\,2 \sqrt{\frac{2}{3}},\, 1 \bigg),\\
    v_3 =& \bigg(1,\, \sqrt{\frac{2}{3}}, \,- \sqrt{\frac{1}{15}},\,- \sqrt{\frac{2}{5}},\,- \sqrt{\frac{1}{15}}, \,\sqrt{\frac{2}{3}},1\bigg). 
\end{align*}
One then computes that $\langle A_tv_i, v_i\rangle_h > 0$ for $i \in \{1,2,3\}$ and any $a_t$ such that $0 \leq a_t \leq \sqrt{\frac{5}{3}}$.
\end{proof} 

Lemma \ref{Devt_Einstein} says that the image of $\dev_t$ is disjoint from $\mathbb{P}(\xi_\rho^3(\zeta))$ for each $\zeta\in \partial \pi_1S$. Hence, the image $\dev_t(\overline{B}_{\Psi_t})$ of the developing map $\dev_t$ lies inside $\Omega_{\rho}$ for $0 \leq t \leq 1$ as claimed. \medskip

By \cite[Theorem 1.5]{Dav25}, we know that for $t=1$, the map $\dev_1: \overline{B}_{\Psi_1} \rightarrow \Omega_{\rho}$ is a diffeomorphism to the domain $\Omega_{\rho}$. 
We will compare this map with $\dev = \dev_0$, the given developing map from Theorem \ref{Thm:DevEin23}.

We now make some useful identifications. Set $\Gamma = \pi_1S$. 
Define $B_{\Psi_t} \subset\V[\Ein^{2,3}]$, which yields smooth compact manifolds for $0 \leq t \leq 1$. We then see the compact manifolds $B_{\Psi_t}$ are all diffeomorphic by Ehresmann's fibration theorem. Now, identifying $B_{\Psi_t}$, the Higgs bundle pencil, with the quotient $\Gamma\backslash \overline{B}_{\Psi_t} $, we conclude that $\Gamma\backslash \overline{B}_{\Psi_t} $ is diffeomorphic to $\Gamma\backslash \overline{B}_{\Psi_0}$ for all $0\leq t\leq1$. 

Since $\image(\dev_t)\subset \Omega_{\rho}$, 
the map $\dev_t$ descends to a map $\overline{\dev}_t: \Gamma \backslash \overline{B}_{\Psi_t}\to \rho(\Gamma) \backslash \Omega_\rho$. By the previous identifications, we may write 
$\overline{\dev}_t: \Gamma\backslash \overline{B}_{\Psi_0} \rightarrow \rho(\Gamma) \backslash \Omega_{\rho} $ for the induced maps. 
Now, for $t =0$, the map $\overline{\dev}_0$ is open as a local diffeomorphism by Theorem \ref{Thm:DevEin23}, and closed as $\Gamma \backslash \overline{B}_{\Psi_0}$ is compact. Note that $\Omega_\rho$ is connected: by Proposition \ref{Prop:OmegaThickComplementEin23}, the complement $K$ is a codimension two submanifold of $\Ein^{2,3}$. 
Hence, $\overline{\dev}_0$ surjects onto $\Omega_{\rho}$. As a proper local diffeomorphism,  $\overline{\dev}_0$ is a finite-sheeted covering map. On the other hand, 
$\overline{\dev}_1$ is a diffeomorphism. 
Since $\overline{\dev}_t$, for $0 \leq t \leq 1$, provides a homotopy between the finite-sheeted covering $\overline{\dev}_0$ and the diffeomorphism $\overline{\dev}_1$, we conclude that $\deg(\overline{\dev}_0)= 1$. This means $\overline{\dev}_0$ and $\dev_0$ are each  diffeomorphisms. 
\end{proof}

\begin{remark}
Note that the proof of Theorem \ref{Thm:BetaTits=Pencil} does not address whether the maps $\mathrm{dev}_t$ are local diffeomorphisms for $0<t<1$. 
\end{remark}

As a consequence of the theorem, we have determined the global topology of the quotient $\rho(\pi_1S) \backslash \Omega_{\rho}$ when $\rho$ is Hitchin, recovering \cite[{Corollary 1.6}]{DE25a}. The result is easier to state for the double cover $\widehat{\Ein}^{2,3} \cong \mathbb{S}^2\times \mathbb{S}^3$, the space of isotropic rays in $\R^{3,4}$. 

\begin{corollary}
Let $\rho: \pi_1S\rightarrow \Gtwosplit$ be Hitchin and $\widehat{\Omega}_{\rho} \subset \widehat{\Ein}^{2,3}$ the pullback domain of discontinuity. Then the quotient $\rho(\pi_1S)\backslash \widehat{\Omega}_{\rho}$ is diffeomorphic to the $(\mathbb{S}^2\times \mathbb{S}^1)$-fiber bundle $\mathbb{S}(\underline{\R^3})\oplus \mathbb{S}(\mathcal{K}^3) \rightarrow S$. \footnote{Here we use that $\underline{\mathbb{R}^3}\cong \mathcal{O}\oplus N = \mathscr{U}$ as smooth real vector bundles over $S$.}
\end{corollary}

Using Theorem \ref{Thm:BetaTits=Pencil}, which handles only the Fuchsian case of $\Gtwosplit$-Hitchin representations, we can easily deduce a stronger result. By a result of Labourie \cite{Lab17}, reproven by Collier-Toulisse in \cite{CT24}, every Hitchin representation is the holonomy of a unique cyclic $\beta$-bundle for some Riemann surface $\Sigma=(S,J)$ on the surface $S$. 

We now explain this point in greater detail. In other words, if $\mathcal{H}_{\beta}(S)$ is the moduli space of (polystable) $\beta$-cyclic Higgs bundles up to gauge, Labourie provides a section $s:\Hit(S,\Gtwosplit)\rightarrow \mathcal{H}_{\beta}(S)$ of the holonomy map $\hol:\mathcal{H}_{\beta}(S)\rightarrow \chi(\pi_1S,\Gtwosplit)$ on the $\Gtwosplit$-Hitchin component. This map $s$ operates as follows. Let $\mathcal{Q}_6\rightarrow T(S)$ be the holomorphic vector bundle over Teichm\"uller space with fiber $\mathcal{Q}_{6}|_{[\Sigma]}\cong H^0(\K^6_{\Sigma})$. By Labourie, the total space of $\mathcal{Q}_6$, diffeomorphic to $\R^{14(6g-6)}$ by Riemann-Roch, is canonically diffeomorphic to $\Hit(S,\Gtwosplit)$. Then Labourie's map takes the form  
$s( [\Sigma, q_6]) = [\,(\V, \Phi(q_6))\,]$, where the holomorphic vector bundle $\V$ and Higgs field $\Phi(q_6)$ are: 
\begin{equation}\label{HitchinHiggsBundle}
    \V =\K^{3} \oplus\K^{3} \oplus\K^{2} \oplus\K^{1} \oplus\mathcal{O} \oplus\K^{-1} \oplus \mathcal{K}^{-2} \oplus \mathcal{K}^{-3},
\end{equation}
\begin{equation}\label{HitchinHiggsField}
    \begin{tikzcd} 
    \mathcal{K}^3 \arrow[r, "1"] & \mathcal{K}\arrow[r,"1"]^2&\mathcal{K}\arrow[r,"-i\sqrt{2}"]&\mathcal{O}\arrow[r,"-i\sqrt{2}"]&\mathcal{K}^{-1}\arrow[r,"1"]&\mathcal{K}^{-2}\arrow[r,"1"]\arrow[lllll, bend left, "q_6"]&\mathcal{K}^{-3} \arrow[lllll, bend left, "q_6"] 
\end{tikzcd}.
\end{equation}
We use the same conventions as in Definition \ref{Defn:G2Cyclic} to define the $\Gtwosplit$-structure here. See \cite[Section 3.6]{Bar10} or \cite[Section 2.2]{Eva25} for further details. 

\begin{theorem}[Differential Geometry = Geometric Group Theory]\label{Thm:BetaTits=PencilV2}
Let $\rho: \pi_1S \rightarrow \Gtwosplit$ be any $\Gtwosplit$-Hitchin representation. The developing map $\dev: \overline{B}_{\Psi_0} \rightarrow \Ein^{2,3}$ from Theorem \ref{Thm:DevEin23}, applied to the $\beta$-bundle $s([\rho])$, is a diffeomorphism onto the domain $\Omega_\rho$ from \eqref{OmegaEin23_GWPerspective}.
\end{theorem}

The idea behind this proof is basically the Thurston-Ehresmann principle, however we can treat the case at hand by a straightforward and direct argument.

\begin{proof}
    Since the Hitchin component is a smooth ball, we can construct $(\rho_t)_{t\in [0,1]}$ be a smooth path of representations from the holonomy $\rho_0$ of a $\Gtwosplit$-Hitchin Hodge bundle (a Fuchsian-Hitchin representation) to $\rho_1=\rho$. We therefore can construct the associated path of cyclic bundles, and the associated path of local diffeomorphism $\dev_t:\overline{B_t} \rightarrow \Ein^{2,3}$, where $\overline{B}_t=(\overline{B}_{\Psi_0})_t$ for the cyclic $\beta$-bundle corresponding to $\rho_t$. Let $\Omega_t \subset \Ein^{2,3}$ be the domain of discontinuity \eqref{OmegaEin23_GWPerspective} for $\rho_t$, well-defined since Hitchin representations are $\lbrace\beta\rbrace$-Anosov, and $K_t=\Ein^{2,3}\backslash \Omega_t$ its complement. Denote also $\Gamma = \pi_1S$.
    
    Let $I\subset [0,1]$ be the set of indices $t$ such that the image of $\dev_t$ lies inside the domain of discontinuity $\Omega_{t}$. By Theorem \ref{Thm:BetaTits=Pencil}, $0 \in I$. We first claim $I$ is open. Indeed, in order to check that $\image(\dev_t)\subset \Omega_t$, it suffices to check $\dev_t(D_t) \subset \Omega_t$, for a compact fundamental domain $D_t$ for the action of $\Gamma$ on $\overline{B}_t$. Next, we show that the complement of $I$ is also open. If $t_0 \notin I$, then there is some point $x \in \overline{B}_t$ such that $\dev_{t_0}(x)=\ell\notin \Omega_{t_0}$. Since $\dev_t$ is a local diffeomorphism for all $t\in [0,1]$, that is varying smoothly, then for all $t$ close enough to $t_0$, the image of $\dev_t$ contains a fixed open neighborhood $O$ of $\ell$. However, for $t$ close enough to $t_0$, we see $O \cap K_t \neq \emptyset$, meaning $t\notin I$. We conclude $I$ is open, closed, and non-empty and thus equal to $[0,1]$. On the quotient, $\overline{\dev_1}:\Gamma \backslash \overline{B}_1 \rightarrow \rho(\Gamma)\backslash \Omega_1$ defines a local diffeomorphism of degree one, since it is homotopic to the diffeomorphism $\overline{\dev}_0$. Hence, $\overline{\dev}_1$ is a diffeomorphism and the same holds for $\dev_1:\overline{B}_1\rightarrow \Omega_1$. 
\end{proof}

\section{\texorpdfstring{$(\Gtwosplit, \Pho^\times)$}{G2',Pho(x)}-Geometric Structures}\label{Sec:PhoXGeometricStructures}

This section is the companion to Section \ref{Sec:G2Ein23GeomStructures}. For any Riemann surface $\Sigma$ on $S$, we consider $\alpha$-bundles from Definition \ref{Defn:G2Cyclic} that are also \emph{Hodge bundles}. For each such Higgs bundle, we construct a $(\Gtwosplit, \Pho^\times)$-structure on an $(\RP^2\times \mathbb{S}^1)$-fiber bundle $M \rightarrow S$ whose holonomy descends to $S$ as the representation $\rho:\pi_1S\rightarrow \Gtwosplit$ associated by the non-abelian Hodge correspondence. 
Here, $\Pho^\times =\Gtwosplit/P_{\alpha}$ is the partial flag manifold associated to the parabolic subgroup $P_{\alpha}$ of the long root $\alpha$ of $\g_2'$. The choice of $\mathbb{RP}^2\times \mathbb{S}^1$ for the fibers is motivated by the following result that we also prove: the fibers of the $(\Gtwosplit, \Pho^\times)$ manifolds for $\Gtwosplit$-Hitchin representations, using the \cite{KLP18} domain of discontinuity, are $\RP^2\times \mathbb{S}^1$. 

\subsection{The Flag Manifold \texorpdfstring{$\Pho^\times$}{Pho(x)}}\label{Sec:PhoxModelSpace}

In this section, we introduce the relevant geometry of the $\Gtwosplit$-partial flag manifold $\Pho^\times$, the flag manifold of the long root. The main purpose of the section is to describe a realization of $\Pho^\times$ as a principal $\mathbb{S}^3$-bundle over $\RP^2$ and to develop some geometric intuition that will be necessary for Section \ref{Sec:PhoxFibersHitchin}. 
 
The space $\Pho^\times$ consists of all  \emph{annihilator photons} in $\imoct$: that is, photons with trivial cross product: 
\[ \mathsf{Pho}^\times:= \{ \omega \in \Pho(\imoct) \; | \; \; \omega \times_{\imoct} \omega = 0 \}. \] 
The space of \emph{all} photons $\Pho(\imoct)=\Pho(\R^{3,4})$ is a seven-dimensional $\SO(3,4)$-flag manifold, and $\Pho^\times$ is a codimension two submanifold. Recall that in $\S$\ref{Subsec:G2SymmetricSpace} we saw there is a natural $\Gtwosplit$-equivariant identification $\Pho^\times \cong \Gtwosplit/P_{\alpha}$, where $P_\alpha$ is the maximal parabolic subgroup associated to $\Theta = \{\alpha\}$, for $\alpha$ the long simple root in the $\g_2'$-root system. 

We are headed towards a description of $\Pho^\times$ as a principal bundle. The following result will help us to understand the fibers of this fibration.
For the statement, given any vector subspace $E < \imoct$, we may write $\Pho^\times(E) := \{\omega \in \Pho^\times \mid \omega \subset E\}$. 

Now, we recall a basic fact about photons. Namely, for any spacelike 3-plane $U \in \Gr_{(3,0)}(\R^{3,4})$, each photon $\omega \in \Pho(\R^{3,4})$ is the graph of a unique map $\phi_{U,\omega}:W \rightarrow U^\bot$ for some 2-plane $W \in \Gr_{2}(U)$, where $\phi_{U,\omega}$ is an anti-isometry onto its image. In fact, if $\pi_U: \R^{3,4}\rightarrow U$ denotes the orthogonal projection, then $W = \pi_U(\omega)$. As a consequence, there is a natural projection $\pi_U:\Pho(\R^{3,4})\rightarrow \Gr_{2}(U)$. We will consider the restriction of this map to $\Pho^\times$ when $U$ is chosen in $\X_{\Gtwosplit}$. 

\begin{proposition}\label{Prop:Pho24}
Let $U \in \X_{\Gtwosplit}$ and let $\pi_U:\Pho^\times \rightarrow \Gr_2(U)$ denote the orthogonal projection. For any $W \in \Gr_{2}(U)$, the fiber $\Pho^\times|_W:=\pi_U^{-1}(W)$ is diffeomorphic to $\mathbb{S}^3$. 
\end{proposition} 

\begin{proof} 
Observe that $\Pho^\times|_W =\Pho^\times(W \oplus U^\bot)$.  

Now, fix any orthonormal basis $(u,v)$ for $W$. Then a point $\omega \in \Pho^\times|_W$ uniquely obtains the form $\omega = \spann\{ u + z, v +(uv) z \}$, for some $z \in Q_-(U^\bot) \cong \mathbb{S}^3$ by Corollary \ref{Cor:AnnihilatorPhotonGraph}.
\end{proof}

By the Ehresmann fibration theorem, the map $\pi_U$ then realizes $\Pho^\times$ as an $\mathbb{S}^3$-fiber bundle over $\Gr_2(U)\cong \RP^2$. We now upgrade `fiber bundle' to `principal bundle'. 
To this end, a relevant subgroup will be the pointwise stabilizer $H_U:=\Stab_{\Gtwosplit}^{pt}(U)$ that fixes $U$ pointwise. By Proposition \ref{Prop:FirstStiefel}, we have $H_U\cong \Sp(1) \cong \mathbb{S}^3$. 

\begin{lemma}\label{AnnihilatorPhotonPrincipalBundle}
Let $U \in \X_{\Gtwosplit}$ and define $H_U = \Stab_{\Gtwosplit}^{pt}(U)$. 
The $\mathbb{S}^3$-fiber bundle $\Pho^\times \rightarrow \Gr_2(U)$ is a principal $H_U$-bundle. 
\end{lemma} 

\begin{proof}
Take any $W \in \Gr_2(U)$. Then $H_U$ fixes $W$ identically and preserves the splitting $U \oplus U^\bot$. Hence, $H_U$ preserves the fiber $\Pho^\times|_W$. 
Now, $H_U$ acts simply transitively on $Q_-(U^\bot) \cong \mathbb{S}^3$ by Proposition \ref{Prop:FirstStiefel}. Writing $\omega \in \Pho^\times|_W$ as in Proposition \ref{Prop:Pho24}, one immediately sees $H_U$ acts simply transitively on $\Pho^\times|_W$. 
Since $H_U$ is compact, it follows that $\Pho^\times \rightarrow \Gr_2(U)$ is a principal $H_U$-bundle.
\end{proof}

In fact, the principal bundle $\pi_U:\Pho^\times \rightarrow \Gr_2(U)$ is trivial. The proof shall naturally occur later, but we record this fact here to clarify the topology of this fiber bundle structure. 

\begin{remark}\label{Thm:PhoxTrivial}
Fix $U \in \mathbb{X}_{\Gtwosplit}$. Lemma \ref{Lemma:SimplifiedDevinPho} shows that the principal bundle $\pi_U:\Pho^\times \rightarrow \Gr_{2}(U)$ admits a section and hence is trivial.
\end{remark}

There is a more complicated yet more insightful perspective on the fibers $\Pho^\times|_W$ that we now pursue, leading to another model for $\Pho^\times$. For $W \in \Gr_2(U)$, define $L := W\times W \in \Gr_1(U)$. Then for $x \in Q_+(L)$, the fiber $\Pho^\times|_W$ is intrinsically identified as follows:
\[ \Pho^\times|_W \cong Q_+(\Hom_{\C}( \, (W, \mathcal{C}_x), (U^\bot, \mathcal{C}_x))). \]
Let us unravel this identification. 

Indeed, for $\phi \in \Hom_{\C}((W, \mathcal{C}_x), (U^\bot, \mathcal{C}_x))$, if the operator norm is $|\phi|_{op}=1$, then for $u \in Q_+(W)$, we have $|q(\phi(u))|=1$. Hence, using an orthonormal basis $(u,x u)$ for $U$, there is a natural map
\begin{align} 
Q_+(\Hom_{\C}( \, (W, \mathcal{C}_x), (U^\bot, \mathcal{C}_x))) &\rightarrow \Pho^\times|_W \\
\phi &\mapsto \graph(\phi)=\spann \{u+\phi(u),xu+x\phi(u)\}
\end{align} 
Setting $v= xu$ and $z=\phi(u)$, we find that $\graph(\phi)=\spann \{u+z,v+(uv)z\} \in \Pho^\times(W\oplus U^\bot)$ by Proposition \ref{Prop:PointingTowardsPhotonPrelim}. \medskip 

We can circumvent the choice of $x \in Q_+(L)$, yielding a natural diffeomorphism between $\Pho^\times$ and a certain model space. 
Let $\mathscr{T}: \mathscr{W} \rightarrow \Gr_2(U)$ denote the tautological vector bundle, with fiber $\mathscr{W}|_W=W$. For any $W \in \Gr_2(U)$, the real vector space $\Hom_{\R}(W, U^\bot)$ admits a natural complex structure as follows. Again, set $L = W^\bot \cap U \in \Gr_{1}(U)$. The key observation is that for any $x \in Q_+(L) \cong \sphere^0$ and any $\phi \in \Hom_{\R}(W, U^{\bot})$, the map $\phi$ commutes with $\mathcal{C}_x$ if and only if it commutes with $\mathcal{C}_{-x}$. 

Hence, we may unambiguously define 
\begin{align}\label{eq:HomC-Natural}
     \Hom_{\C}(W,U^\bot):=\{ \phi \in \Hom_{\R}(W,U^\bot) \mid \phi \circ \mathcal{C}_x = \mathcal{C}_x \circ \phi, \, \forall x \in Q_+(W\times W)\}.
\end{align}
This definition leads to a complex rank two vector bundle $\Hom_{\C}(\mathscr{W}, \underline{U^\bot}) \rightarrow \Gr_2(U)$, with fiber at $W \in \Gr_2(U)$ as in \eqref{eq:HomC-Natural}. 
Out of this vector bundle, we obtain a canonical model space for $\Pho^\times$, relative to each choice of $U \in \mathbb{X}_{\Gtwosplit}$.

\begin{corollary}[$\Pho^\times$ Model]\label{Cor:PhoxModel}
Let $U \in \mathbb{X}_{\Gtwosplit}$ and set $K_U = \Stab_{\Gtwosplit}(U)$. There is a $K_U$-equivariant diffeomorphism 
$G:Q_+(\Hom_\C(\mathscr{W}, \underline{U^\bot})) \rightarrow \Pho^\times$ given by 
$(W, \phi) \mapsto \graph(\phi)$.
\end{corollary}

This model is especially useful because it describes the geometry of the $\Gtwosplit$-symmetric space $\mathbb{X}$ and how $\Pho^\times$ is embedded in $\vis\mathbb{X}$. The following proposition clarifies this point and justifies our labor to define this model. Recall the notation from Section \ref{Sec:SymmSpace} that for $f \in \vis\X$ and $x \in \X$, we set $v_{x,f} \in\T^1_x\X$ to be the unique unit tangent vector pointing towards $f$. 

To state the next lemma, we define an embedding $ Q_+(\Hom_\C(\mathscr{W}, \underline{U^{\bot}}) \hookrightarrow \T^1_U\X_{\Gtwosplit} \subseteq \Hom^\times(U,U^\bot)$ as follows: up to a multiplicative constant for renormalization, we extend an element $\phi:W \rightarrow U^\bot$ in the domain to $\hat{\phi}:U \rightarrow U^\bot$, where  
\begin{align*}
    \hat{\phi}|_W &= \phi \\
    \hat{\phi}|_{W\times W} &= 0. 
    \end{align*}
We now use this inclusion to identify the tangent vectors pointing towards $\Pho^\times$. 
\begin{proposition}\label{Prop:PointingTowardsPhoXAgain}
Let $U \in \X_{\Gtwosplit}$. Then the map $\Pho^\times \rightarrow \mathrm{T}^1_U\X_{\Gtwosplit}$ given by 
$\omega \mapsto v_{U, \omega}$ is given by the following composition: 
$\Pho^\times \stackrel{\cong}{\longrightarrow} Q_+(\Hom_\C(\mathscr{W}, \underline{U^\bot})) {\hookrightarrow} \T^1_U\X_{\Gtwosplit}$.
\end{proposition}

Given a \emph{Frenet frame}, $(\mathscr{L},T,N,B)$  (recall Definition \ref{defn:FrenetFrame}), one can trivialize the vector bundle $\Hom_\C(\mathscr{W}, \underline{U^\bot}))$ to obtain a simpler geometric parametrization of $\Pho^\times$. Unfortunately, this map $F$ is only continuous. 

\begin{lemma}[Trivializing $\Pho^\times$]\label{Lem:PhoXtrivialization}
Fix a Frenet frame $(\mathscr{L}, T,N,B)$ and define $U = \mathscr{L} \oplus N$. Then there is a homeomorphism $F:  \Gr_2(U) \times Q_+(\Hom_{\C}(N, U^\bot)) \rightarrow \Pho^\times$ given as follows: 
$\omega = F(W,\delta)$ is the unique annihilator photon that is a graph of a map $\theta:W \rightarrow U^\bot$ such that $\theta|_{W \cap N } = \delta |_{W \cap N}$. 
\end{lemma}

\begin{proof}
\textbf{Step 1: Setup.} By Corollary \ref{Cor:PhoxModel}, $\Pho^\times$ is the sphere bundle of the rank two complex vector bundle $\Hom_{\C}(\mathscr{W}, \underline{U}^\bot)$ on $\Gr_2(U)$. We shall produce a $(C^{0})-$ complex vector bundle isomorphism \[ L:\Gr_2(U)\times \Hom_{\C}(N, U^\bot) \rightarrow \Hom_{\C}(\mathscr{W}, \underline{U}^\bot),\]
which induces the map $F$. In particular, the existence of such a map $L$ will easily imply the result. \medskip 

The map $L$ shall obtain the form $L(W,\delta) = (W, \phi(W,\delta))$, where 
$\phi(W,\delta) \in \Hom_{\C}( W, U^\bot)$. 
We can define $\phi(W,\delta)$ as follows: 
choose any nonzero vectors $u \in W \cap N$, $x \in W^{\bot} \cap N$, and set 
\[ \phi(W,\delta) =\begin{cases} u &\mapsto \delta(u) \\
x \times u &\mapsto x \times \delta(u) \end{cases}\]
Well-definedness of $\phi(W,\delta)$ is easily verified in two cases: \medskip 

\emph{Case 1: } ($W=N$). In this case, $x\times u \in N$ and $\delta(x \times u) = x \times \delta(u)$, so we merely have $\phi(W,\delta) = \delta$. Hence, $\phi$ is well-defined in this case. \medskip 

\emph{Case 2:} ($W \neq N)$. In this case, $W \cap N$ is a line. Thus, $u \in W \cap N$ is well-defined up to real scalars. Thus, $\phi(W,\delta)$ is clearly well-defined independent of choice of $x$ and $u$. \medskip 

We conclude that $\phi$, and hence $L$, is well-defined. By definition, $L$ is a bundle map, which is easily seen to be $\C$-linear. Note also that $L$ is injective. \medskip

\textbf{Step 2: Continuity.} 
Note that we have not addressed the continuity of $L$ yet.
There are naturally two cases to consider: $W \neq N$ and $W = N$. 
In fact, the map $F$ is smooth at points of the former type. To verify smoothness, the key idea is to replace the Grassmannian $\Gr_2(U)$ with the Stiefel manifold $V_2^N(U)$ of orthonormal pairs $(u,v) \in V_2(U)$ such that $u \in N$. In doing so, we obtain an explicit coordinate formula for $L$. That is, we can define an evidently smooth map $\hat{L}$ upstairs by 
\begin{align*} \hat{L}: V_2^N(U) \times \Hom_{\C}(N,U^\bot) &\rightarrow \Hom_{\C}(N, U^\bot)\\
    \hat{L}(u,v,\delta) &= \begin{cases} u\mapsto \delta(u) \\ v \mapsto uv\,\delta(u)\end{cases}.
\end{align*}
There is a projection map $\pi: V_2^N(U)\rightarrow \Gr_2(U)$ by $(u,v)\mapsto \spann \{ u,v\}$, and by definition $\hat{L}$ satisfies $\hat{L} = L \circ (\pi, \id)$. 
Of course, $\pi$ is \emph{not} a submersion, and fails to be so exactly when $(u,v) \in \pi^{-1}(N)$. However, by smoothness of $\hat{L}$, we can deduce smoothness of $L$ downstairs for $W \neq N$. \medskip 

Now, we address continuity of $L$ when $W=N$. In this case, it seems a descent to local coordinates cannot be avoided. Moreover, the resulting formula will make non-smoothness evident. Let us undertake the task now.

Fix an orthonormal basis $(n_1,n_2,e)$ for $U$ once-and-for-all such that $N =\spann \{ n_1,n_2\}$ and $n_2= en_1$. 
We place a smooth local chart for $\Gr_2(U)$ on a neighborhood of $N$ in standard fashion: let $\mathbf{c} =(c_1,c_2)$ correspond to the 2-plane
$ W(\mathbf{c})= \spann \langle n_1+c_1e, n_2+c_2e\rangle.$ 
With some work, one can then verify the following local coordinate formula for $L$, when $W= W(\mathbf{c})$: 
\begin{align}\label{F_incoordinates}
    L(\mathbf{c},\delta) = \begin{cases}
         n_1+c_1e\mapsto \delta(n_1)-\frac{c_1}{c_1^2+c_2^2}\big(c_1n_1+c_2n_2\big)\big(c_2\delta(n_1)-c_1\delta(n_2)\big) ,\\
         n_2+c_2e\mapsto \delta(n_2)-\frac{c_2}{c_1^2+c_2^2}\big(c_1n_1+c_2n_2\big)\big(c_2\delta(n_1)-c_1\delta(n_2)\big) .
    \end{cases}
\end{align}
Here are a few details on the derivation of this formula. First, note that when $\mathbf{c}=0$, we do have $L(\mathbf{0},\delta)=\delta$. Now, suppose $\mathbf{c}\neq \mathbf{0}$. In this case, we must change basis. Define 
\begin{align*}
    w_1&=w_1(\mathbf{c}) := c_2n_1-c_1n_2, \\  w_2&=w_2(\mathbf{c}):= c_1n_1+c_2n_2+(c_  1^2+c_2^2)e.
\end{align*} 
Note that 
$W(\mathbf{c})=\spann \{ w_1,w_2\}$. By definition, $L(\mathbf{c},\delta)$ satisfies $L(\mathbf{c},\delta)(w_1)=\delta(w_1)$. Hence,  $L(\mathbf{c},\delta)(w_2)= \frac{1}{c_1^2+c_2^2}(w_1w_2)\delta(w_1)$ by $\C$-linearity. Here, the multiplicative factor appears because $\frac{1}{||w_1||^2}(w_1w_2)w_1=w_2$.  
A calculation changing basis from $(w_1,w_2)$ to $(u_1,u_2)$, where $u_1=n_1+c_1e$ and $u_2=n_2+c_2e$, and using $\delta(n_2)=e\delta(n_1)$, leads to the formula for $L$ in \eqref{F_incoordinates}. That is, write $A\begin{pmatrix} u_1 \\ u_2 \end{pmatrix} =\begin{pmatrix} w_1 \\ w_2 \end{pmatrix}$ and then one finds $\begin{pmatrix} L(u_1) \\ L(u_2) \end{pmatrix} =A^{-1} \begin{pmatrix}L(w_1) \\ L(w_2) \end{pmatrix}$.

From the expression \eqref{F_incoordinates} for $L$ in local coordinates, we conclude $L$ is 
continuous on an open neighborhood of $\{N\}\times \Hom_{\C}(N,U^\bot)$. Thus, $L$ is a continuous map globally. \medskip 

\textbf{Step 3: Conclude.} 
We have shown $L$ is a continuous isomorphism of complex vector bundles.
Let us write $\underline{\Hom_{\C}(N,U^\bot)}:=\Gr_2(U)\times \Hom_{\C}(N,U^\bot)$ for the trivial vector bundle. 
Consequently, $L$ induces a homeomorphism $\sphere(L)$ between associated sphere bundles:
\[ \sphere(L): \sphere(\underline{\Hom_{\C}(N,U^\bot)}) \rightarrow \sphere(\Hom_{\C}(\mathscr{W}, U^\bot)), \] 
by restriction to the sphere bundle of the domain and renormalization of the output as necessary. 
By definition, $F = G\circ \sphere(L)$, where $G$ is the diffeomorphism from Corollary \ref{Cor:PhoxModel}, and thus $F$ is a homeomorphism.  
\end{proof}

\begin{remark}[Non-smoothness]
The $C^0$-vector bundle isomorphism $L$ in Lemma \ref{Lem:PhoXtrivialization} is not $C^1$ at points $(W,\delta)$ in the domain such that $W = N$, as demonstrated by formula \eqref{F_incoordinates}. Thus, $F$ is also not $C^1$ at such points.
\end{remark}

\subsection{Tits Metric Thickening Domain in \texorpdfstring{$\Pho^\times$}{Pho(x)}}\label{Sec:PhoxKLPComplement}  

In this section, we study the Tits metric thickening domain $\Omega_{\rho}^{\Thick}$ \eqref{Omega_Thick}, for a $P_{\alpha}$-Anosov representation $\rho: \pi_1S \rightarrow \Gtwosplit$.
In particular, we study the topology of the thickening $K_{\Lambda}$ of the limit set $\Lambda = \im(\xi^{\alpha})$, from which $\Omega_{\rho}^{\Thick} = \Pho^\times \setminus K_{\Lambda}$ is built. 
The precise topology of the thickening $K_{\omega}$ of a single annihilator photon $\omega$ is not relevant to the construction of geometric structures, however, the qualitative description of $K_{\omega}$ provided is needed in Theorem \ref{Thm:AlphaHodgeTits=Pencil}. \medskip 

Let $\rho: \pi_1S \rightarrow \Gtwosplit$ be a $P_{\alpha}$-Anosov representation and $\xi^2:\partial_{\infty}\pi_1S \rightarrow \Pho^\times$ the associated limit map. For each $\omega \in \Pho^\times$, define the following thickening: 
\[ K_{\omega}:= \left \{ \omega' \in \Pho^\times \; |\ \titsangle(\omega, \omega') \leq \frac{\pi}{2}  \right \}.\]
The Tits metric thickening domain is again the complement of the thickenings across all elements in the limit set:
\begin{align}\label{Omega_Thick_PhoX}
    \Omega_{\rho}^{\mathrm{Thick}}:= \Pho^\times \backslash \bigcup_{x \in \partial_{\infty}\pi_1S}K_{\xi^2(x)} 
\end{align}
We will see a simpler description of $\Omega^{\mathrm{Thick}}_{\rho}$ shortly. \medskip 

Given $\omega \in \Pho^\times$, we can form the orthogonal complement $\omega^{\bot} \subset \imoct$, then consider $\Pho^\times(\omega^{\bot})$, the space of all annihilator photons $\omega'$ orthogonal to $\omega$. The associated space $\Pho^\times(\omega^\bot)$ is actually just the thickening $K_{\omega}$, as we will show in a technical lemma momentarily. 

We now describe the orbits of the $\Gtwosplit$ on $\Pho^\times \times \Pho^\times$, including model representatives, which are arranged neatly in Figure \ref{Fig:Apartment} in an apartment in $\vis\X$. Here, we write $\sig(V)=(p_+,p_-,p_0)$ for $p_+$, $p_-$, $p_0$ the positive, negative, and null part of the signature of a subspace $V \subset \imoct$. 

\begin{proposition}[$\Gtwosplit$-orbits in $\Pho^\times \times \Pho^\times$]\label{Prop:PhoXOrbits} 
There are four orbits $(\mathcal{O}_i)_{i=1}^4$ of the diagonal $\Gtwosplit$ action on $\Pho^\times \times \Pho^\times$. Fix an $\R$-cross product basis $(x_i)_{i=3}^{-3}$ for $\imoct$. The orbits $\mathcal{O}_i$ admit representatives $\boldsymbol{\omega}_i$, written $\boldsymbol{\omega}_i = (\omega^1_i, \omega^2_i)$, as follows:  
\begin{itemize}[noitemsep]
	\item $\mathcal{O}_0 =  \{ (\omega_1, \omega_2) \in (\Pho^\times)^2 \mid \dim(\omega_1 + \omega_2) = 2 \} $, represented by $\boldsymbol{\omega}_1 = ( \langle x_3, x_2 \rangle, \langle x_3, x_2 \rangle) $ 
	\item $\mathcal{O}_1 = \{ (\omega_1, \omega_2) \in (\Pho^\times)^2 \mid \dim(\omega_1 + \omega_2) = 3 \} $, represented by $\boldsymbol{\omega}_2 = (\langle x_3, x_2 \rangle, \langle x_{3}, x_1 \rangle )$
	\item $\mathcal{O}_2 = \{ (\omega_1, \omega_2) \in (\Pho^\times)^2 \mid \sig(\omega_1 \oplus \omega_2) = (1,1,2) \} $, represented by $\boldsymbol{\omega}_3 = (\langle x_3, x_2 \rangle, \langle x_{-3}, x_{-1} \rangle )$.  
	\item $\mathcal{O}_3 = \{ (\omega_1, \omega_2) \in (\Pho^\times)^2 \mid  \sig(\omega_1 \oplus \omega_2) = (2,2,0)   \} $, represented by $\boldsymbol{\omega}_4 = (\langle x_3, x_2 \rangle, \langle x_{-2}, x_{-3} \rangle )$.
\end{itemize} 
Moreover, if $(\omega, \omega') \in \mathcal{O}_k$, then  $\titsangle(\omega, \omega') = \frac{k\pi}{3}$. 
\end{proposition} 

\begin{proof}
Note that $\dim(\boldsymbol{\omega}_1 + \boldsymbol{\omega}_2) = 2$ if and only if $\boldsymbol{\omega}_1 = \boldsymbol{\omega}_2$. Next, suppose $\dim(\boldsymbol{\omega}_1 + \boldsymbol{\omega}_2) = 3$. Then $\boldsymbol{\omega}_1 \cap \boldsymbol{\omega}_1 \neq \{0\}$, so choose $0\neq w \in \boldsymbol{\omega}_1 \cap \boldsymbol{\omega}_2$. By dimension count, $\boldsymbol{\omega}_1 + \boldsymbol{\omega}_2 = \Ann(w)$. Such a 3-plane $\Ann(w)$ obtains the form $\Ann(w) = \langle x_3, x_2, x_1 \rangle$ for an appropriate $\R$-cross product basis. 

Otherwise, $\dim(\boldsymbol{\omega}_1 + \boldsymbol{\omega}_2) = 4$. Since any maximal isotropic plane is of dimension three in $\R^{3,4}$, there are $y_i \in \omega_i$ such that $q(y_1, y_2) \neq 0$. 
Choose nonzero vectors $z_1 \in \omega_1 \cap y_2^\bot$ and $z_2 \in \omega_2 \cap y_1^{\bot}$ 
so that $(y_i, z_i)$ is a basis for $\boldsymbol{\omega}_i$. Then there are two cases.\medskip

\textbf{Case 1}: $z_1 \cdot z_2 = 0$. Then $\sig(\boldsymbol{\omega}_1 \oplus \boldsymbol{\omega}_2 ) = (1,1,2)$. We can assume $\boldsymbol{\omega}_1 = \langle x_3, x_2 \rangle$. Hence, $\boldsymbol{\omega}_2$ obtains the form $\boldsymbol{\omega}= \langle x_{-3}, z_2 \rangle $
for $z_2 \in \Ann(x_{-3})$. Then $z_2 \, \bot \, \boldsymbol{\omega}_1$ forces $z_2 \in \R \{ x_{-1} \}$. \medskip 

\textbf{Case 2}: $z_1 \cdot z_2 \neq 0$. Up to the $\Gtwosplit$-action, we can assume  $\boldsymbol{\omega}_1 = \langle x_3, x_2 \rangle$. We can also assume $\boldsymbol{\omega}_2= \langle x_{-3}, y \rangle $
for some $y \in \Ann(x_{-3}) $. Since $y \cdot x_{2} \neq 0 $, we have $y = ax_{-2} + bx_{-1}$ for $ a \neq 0$. Using Lemma \ref{Lem:NullStiefel}
find $\varphi \in \Gtwosplit$ such that 
$\varphi \cdot \boldsymbol{\omega}_1= \boldsymbol{\omega}_1$ and 
$\varphi \cdot \boldsymbol{\omega}_2 = \langle x_{-3}, x_{-2} \rangle.$

The `moreover' statement follows by considering the model representatives in an apartment in $\vis\X$. 
\end{proof}

We now see a multiplicity of descriptions of the thickening $K_{\omega}$.
\begin{proposition}[Characterizing $\Pho^\times$-thickenings]\label{Prop:OrthgonalAnnihilatorPhotons}
Let $\omega \in \Pho^\times$. Then for any $\omega' \in \Pho^\times$, the following are equivalent: 
\begin{enumerate}[noitemsep]
    \item $\omega' \in K_{\omega}$, 
    \item $\omega' \in \Pho^\times(\omega^\bot)$, 
    \item $\omega+ \omega' \subseteq \Ann(x)$ for some $x \in \Ein^{2,3}$, 
    \item $\omega \cap \omega' \neq \{0\}$. 
\end{enumerate}
In particular, $K_{\omega} = \Pho^\times(\omega^\bot)$. 
\end{proposition}

\begin{proof}
By Proposition \ref{Prop:PhoXOrbits}, we have the equivalence of (1), (2), and (4). Next, we show (2) implies (3). If $\omega' \in \Pho^\times(\omega^\bot)$, then we must have $\dim(\omega+ \omega') \leq 3$ since a maximal isotropic subspace of $\imoct$ is three dimensional. Thus, take any nonzero vector 
$x \in \omega \cap \omega'$. Then $\omega+\omega'\subset \Ann(x)$.  Hence, (2) implies (3). Conversely, since $\Ann(x)$ is isotropic, it is evident that (3) implies (2).\end{proof}

We now show $K_{\omega}$ is a projective plane, just like in the $\Ein^{2,3}$ case. 

\begin{proposition}[Thickenings in $\Pho^\times$]\label{Prop:PhoXThickeningTopology}
For any annihilator photon $\omega \in \Pho^\times$, its thickening $K_{\omega}$ is homeomorphic to $\RP^2$.
\end{proposition}

\begin{proof}
Let $\omega \in \Pho^\times$ be arbitrary. 
A carefully selected version of the model of $\Pho^\times$ from Lemma \ref{Lem:PhoXtrivialization} makes the topology of the thickening $K_{\omega}$ become transparent. 

Fix any point $U \in \X_{\Gtwosplit}$. Using $\omega$, we refine the pair $(U,U^\bot)$ to a Frenet frame. To this end, define $N := \pi_U(\omega)$. 
Then $\omega$ is uniquely the graph of a map $\delta: N \rightarrow U^\bot$. Define $B$ to be the image of $\delta$. Finally, set $\mathscr{L} := N^\bot \cap U$
and $T:= B^\bot \cap U^\bot$. 
Hence, $\mathscr{F}=(\mathscr{L},T,N,B)$ is a Frenet frame. Let $F: \Gr_2(U) \times Q_+(\Hom_{\C}(N,U^\bot)) \rightarrow \Pho^\times$ be the homeomorphism from Lemma \ref{Lem:PhoXtrivialization} relative to $\mathscr{F}$. Note that $F(N,\delta)= \omega$. We will show 
\begin{align}\label{eq:ThickeningInModel}
     F^{-1}(K_{\omega} ) = \Gr_2(U) \times \{\delta\}. 
\end{align}

For any $W' \in \Gr_2(U)$, the annihilator photon $F(W', \delta)$ by definition intersects $\omega$ non-trivially. Thus, by Proposition \ref{Prop:OrthgonalAnnihilatorPhotons} part (5), we conclude 
$\Gr_2(U) \times \{\delta\} \subseteq F^{-1}(K_{\omega})$. 
Conversely, for any annihilator photon $\omega' = F(W',\delta')$, suppose $\omega'\cap \omega \neq \{0\}$. We will show that $\delta=\delta'$. 
Select a nonzero vector $u \in W' \cap N$, so that  $u+\delta(u) \in \omega \cap \omega'$. Note that $\delta(u)=\delta'(u)$ is forced by definition of $F$, which implies $\delta = \delta'$. Hence, 
$ F^{-1}(K_{\omega}) \subseteq \Gr_2(U) \times \{\delta\}$. Thus, \eqref{eq:ThickeningInModel} holds. Since $F$ is a homeomorphism, we conclude $K_{\omega}$ is homeomorphic to $\Gr_2(U) \cong \RP^2$. 
\end{proof}

Using these observations, we can now describe the union of thickenings $K=\bigcup_{\omega\in \text{im}(\xi^\alpha)}K_\omega$, which is used to define the domain of discontinuity $\Omega_\rho^{\mathrm{Thick}}=\Pho^\times \setminus K$.

\begin{proposition}[Union of Thickenings] \label{Prop:PhoXThickening}
Let $\rho:\pi_1S \rightarrow \Gtwosplit$ be $P_{\alpha}$-Anosov. The thickening of the limit set $K = \bigcup_{\omega\in \text{im}(\xi^\alpha)}K_\omega$ is homeomorphic to the total space of an $\RP^2$-fiber bundle over $\sphere^1$. \end{proposition}

\begin{proof}
By Proposition \ref{Prop:PhoXThickeningTopology}, for $\omega \in \Pho^\times$, the thickening $K_{\omega}$ is diffeomorphic to $\RP^2$. 
By Proposition \ref{Prop:PhoXOrbits}, $K_{\omega} \cap K_{\omega'}=\emptyset$ if $\omega\pitchfork \omega'$, since transversality is equivalent to $\titsangle(\omega,\omega')=\pi$. Thus, there is a continuous projection map $\pi: K \rightarrow \im(\xi^{\alpha})\cong \sphere^1 $, each fiber of which is a copy of $\RP^2$, defining a fiber bundle structure on $K$. 
\end{proof}

\begin{remark}
In fact, in \cite{DE26Anosov}, we prove that $K$ is homeomorphic to $\sphere^1 \times \RP^2$, but we do not need this extra information presently. 
\end{remark}

\subsection{\texorpdfstring{$\Pho^\times$}{Phox}-Fibers for \texorpdfstring{$\alpha$}{a}-Fuchsian Representations}  \label{Subsec:Pho^XSimplePencil}

Let $\rho: \pi_1S\rightarrow \SL(2,\R)_{\alpha}\rightarrow \Gtwosplit$ be an  $\alpha$-Fuchsian representation, factoring through the $\SL(2,\R)$-subgroup associated to the $\mathfrak{sl}_2\R$-subalgebra $\mathfrak{s}_\alpha=\langle E_{\alpha},E_{-\alpha},T_{\alpha}\rangle$ of the simple long $\alpha$ in $\g_2'$; see Appendix \ref{Appendix:SL2} for details on $\mathfrak{s}_{\alpha}$. 
In this subsection we consider a pencil $\mathcal{P}_0\subset \T_U \X_{\Gtwosplit}$ that is tangent to the associated sub-symmetric space $\Ha^2_{\alpha}$. 
We compute the topology of the $\alpha$-base of pencil $\mathcal{B}_{\alpha}(\mathcal{P}_0) \subset \Pho^\times$, which is the fiber of the domain $\Omega^{\Thick}_{\rho}$ over $\tilde{S}$ by Lemma \ref{Lem:NearestPointProjection}. Recall that bases of pencils were defined in Definition \ref{defn:Pencil}. \medskip 

Up to the $\Gtwosplit$-action, we can give an explicit description of $\mathcal{P}_0$ in a model basis. At the end of the section we will verify  this description is compatible with the geometric interpretation from the above paragraph.

Now, let consider a Frenet frame $(\mathscr{L},T,N,B)$ splitting of $\imoct$ as in Definition \ref{defn:FrenetFrame}. 
The Frenet splitting gives an associated point $U\in \X_{\Gtwosplit}$ by $U=\mathscr{L}\oplus N$.
We shall actually want a refined Frenet frame splitting $(x,T,N,B)$, where $x \in Q_+(\mathscr{L})$. 
Recall that the cross product endomorphism $\mathcal{C}_x$ of $x\in \mathscr{L}$ defines a complex structure on $T$, $N$ and $B$. We now obtain a pencil $\mathcal{P}_0$ by defining:
\begin{align}\label{ModelAlphaPencil}
    \mathcal{P}_0=\Hom_\C(N,T)\subset \Hom_{\R}^{\times}(U, U^\perp)\cong \T_U\X_{\Gtwosplit}. 
\end{align}
Here, the map $\Hom_{\R}(N,T)\hookrightarrow \Hom_{\R}^{\times}(U,U^\bot)$ is by extension: namely each $\psi \in \Hom_{\R}(N,T)$ extends uniquely to $\hat{\psi}$ with $\hat{\psi}|_N=\psi$ and $\hat{\psi}|_{\mathscr{L}}=0$. One directly verifies $\hat{\psi}$ is a derivation of $\times$, so that $\hat{\psi} \in \Hom_{\R}^\times(U,U^\bot)$.

We now define an auxiliary construction that will be central to the rest of this section. 
\begin{proposition}[$\Pho^{\times}$-Base of Pencil]\label{Prop:PhoXBase}
For any pencil $\mathcal{P}\subset \T_U\X_{\Gtwosplit}$ and orthonormal vectors $(u,v) \in V_2(U)$, define the following subspace:
\begin{align}\label{RbundleTau2Base}
 \hat{\mathcal{R}}^\perp|_{(u,v)} = \{ z \in U^{\bot} \; |\; \psi(u) \cdot z + ((uv)z)\cdot \psi(v) = 0, \, \forall \psi \in \mathcal{P}  \}.
\end{align} 
Let $z \in Q_-(U^\bot)$. Then  
$\omega = \spann \{ u+z, v+(uv)z\} \in \Pho^\times$ satisfies 
$\omega \in \mathcal{B}_{\alpha}(\mathcal{P})$ exactly when $z \in \hat{\mathcal{R}}^\perp_{u,v}$. 
\end{proposition}

Thus, the $\alpha$-base is explicitly described in terms of $\hat{\mathcal{R}}^{\bot}$. This construction is the $\Pho^\times$-analogue of $\mathcal{R}_{\Psi_0}^\bot$, defined in Proposition \ref{Prop:Ein23Psi0BaseViaR}, which described corresponding base of pencil in $\Ein^{2,3}$.

\begin{proof}
Consider the following space of Stiefel triplets: 
\[ V_{(+,+,-)}(U):= \{(u,v,z) \in V_{(+,+,-)}(\imoct) \mid u,v \in U, z \in U^\bot \}.\]
Recall from Proposition \ref{Prop:PointingTowardsPhotonPrelim}, for any triplet $(u,v,z) \in V_{(+,+,-)}(U)$, we obtain a tangent vector $\phi_{u,v,z} \in \T_{U}\X$ given by \eqref{ModelAlphaGeodesic}, and $\phi_{u,v,z}$ points towards the annihilator photon $\omega = \spann \langle u+z, v+(uv)z\rangle$. Moreover, every tangent vector pointing towards $\Pho^\times$ obtains this form. Thus, the claim follows by the fact that for any tangent vector $\psi \in \T_{U}\X$, we have 
\[ \langle \phi_{u,v,z}, \psi\rangle_{\X} =0 \iff \langle \psi(u), z\rangle_{q}+\langle (uv)z, \psi(v)\rangle_{q} =0,\] 
by Corollary \ref{Cor:Orthogonality}. 
\end{proof}

Let us make a small observation regarding the subspaces $\hat{\mathcal{R}}^\perp_{u,v}$.

\begin{proposition}\label{Prop:R_Invariance}
Let $\mathcal{P} \subset \T_U\X_{\Gtwosplit}$ be a pencil. Then for any orthonormal basis $(u,v)$ of $U$, we have $\hat{\mathcal{R}}^\perp_{-u,v} = \hat{\mathcal{R}}^\perp_{u,v} = \hat{\mathcal{R}}^\perp_{u,-v} $. 
\end{proposition} 

\begin{proof}
The equations for $z \in \hat{\mathcal{R}}^\perp$ are linear $z$ and $u$ and quadratic in $v$. The claims easily follow.
\end{proof} 

\begin{remark}
In fact, we will see the subspaces $\hat{\mathcal{R}}^\perp_{u,v}$ define a rank two vector bundle $\hat{\mathcal{R}}^\perp\rightarrow V_2(U)$. Due to Proposition \ref{Prop:R_Invariance}, $\hat{\mathcal{R}}^\perp$ descends to a vector bundle over $\Flag(U)$, the space of full flags of $U$. 
\end{remark}

Given a unit vector $u\in Q_+(U)$ and a two-plane $W\in \Gr_2(U)$ containing $u$, we may then denote by $\hat{\mathcal{R}}^\perp_{[u],W}$ the plane $\hat{\mathcal{R}}^\perp_{u,v}$ where $v\in W$ is any unit element orthogonal to $u$.
We now describe the subspaces $\hat{\mathcal{R}}^\perp_{[u],W}$ for the pencil $\mathcal{P}_0$ in the case that $u \in N$. 

\begin{lemma}[$\hat{\mathcal{R}}^\perp$, Essential Case]\label{Lemma:ProjectionOntoB}
Let $\mathcal{P}_0 \subset \T_U\X_{\Gtwosplit}$ be the pencil \eqref{ModelAlphaPencil}, $W \in \Gr_{2}(U)$ and $u \in Q_+(N \cap W)$. Then the orthogonal projection map $\pi_B: \hat{\mathcal{R}}^\perp_{[u],W} \rightarrow B$ is a linear isomorphism.
\end{lemma}

Before the next technical proof, recall that in $\imoct$, we have $u\times v = uv $ if $u \, \bot \,v$. We will also frequently use cross product relations among $(\mathscr{L},T,N,B)$ from \eqref{ModelFrenetFrame}. For example, $N \times_{\imoct} N=\mathscr{L}$ and any pair of distinct subspaces from $(T,N,B)$ multiply to the third subspace.

\begin{proof} We consider the standard multiplication basis $(\mathbf{i}, \mathbf{j},\mathbf{k},\mathbf{li},\mathbf{lj},\mathbf{lk})$ for $\imoct$ again, as well as its Frenet frame splitting \eqref{ModelFrenetFrame}. Extend $u$ to an orthonormal basis $(u,v)$ for $W$. 

Let us set $\psi_1 \in \mathcal{P}_0$ as the unique element in $\Hom_\C(N,T)$ that maps $\mathbf{j}$ to $\mathbf{l}$. We then define $\psi_2 \in \mathcal{P}_0$ as $\psi_2 := \mathcal{C}_{\mathbf{i}} \circ \psi_1$, which maps $\mathbf{j}$ to $\mathbf{il}=-\mathbf{li}$. A candidate element $z \in T\oplus B$ satisfies $z \in \mathcal{R}_{u,v}$ if and only if the following equations both hold:
\begin{align} 
	\psi_1(u) \cdot z + \psi_1(v) \cdot ((uv)z) &= 0 \label{Orthogonality1} \\
	\psi_2(u) \cdot z + \psi_2(v) \cdot ((uv)z) &= 0 \label{Orthogonality2}. 
\end{align}

\medskip
 
Since $u \in Q_+(N)$, without any loss of generality we can assume $u=\mathbf{j}$. Indeed, we can re-gauge by $\mathscr{T}=\Stab_{\Gtwosplit}(\mathbf{i},T,N,B)$, which acts transitively on $Q_+(N)$ by Proposition \ref{Prop:FirstStiefel}. 

Hence, any unit element $v \in Q_+(U)$ orthogonal to $u$ obtains the form 
\[v = a_1 \mathbf{i}+ a_2 \mathbf{k} \] 
for some $a_1,a_2\in \R$ such that $a_1^2+a_2^2=1$. 
By definition of $\mathcal{P}_0$, for all $\psi\in \mathcal{P}_0$ one has $\psi(\mathbf{i})=0$. In particular:
\[ \psi_1(v) = \psi_1(a_2 \mathbf{k}) = a_2 \psi_1(\mathbf{ij}) = a_2 \, \mathbf{il} \] 
Note that $\psi_2(\mathbf{j})=\mathbf{il}$ and  $\psi_2(\mathbf{k}) = -\mathbf{l}$. 
 
 Here is an auxiliary calculation that is useful: take any $z \in U^{\bot}$ and write $z = z^T + z^N$ for the orthogonal projections onto $T, N$, respectively. Then 
 \begin{align}\label{zuv_calculation}
 (uv)z = \underbrace{( (uv)^{\mathscr{L}} z^B + (uv)^N z^T)}_{B-\text{projection}} + \underbrace{( (uv)^{\mathscr{L}}z ^T + (uv)^N z^B)}_{T-\text{projection}}.
\end{align}
We shall denote $\pi_T, \pi_N, \pi_B$ as the orthogonal projections from $\imoct$ to the given subspace. 
Now, since $u \in N$, writing $v = v^{\mathscr{L}} + v^{N}$, then $(uv)^{\mathscr{L}} = uv^N$ and $(uv)^N = uv^{\mathscr{L}}$. Hence,  
\[ \pi_{T}( (uv)z)  = (u v^N)z^T + (uv^{\mathscr{L}}) z^B = a_2(\mathbf{i}z^T) - a_1(\mathbf{k}z^B). \]
Hence, equations \eqref{Orthogonality1} and \eqref{Orthogonality2} take the following form: 
\begin{align}
	\mathbf{l} \cdot z + (a_2 \mathbf{il}) \cdot \big( (a_2 \mathbf{i})z^T - (a_1\mathbf{k})z^B \big) &= 0 \label{Ortho1}, \\
	\mathbf{il} \cdot z + (-a_2 \mathbf{l}) \cdot \big((a_2 \mathbf{i})z^T - (a_1\mathbf{k})z^B \big) &= 0. \label{Ortho2}
\end{align}
We can simplify further by decomposing $(uv)z$. Let us write 
\begin{align}\label{General_z}
	z= z^T + z^{B} = c_1 \mathbf{l} + c_{2} \mathbf{il} + c_3 \mathbf{lj}  +c_4 \mathbf{lk}.
\end{align} 
Then one finds 
\begin{align}
	 \pi_{[\mathbf{l}]}( (uv)z) &= a_2c_2 (\mathbf{i} \times \mathbf{il})  -a_1c_4 (\mathbf{k} \times \mathbf{lk}) = (-a_2 c_2 -a_1c_4) \mathbf{l} \label{Projection_t}\\
	 \pi_{[\mathbf{il}]}( (uv)z) &= a_2c_1 (\mathbf{i} \times \mathbf{l}) - a_1c_3 (\mathbf{k} \times \mathbf{lj}) = (a_2 c_1 + a_1c_3)\mathbf{il}. \label{Projection_t'}\\
	 \pi_{[\mathbf{lj}]}( (uv)z) &= a_2c_4(\mathbf{i}\times \mathbf{lk}) - a_1c_2(\mathbf{k} \times \mathbf{il}) = (a_2c_4-a_1c_2)\mathbf{lj} \label{Projection_b}\\
	 \pi_{[\mathbf{lk}]}((uv) z) &= a_2c_3(\mathbf{i} \times \mathbf{lj}) - a_1c_1(\mathbf{k} \times \mathbf{l}) = (-a_2c_3+a_1c_1)\mathbf{lk} \label{Projection_b'}
\end{align}
Combining \eqref{Ortho1}-\eqref{Ortho2} as well as \eqref{Projection_t}-\eqref{Projection_t'}, we find that $z\in\hat{\mathcal{R}}^\perp_{u,v}$ exactly when:
\begin{align}
	c_1 + a_2^2c_1+ a_1a_2c_3 &=0 \label{c1Eqn}\\
	c_2 +a_2^2c_2+a_1a_2c_4&=0 \label{c2Eqn}. 
\end{align} 
Clearly, we may solve for $c_1, c_2$ in terms of $c_{3}, c_4$. Hence the subspace $\hat{\mathcal{R}}^\perp_{u,v}$ defined by these equation is a plane and the orthogonal projection to $B=\langle\mathbf{lj},\mathbf{lk} \rangle$ is an isomorphism. \end{proof}  

As a corollary, we can define the following construction. 

\begin{definition}
Let $W \in \Gr_2(U)$ and $u \in Q_+(N \cap W)$. Denote by $\Gamma_{[u],W}$ the unique linear map 
$\Gamma_{[u],W}: B \rightarrow T$ whose graph is $\hat{\mathcal{R}}^\perp_{[u],W}$.
\end{definition}

We can write down the map $\Gamma_{[u],W}$ in coordinates as a consequence of the previous lemma. 

\begin{corollary}\label{Cor:RGraphEquation}
For any $W \in \Gr_2(U)$, write $W= \langle u, \,a_1w_1+a_2w_1\times u\rangle$, for $u \in Q_+(N \cap W)$ where $w_1\in \mathscr{L}$ is a norm one vector and $a_1^2+a_2^2=1$. Then 
$\Gamma_{[u],W}: B \rightarrow T$ is given as follows:
\[ \Gamma_{[u],W}(b) = \frac{a_1a_2}{1+a_2^2}(b \times u). \]
\end{corollary}

Note that the above equation is well-defined independent of choice of $u$. 

\begin{proof}
Using previous notation, assume without any loss of generality that $u=\mathbf{j}$. The equations \eqref{c1Eqn}, \eqref{c2Eqn} imply that elements $z\in\hat{\mathcal{R}}^\perp_{[u],W}$ are of the following form for $c_3,c_4\in \R$:
\begin{equation}z=-c_3\frac{a_1a_2}{1+a_2^2}\mathbf{l}-c_4\frac{a_1a_2}{1+a_2^2}\mathbf{il}+c_3\mathbf{lj}+c_4\mathbf{lk}.\label{eq:ValueZinRuv}\end{equation}
The result then follows from the fact that $\Gamma_{[u],W}(\mathbf{lj})= -\frac{a_1a_2}{1+a_2^2}\mathbf{l}$ and $\Gamma_{[u],W}(\mathbf{lk})= -\frac{a_1a_2}{1+a_2^2}\mathbf{il}$.
\end{proof}

Finally, we can show the base of pencil $\mathcal{B}_{\alpha}(\mathcal{P}_0)$ is an explicit circle bundle over $\Gr_2(U)$. 
For the statement and the consequences, it is useful to change perspective the bundle $\hat{\mathcal{R}}^\perp$ and as a vector sub-bundle of 
the bundle $\Hom_{\C}(\mathscr{W},U^\bot) \rightarrow \Gr_2(U)$ from $\S$\ref{Sec:PhoxModelSpace}. 
We define the real rank two vector bundle $\mathcal{R}^\perp\rightarrow \Gr_2(U)$ with fiber at $W \in \Gr_2(U)$ given by  
\[ \mathcal{R}^\perp|_{W} = \{ \phi \in \Hom_{\C}(W,U^\bot) \mid \phi \in \mathcal{P}_0^{\bot} \}. \]
The fiber $\mathcal{R}^{\bot}|_W$ carries the information of all of the subspaces $\hat{\mathcal{R}}^\perp_{[u],W}$ as follows: 
\begin{align}\label{R-equivalence}
    \hat{\mathcal{R}}^\perp_{[u],W} = \{ \phi(u) \mid \phi \in {\mathcal{R}}^\perp|_{W} \}. 
\end{align}

Via the construction of ${\mathcal{R}}^\perp$, in the model $G=G_U$ from Corollary \ref{Cor:PhoxModel} we see that $\mathcal{B}_{\alpha}(\mathcal{P}_0)$ is an $\sphere^1$-fiber sub-bundle of $\Pho^\times \rightarrow \Gr_2(U)$. In particular, $\mathcal{B}_{\alpha}(\mathcal{P}_0)$ is the sphere bundle of ${\mathcal{R}}^\perp$. 

\begin{corollary}\label{Cor:PhoXBaseCircleBundle}
The orthogonal projection $\pi_U:\mathcal{B}_{\alpha}(\mathcal{P}_0) \rightarrow \Gr_{2}(U)$ defines a surjective submersion and an $\mathbb{S}^1$-fiber bundle,
with fiber at $W$ given by 
\[ \mathcal{B}_{\alpha}(\mathcal{P}_0)|_{W} = \{ \graph(\phi) \in \Pho^\times \mid \phi \in {\mathcal{R}}^\perp|_{W} ,\, ||\phi||_{op}=1\}. \]
\end{corollary} 

\begin{proof} 
Fix any $W \in \Gr_{2}(U)$. 
Note that $W \cap N \neq \{0\}$. Thus, we may choose a unit spacelike vector $u \in Q_+(N\cap W)$. Let $v=wu$, so $(u,v)$ is an orthonormal basis of $W$. 
By Proposition \ref{Prop:PhoXBase}, 
if $z \in Q_-(U^\bot)$, the annihilator photon $\omega = \langle u+z, v+(uv)z\rangle \in \Pho^\times$ satisfies $\omega \in \mathcal{B}_{\alpha}(\mathcal{P}_0)$ if and only if $z \in \hat{\mathcal{R}}^\perp_{u,v}$. In particular, Lemma \ref{Lemma:ProjectionOntoB} implies $\mathcal{B}_{\alpha}(\mathcal{P}_0)|_{W} \cong Q_-(\hat{\mathcal{R}}^\perp_{[u],W}) \cong \mathbb{S}^1$. 
Equivalently, $\omega= \graph(\phi)$ for $\phi \in {\mathcal{R}}^\perp{|_W}$ the unique $\C$-linear map
$\phi: W \rightarrow U^\bot$ satisfying $\phi(u)=z$. The description of $\mathcal{B}_{\alpha}(\mathcal{P}_0)|_W$ then holds by the equivalence in \eqref{R-equivalence}. 

The projection $\pi_U: \mathcal{B}_{\alpha}(\mathcal{P}_0) \rightarrow \Gr_2(U)$ is proper by compactness, smooth, and surjective. Since $\pi_U$ is also a submersion, the Ehresmann fibration theorem says that $\pi_U$ defines an $\mathbb{S}^1$-fiber bundle.  
\end{proof}

In fact, we now show the fiber bundle $\pi_U:\mathcal{B}_{\alpha}(\mathcal{P}_0) \rightarrow \Gr_2(U)$ is trivial. 
By the previous corollary, this amounts to trivializing the vector bundle ${\mathcal{R}}^\perp \rightarrow \Gr_2(U)$. We introduce a definition that will aid in this endeavor. 

\begin{definition}
Denote 
$\underline{\Hom_{\C}(N,B)}:=\Gr_2(U) \times \Hom_{\C}(N,B)$ and consider the vector bundle 
morphism 
\[ \hat{\sigma}: \underline{\Hom_{\C}(N,B)}\rightarrow {\mathcal{R}}^\perp, \]
given as follows. For any $W \in \Gr_2(U)$,
take any $u \in Q_+(W \cap N)$ and define 
$\hat{\sigma}(W,\delta)$ to be the unique $\C$-linear map $\hat{\sigma}(W,\delta):W \rightarrow T \oplus B$ satisfying
\[ \hat{\sigma}(W,\delta)(u) = (\Gamma_{[u],W}\circ \delta)(u) +\delta(u).  \]
One can verify that $\hat{\sigma}$ is well-defined independent of choice of $u$ and $\hat{\sigma}(W,\delta) \in \mathcal{R}^\perp$ by Corollary \ref{Cor:RGraphEquation}. 
\end{definition}

We now finally give an explicit geometric and topological description of the base of pencil $\mathcal{B}_{\alpha}(\mathcal{P}_0)$. 

\begin{lemma}[Simplified $\alpha$-Base Topology]\label{Lemma:SimplifiedDevinPho}

The morphism $\hat{\sigma}$ is a ($C^0$) bundle isomorphism. 
Thus, there is a homeomorphism
\[ \psi: \Gr_2(U) \times Q_+(\Hom_\C(N,B)) \rightarrow \mathcal{B}_{\alpha}(\mathcal{P}_0),\]
which is given by $\psi=F\circ \sphere(\hat{\sigma})$. Equivalently, $\psi$ can also be written as follows: given any arbitrary $u \in Q_+(N \cap W)$ and $x \in Q_+(W \times W)$, then:
\[\psi( W, \delta) = \spann \big\langle u + \sphere(\hat{\sigma})(W,\delta)(u),\, xu + x \sphere(\hat{\sigma})(W,\delta)(u) \,\big\rangle.\]

Moreover, $\mathcal{B}_{\alpha}(\mathcal{P}_0)$ is diffeomorphic to $ \RP^2 \times \mathbb{S}^1$. 
\end{lemma}

Here $\sphere(\hat{\sigma})$ denotes the induced map of sphere bundles, meaning that on each fiber we restrict $\hat{\sigma}$ to input unit vectors and renormalize the output of $\hat{\sigma}$ so that it has operator norm $1$, and $F$ is the homeomorphism from Lemma \ref{Lem:PhoXtrivialization}.

\begin{proof}
The coordinate formula in Corollary \ref{Cor:RGraphEquation} for $\Gamma_{[u],W}$ determines the formula for $\hat{\sigma}$, and validates that $\hat{\sigma}$ is continuous. It is evident that $\hat{\sigma}$ is a vector bundle morphism lifting the identify on $\Gr_2(U)$. 
The map $\hat{\sigma}$ is injective on fibers and thus an isomorphism fiberwise. Hence, $\hat{\sigma}$ is a continuous isomorphism of vector bundles over $\Gr_2(U)$. 

The fact that the explicit formula for $\psi$ is correct comes from definition unraveling. First, since $u\in Q_+(N\cap W)$, then $\sphere(\hat{\sigma})(W,\delta)(u)\in Q_-(\hat{\mathcal{R}}^{\bot}_{u,W})$. Here, recall the relation between $\mathcal{R}^\bot$ and $\hat{\mathcal{R}}^\perp$ from \eqref{R-equivalence}. Then $W \in \Gr_2(U)$ and $[u+\sphere(\hat{\sigma})(W,\delta)(u)]\in \Ein^{2,3}$ uniquely determine the annihilator photon $\psi(W,\delta)$ using Corollary \ref{Cor:AnnihilatorPhotonGraph}.

 We have verified that the smooth fiber bundle $\pi_U: \mathcal{B}_{\alpha}(\mathcal{P}_0)\rightarrow \Gr_2(U)$ is topologically trivial. However, this implies $\pi_U$ is also smoothly trivial and thus diffeomorphic to $\RP^2\times \sphere^1$. 
\end{proof}

We finally reach the desired conclusion of this subsection.
\begin{corollary}[$\Pho^\times$-fibers for $\alpha$-Fuchsians]\label{Cor:PhoXSimpleFibers}
Let $\rho:\pi_1S \rightarrow \SL(2,\R) \stackrel{\iota_{\alpha}}{\hookrightarrow} \Gtwosplit$ be an $\alpha$-Fuchsian representation in $\Gtwosplit$ and $\Omega_{\rho}^{\Thick}$ the domain \eqref{Omega_Thick_PhoX}. Then there is a smooth fiber bundle $\Omega^{\Thick}_{\rho} \rightarrow \Ha^2$ with fiber diffeomorphic to $\RP^2\times \mathbb{S}^1$. 
\end{corollary}

\begin{proof}
The representation $\rho$ is $\alpha$-Anosov and admits a $\rho$-equivariant totally geodesic minimal surface $f: \tilde{S} \rightarrow \X_{\Gtwosplit}$ whose image is the sub-symmetric space $\mathbb{H}^2_{\alpha}$ associated to the $\SL(2,\R)$-subgroup $\SL(2,\R)_{\alpha}$ of the long root $\alpha$. Now, the fiber of $\Omega^{\Thick}$ is diffeomorphic to  $\mathcal{B}_{\alpha}(\mathcal{P}_\alpha)$, for $\mathcal{P}_\alpha=df(\mathrm{T}_p\tilde{S})$, by Lemma \ref{Lem:NearestPointProjection}. The result then follows from Lemma \ref{Lemma:SimplifiedDevinPho} if we know the pencils $\mathcal{P}_\alpha$ and the model pencil $\mathcal{P}_0$ in \eqref{ModelAlphaPencil} are equivalent up to the $\Gtwosplit$-action. This holds because, in an appropriate $\R$-cross product basis $(x_i)_{i=3}^{-3}$, the group $\SL(2,\R)_{\alpha}$ fixes $x_3, x_0, x_{-3}$ identically, and acts faithfully on the 2-planes  $\spann \{x_2,x_1\}$ and $\spann \{x_{-1},x_{-2}\}$ by Appendix \ref{Subsec:SL2LongRoot}. In particular, every $\psi \in \mathcal{P}_\alpha$ obtains the form of a rank two map $\psi: U \rightarrow U^\bot$, for $U = f(p)$, with the same image and kernel. We conclude $\mathcal{P}_0$ and $\mathcal{P}_{\alpha}$ are $\Gtwosplit$-equivalent.
\end{proof}

\subsection{\texorpdfstring{$\Pho^\times$}{Phox}-Fibers for Hitchin Representations}\label{Sec:PhoxFibersHitchin}

In the previous subsection, we described the fiber of the quotient of the domain of discontinuity in $\Pho^\times$ for $\alpha$-Fuchsian representations. In particular, these fibers are diffeomorphic to $\RP^2\times \mathbb{S}^1$. 
We now briefly explain that why the corresponding fibers for $\Gtwosplit$-Hitchin representations are diffeomorphic. \medskip 

Let $\rho_0: \pi_1S \rightarrow \Gtwosplit$ be a fixed Fuchsian-Hitchin representation. The $\alpha$-bundle $(\V,\Phi)$ associated to $\rho_0$ via Labourie's map, defined in Section \ref{Subsec:G2Ein23HitchinCase}, obtains the form \eqref{HitchinHiggsBundle},\eqref{HitchinHiggsField}, for some Riemann surface $\Sigma = (S,J)$, where $q_6=0$. In particular, for the following discussion, this Higgs bundle is fixed once-and-for-all. 

Using $\rho_0$ to uniformize $\tilde{S}$, the associated $\rho_0$-equivariant harmonic map is a totally geodesic embedding $\iota_{\Delta}:\Ha^2\rightarrow \X_{\Gtwosplit}$. Let $\Ha^2_{\Delta} = \image(\iota_{\Delta})$ be the sub-symmetric space of the principal $\PSL(2,\R)$-subgroup and take $\mathcal{P} \subset \T_x\X$ be the tangent pencil to $\Ha^2_{\Delta}$ at $x \in \Ha_{\Delta}$. By Lemma \ref{Lem:NearestPointProjection}, the fiber of the compact quotient $M_{\rho_0} = \rho_0(\pi_1S)\backslash \Omega_{\rho_0}$ over $S$ is the base of pencil $\mathcal{B}_{\alpha}(\mathcal{P})\subset \Pho^\times$. We now explain why this base of pencil $\mathcal{B}_{\alpha}(\mathcal{P})\subset \Pho^\times$ is diffeomorphic to the base of pencil $\mathcal{B}_{\alpha}(\mathcal{P}_0)$, where $\mathcal{P}_0$ is the pencil studied previously: a tangent space to the sub-symmetric space $\Ha^2_{\alpha}$.

We use the same strategy as in \cite{DE26}. Viewing the two pencils $\mathcal{P},\mathcal{P}_0$ as living in a common tangent space $\T_x\X$, we shall deform $\mathcal{P}$ to $\mathcal{P}_0$. It is easiest to see this deformation by describing the pencils in terms of the associated Higgs bundle $(\V,\Phi)$ to $\rho_0$. In particular, we can describe this deformation via a one parameter family \emph{auxiliary} endomorphism-valued one-forms on the \emph{fixed} Higgs bundle $(\V,\Phi)$. 
Now, define following one parameter family $\Phi_t \in \Omega^{1,0}(\Sigma, \End^\times(\V))$, for $0\leq t\leq 1$,  
\[ \Phi_t = \left [ \mathcal{K}^3 \stackrel{t}{\longrightarrow}\mathcal{K}^2 \stackrel{1}{\longrightarrow}\mathcal{K}^1 \stackrel{-t\sqrt{2}i}{\longrightarrow} \mathcal{O} \stackrel{-t\sqrt{2}i}{\longrightarrow} \mathcal{K}^{-1} \stackrel{1}{\longrightarrow}\mathcal{K}^{-2} \stackrel{t}{\longrightarrow}\mathcal{K}^{-3}  \right ].\]
Note that $\Phi = \Phi_1$ is the Fuchsian-Hitchin Higgs field. 
Let $h_0$ be the harmonic metric on $(\V,\Phi)$ and $f_0: \tilde{\Sigma} \rightarrow \X_{\Gtwosplit}$ the associated harmonic map.   
Fix a basepoint $p_0 \in \tilde{S}$ and write $x_0 = \pi(p_0)$ for $\pi : \tilde{S} \rightarrow S$. Under this setup, we may define a family of pencils as follows: 
\begin{align}\label{PhoxPencils}
 \mathcal{P}_t = \{ \;(\Phi_t+\Phi_t^{*h_0})(X) \in \End^\times(\V^\R|_{x_0}) \; | \; X \in \T_{x_0}S .\} 
\end{align} 
Using the same identifications as in Section \ref{Sec:BasesPencilEin23}, we may regard $\mathcal{P}_t$ as a pencil in  $\T_{f(p_0)}\mathbb{X}_{\Gtwosplit}$. 

We now conclude that the simplified $\alpha$-base of pencil has the same topology as the original. 

\begin{lemma}\label{Lem:PhoXPencilDeformation}
The bases of pencil $\mathcal{B}_{\alpha}(\mathcal{P}_0)$ and $ \mathcal{B}_{\alpha}(\mathcal{P}_1)$ are diffeomorphic. 
\end{lemma}

\begin{proof}
By Proposition \ref{Prop:Alpha2RegularPencil}, the family $\mathcal{P}_t$ of pencils is $\alpha$-regular for $0 \leq t \leq 1$. Thus, by \cite[Corollary 6.8]{Dav25}, the result holds. 
\end{proof} 

\begin{corollary}[$\Pho^\times$-fibers for Hitchin Representations]\label{Cor:PhoxFibersHitchin}
Let $\rho:\pi_1S\rightarrow \Gtwosplit$ be Hitchin. The quotient  $\rho(\pi_1S)\backslash \Omega^{\Thick}_{\rho}$ of the domain \eqref{Omega_Thick_PhoX} is diffeomorphic to an $\RP^2 \times \mathbb{S}^1$-fiber bundle over $S$.
\end{corollary}

\begin{proof}
Fix $\rho: \pi_1S\rightarrow \Gtwosplit$ a Fuchsian-Hitchin representation. 
By Lemma \ref{Lem:NearestPointProjection}, the quotient $M_{\rho}$ is a fiber bundle over $S$ whose fiber is diffeomorphic to $\mathcal{B}_{\alpha}(\mathcal{P}_1)$. Then Lemma  \ref{Lem:PhoXPencilDeformation} and Corollary  \ref{Cor:PhoXSimpleFibers} give the desired result for $\rho$. For any Hitchin representation $\rho':\pi_1S\rightarrow \Gtwosplit$, by \cite[Corollary 7.15]{Dav25} we have a diffeomorphism $\rho'(\pi_1S)\backslash \Omega_{\rho'} \cong \rho(\pi_1S)\backslash \Omega_{\rho}$. 
\end{proof}

\subsection{Stability for \texorpdfstring{$\alpha$}{a}-cyclic Bundles}\label{Sec:StabilityAlplhaBundles}

In this subsection, we describe the stability conditions for  $\alpha$-cyclic Higgs bundles. 

Recall that $\alpha$-cyclic Higgs bundles, defined in Definition \ref{Defn:G2Cyclic}, obtain the following form:
\begin{equation}\label{AlphaCyclic}
    \begin{tikzcd}
\mathcal{T}^2\mathcal{K} \arrow[r, "\beta"] & \mathcal{T}\mathcal{K}\arrow[r, "1"]&\mathcal{T}\arrow[r, "-i\sqrt{2}\beta"]&\mathcal{O} \arrow[r, "-i\sqrt{2}\beta"]&\mathcal{T}^{-1} \arrow[r, "1"]&\mathcal{T}^{-1}\mathcal{K}^{-1} \arrow[r, "\beta"]\arrow[bend left, lllll, "\delta"]&\mathcal{T}^{-2}\mathcal{K}^{-1} \arrow[bend left, lllll, "\delta"]
\end{tikzcd}.
\end{equation}

 In this construction:
 \begin{itemize}[noitemsep]
     \item $\mathcal{T}$ is a holomorphic line bundle on $\Sigma$. 
     \item $\beta \in H^0(\Sigma, \mathcal{K}\mathcal{T}^{-1} ).$
     \item $\delta \in H^0(\Sigma, \mathcal{T}^3\mathcal{K}^3)$. 
 \end{itemize}

\begin{remark}
Just as in Section \ref{Subsec:Ein23Structures}, we must demand $\beta \neq 0$ for such a Higgs bundle to yield an associated almost-complex curve in $\quadric$. The resulting curve $\nu:\tilde{\Sigma} \rightarrow \quadric$ is immersed if and only if the associated representation is Hitchin (when $\mathcal{T} \cong \mathcal{K}$), but always has pointwise non-vanishing $\sff$, evidenced by $\alpha$. See \cite[Theorem 3.24]{CT24} for further details on the complex Frenet framing of such a curve $\nu$ and the relation with the holomorphic differentials $\alpha=1, \beta,\delta$.\footnote{We caution the reader that $\alpha,\beta$ have roles reversed in \cite{CT24} as they do here.}
\end{remark}

We now describe the stability conditions on such Higgs bundles in the sense of Definition \ref{Defn:StabilitySL}. 
The argument here is similar to that of \cite[Proposition 5.6]{CT24}, with one additional technical complication. 

\begin{proposition}[Stability of $\alpha$-bundles]\label{Prop:AlphaStability}
Let $\mathcal{H} =(\V, \Phi)$ be an $\alpha$-cyclic Higgs bundle.
\begin{enumerate}[noitemsep]
    \item If $\mathcal{H}$ is polystable, then $-2g+2\leq \deg(\mathcal{T}) \leq 2g-2$. 
    \item (Generic locus). If $\delta\neq 0$ and $\beta \neq 0$, then $\mathcal{H}$ is polystable if and only if $-2g+2\leq \deg(\mathcal{T}) \leq 2g-2$. Moreover, $\mathcal{H}$ is stable except if $\deg(\mathcal{T}) \neq -g+1$ and the sections $\beta,\delta$ have the same divisor, in which case it is strictly polystable. 
    \item ($\delta=0$ locus). If $\delta =0$ and $\beta\neq 0$, then $\mathcal{H}$ is polystable if and only if and $-g+1<\deg(\mathcal{T}) \leq 2g-2$, in which case $\mathcal{H}$ is stable.  
    \item ($\beta=0$ locus). If $\beta=0$ and $\delta\neq 0$, then $\mathcal{H}$ is polystable if and only if $-2g+2\leq \deg(\mathcal{T}) < -g+1$, in which case $\mathcal{H}$ is strictly polystable. 
    \item ($\beta,\delta=0$ locus). If $\beta,\delta=0$, then $\mathcal{H}$ is polystable if and only if $\deg(\mathcal{T})=-g+1$, in which case it is strictly polystable.
\end{enumerate}
\end{proposition}

\begin{remark} 
In every connected component of the space of polystable $\alpha$-cyclic Higgs bundles, there is a Hodge bundle with $\delta=0$ or $\beta=0$.
\end{remark}

\begin{proof}
We first make two observations, useful in all cases. Note that $0 \neq \beta \in H^0(\mathcal{K}\mathcal{T} ^{-1}) $ implies $\deg(\mathcal{K}\mathcal{T}^{-1}) \geq 0$, so that  $\deg(\mathcal{T}) \leq \deg(\mathcal{K}) = 2g-2$. Similarly, $0 \neq \delta \in H^0(\mathcal{T}^3\mathcal{K}^3)$ implies that $\deg(\mathcal{T}) \geq -2g+2$.

Let us consider case (3). Now, in this case, the Higgs bundle is 7-cyclic. For stability, the condition $\mathcal{T}^{-2}\mathcal{K}^{-1} \in \ker(\Phi)$ implies $\deg(\mathcal{T}) > -\frac{1}{2}\deg(\mathcal{K})=-g+1$. The other stability considerations, namely $\deg(\bigoplus_{i=-3}^{j}\mathcal{L}_i)<0$, for $j <3$, yield only weaker demands. Claim (3) follows. \medskip 

Let us now consider case (4), which is similar to (3). We observe that $\V$ decomposes as the following sum of three degree zero Higgs sub-bundles:
\[ \big(\begin{tikzcd}
\mathcal{T}^{-2}\mathcal{K}^{-1} \arrow[r, "\delta"]&\mathcal{T}\mathcal{K}\arrow[r, "1"]&\mathcal{T}
\end{tikzcd}\big) \oplus
\big(\begin{tikzcd} \mathcal{T}^{-1} \arrow[r, "1"]&\mathcal{T}^{-1}\mathcal{K}^{-1}\arrow[r, "\delta"] &\mathcal{T}^2\mathcal{K}
\end{tikzcd}\big) \oplus \begin{tikzcd}\mathcal{O}
\end{tikzcd}\]
The stability condition $\deg(\mathcal{T}^2\mathcal{K}) <0$ for the middle Higgs sub-bundle implies the stability of all three Higgs sub-bundles. We conclude point (4). 
\medskip 

Let us consider case (2). Now, suppose $\beta, \delta \neq 0$. 
In this case, $\mathcal{H}$ is 6-cyclic, with holomorphic sub-bundles $(V_i)_{i=1}^6$, where
\[V_1 = (\mathcal{L}_3\oplus \mathcal{L}_{-3}),  \;V_2= \mathcal{L}_2, \;V_3= \mathcal{L}_1, \; V_4=  \mathcal{O}, \; V_5= \mathcal{L}_{-1},\; V_6 =  \mathcal{L}_{-2}.  \]
By \cite{Sim09}, it suffices to check the stability condition only for $\Phi$-invariant holomorphic sub-bundles compatible with the cyclic splitting. Now, such a $\Phi$-invariant, holomorphic, proper sub-bundle $\mathcal{V} = \bigoplus_{i=3}^{-3} \mathcal{V}_i$, such that $\mathcal{V}_i \subset V_i$ must satisfy $\mathcal{V} \subset \ker (\Phi|_{V_1})$. The remainder of the argument examines this possibility. 
Observe that $\Phi$ restricts to a map $\eta:=\Phi|_{V_1}: V_1 \rightarrow \mathcal{T}\mathcal{K}\otimes \mathcal{K}$ given by $\eta = (\beta, \delta)$. 

The bundle $\mathcal{V}$ is either trivial or a line subbundle $\mathcal{V}\subset \mathcal{L}_3\oplus \mathcal{L}_{-3}$ such that $\delta_{\mid \mathcal{V}}=-\beta_{\mid \mathcal{V}}$. The projection from $\mathcal{V}$ to $\mathcal{L}_{3}$ and $\mathcal{L}_{-3}$ define holomorphic sections $s_+ \in  H^0(\mathcal{V}^{-1}\mathcal{L}_{3})$ and $s_- \in H^0(\mathcal{V}^{-1}\mathcal{L}_{-3})$ that are each  non-zero since $\delta,\beta\neq 0$. Let $D_+$ and $D_-$ be their respective divisors. Observe that the degree of $\mathcal{V}$ is given simultaneously by: 
\[ \deg(\mathcal{L}_3)-\deg(D_+)= \deg(\mathcal{V})=\deg(\mathcal{L}_{-3})-\deg(D_-). \] 
For example, to see the first equality, if $s_{\mathcal{V}}$ a nonzero meromorphic section of $\mathcal{V}$, then $s_+\circ s_{\mathcal{V}}$ is a nonzero meromorphic section of $\mathcal{L}_3$ and hence 
$ \deg(\mathcal{L}_3) = \deg\big(\mathrm{div}(s_+\circ s_{\mathcal{V}})\big)=\deg(D_+)+\deg(\mathcal{V})$. 

Now, since $s_+, s_-$ are holomorphic,  $\deg(D_+),\deg(D_-)\geq 0$. We now consider some cases based on $\deg(\mathcal{T})$, keeping in mind the bounds $-2g+2 \leq \deg(\mathcal{T}) \leq 2g-2$ found earlier.\medskip  

\textbf{Case 2(a)}. If $\deg(\mathcal{T})>-g+1$, then $\deg(\mathcal{L}_{-3})<0$, so $\deg(\mathcal{V})<0$. Hence, $\mathcal{H}$ is stable. 

\textbf{Case 2(b)}. If $\deg(\mathcal{T})<-g+1$, then $\deg(\mathcal{L}_3)<0$ and hence $\deg(\mathcal{V})<0$. Hence, $\mathcal{H}$ is stable. 

\textbf{Case 2(c)}. Now assume that $\deg(\mathcal{T})=-g+1$, hence $\deg(\mathcal{L}_3)=\deg(\mathcal{L}_{-3})=0$. In this case, $\deg(\mathcal{V})\leq 0$, and we have two further possibilities.

\emph{2c(i)} $D_+$ or $D_-$ is non-trivial, which implies $\deg(\mathcal{V})<0$ and stability. 

\emph{2c(ii)} $D_+=D_-=0$ are trivial. In this case, one can identify $\mathcal{L}_3$ and $\mathcal{L}_{-3}$ by $s_-\circ s_+^{-1}$ and under this identification $\beta+\delta=0$ so $\beta=-\delta$. 
Note that $\deg(\mathcal{V})= \deg(\mathcal{L}_3)=0$.  
Let us construct another line bundle $\mathcal{V}^\bot$ by :
\[ \mathcal{V}^\perp=\lbrace s_+(z)-s_-(z)\mid z\in \mathcal{V}\rbrace\subset \mathcal{L}_{3}\oplus \mathcal{L}_{-3}.\]
We have decomposed $\mathcal{L}_3\oplus \mathcal{L}_{-3} = \mathcal{V} \oplus \mathcal{V}^\bot$ as a holomorphic direct sum. 

As in the end of the proof of \cite[Proposition 5.6]{CT24}, the Higgs bundle decomposes as follows:
\[
\begin{tikzcd} 
\mathcal{V}^\perp \arrow[r, "\eta|_{\mathcal{V}^\bot}"] &\mathcal{T}\mathcal{K}\arrow[r, "1"]&\mathcal{T}\arrow[r, "-i\sqrt{2}\beta"]&\mathcal{O} \arrow[r, "-i\sqrt{2}\beta"]&\mathcal{T}^{-1} \arrow[r, "1"]&\mathcal{T}^{-1}\mathcal{K}^{-1}\arrow[bend left, lllll, "(\beta\text{,}\delta)"]
\end{tikzcd} \;\; \oplus
\begin{tikzcd}
\mathcal{V}
\end{tikzcd}\]
Hence, $\mathcal{H}$ is strictly polystable in case 2(c)(ii).

\medskip
Point (1) is a consequence of the argument as a whole. \end{proof}

\begin{remark}[Special $\alpha$-bundles]\label{Remk:SpecialAlphaBundles}
We record here some noteworthy cases of polystable $\alpha$-cyclic Higgs bundles:
\begin{itemize}[noitemsep]
     \item If $\mathcal{T} \cong \mathcal{K}^{-1/2}$, then $\beta, \delta \in H^0(\mathcal{K}^{3/2})$. 
        \begin{itemize}
            \item If $\beta =\delta=0$, then the corresponding representation factors through the $\SL(2,\R)$-subgroup of the long root $\alpha$; see Appendix \ref{Subsec:SL2LongRoot}. The uniformizing Higgs bundle is 
            \[\mathcal{O} \;\;\;\;\;\;\;\mathcal{K}^{1/2} \stackrel{1}{\longrightarrow} \K^{-1/2} \;\;\;\;\;\;\; \mathcal{O} \;\;\;\;\;\;\;\;\;\mathcal{K}^{1/2}\stackrel{1}{\longrightarrow}\mathcal{K}^{-1/2}\;\;\;\;\;\;\;\mathcal{O}. \]
            \item If $\beta = \delta \neq 0$, 
            then the corresponding representation factors through $\SL(3,\R)$ as the deformations of Barbot representations studied in \cite{DB25}.
        \end{itemize}
    \item $\deg(\mathcal{T}) =2g-2$ if and only if  $\mathcal{T}\cong \mathcal{K}$ if and only if the corresponding representation is $\Gtwosplit$-Hitchin.
    The `if' statement come from Hitchin's parametrization. The `only if' statement is seen as follows:
    $\deg(\mathcal{T})=2g-2$ forces $\beta\neq 0$ by Proposition \ref{Prop:AlphaStability}. Then $0\neq \beta \in H^0(\K \mathcal{T}^{-1})$ causes the degree zero bundle $\K\mathcal{T}^{-1}$ to be trivial and $\beta$ to be non-vanishing. 
    \item If $\mathcal{T}=\K$, $\beta =1$, and $\delta=0$, the corresponding representation is $\Gtwosplit$-Fuchsian-Hitchin. 
    \item If $\beta=0$, then the corresponding representations factor through $\SU(2,1)$.
    \item If $\mathcal{T} \cong \mathcal{K}^{-1}$, then $\delta \in H^0(\mathcal{O})$ and $\beta \in H^0(\mathcal{K}^2)$. 
        \begin{itemize}
            \item If $\delta =1$ and $\beta =0$, then the corresponding representation is a Fuchsian representation factoring through the $\PSL(2,\R)$-subgroup in Appendix \ref{Appendix:SL2Number5}. The uniformizing Higgs bundle is:
    \[  \begin{tikzcd}
    \mathcal{K}^{-1} & \mathcal{O}\arrow[r, "1"]&\mathcal{K}^{-1} &\mathcal{O} &\mathcal{K} \arrow[r, "1"]&\mathcal{O} \arrow[bend left, lllll, "1"]&\mathcal{K}\arrow[bend left, lllll, "1"]
    \end{tikzcd}. \]
        \end{itemize}
\end{itemize}
\end{remark}

\subsection{\texorpdfstring{$\Pho^\times$}{Phox}-Structures for \texorpdfstring{$\alpha$}{a}-Hodge Bundles}\label{Sec:PhoxStructuresAlphaHodge}  

Recall that in $\S$\ref{Sec:PhoxFibersHitchin}, we showed the fibers of the Tits metric thickening domain in $\Pho^\times$ are $\mathbb{RP}^2\times \mathbb{S}^1$ for $\Gtwosplit$-Hitchin representations. This means that there is an $(\RP^2\times \sphere^1)$-fiber bundle over $S$ that carries $(\Gtwosplit, \Pho^\times)$-geometries with Hitchin holonomy. In this section, we shall build this manifold explicitly, and find such fiber bundles in greater generality.
In particular, we build fibered geometric structures on $(\RP^2\times \sphere^1)$-fiber bundle over $S$ from stable $\alpha$-cyclic $\Gtwosplit$-Higgs bundles that are also \emph{Hodge bundles}. 

Fix $\Sigma = (S,J)$ a Riemann surface on $S$. 
Now $\alpha$-Hodge bundles obtain the following form, where either $\delta=0$ or $\beta=0$: 
 \begin{equation}\label{alpha2HodgeBundle}
\begin{tikzcd}
\mathcal{T}^2\mathcal{K} \arrow[r, "\beta"] & \mathcal{T}\mathcal{K}\arrow[r, "1"]&\mathcal{T}\arrow[r, "-\sqrt{2}i \beta"]&\mathcal{O} \arrow[r, "-\sqrt{2}i\beta"]&\mathcal{T}^{-1} \arrow[r, "1"]&\mathcal{T}^{-1}\mathcal{K}^{-1}\arrow[lllll,bend left, "\delta"] \arrow[r, "\beta"]&\mathcal{T}^{-2}\mathcal{K}^{-1}\arrow[lllll,bend left, "\delta"]
\end{tikzcd}.
\end{equation}
The stability of such Higgs bundles is described by (3) of Proposition \ref{Prop:AlphaStability}. 

The Higgs field $\Phi$, and the section $\Phi_0$ of $\End(\mathcal{E})\otimes \mathcal{K}$ that we will call the \emph{$\alpha$-Fuchsian} part of the Higgs field are described by the following matrices:

$$\Phi=\begin{pmatrix}
    0& 0 &0&0&0&\delta&0\\
     \beta& 0 &0&0&0&0&\delta\\
    0& 1&0&0&0&0&0\\
    0& 0 &-i\sqrt{2} \beta&0&0&0&0\\
    0& 0 &0&-i\sqrt{2} \beta&0&0&0\\
    0& 0 &0&0&1&0& 0\\
    0& 0 &0&0&0& \beta&0
       \end{pmatrix},\;
       \Phi_0=\begin{pmatrix}
    0& 0 &0&0&0&0&0\\
     0& 0 &0&0&0&0&0\\
    0& 1&0&0&0&0&0\\
    0& 0 &0&0&0&0&0\\
    0& 0 &0&0&0&0&0\\
    0& 0 &0&0&1&0& 0\\
    0& 0 &0&0&0& 0&0
       \end{pmatrix}.$$
In particular, the sub-endomorphism $\Phi_0$ is the component of $\Phi$ in the simple root space $e_{-\alpha}$.

\subsubsection{Constructing the geometric structures} \label{Subsec:PhotoStructuresHodge} 

We keep the running notation from the previous section. 
Let $\rho: \pi_1S\rightarrow \Gtwosplit$ denote the associated representation to an $\alpha$-cyclic Hodge bundle. 
We construct an $(\RP^2 \times \mathbb{S}^1)$-fiber bundle $M_{\rho} \rightarrow S$, upon which there will be a fibered $(\Gtwosplit, \Pho^\times)$-structure with holonomy that descends to $\pi_1S$ as $\rho$. 
This geometric structure is constructed as in Section \ref{Subsec:Ein23Structures}.

We first describe in two ways the 5-manifold which carries the fibered $(\Gtwosplit,\Pho^\times)$-structures. 
We write $B_{\Psi_0}$, with $B$ to stand `base of pencil', for the more complicated space, which carries a more natural developing section to the flat $\Pho^\times$-bundle. 
On the other hand, we denote $M_{\Psi_0}$, with $M$ to stand for `model', for a geometrically simpler space, which is homeomorphic to $B_{\Psi_0}$. 

Let us first define the key object $\Psi_0$. As in $\S$\ref{Subsec:Ein23Structures}, this object can be built either from the Higgs bundle or the associated harmonic map. We describe from the latter perspective first. 
Recall the $\rho$-equivariant maps $f:\tilde{\Sigma} \rightarrow \X$, and its refinement, the Frenet frame $\mathscr{F} = (\mathscr{L},T,N,B)$. 
In the case, $\beta \neq 0$, the map $\mathscr{F}$ is the Frenet frame of a $J$-holomorphic curve
$\nu: \tilde{\Sigma} \rightarrow \quadric$, which is not necessarily immersed, but second fundamental form $\sff \in \Omega^1(\tilde{\Sigma},\Hom(T,N))$ that is non-vanishing. 

First, we reconsider the trivial flat bundle $\underline{\R^{3,4}} = \tilde{\Sigma}\times \R^{3,4}$ with its canonical trivial connection $D$. Under the Frenet frame splitting $(\mathscr{L}, T,N,B)$, this connection decomposes as 
\[ D = \begin{pmatrix} \nabla_{\mathscr{L}} & -\mathrm{I}^* & & \\
\mathrm{I} & \nabla_T & -\sff^* & \\
 &  \sff & \nabla_{N}& -\tff^* \\
 &  &  \tff&\nabla_{B} \end{pmatrix}. \]
The object entitled $\Psi_0 \in \Omega^1(\tilde{\Sigma}, \End(\underline{\R^{3,4}}))$, from this point of view, is given by  
\[ \Psi_0 = \begin{pmatrix} 0&  & & \\
 &0 & -\sff^* & \\
 &  \sff & 0&  \\
 &  &  &0 \end{pmatrix}. \]
Observe that for $p \in \tilde{\Sigma}$, setting $U=f(p)$, then $\Psi_0(X)$ restricts pointwise a map $\Psi_0(X)|_p:U\rightarrow U^\bot$. Hence, we may equivalently view $\Psi_0 \in \Omega^1(\tilde{\Sigma}, f^*\T\X)$. 
Note that for every $X \in \T_p\tilde{\Sigma}$, the tangent vector $\Psi_0(X)|_p:U\rightarrow U^\bot \in \T_U\X$ has the same one-dimensional kernel and two-dimensional image. By $\S$\ref{Subsec:Pho^XSimplePencil}, we see that the pencil parametrized by $\Psi_0$ is pointwise $\Gtwosplit$-equivalent to the tangent pencil of the sub-symmetric space $\Ha^2_{\alpha}$. 

With all of this said, the $\pi_1S$-cover of our $(\Gtwosplit,\Pho^\times)$-manifold, denoted $\overline{B}_{\Psi_0}$, is the $\alpha$-base of $\Psi_0$: 
\begin{align}\label{BPsi0Cover}
     \overline{B}_{\Psi_0}= \big\lbrace \,(p,\omega) \in \tilde{\Sigma} \times \Pho^\times \mid  \, \omega\in \mathcal{B}_\alpha(\Psi_0|_p) \big\rbrace 
\end{align}
Note that there is a tautological map $\dev: \overline{B}_{\Psi_0} \rightarrow \Pho^\times$ by $(p,\omega)\mapsto \omega$. We will study this map in Theorem \ref{Thm:DevMapPhoX}, verifying that it is a local diffeomorphism, which furnishes a $\Pho^\times$-structure on $B_{\Psi_0}$. 

\begin{remark}
As in $\S$\ref{Subsec:Ein23Structures}, the description of $\Psi_0$ depends on which model for $\T_U\X_{\Gtwosplit}$ used. The two models are $\T_U\X_{\Gtwosplit} \cong \mathfrak{p}(U)$, where $\g_2'=\mathfrak{k}(U) \oplus \mathfrak{p}(U)$ is the associated Cartan decomposition, and $\T_U\X_{\Gtwosplit} \cong \Hom^\times(U,U^\bot)$. The map $A \mapsto A-A^{*q}$ identifies $\Hom^\times(U,U^\bot)\cong \mathfrak{p}(U)$. In particular, the pencil $\Psi_0$ in the latter model is viewed as $\sff-\sff^*$, while in the former it is just $-\sff^*$, by slight abuse, where we extend $\sff^*|_p$ to a map on $\mathscr{L}_p \oplus N_p$ by declaring $\mathscr{L}_p$ to be in the kernel.
\end{remark}

We can now translate this whole picture to the Higgs bundle. 
The object $\Psi_0 \in \Omega^1(\tilde{\Sigma}, f^*\T\X)$ descends by $\rho$-equivariance to give an object of the form $\Psi_0 \in \Omega^1(\Sigma, \End(\V^\R))$. We can express this object in terms of the Higgs bundle: it is written $\Psi_0=\Phi_0+\Phi_0^*$ where $\Phi_0$ is again the $\alpha$-Fuchsian part of the Higgs field shown below.
\begin{equation}
\begin{tikzcd}
\mathcal{T}^2\mathcal{K}  & \mathcal{T}\mathcal{K}\arrow[r, "1"]&\mathcal{T}&\mathcal{O} &\mathcal{T}^{-1} \arrow[r, "1"]&\mathcal{T}^{-1}\mathcal{K}^{-1} &\mathcal{T}^{-2}\mathcal{K}^{-1}
\end{tikzcd}.
\end{equation}
We can form a moving base of pencil of $\Psi_0$ then in the Higgs bundle picture, which gives a geometric description of the $\pi_1S$-quotient of $B_{\Psi_0} = \pi_1S \backslash \overline{B}_{\Psi_0}$. Indeed, for $\Psi_0 \in \Omega^1(\Sigma, \End(\V^\R))$, we define $B_{\Psi_0}$ as a fiber sub-bundle of the associated fiber bundle $\Pho^\times(\V^\R) \rightarrow S$:  
\[ B_{\Psi_0} = \{ (p, \omega) \in \Pho^\times(\V^\R) \mid \omega \in \mathcal{B}_{\alpha}(\Psi_0|_p) \}.\]

\vspace{1ex}

We will now describe $B_{\Psi_0}$ as an $(\RP^2 \times \mathbb{S}^1)$-bundle over $S$, using the explicit description of bases of pencils from Proposition \ref{Prop:PointingTowardsPhotonPrelim}. Here, we must introduce quite a bit of notation.
We will then produce a homeomorphism $\mathcal{D}: M_{\Psi_0} \rightarrow B_{\Psi_0}$ of fiber bundles that elucidates the topology of $B_{\Psi_0}$.  

In terms of the Higgs bundle Frenet frame $(\mathscr{L},T,N,B)$ from \eqref{LTNB}, recall the real sub-bundle $\mathscr{U} = \mathscr{L} \oplus N$ of $\V^\R$. Consider the $h$-orthogonal splitting $\V^\R = \mathscr{U} \oplus \mathscr{U}^\bot$, whereby $h|_{\V^\R} = Q|_{\mathscr{U}} -Q|_{\mathscr{U}^\bot}$. 
We consider the restriction of $\Psi_0$ to $\mathscr{U}$: namely, $\Psi_{0}|_{\mathscr{U}} \in \Omega^1(\Sigma, \Hom(\mathscr{U}, \mathscr{U}^\bot))$. For $p \in \Sigma$, $X \in \mathrm{T}_U\Sigma$, and any two-dimensional sub-space $W \in \mathscr{U}|_p$, we may consider further the restriction $\Psi_0(p,X)|_W$ as a linear map $\Psi_0(p,X)|_W:W \rightarrow \mathscr{U}^\bot|_p$. As explained before Corollary \ref{Cor:PhoxModel}, the space $\Hom_\R(W, \mathscr{U}^\bot|_p)$ is naturally a complex vector space. 
We shall regard $\Hom_\C(W, \mathscr{U}^\bot|_p)$ as a four-dimensional real sub-bundle of $\Hom_{\R}(W, \mathscr{U}^\bot|_p)$ as we shall be interested in non-complex sub-bundles of the former.  
Next, we place an inner product on $\Hom_\R(W,\mathscr{U}^\bot|_p)$. Define 
$\langle \eta, \mu\rangle := \tr(\eta \circ \mu^{*h})$. It then makes sense to consider then $\Psi_0(p,X)|_W^\bot \subset \Hom(W, \mathscr{U}_p^\bot)$ in the inner product space $\Hom_\R(W, \mathscr{U}_p^\bot)$.

We now describe the fibers of $B_{\Psi_0}$ a bit more concretely: 
\begin{align}\label{BPsi0Explicit}
     B_{\Psi_0}|_{p} \simeq \{ \,(W,L) \in \Gr_2(\mathscr{U}|_p)\times Q_+(\Hom_\C(W, \mathscr{U}^\bot|_p)) \; | \; L \in (\Psi_0(p,X)|_W)^\bot , \, \forall X \in \mathrm{T}_U\Sigma\}.
\end{align}
The above smooth description of $B_{\Psi_0}$ is a bundle analogue of Corollary \ref{Cor:PhoXBaseCircleBundle}, which now shows that $B_{\Psi_0}$ is a smooth (orientable) circle bundle over the total space of the bundle $\Gr_2(\mathscr{U}) \rightarrow S$. More precisely, we can define a complex line bundle 
$\mathscr{R}^{\bot} \rightarrow \Gr_2(\mathscr{U})$ with 
fiber 
\[ \mathscr{R}^{\bot}|_{(p,W)}= \{L \in \Hom_{\C}(W, \mathscr{U}_{\mid p}^\bot) \mid L \in \Psi_0(p)|_{W}^\bot \} ,\]
and then \eqref{BPsi0Explicit} reads 
as a diffeomorphism as follows: 
\[ Q_+(\mathscr{R}^\bot) \rightarrow B_{\Psi_0}, \qquad (W,L)\mapsto \graph(L). \]

We can give another description of the bundle $B_{\Psi_0}$, using ideas from $\S$\ref{Sec:PhoxFibersHitchin}. Applying the identification in Lemma \ref{Lemma:SimplifiedDevinPho} fiberwise, we obtain the following canonical fiber bundle homeomorphism:
\begin{align}\label{eq:D-Identification}
 \mathcal{D}: M_{\Psi_0}:=\Gr_2(\mathscr{U}) \oplus Q_+(\Hom_\C(N,B)) \rightarrow B_{\Psi_0},
\end{align}

which describes the global topology of $B_{\Psi_0}$ as a fiber bundle. Indeed, since $\mathscr{U} \rightarrow S$ is smoothly trivial, then $\Gr_2(\mathscr{U})\cong S \times \RP^2$. Moreover, the identifications of Lemma \ref{Lemma:SimplifiedDevinPho} show that the bundle $\mathscr{R}^\bot \rightarrow \Gr_2(\mathscr{U})$ is $(C^0)$ isomorphic to the pullback of the complex line bundle $\Hom_{\C}(N, B) \rightarrow \Sigma$. 
Hence, the topology of $B_{\Psi_0}$ is directly controlled by $\deg(\Hom_{\C}(N,B))=\deg(\Hom(\mathcal{T}^{3}\mathcal{K}^{3}))$, which is controlled by $\deg(\mathcal{T})$.  
We explore the consequences of the topology of $B_{\Psi_0}$ further in \cite{DE26Anosov}, by relating the varying topologies to the connected components of $\alpha$-Anosov representations in the $\Gtwosplit$-character variety. 

\begin{remark}\label{Remk:TopologyAlpha}
    The map $\mathcal{D}$ is only a homeomorphism of circle bundles. However, both of $M_{\Psi_0},B_{\Psi_0}$ are smooth orientable circle bundles over $\Gr_2(\mathscr{U}|_p)$, which means that their diffeomorphism type is characterized by their Euler number (orientable circle bundles can always be given the structure of principal circle bundles). These numbers agree because of the homeomorphism, hence these spaces are diffeomorphic as fiber bundles over $S$.
\end{remark}

With our geometric descriptions of $B_{\Psi_0}$ complete, we now proceed towards the proof that the map $\dev: B_{\Psi_0}\rightarrow \Pho^\times$ is a local diffeomorphism. For the proof, we will use the fact that $\Psi$ is not too far from $\Psi_0$, or in other words that the norms of $\beta$ and $\delta$ are controlled. This is a consequence of a maximum principle very similar to the one for $\beta$-cyclic bundles in Lemma \ref{Lem:MaximumPrincipleBeta}.

\begin{lemma}[Maximum Principle for $\alpha$-Hodge Bundles]\label{Lem:MaximumPrincipleAlpha}
    Let $(\V,\Phi)$ be an $\alpha$-cyclic bundle on $\Sigma$ with $\delta=0$. Then we have the global inequality \[\lVert \beta\rVert_h \leq  \sqrt{\frac{3}{5}} \lVert \alpha \rVert_h. \]
Let $(\V,\Phi)$ be an $\alpha$-cyclic bundle on $\Sigma$ with $\beta=0$. Then we have the global inequality \[\lVert \delta\rVert_h \leq  \lVert \alpha \rVert_h. \]
\end{lemma}

\begin{proof}
Let $g$ be a conformal metric on $\Sigma$. 
First, we apply the general equations from  Lemma \ref{Lem:HitchinsEquationsGeneral} in the case of $\alpha$-cyclic bundles with $\delta=0$, obtaining the following:
\begin{align}
     \Delta_g \log\lVert \alpha \rVert^2&=2\lVert \alpha \rVert ^2-3\lVert \beta \rVert^2+ \kappa_g \\
    \Delta_g \log\lVert \beta \rVert^2&=2\lVert \beta\rVert ^2-\lVert \alpha\rVert^2+ \kappa_g.
\end{align}

We now consider a point $x_0 \in S$ where the ratio $\left(\frac{\lVert \beta\rVert}{\lVert \alpha \rvert}\right)$ attains its maximum. Such a point exists since the tautological section $\alpha$ does not vanish. At a local maximum, the Laplacian is non-positive. Now, $\left(\frac{\lVert \alpha \rVert}{\lVert \beta \rVert}\right)(x_0) >0$ since $\beta$ is not identically zero.  Thus,
\[0 \geq2 \Delta_g \log\left(\frac{\lVert \beta \rVert}{\lVert \alpha \lVert}\right)(x_0) =5\lVert \beta \rVert^2(x_0)-3\lVert \alpha \rVert^2(x_0).\]
In particular $\left(\frac{\lVert \beta\rVert}{\lVert \alpha \rVert}\right) (x_0)\leq\sqrt{\frac{3}{5}}$. As $x_0$ is a maximum, this inequality holds everywhere.\medskip

Similarly, we apply the general equations from  Lemma \ref{Lem:HitchinsEquationsGeneral} in the case of $\alpha$-cyclic bundles with $\beta=0$, obtaining the following:
\begin{align}
     \Delta_g \log\lVert \alpha \rVert^2&=2\lVert \alpha \rVert ^2-\lVert \delta \rVert^2+ \kappa_g \\
    \Delta_g \log\lVert \delta \rVert^2&=2\lVert \delta\rVert ^2-\lVert \alpha\rVert^2+ \kappa_g.
\end{align}
In particular, $\Delta_g\log(\frac{\lVert \delta \rVert^2}{\lVert \alpha \rVert^2})=2\lVert \delta \rVert^2-2\lVert \alpha \rVert^2$. Hence, by the same argument considering a global maximum of $\log(\frac{\lVert \delta \rVert^2}{\lVert \alpha \rVert^2})$, we obtain that $\lVert \delta\rVert_h \leq  \lVert \alpha \rVert_h$.
\end{proof}

Let us write $h = \diag(h_i)_{i=3}^{-3}$ on $\V = \bigoplus_{i=3}^{-3} \mathcal{L}_i$. Just as with $\beta$-cyclic bundles, a certain conformal metric $g$ on $\Sigma$ will play a key role. 

\begin{definition}[$\alpha$-Projected Metric]\label{Defn:AlphaMetric}
Define the metric $g: = h_{1}h_{2}^{-1}$, which is a hermitian metric on $\K^{-1} \cong \mathcal{L}_{2}^{-1}\otimes \mathcal{L}_{1}$ by \eqref{alpha2HodgeBundle} and hence equivalently a Riemannian metric on $S$. 
\end{definition}
The metric $g$ is a conformal metric on $\Sigma$. 
The following is the analogue of Proposition \ref{Prop:ParallelismPsi0Alpha1}. 

\begin{proposition}[$\Psi_0$-parallelism for $\alpha$-cyclic Bundles]\label{Prop:ParallelismPsi0Alpha2}
Let $\nabla^h,\nabla^g$ be the Chern connections on $\End(\V)$ and $\K^{-1}\cong \T\Sigma$, respectively, where $g=g_{\alpha}$. 
Take $X \in \Omega^0(\Sigma, \T\Sigma)$ to be any vector field.
Then $\Psi_0$ has the following parallelism property:  
\[ (\nabla^h \circ \Psi_0)(X) = (\Psi_0\circ \nabla^g)( X). \]  
\end{proposition}

The proof of Proposition \ref{Prop:ParallelismPsi0Alpha2} is directly analogous to that of Proposition \ref{Prop:ParallelismPsi0Alpha1}, now using that $h$ induces the same metric $g$ on $\Hom(\mathcal{L}_2,\mathcal{L}_1)$ and on $\Hom(\mathcal{L}_{-1}, \mathcal{L}_{-2})$.\medskip

The main result of this section is the verification that the tautological map, suggestively denoted $\dev$, is a developing map for $B_{\Psi_0}$. 
Since the proof is rather long, we provide a summary here:
\begin{itemize}
    \item \textbf{Step 0}: Fix $p_0 \in \tilde{\Sigma}$ as well as a point $x_0:=(p_0,\omega_0)\in{\overline{B}_{\Psi_0}}$. We recall that $\dev$ is an injective immersion on the fiber over $p_0$ by the base of pencil construction. 
    Now, fix a tangent vector $X_0 \in \mathrm{T}_{p_0}\tilde{\Sigma}$, and we consider the quantity $A:=\frac{d}{dt}\big |_{t=0} \langle \Psi_0(X_0), \varphi_{\omega_t, U_0}\rangle_{\X}$, where $U_0 \in \X$, $(\omega_t)_{t \in (-\varepsilon, \varepsilon)}$ is a horizontal and parallel curve in $\overline{B}_{\Psi_0}$ with initial velocity that projects down to $X_0$, and $\varphi_{\omega_t,U_0}$ expresses the annihilator photon $\omega_t$ as a graph over $W_0:=\pi_{U_0}(\omega_0)$ as in Corollary \ref{Cor:AnnihilatorPhotonGraph}. In other words, $\varphi_{\omega_t,U_0} \in \T_{U_0}\X_{\Gtwosplit}$ points towards $\omega_t$ in the visual boundary. 
    The condition $A\neq 0$ for all $X_0 \in \T_p\Sigma$ allows us to conclude $\dev$ is an immersion at $x_0$. 
    \item \textbf{Step 1}: Express the quantity $A$ in terms of $\Psi_0(X_0)$, $\Psi(X_0)$, and $\omega_0$ only. Here, parallelism of $\Psi_0$ from Proposition \ref{Prop:ParallelismPsi0Alpha2} is critical.
    \item \textbf{Step 2}: Express the quantity $A$ in terms of local coordinates for the fiber $(B_{\Psi_0})|_p$. Here, we need the understanding of the $\alpha$-base $\mathcal{B}_{\alpha}(\mathcal{P}_0)$ of the tangent pencil $\mathcal{P}_0$ to $\Ha^2_{\alpha}$ from $\S$\ref{Subsec:Pho^XSimplePencil}. 
    \item \textbf{Step 3}: Prove that $A> 0$. Here, we perform elementary algebra and inequalities and eventually win by verifying certain polynomials in $\Z[x]$ have no real roots. 
\end{itemize}

Now, we proceed with the theorem. 

\begin{theorem}[$\Pho^\times$-structures for $\alpha$-Hodge bundles]\label{Thm:DevMapPhoX}
Let $(\V, \Phi)$ be a polystable $\alpha$-cyclic bundle with either $\beta=0$ or $\delta=0$. 
Then the tautological map $\dev: \overline{B}_{\Psi_0} \rightarrow \Pho^\times$ is a local diffeomorphism. 
\end{theorem}

\begin{proof}
\textbf{Step 0: Setup.} Fix any point $x_0=(p_0,\omega_0 )\in \overline{B}_{\Psi_0}$. 
Denote by $\mathcal{F}_{p_0}$ the developed image of the fiber $\overline{B}_{\Psi_0}|_{p_0}$. Now, to prove $\dev$ is an immersion, it suffices to show the differential $d\dev_{x_0}$ surjects.  By construction, $\dev$ is an immersion when restricted to the fiber $\overline{B}_{\Psi_0}|_{p_0}$. 
Thus, 
we need only prove $d\dev_{x_0}$ surjects onto the quotient space $H_{\omega_0} :=\T_{\omega_0}\Pho^\times/\T_{\omega_0}\mathcal{F}_{p_0}$. 
The following procedure defined next will allow us to prove this. 

Recall the metric $g$ from Definition \ref{Defn:AlphaMetric}. 
For any tangent vector $X_0 \in \T_{p_0}\tilde{\Sigma}$, take $\gamma_b:(-\varepsilon, \varepsilon) \rightarrow \tilde{\Sigma}$ to be the $g$-geodesic with $\dot{\gamma}_b(0) = X_0$. Let us express $x_0 \in\Pho^\times$ in the usual Stiefel triplet description. 
Fix a background choice of $u_0 \in Q_+(N), v_0 \in Q_+(U \cap u_0^\bot), z_0 \in Q_-(\hat{\mathcal{R}}_{u_0, v_0}^{\bot})$ such that $\omega_0 = \spann\langle u_0+z_0, v_0+(u_0v_0)z_0 \rangle \in \Pho^\times$. 
Next, we choose the unique $ \nabla^h$-parallel curve $\gamma_{t}: (-\varepsilon , \varepsilon) \rightarrow (\V^{\R})^3$ along $\gamma_b$ such that $\gamma_0=(u_0, v_0, z_0)$. 
Write $\gamma_t=(u_t, v_t,z_t)$. Since $\nabla^h$ preserves the Frenet frame splitting, 
$u_t \in N_{\mid\gamma_b(t)}$, $v_t \in \mathscr{U}_{\mid\gamma_b(t)}$, and $z_t \in \mathscr{U}^\bot_{\mid\gamma_b(t)}$ for all $t$. 
Define $\omega_t =\spann_{\R} \langle u_t+z_t, v_t+(u_tv_t)z_t\rangle $. Since $\nabla^h$ is a $\Gtwosplit$-connection, $(u_t,v_t,z_t) \in V_{(+,+,-)}(\V^{\R})$ and hence $\omega \in \Pho^\times $ by Corollary \ref{Cor:AnnihilatorPhotonGraph}. We claim that, in fact, $\omega_t $ is a section of $\overline{B}_{\Psi_0}$ along $\gamma_b$. 
Here, we recall that any annihilator photon $\omega \in \Pho^\times $ satisfies $\omega \in \overline{B}_{\Psi_0}|_{p}$ if and only if 
\[ \langle \Psi_0(X), \varphi_{\omega, U_0} \rangle_{\X} = 0 , \; \forall X \in \mathrm{T}_U\tilde{\Sigma}, \] 
using the notation $\varphi_{\omega,U_0}$ from Corollary \ref{Cor:AnnihilatorPhotonGraph}, and writing $U_0=f(p_0)$ for $f: \tilde{\Sigma} \rightarrow \X_{\Gtwosplit}$ the associated harmonic map. That is, $U_0=\mathscr{U}_{\mid {p_0}}$. However, using Lemma \ref{Lem:PhoxFiberConditionSimple}, we see $\omega_t \in \overline{B}_{\Psi_0}|_{\gamma_b(t)}$ if and only if 
\begin{align}\label{ParallelCheck}
    \langle \Psi_0(Y)(u_t), z_t\rangle_{h} +\langle \Psi_0(Y)(v_t), w_t\rangle_h=0, \, \forall Y \in \T_{\gamma_b(t)}\tilde{\Sigma},
\end{align}
where we also denote $w_t := (u_tv_t)z_t$. 
To prove \eqref{ParallelCheck}, we will prove that for $\mathscr{Y}$ any $\nabla^g$-parallel vector field along $\gamma_b$ that the function $f_{\mathscr{Y}}:(-\varepsilon, \varepsilon) \rightarrow \R$ given by  
\[ f_{\mathscr{Y}}(t) = \langle \Psi_0(\mathscr{Y}_t)u_t, z_t\rangle_h + \langle \Psi_0(\mathscr{Y}_t)v_t, w_t\rangle_h, \]
is identically vanishing, which then proves \eqref{ParallelCheck}. To check $f_{\mathscr{Y}} \equiv 0$, we need only prove $f_{\mathscr{Y}}' \equiv 0$, since $f_{\mathscr{Y}} (0) = 0$ by $\omega_0\in B_{\Psi_0}$. On the other hand, using $ \nabla^hh=0$, the parallelism property of $\Psi_0$ from Proposition \ref{Prop:ParallelismPsi0Alpha2}, and the $ \nabla^h$-parallelism of $u_t,v_t,z_t$, one finds $f_{\mathscr{Y}}'\equiv 0$ directly. We conclude that $\omega_t \in B_{\Psi_0}$ for all $t \in (-\varepsilon, \varepsilon)$. 

The key step to this proof will be to show the following:
\begin{align}\label{KeyDerivative}
    \frac{d}{dt}\bigg|_{t=0} \langle \Psi_0(X_0), \varphi_{\omega_t,U_0} \rangle_{\X} \neq 0.
\end{align} 
The equation \eqref{KeyDerivative} implies that $\omega_t$ is moving off the developed fiber $\mathcal{F}_{p_0}$. Now, the aforementioned process defines for each $X_0 \in \T_{p_0}\tilde{\Sigma}$ a tangent vector $\hat{X}_0 := \frac{d}{dt}\big|_{t=0}\,\omega_t \in \T_{x_0}\overline{B}_{\Psi_0}$ such that $d\pi_{\tilde{\Sigma}}(\hat{X}_0) = X_0$ and $\pi_{H_{\omega_0}}(d \dev_{x_{0}}(\hat{X})) \neq 0$. By dimension count, this means $d\dev_{x_0}$ surjects onto $H_{\omega_0}$, so $\dev$ is an immersion at $x_0$ as desired. 
The remainder of the proof is to prove  \eqref{KeyDerivative}. 

\textbf{Step 1: Differentiate.} Let us write $\psi_0: = \Psi_0(X_0)$ as well as $\psi := \Psi(X_0)$. Recall by Proposition \ref{Prop:ModelSpaceMetric} that $c\langle\eta,\gamma\rangle_{\X}=-\tr(\eta^{*q}\circ \gamma)$ for some constant $c >0$. The remainder of the proof is to compute 
\[A := \frac{1}{2c}\frac{d}{dt} \big|_{t=0}\langle \varphi_{\omega_t, U_0}, \psi_0\rangle_{\mathbb{X}} \]
We shall need the following formulas.

\begin{lemma}[Fiber Equation, Simple]\label{Lem:PhoxFiberConditionSimple}
Let $\omega \in \Pho^\times$. For $U \in \mathbb{X}_{\Gtwosplit}$ let $\varphi_{\omega, U}: U \rightarrow U^\bot$ be the unique rank two linear map such that $\graph^*(\varphi_{\omega, U})= \omega$ by Proposition \ref{Prop:PointingTowardsPhotonPrelim}. Choose any basis $(x_1, x_2)$ for $\omega$ such that $(\pi_U(x_1), \pi_U(x_2))$ are orthogonal and $q(\pi_U(x_1))= q(\pi_U(x_2))$. Define the Euclidean quadratic form $g_U := q|_U - q|_{U^\bot}$. Then for any 
tangent vector $\psi \in \T_{U} \mathbb{X}_{\Gtwosplit}$,  
\[ \langle \varphi_{\omega, U}, \psi\rangle_{\mathbb{X}} = 0 \Leftrightarrow  \langle x_1, \psi(x_1) \rangle_{g_U} + \langle x_2, \psi(x_2)\rangle_{g_U} = 0.\]
\end{lemma}
Recall that $\graph^*(L:V \rightarrow W)$ denotes $\graph(L|_{\ker(L)^\bot})$. 

\begin{proof}
Let $(x_1, x_2)$ be a basis for $\omega$ satisfying the hypotheses. Write $x_i = u_i +z_i$ uniquely for $u_i \in U$ and $z_i \in U^\bot$. 
By Corollary \ref{Cor:Orthogonality}, 
\[ \langle \varphi_{\omega, U}, \psi \rangle_{\X} =0 \; \Leftrightarrow
\langle u_1, \psi(u_1)\rangle_{q_{3,4}} + \langle u_2 , \psi(u_2) \rangle_{q_{3,4}} =0 \Leftrightarrow \langle u_1, \psi(u_1)\rangle_{g_U} + \langle u_2 , \psi(u_2) \rangle_{g_U}=0 .\] 
The claim follows. 
\end{proof}

Next, we describe a more general version of the lemma using Gram-Schmidt.

\begin{lemma}[Fiber Condition, General]\label{Lem:PhoxFiberConditionGeneral}
Let $\omega \in \Pho^\times$, $U \in \X_{\Gtwosplit}$, and $\psi \in \T_{U}\X_{\Gtwosplit}$. Take $\varphi_{\omega, U}: U \rightarrow U^\bot$ such that $\omega= \graph^*(\varphi_{\omega, U})$. Now, for any basis $(x_1, x_2)$ of $\omega$, we denote $x_i = u_i +z_i$ for $u_i \in U$, $z_i \in U^\bot$. Then we compute that 
\begin{align*}
   -\frac{1}{2} \tr(\psi\circ \varphi_{\omega,U})= f_1^2\langle\psi(u_1), z_1 \rangle_{q}+f_2^2 \bigg( \langle \psi(u_2), z_2\rangle_q + f_3^2 \langle\psi(u_1), z_1\rangle_q -f_3 \langle \psi(u_1), z_2 \rangle_q - f_3 \langle \psi(u_2), z_1 \rangle_q\bigg),
\end{align*}
where $f_1 =q(u_1)^{-1/2},\, \, f_2= q(u_2-f_3u_1)^{-1/2}$, \,$f_3 =\frac{ \langle u_2, u_1 \rangle}{\langle u_1, u_1\rangle}$. 
\end{lemma}

\begin{proof}
Apply Gram-Schmidt to the ordered basis $(u_1,u_2)$ to obtain $(u_1', u_2')$ given by $u_1' = f_1u_1$ and $u_2 ' = f_2(u_2-f_3u_1)$. Write
$z_i'=\varphi_{\omega, U}(u_i').$ This yields another basis $(x_1',x_2')$ for $\omega$ by $x_i'=u_i'+z_i'$. Then Corollary \ref{Cor:Orthogonality} shows that for some $c >0$,  
\[ c\langle\varphi_{\omega,U},\psi\rangle_{\X}= \tr(\psi\circ  \varphi_{\omega, U})= -2\langle \psi (u_1'), z_1'\rangle_{q} -2\langle \psi(u_2'), z_2'\rangle_q. \]
Expanding this equation yields the result. 
\end{proof}

\begin{remark}Lemmas \ref{Lem:PhoxFiberConditionSimple} and \ref{Lem:PhoxFiberConditionGeneral}, while elementary, are vital. 
Testing if an annihilator photon $\omega \in \Pho^\times$ lies in $\in \overline{B}_{\Psi_0}|_p$ requires viewing $\omega$ as a graph over its cokernel in $U = f(p)$. To check that the candidate developing map is an immersion, we must confirm the $\Psi_0$-bases of pencil provide a local fibration over $\tilde{\Sigma}$. Such a question involves comparing nearby fibers $\overline{B}_{\Psi_0}|_{p}, \;\overline{B}_{\Psi_0}|_{p'}$, for $p, p' \in \tilde{\Sigma}$. However, these fibers are each understood with respect to different splittings $\imoct = U\oplus U^\bot$ and $\imoct = U' \oplus (U')^\bot$, where $f(p) = U, f(p') = U'$. Lemma \ref{Lem:PhoxFiberConditionGeneral} allows us to test if these fibers are actually locally disjoint. 
\end{remark}

We now resume the proof of Step 1 of Theorem \ref{Thm:DevMapPhoX}. Recalling $w_t=(u_tv_t)z_t$, let us write $x_1^t := u_t+z_t$ and $x_2^t := v_t+w_t$, so that $\omega_t=\spann \{x_1^t, x_2^t\}$. We can then further decompose 
$x_i^t = u_i^t+z_i^t$ for $u_i^t \in U_0$ and $z_i^t \in U_0^\bot$.
We emphasize: in general, the two splittings of $x_i^t$ are not the same, however, we do have $u_1^0=u_0, \,z_1^0=z_0$ and $u_2^0=v_0, \, z_2^0=w_0$. \medskip 

Here, we make a background identification of the pullback bundle $\pi^*\mathcal{E}^\R$ over $\widetilde{\Sigma}$ with $\widetilde{\Sigma}\times \R^{3,4}$ using the flat connection $\nabla$. After this background identification, the covariant derivative $D_t$ with respect to $\nabla$ is just $\frac{d}{dt}$. Note also that since the triple $(u_t,v_t,z_t)$ is $\nabla^h$-parallel, we then have $\frac{d}{dt}\big|_{t=0}u_t=\psi u_0$, and similarly for $v_t,z_t$.

Applying Lemma  \ref{Lem:PhoxFiberConditionGeneral} to $\omega_t$ with the basis $(x_1^t,x_2^t)$, but replacing $q$ with $h$ yields
\begin{align*}
    A = \frac{d}{dt}\bigg|_{t=0} \bigg (
   \underbrace{ f_1^2(t) \langle \psi_0(u_1^t), z_1^t\rangle_h}_{\text{(i)}}+ \underbrace{f_2^2(t) \langle \psi_0(u_2^t), z_2^t \rangle_h}_{\text{(ii)}} +\underbrace{f_2^2(t)f_3^2(t) \langle \psi_0(u_1^t), z_1^t\rangle_h}_{\text{(iii)}}  \\
   -\underbrace{f_2^2(t)f_3(t)\langle\psi_0(u_1^t), z_2^t\rangle_h}_{\text{(iv)}} -\underbrace{f_2^2(t)f_3(t) \langle \psi_0(u_2^t), z_1^t\rangle_h}_{\text{(v)}} \bigg ).
\end{align*}
First, observe that since $f_3(0) = 0$, term (iii) vanishes under differentiation. Next, we compute $f_1'(0), f_2'(0), f_3'(0)$. Recall $\psi = \Psi(X_0)$. It is critical going forwards that $\psi \in \T_{U_0}\X_{\Gtwosplit}$, so $\psi$ exchanges $U_0$ and $U_0^\bot$. Since the orthogonal projection map $\pi_{U_0}$ to $U_0$ and differentiation commute, 
\begin{align}\label{u1dot}
    \frac{d}{dt}\bigg|_{t=0} u_1^t =\pi_{U_0}\left( \frac{d}{dt}\bigg|_{t=0} x_1^t \right)= \pi_{U_0}\left( \psi(x_1^0) \right)= \psi(z_0).
\end{align}
Similar reasoning leads to the following identities: 
\begin{align}
     \frac{d}{dt}\big|_{t=0}z_1^t &= \psi(u_0) \\
      \frac{d}{dt}\big|_{t=0}u_2^t &= \psi(w_0) \\
       \frac{d}{dt}\big|_{t=0}z_2^t &= \psi(v_0) .
\end{align}
Thus, we have:
\[ \frac{d}{dt}\bigg|_{t=0} \langle u_1^t, u_1^t\rangle_h  = 2 \big\langle \frac{d}{dt}\big|_{t=0} u_1^t, u_1^0 \big\rangle_h = 2\langle \psi(z_0), u_0\rangle_h.  \]
As $f_1^2(t) = \langle u_1^t, u_1^t \rangle_h^{-1}$, we find $\frac{d}{dt}\big|_{t=0}f_1^2(t) = -2\langle \psi(u_0), z_0\rangle_h $. Using $f_1(0) = 1$, the derivative of term (i) is 
\[ \frac{d}{dt}\bigg|_{t=0}\mathrm{(i)}= -2\langle \psi(u_0),z_0\rangle_h \langle \psi_0(u_0), z_0\rangle_h + \langle \psi_0 \psi(u_0), u_0 \rangle_h + \langle \psi_0 \psi(z_0), z_0\rangle_h.\]
Expanding $f_2^2(t)=\langle u_2^t-f_3(t)u_1^t,u_2^t-f_3(t)u_1^t\rangle$, using that $f_3(0)=0, \langle u_2^0,u_1^0\rangle =0$, one finds
\[ \frac{d}{dt}\bigg|_{t=0}f_2^2(t) = -2\langle \psi(v_0), w_0\rangle.\]
Then one obtains 
\[ \frac{d}{dt}\bigg|_{t=0}\mathrm{(ii)} = -2\langle \psi(v_0), w_0\rangle_h \langle \psi_0(v_0), w_0\rangle_h+\langle \psi_0\psi(v_0), v_0\rangle_h+\langle\psi_0\psi(w_0),w_0\rangle_h.\]

Next, we find a contribution from the cross-terms: 
\[ f_3'(0) = \frac{d}{dt}\bigg|_{t=0}\langle u_1^t, u_2^t\rangle =\langle \psi(z_0), v_0\rangle_h +\langle u_0, \psi(w_0)\rangle_h = \langle \psi(v_0), z_0\rangle_h+\langle \psi(u_0), w_0\rangle_h. \]
Recall $f_2(0) = 1$. Using $f_3(0)$ again, we then find the derivatives of (iv), (v) as follows:
\begin{align}
    \frac{d}{dt}\bigg|_{t=0} \mathrm{(iv)}&= f_3'(0) \langle\psi_0(u_0), w_0\rangle_h= \langle \psi(v_0), z_0\rangle_h \langle\psi_0(u_0), w_0\rangle_h+\langle \psi(u_0), w_0\rangle_h\langle\psi_0(u_0), w_0\rangle_h\\
    \frac{d}{dt}\bigg|_{t=0} \mathrm{(v)}&= f_3'(0) \langle\psi_0(v_0), z_0\rangle_h= \langle \psi(v_0), z_0\rangle_h\langle\psi_0(v_0), z_0\rangle_h+\langle \psi(u_0), w_0\rangle_h\langle\psi_0(v_0), z_0\rangle_h.
\end{align}
We now suppress subscripts, and all parings are with $h$ unless otherwise specified. In total, we have computed $A$ to be the following. 
\begin{equation}
    \label{eq:A-coordinates}
\begin{aligned}
    A &= \left\langle \psi_0 \psi(u_0), u_0 \right\rangle + \left\langle \psi_0 \psi(v_0), v_0 \right\rangle \\
  &+ \left\langle \psi_0 \psi(z_0), z_0 \right\rangle + \left\langle \psi_0 \psi(w_0), w_0 \right\rangle \\
  &-2 \left\langle \psi_0(u_0), z_0 \right\rangle \left\langle \psi(u_0), z_0 \right\rangle \\
  &-2 \left\langle \psi_0(v_0), w_0 \right\rangle \langle \psi(v_0), w_0\rangle \\
  &- \left\langle \psi_0(u_0), w_0 \right\rangle \left\langle \psi(u_0), w_0 \right\rangle \\
  &- \left\langle \psi_0(u_0), w_0 \right\rangle \left\langle \psi(v_0), z_0 \right\rangle\\
   & - \left\langle \psi_0(v_0), z_0 \right\rangle \left\langle \psi(u_0), w_0 \right\rangle \\
  &- \left\langle \psi_0(v_0), z_0 \right\rangle \left\langle \psi(v_0), z_0 \right\rangle.
\end{aligned}
\end{equation}
\textbf{Step 2: Express $A$ in Local Coordinates of the fiber.}
 
By Proposition \ref{Prop:HiggsBundleCrossProduct}, we can fix a unitary basis $e_j$ of $\mathcal{L}_j$ for $-3\leq j \leq 3$ for the harmonic metric at the point $p_0\in \widetilde{\Sigma}$ that is a complex cross product basis for $\times_{\V}$ and satisfies multiplication table \ref{Table:ComplexCrossProductBasis}. 
There is still some freedom left in the unitary basis $(e_k)_{k=3}^{-3}$. 
By choosing $e_2,e_1$ appropriately, which determine the whole basis $(e_i)_{i=3}^{-3}$, we can assume that the matrix of $\psi$ in this basis is equal, for some $\beta_0,\delta_0\in \C$, to the following: 

\[ \psi=\frac{1}{||\alpha||}\begin{pmatrix}
    0& \overline{\beta_0} &0&0&0&\delta_0&0\\
     \beta_0& 0 &1&0&0&0&\delta_0\\
    0& 1&0&i\sqrt{2} \,\overline{\beta_0}&0&0&0\\
    0& 0 &-i\sqrt{2} \beta_0&0&i\sqrt{2}\,  \overline{\beta_0}&0&0\\
    0& 0 &0&-i\sqrt{2} \beta_0&0&1&0\\
    \overline{\delta_0}& 0 &0&0&1&0& \overline{\beta_0}\\
    0& \overline{\delta_0} &0&0&0& \beta_0&0
       \end{pmatrix},\;
       \psi_0=\frac{1}{||\alpha||}\begin{pmatrix}
    0& 0 &0&0&0&0&0\\
     0& 0 &1&0&0&0&0\\
    0& 1&0&0&0&0&0\\
    0& 0 &0&0&0&0&0\\
    0& 0 &0&0&0&1&0\\
    0& 0 &0&0&1&0& 0\\
    0& 0 &0&0&0& 0&0
       \end{pmatrix} \]
We shall not apply the condition $\delta_0=0$ or $\beta_0=0$, corresponding to the Hodge bundle hypothesis, until Step 3.

Using our work from $\S$\ref{Subsec:Pho^XSimplePencil}, we can now parametrize the fiber $\overline{B}_{\Psi_0}|_{|p_0}$ in coordinates. 
For every annihilator photon $\omega_0\in \mathcal{F}_{p_0}$, we shall write $\omega_0=\spann \{u_0+z_0, v_0+(u_0v_0)z_0\}$ as before, where $(u_0,v_0,z_0) \in V_{(+,+,-)}(\V^\R)$ and $u_0 \in N$. In particular, $u_0\in \mathcal{L}_2\oplus\mathcal{L}_{-2}$ can be written in the unitary basis for some $z\in \C$ with $\abs{z}=1$ as:
$$u_0=\frac{1}{\sqrt{2}}(z,0,\overline{z})\in \mathcal{L}_2\oplus\mathcal{L}_0\oplus \mathcal{L}_{-2}.$$
For some $a,b\in \R$ such that $a^2+b^2=1$, one can therefore write $v_0$ as:
$$v_0=\left( \frac{biz}{\sqrt{2}},\,a,\,-\frac{bi\overline{z}}{\sqrt{2}} \right)\in \mathcal{L}_2\oplus\mathcal{L}_0\oplus \mathcal{L}_{-2}.$$
By Corollary \ref{Cor:RGraphEquation}, we can write $z_0$ as 
\[ z_0=\lambda(1+b^2)Z+ \lambda ab Z\times u_0=\lambda(a^2+2b^2)Z+ \lambda ab Z\times u_0,\]
where $Z\in \mathcal{L}_3\oplus \mathcal{L}_{-3}$ is a norm one real vector and $\lambda^{-1}=\sqrt{(ab)^2+(1+b^2)^2}= \sqrt{a^2+4b^2}.$ Let us write $Z= \frac{1}{\sqrt{2}}(w,\overline{w})\in \mathcal{L}_3\oplus \mathcal{L}_{-3}$ in the basis $(e_3,e_{-3})$, where $w\in \C$ satisfies $\abs{w}=1$. We obtain:

\[z_0=\frac{\lambda}{\sqrt{2}}\left((a^2+2b^2)w,ab w\overline{z},ab\overline{w}z,(a^2+2b^2)\overline{w}\right)\in \mathcal{L}_3\oplus\mathcal{L}_1\oplus \mathcal{L}_{-1}\oplus \mathcal{L}_{-3}.\]

Finally, we compute similarly $w_0=(u_0v_0)z_0=(u_0\times v_0)\times z_0$: 
\begin{align*}
    u_0v_0&=\left(\frac{aiz}{\sqrt{2}},-b, -\frac{ai \,\overline{z}}{\sqrt{2}} \right) \in \mathcal{L}_2\oplus \mathcal{L}_0\oplus \mathcal{L}_{-2}, \\
    w_0&=\frac{\lambda}{\sqrt{2}}\left(-2b^3iw,+(3ab^2+a^3)iw\overline{z},-(3ab^2+a^3)i\overline{w}z,2b^3i\overline{w}\right)\in \mathcal{L}_3\oplus\mathcal{L}_1\oplus \mathcal{L}_{-1}\oplus \mathcal{L}_{-3}.
\end{align*}

In order to express $A$ in these coordinates, we  need to record the matrices associated with $\psi$, $\psi_0$, $\psi\psi_0$. First, we write the matrix of the linear map $\psi$ restricted to $\mathcal{L}_2\oplus \mathcal{L}_0\oplus \mathcal{L}_{-2}$ and co-restricted to $\mathcal{L}_3\oplus \mathcal{L}_1\oplus \mathcal{L}_{-1}\oplus \mathcal{L}_{-3}$:
 \[ \psi =
 \begin{pmatrix}
 {\color{blue}\overline{\beta_0}} & 0 & {\color{olive}{\delta_0}}  \\
  1 & {\color{blue}\sqrt{2}i\overline{\beta_0}} & 0 \\
  0 & {-\color{blue}\sqrt{2}i\beta_0} & 1 \\
 {\color{olive}{\overline{\delta_0}} }& 0 & {\color{blue}\beta_0}
 \end{pmatrix} .\] 
Here, the terms different from $1$, that are in color, are the ones that are to be replaced by $0$ to obtain the corresponding matrix representative of $\psi_0$. 

Next, we write the matrix associated to $\psi_0\psi$ restricted to the subspace $U = \mathcal{L}_2\oplus \mathcal{L}_{0}\oplus \mathcal{L}_{-2}$:
\[
\psi_0\psi |_U=
\begin{pmatrix}
1 & {\color{blue}\sqrt{2}i\overline{\beta_0}} & 0  \\
 0 & 0 & 0 \\
 0 & {-\color{blue}\sqrt{2}i\beta_0} & 1 
\end{pmatrix}
\]

Finally, we write the matrix associated to $\psi_0\psi$ restricted to the subspace $U^\bot = \mathcal{L}_3\oplus \mathcal{L}_{1}\oplus \mathcal{L}_{-1}\oplus \mathcal{L}_{-3}$:
 \[
 \psi_0\psi|_{U^{\bot}} =
 \begin{pmatrix}
 0 & 0 & 0 & 0  \\
  {\color{blue}\beta_0} & 1 & 0 & \color{olive}{\delta_0}\\
  \color{olive}{\overline{\delta_0}} & 0 & 1 & {\color{blue}\overline{\beta_0}}\\
  0 & 0 & 0 & 0\\
 \end{pmatrix} .\] 

Let us compute the 8 terms in the expression \eqref{eq:A-coordinates} for $A$ explicitly in terms of $z,w,a,b,\beta_0,\delta_0,\lambda$. Once again, we keep in color the parts containing $\beta_0$ and $\delta_0$:
\begin{align*}
\left\langle \psi_0\psi(u_0), u_0 \right\rangle&=1,\\
\left\langle \psi_0\psi(v_0), v_0 \right\rangle&=b^2+{\color{blue}\left(2ab\right)\text{Re}\left(\beta_0z\right)},\\
2\lambda^{-2}\left\langle \psi_0\psi(z_0), z_0 \right\rangle&=2a^2b^2+{\color{blue}\left(2a^3b+4ab^3\right)\text{Re}\left(\beta_0z\right)}+\color{olive}{2ab(a^2+2b^2)\text{Re}(\delta_0\overline{w}^2z)},\\
2\lambda^{-2}\left\langle \psi_0\psi(w_0), w_0 \right\rangle&=2(a^3+3ab^2)^2-{\color{blue}\left(4a^3b^3+12ab^5\right)\text{Re}\left(\beta_0z\right)}+\color{olive}{4ab^3(a^2+3b^2)\text{Re}(\delta_0\overline{w}^2z)},\\
\sqrt{2}\lambda^{-1}\left\langle \psi(u_0), z_0 \right\rangle&=\sqrt{2}ab \text{Re}(z^2\overline{w})+{\color{blue}{\sqrt{2}}\left(a^2+2b^2\right)\text{Re}\left(\beta_0\overline{z}w\right)}+\color{olive}{\sqrt{2}(a^2+2b^2)\text{Re}(\delta_0\overline{wz})},\\
\sqrt{2}\lambda^{-1}\left\langle \psi(v_0), w_0 \right\rangle&=\sqrt{2}\left(a^3b+3ab^3\right) \text{Re}(z^2\overline{w})+{\color{blue}2{\sqrt{2}}\left(a^4+3a^2b^2-b^4\right)\text{Re}\left(\beta_0\overline{z}w\right)}+\color{olive}{2\sqrt{2}b^4\text{Re}(\delta_0\overline{wz})}.\\
\sqrt{2}\lambda^{-1}\left\langle \psi(u_0), w_0 \right\rangle &=-\sqrt{2}\left(a^3+3ab^2\right) \text{Re}(iz^2\overline{w})-{\color{blue}{2\sqrt{2}}b^3\text{Re}\left(i\beta_0\overline{z}w\right)}+\color{olive}{2\sqrt{2}b^3\text{Re}(i\delta_0\overline{wz})},\\
\sqrt{2}\lambda^{-1}\left\langle \psi(v_0), z_0 \right\rangle &=\sqrt{2}\left(ab^2\right) \text{Re}(iz^2\overline{w})-{\color{blue}{\sqrt{2}}\left(3a^2b+2b^3\right)\text{Re}\left(i\beta_0\overline{z}w\right)}-\color{olive}{\sqrt{2}b(a^2+2b^2)\text{Re}(i\delta_0\overline{wz})},
\end{align*}

By expanding the individual terms of $A$ and grouping them, we write out $2\lambda^{-2} A $ as follows: 
\begin{equation}
    \label{PolynomialEquationA} 
\begin{aligned}
    2\lambda^{-2} A &= 2a^6+12a^4b^2+18a^2b^4+8b^4+4a^2b^2+2a^2+8b^2\\
    &+\textcolor{blue}{\mathrm{Re}(\beta_0z)}\bigg(-4a^3b^3-12ab^5+6a^3b+20ab^3\bigg) \\
    &-4\mathrm{Re}(z^2\overline{w})^2 \big( a^6 b^2 + 6 a^4 b^4 + 9a^2 b^6 + a^2 b^2\big ) \\
    &-4 \mathrm{Re}(z^2\overline{w}) \textcolor{blue}{\mathrm{Re}(\beta_0\overline{z}w)}\bigg( 2 a^7b + 12 a^5 b^3 + 16 a^3 b^5 - 6 a b^7 +a^3 b + 2 ab^3\bigg) \\
    &-2 \mathrm{Re}(iz^2\overline{w})^2\bigg(a^6+4a^4b^2+4a^2b^4\bigg)\\
    &-\mathrm{Re}(iz^2\overline{w})\textcolor{blue}{\mathrm{Re}(i\beta_0\overline{z}w)} \bigg( 6 a^5 b + 20 a^3 b^3 + 16 a b^5 \bigg).\\
    & +\textcolor{olive}{\text{Re}(\delta_0\overline{w}^2z)} \bigg(4a^3b^3+2a^3b+12ab^5
    +4ab^3 \bigg)\\
    & -4\mathrm{Re}(z^2\overline{w})\textcolor{olive}{\mathrm{Re}(\delta_0 \overline{z}\overline{w})}\bigg(2a^3b^5+a^3b+6ab^7+2ab^3\bigg) \; \\ 
    &-\mathrm{Re}(iz^2\overline{w})\textcolor{olive}{\mathrm{Re}(i\delta_0 \overline{z}\overline{w})}\bigg(2 a^5 b + 4 a^3 b^3\bigg) \\
\end{aligned} 
\end{equation}

\begin{remark}
As in Remark \ref{rem:explainationComputation}, we know a priori that the sum of uncolored terms is positive. All that remains is to control the colored terms via maximum principles.
\end{remark}

\textbf{Step 3: Prove $A > 0$.} 

First note that $z^2\overline{w} \times \beta_0\overline{z}w=\beta_0z$ and $z^2\overline{w} \times \delta_0\overline{zw} =\delta_0\overline{w}^2z$ and hence 
\begin{align}
    \textcolor{blue}{\mathrm{Re}(\beta_0z)}&=\mathrm{Re}(z^2\overline{w}) \textcolor{blue}{\mathrm{Re}(\beta_0\overline{z}w)}-\mathrm{Re}(iz^2\overline{w})\textcolor{blue}{\mathrm{Re}(i\beta_0\overline{z}w)}. \\
    \textcolor{olive}{\mathrm{Re}(\delta_0\overline{w}^2z)}&=\mathrm{Re}(z^2\overline{w}) \textcolor{olive}{\mathrm{Re}(\delta_0\overline{zw})}-\mathrm{Re}(iz^2\overline{w})\textcolor{olive}{\mathrm{Re}(i\delta_0\overline{zw})} 
\end{align}Let us write $x:=\mathrm{Re}(z^2\overline{w})$ and $y:=\mathrm{Re}(iz^2\overline{w})$, so that  $x^2+y^2=1$. We can rewrite $ 2\lambda^{-2} A$ as follows

\begin{equation}
    \label{PolynomialEquationA2}
\begin{aligned}
    2\lambda^{-2} A &= 2a^6+12a^4b^2+18a^2b^4+8b^4+4a^2b^2+2a^2+8b^2\\
    &-4x^2 \bigg( a^6 b^2 + 6 a^4 b^4 + 9a^2 b^6 + a^2 b^2\bigg) \\
    &-2 y^2\bigg(a^6+4a^4b^2+4a^2b^4\bigg)\\
     &-x \textcolor{blue}{\mathrm{Re}(\beta_0\overline{z}w)}\bigg( 8 a^7b + 48 a^5 b^3 + 64 a^3 b^5 - 24 a b^7 + 4a^3b^3+12ab^5 -2a^3b-12ab^3\bigg)\\
    &- y\textcolor{blue}{\mathrm{Re}
    (i\beta_0\overline{z}w)} \bigg( 6 a^5 b + 16 a^3 b^3 + 4 a b^5 +6a^3b+20ab^3\bigg) \\
    & -x \textcolor{olive}{\mathrm{Re}(\delta_0\overline{z}\overline{w})}\bigg( 2 a b (a^2 + 2 b^2 - 2 a^2 b^2 - 6 b^4 + 4 a^2 b^4 + 12 b^6) \bigg) \\
    &-y\textcolor{olive}{\mathrm{Re}(i\delta_0 \overline{z}\overline{w})}\bigg(2ab(a^4+4a^2b^2+a^2+2b^2+6b^4)\bigg). \\
\end{aligned}
\end{equation}
We can simplify expressions in $a,b$ using the relation $a^2+b^2=1$. Doing so, one finds the following:  
\[\begin{aligned}
2a^6+12a^4b^2+18a^2b^4+8b^4+4a^2b^2+2a^2+8b^2 &= 4 + 16 b^2 + 4 b^4 - 8 b^6,\\
a^6 b^2 + 6 a^4 b^4 + 9a^2 b^6 + a^2 b^2&=2b^2(1+b^2-2b^6), \\
a^6+4a^4b^2+4a^2b^4 &= 1 + b^2 - b^4 - b^6, \\
 \textcolor{blue}{8 a^7b + 48 a^5 b^3 + 64 a^3 b^5 - 24 a b^7 + 4a^3b^3+12ab^5 -2a^3b-12ab^3} &= \textcolor{blue}{6ab(1+3b^2-8b^6)},\\
  \textcolor{blue}{6 a^5 b + 16 a^3 b^3 + 4 a b^5 +6a^3b+20ab^3} &= \textcolor{blue}{6 a b (2 + 3 b^2 - b^4)},\\
\textcolor{olive}{2 a b (a^2 + 2 b^2 - 2 a^2 b^2 - 6 b^4 + 4 a^2 b^4 + 12 b^6)}&= \textcolor{olive}{2 a b (8 b^6 - b^2 + 1)}\\
  \textcolor{olive}{2ab(a^4+4a^2b^2+a^2+6b^4+2b^2)} &= \textcolor{olive}{2ab(3b^4+3b^2+2)} .
  \end{aligned} \] 

Using these relations, we can re-write $2\lambda^{-2}A$ in much simpler form.
\begin{equation}
    \label{PolynomialEquationA3}
\begin{aligned}
    2\lambda^{-2} A &= 4 + 16 b^2 + 4 b^4 - 8 b^6\\
    &-4x^2 \bigg( 2b^2+2b^4-4b^8\bigg) \\
    &-2 y^2\bigg(1 + b^2 - b^4 - b^6\bigg)\\
    &-x \textcolor{blue}{\mathrm{Re}(\beta_0\overline{z}w)}\bigg( 6ab(1+3b^2-8b^6)\bigg) \\
    &- y\textcolor{blue}{\mathrm{Re}(i\beta_0\overline{z}w)} \bigg(6 a b (2 + 3 b^2 - b^4)\bigg).\\
    &-\, x\textcolor{olive}{\mathrm{Re}(\delta_0 \overline{z}\overline{w})}\bigg(2 a b (8 b^6 - b^2 + 1)\bigg)\\
    &-y\textcolor{olive}{\mathrm{Re}(i\delta_0 \overline{z}\overline{w})}\bigg(2ab(3b^4+3b^2+2))\bigg) \\
\end{aligned}
\end{equation}
We now gather the terms in \eqref{PolynomialEquationA3} as follows: 
\[ 2\lambda^{-2} A = A_1
    -x^2 A_2 -y^2 A_3-x \textcolor{blue}{\mathrm{Re}(\beta_0\overline{z}w)}A_4 - y\textcolor{blue}{\mathrm{Re}(i\beta_0\overline{z}w)} A_5-x\textcolor{olive}{\mathrm{Re}(\delta_0 \overline{z}\overline{w})}A_6-y\textcolor{olive}{\mathrm{Re}(i\delta_0 \overline{z}\overline{w})}A_7 ,\]
where the coefficients in $a,b$ are given by:
\[\begin{aligned}
A_1&=4+16b^2+4b^4-8b^6, \\
A_2&=8b^2+8b^4-16b^8, \\
A_3&=2+2b^2-2b^4-2b^6, \\
A_4&={6ab(1+3b^2-8b^6)},\\
A_5&=6ab(2+3b^2-b^4),\\
 A_6&= 2 a b (8 b^6 - b^2 + 1),\\
 A_7&= 2ab(3b^4+3b^2+2)
\end{aligned} \]
Recall that $|z|=|w|=1$. 
The vectors $\textcolor{blue}{\mathbf{v}_{\beta}}: = (\textcolor{blue}{\mathrm{Re}(\beta_0\overline{z}w)},\textcolor{blue}{\mathrm{Re}(i\beta_0\overline{z}w)})$ and $\textcolor{olive}{\mathbf{v}_{\delta}}=(\textcolor{olive}{\mathrm{Re}(\delta_0\overline{zw})},\textcolor{olive}{\mathrm{Re}(i\delta_0\overline{zw})})$ in $\R^2$ satisfy $|| \mathbf{v}_{\beta} || \leq |\beta_0|$ and $||\mathbf{v}_{\delta}|| \leq |\delta_0|$. 
Cauchy-Schwarz then gives the following inequality: 
\begin{align}\label{ALowerBound1}
2\lambda^{-2} A \geq A_1-x^2A_2-y^2A_3-|\beta_0|\sqrt{x^2A_4^2+y^2A_5^2} -|\delta_0|\sqrt{x^2A_6^2+y^2A_7^2}.
\end{align}
We now finally consider cases corresponding to whether $\delta$ or $\beta$ vanishes. \medskip 

\textbf{Case 1}: $\delta=0$. Observe that 
$A_1-x^2A_2-y^2A_3\geq A_1-A_2-A_3>0$. 
By Lemma \ref{Lem:MaximumPrincipleAlpha}, we have $|\beta_0| =\frac{||\beta||_h}{||\alpha||_h}\leq \sqrt{\frac{3}{5}}$.
Thus, in order to check that the lower bound of \eqref{ALowerBound1} is positive, it suffices to check that \eqref{ALowerBound2} is positive:
\begin{align}\label{ALowerBound2}
 \left((A_1-A_2)x^2+(A_1-A_3)y^2\right)^2-\frac{3}{5}x^2A_4^2(x^2+y^2)-\frac{3}{5}y^2A_5^2(x^2+y^2).
\end{align}
We homogenized the equation, multiplying some terms with $x^2+y^2=1$. The new polynomial lower bound is now homogeneous of degree 4. 
Define $X=x^2$ and $Y=y^2$. We rewrite this expression in the form $X^2C_{XX}+XYC_{XY}+Y^2C_{YY}$ below:
$$X^2\left((A_1-A_2)^2-\frac{3}{5}A_4^2\right)+XY\left(2(A_1-A_2)(A_1-A_3)-\frac{3}{5}A_4^2-\frac{3}{5}A_5^2\right)+Y^2\left((A_1-A_3)^2-\frac{3}{5}A_5^2\right).$$
The three coefficients $C_{XX},C_{XY},C_{YY}$ are the following polynomials in $b$. 
\begin{align*}
    C_{XX}(b) &= \frac{1}{5}\left(80 + 212 b^2 - 380 b^4 - 964 b^6 + 2780 b^8 + 5056 b^{10} - 5504 b^{12}- 
 8192 b^{14} + 8192 b^{16} \right), \\
    C_{XY}(b) &= \frac{1}{5}\left( 80 + 180 b^2 - 124 b^4 - 48 b^6 + 2368 b^8 + 4700 b^{10} - 3636 b^{12} - 
 7872 b^{14} + 6912 b^{16}\right),\\
 C_{YY}(b) &= \frac{1}{5}\left(20 - 152 b^2 + 236 b^4 + 1476 b^6 + 528 b^8 - 1116 b^{10} + 288 b^{12}\right).
\end{align*}
These three polynomials are positive on all of $\R$, which concludes the proof of \emph{Case 1}.\medskip 

\textbf{Case 2}: $\beta=0$. We obtain from Lemma \ref{Lem:MaximumPrincipleAlpha} that $\abs{\delta_0}=\frac{||\delta||_h}{||\alpha||_h}\leq 1$. We proceed similarly: we rewrite the lower bound for $\lambda^{-2}A$ in \eqref{ALowerBound1} in the form $X^2C_{XX}+XYC_{XY}+Y^2C_{YY}$ below:
\[X^2\left((A_1-A_2)^2-A_6^2\right)+XY\left(2(A_1-A_2)(A_1-A_3)-A_6^2-A_7^2\right)+Y^2\left((A_1-A_3)^2-A_7^2\right).\]

We obtain the following: 
\begin{align*}
    C_{XX}(b) &= 16 + 60b^2 + 44b^4 - 140b^6 - 44b^8 + 448b^{10} - 128b^{12} - 512b^{14} + 512b^{16} , \\
    C_{XY}(b) &= 16 + 124b^2 + 236b^4 - 144b^6 - 352b^8 + 564b^{10} + 260b^{12} - 448b^{14} + 256b^{16} , \\
    C_{YY}(b) &= 4 + 40b^2 + 188b^4 + 108b^6 - 120b^8 - 36b^{10} + 72b^{12} .
\end{align*}

These three polynomials are positive on $\R$, which concludes the proof of \emph{Case 2}. 
\end{proof}

 \begin{remark}
    Note that in both Cases 1 and 2 the polynomials $C_{XX},C_{XY}, C_{YY}$ lie in $\Z[x]$, up to a real scalar, and hence Sturm's method enables algorithmic verification that they are positive on $\R$. 
\end{remark}

\subsection{The \texorpdfstring{$\Gtwosplit$}{G2'}-Hitchin Case: Comparison with Tits Metric Thickening}

We will now prove that our developing map is a finite covering of the domain $\Omega^{\Thick}$ in the case of a $\Gtwosplit$-Hitchin Hodge bundle. Thus, the manifold $B_{\Psi_0}$ coincides with the  quotient $\rho(\pi_1S)\backslash \Omega_{\rho}^{\Thick}$ up to finite cover. 

\begin{theorem}\label{Thm:AlphaHodgeTits=Pencil}
Let $\rho:\pi_1S \rightarrow \Gtwosplit$ be Fuchsian-Hitchin and $\Sigma$ a Riemann surface such that $\rho$ corresponds to a Hodge bundle. Then $\dev:\overline{B}_{\Psi_0}\rightarrow \Pho^\times$ has image the domain $\Omega^{\Thick}_{\rho}$ from \eqref{Omega_Thick_PhoX}. Moreover, $\dev$ is a finite covering map onto $\Omega^{\Thick}_{\rho}$. 
\end{theorem}

\begin{proof}
Take any point $p_0 \in \tilde{\Sigma}$ and fix the unique unit tangent vector $v \in \T_{p_0}\tilde{\Sigma}$ pointing towards $x \in \partial \Gamma$, so that $\xi^2(x) =\omega$, where $\xi^2:\partial \Gamma\to \Pho^\times$ is the limit map of $\rho$. 
We begin by normalizing the matrices representing $\Psi_0(v)$ and $\Psi(v)$ in an $h$-unitary cross product basis $(e_k)_{k=3}^{-3}$. In particular, re-gauging by $g= \diag(ab,a,b,1, \overline{b}, \overline{a}, \overline{ab}) \in \Gtwosplit$, where $a,b \in \mathbb{S}^1$, we can guarantee $\Psi_0(v)$ and $\Psi(v)$ may be represented by some matrices $A_0, A$ of the following form: 
\[ \frac{1}{||\alpha||}A_0=\begin{pmatrix}
    0& 0 &0&0&0&0&0\\
    0& 0 &1&0&0&0&0\\
    0& 1&0&0 &0&0&0\\
    0& 0 &0&0&0&0&0\\
    0& 0 &0&0&0&1&0\\
    0& 0 &0&0&1&0&0\\
    0& 0 &0&0&0&0&0
       \end{pmatrix}, \; \; \frac{1}{||\alpha||}A=\begin{pmatrix}
    0& \sqrt{\frac{3}{5}} &0&0&0&0&0\\
    \sqrt{\frac{3}{5}}& 0 &1&0&0&0&0\\
    0& 1&0&i\sqrt{\frac{6}{5}}&0&0&0\\
    0& 0 &-i\sqrt{\frac{6}{5}} &0&i\sqrt{\frac{6}{5}} &0&0\\
    0& 0 &0&-i\sqrt{\frac{6}{5}}&0&1&0\\
    0& 0 &0&0&1&0& \sqrt{\frac{3}{5}}\\
    0& 0 &0&0&0&\sqrt{\frac{3}{5}}&0
       \end{pmatrix} \]
Since the minimal surface $f: \tilde{\Sigma} \rightarrow \X$ is tangent to a flat and $\Delta$-regular, we know that $\Psi(v)$ points towards a pointed annihilator photon $(x, \omega) \in \mathcal{F}_{1,2}^\times$. The annihilator photon $\omega$ is the span of the top two eigenvectors of $A$.
Thus, in our chosen $h$-unitary complex cross product basis $(e_k)_{k=3}^{-3}$ of $\imoct^\C$, with real locus $(z,w,\zeta,r, \overline{\zeta}, \overline{w}, \overline{z})$ for $z,w,\zeta \in \C$ and $r \in \R$, we have  
\begin{align}\label{LimitSetAnnihilatorPhoton}
 \omega = \spann_{\R}\{(-i, \,-\sqrt{6}i, \,-\sqrt{15}i, \,-2 \sqrt{5}, \,\sqrt{15}i, \,\sqrt{6}i, \,i), \, (\sqrt{3}, \, 2\sqrt{2}, \,\sqrt{5}, \,0, \sqrt{5}, \,2\sqrt{2}, \,\sqrt{3})\} 
\end{align}

We fix the point $U \in \X$ corresponding to $f(p_0)$. In the given coordinates, the 3-plane $U$ is given by $U = \spann_{\R} \langle \mathbf{e}_2+\mathbf{e}_{-2},\mathbf{e}_0, i(\mathbf{e}_2-\mathbf{e}_{-2}) \rangle$.\medskip 

To prove the fiber $(B_{\Psi_0})|_p$ is disjoint from the thickening $K_{\omega}$, by Proposition \ref{Prop:PhoXThickening} it suffices to show the following: every annihilator photon $\omega' \in (B_{\Psi_0})|_p$ in the fiber is not orthogonal to $\omega$. Thus, we need only exhibit elements $x' \in \omega'$ and $x \in \omega$ such that $x\cdot x'\neq 0$. We shall do this now.

Take any $\omega' \in (B_{\Psi_0})|_p$. Define $W:= \pi_{U}(\omega')$ as the orthogonal projection of $\omega'$ onto $U$. Of course, $W \cap N \neq\emptyset$. Now, by Lemma \ref{Lemma:ProjectionOntoB}, there is some element $x' \in \omega'$ of the form $x' = u'+z'$ for $u'\in Q_+(N\cap W)$ and $z' \in Q_-(B)$. In other words, $\omega'$ contains some line $\ell \in \Ein(N \oplus B)$. To finish the proof, it then suffices to verify that every element 
$\ell \in\Ein(N\oplus B)$ is not orthogonal to $\omega$. 
We now prove this claim. Let $0 \neq y \in Q_0(N\oplus B)$ be an arbitrary element spanning a line $\ell = [y] \in \Ein(N \oplus B)$.  
Then $y$ obtains the form 
\[ y =(\cos(\theta) +i\sin(\theta),\,\cos(\alpha)+i\sin(\alpha),0,0,0,\cos(\alpha)-i\sin(\alpha),\cos(\theta)-i\sin(\theta)). \]

Denote $E_1, E_2$ as the eigenvectors of $A$ spanning $\omega$ in \eqref{LimitSetAnnihilatorPhoton}. One can then explicitly compute $y\cdot E_1$ and $y\cdot E_2$ and see that the system 
\[ \begin{cases}
    y \cdot E_1=0\\
    y \cdot E_2 =0 .
\end{cases} \] 
has no solutions.

We conclude that the fiber $(B_{\Psi_0})_{p_0}$ does not intersect the thickening $K_{\omega}$ of the particular annihilator photon $\omega \in \image(\xi^2)$. Now, we may identify $\tilde{\Sigma}$ with $\Ha^2$, allowing us to define $H$ as the stabilizer subgroup in $\PSL(2,\R)$ of $p_0 \in \Ha^2$. Since $H$ acts transitively on $\RP^1 =\vis\Ha^2$, by $H$-equivariance we have $\dev(B_{\Psi_0}|_{p_0}) \subset \Omega^{\Thick}_{\rho}$. That is, one fiber is mapped inside the domain $\Omega_{\rho}^{\Thick}$. 
By $\PSL(2,\R)$-invariance of $\dev$, we then conclude the entire image of $\dev$ is contained inside $\Omega^{\Thick}_{\rho}$. 

Using the same reasoning as in the end of Theorem \ref{Thm:BetaTits=Pencil}, we conclude that $\dev$ surjects onto $\Omega^{\Thick}_{\rho}$ and is a finite covering map. In particular, we need Proposition \ref{Prop:PhoXThickening} to see that $\Omega^{\Thick}_{\rho}$ is connected. 
\end{proof}  

We can then describe the quotient of the \cite{KLP18}-domain of discontinuity up to finite cover. 
\begin{corollary}[Topology of the Quotient]\label{Cor:QuotientTopologyAlpha}
Let $\rho: \pi_1S\rightarrow \Gtwosplit$ be Hitchin and $\mathcal{M} =\rho(\pi_1S)\backslash \Omega_{\rho}$ the quotient of the  $\frac{\pi}{2}$-Tits metric thickening domain in $\Pho^\times$. Then $\mathcal{M}$ admits a (possibly trivial) smooth covering by the total space of the $(\RP^2\times \sphere^1)$-fiber bundle $\underline{\RP^2} \oplus \sphere(\mathcal{T}^3\K^3) \rightarrow S$. 
\end{corollary}

\begin{proof}
This follows from the description of $B_{\Psi_0}$ in Remark \ref{Remk:TopologyAlpha}, the fact that $\Gr_2(\mathscr{U}) \rightarrow S$ is a smoothly trivial fiber bundle, and that $\Hom_{\C}(N,B)\cong \mathcal{T}^3\mathcal{K}^3$ as a complex line bundle over $\Sigma$. 
\end{proof}

\appendix

\section{Unified Constructions of Geometric Structures}\label{Appendix:Unified}

In this appendix, we briefly summarize a unified approach for analytically building geometric structures, which we have used presently to construct $(\Gtwosplit, \Ein^{2,3})$ and $(\Gtwosplit, \Pho^\times)$-structures. We will show this same strategy also gives a reinterpretation of other known geometric structures constructions. The main steps are as follows.  

\begin{enumerate}
\item Fix a Lie group $G$ and a representation of $\SL(2,\R)$ in $G$. Consider the associated totally geodesic copy $\mathcal{H}$ of $\mathbb{H}^2$ in the symmetric space $\X = \X_{G}$. Let $\hat{H}$ be the stabilizer in $G$ of $(x,\mathcal{H})$, for some $x\in\mathcal{H}$, and $G/\hat{H}$ the associated homogeneous space .

\item Consider a representation $\rho: \pi_1S\rightarrow G$ equipped with a pair of $\rho$-equivariant objects: 
\begin{itemize}
    \item a map $u:\tilde{S}\to\X$,
    \item and distribution $\mathcal{P}_x$ of smoothly varying 2-planes, or \emph{pencils}, such that $\mathcal{P}_x\subset \T_{u(x)}\X$ and moreover $\mathcal{P}_x$ is tangent to a totally geodesic copy $\mathcal{H}_x$ of $\mathbb{H}^2$ in $\X_G$.
\end{itemize} 
One can view this distribution of planes as the image of a $1$-form $\Psi_0\in \Omega^1(\tilde{S}, u^*\T\X)$. Such a pair $(u,\Psi_0)$ describes, and is described by, a $\rho$-equivariant map $\hat{\sigma}:\tilde{S}\to G/\hat{H}$. 
We demand moreover that these planes $\mathcal{P}_x$ are transported parallely by the Levi-Civita connection.  

\item Check that the $\mathcal{F}$-base of pencil of $\Psi_0$, for well-chosen flag manifold $\mathcal{F}$, locally define a fibration, to construct a $(G, \mathcal{F})$-structure on a fiber bundle over the surface.

\end{enumerate}

The proof of point (3) in each case is a different computation that uses the parallelism of the distribution of planes and leverages maximum principles for the Hitchin system. In particular, to proceed from (1) to (2), one must \emph{find} such special representations $\rho$. In the present case $G=\Gtwosplit$, this is exactly the purpose of the $\alpha$ and $\beta$-bundles studied. 

\medskip

The map $\hat{\sigma}$ induces a \emph{projected metric} $g$ on $\tilde{S}$, such that the norm of $v\in \T_x\tilde{S}$ is the norm of the orthogonal projection of $\mathrm{d}u(v)$ to $\mathcal{P}_x$. The parallelism condition on the map $\hat{\sigma}$ can be re-written as $\nabla^{\X} \circ \Psi_0=\Psi_0\circ \nabla^{g} $.

\begin{question}
\label{que:Geometric structures appendix}
For which other Lie groups $G$ and representations of $\SL(2,\R)$ can one construct slices of representations admitting parallel equivariant maps $\hat{\sigma}:\tilde{S}\to G/\hat{H}$, and for which the moving bases of pencils construction describes a geometric structure? 
\end{question}

One way to construct such parallel equivariant maps is to start with special cyclic Higgs bundles, such as Higgs bundles in certain Slodowy slices.\medskip

In the cases we describe now, the map $\hat{\sigma}$ is induced by a simpler map $\sigma:\tilde{S}\to G/H$ together with its derivatives, for some $\hat{H}\subset H\subset G$. For example, $\sigma = \nu: \tilde{\Sigma} \rightarrow \quadric$ in the present case. In the following examples we do not discuss $\hat{\sigma}$ directly, but rather describe $u$ and $\Psi_0$ that characterize $\hat{\sigma}$. 

\subsection{\texorpdfstring{$G=\SL(3,\R)$}{G=PSL(3,R)} and \texorpdfstring{$X = \Flag(\R^3)$}{X=Flag(R3)}}
Suppose $\rho: \pi_1S \rightarrow \SL(3,\R)$ is $P_{\Delta}$-Anosov, with limit map $\xi: \partial_{\infty}\pi_1S \rightarrow \Flag(\R^3)$, written $\xi = (\xi^1,\xi^2)$. Given $F = (\ell, H) \in \Flag(\R^3)$, consider the thickening:
\[ K_{F} = \{ (\ell',H') \in \Flag(\R^3) \mid \ell' \subset H \; \text{or} \,\;\ell\subset H' \} \]
Note that $K_{F} \cong \mathbb{S}^1 \vee \mathbb{S}^1$ is a wedge of two circles. Guichard-Wienhard \cite{GW12} introduced the following co-compact domain of discontinuity for $\rho$: 
\begin{align}
\Omega_{\rho}:= \Flag(\R^3) \setminus \bigcup_{x \in \partial_{\infty}\pi_1S} K_{\xi(x)}, 
\end{align}
obtained by removing the thickening of the limit set $\Lambda := \image(\xi)$. A special feature of this example is that $\Omega_{\rho}$ is disconnected. In particular, $\Omega_{\rho}$ is well-known to have 3 components. We label the components as in \cite{NR25}, where Nolte and Riestenberg extensively study the construction of foliated $(\SL(3,\R), \Flag(\R^3))$-structures for $\SL(3,\R)$-Hitchin representations in analogy to the foliated $(\PSL(4,\R), \RP^3)$-structures for $\PSL(4,\R)$-Hitchin representations in \cite{GW08}. To describe the components, we recall the image of $\xi^1$ bounds an open convex domain $\mathcal{C}_{\rho} \subset \RP^2$, which is exactly the convex domain associated to the Choi-Goldman $(\SL(3,\R), \RP^2)$-structures \cite{Gol90, CG93}. Viewing $H$ as a projective line in $\RP^2$, the components of $\Omega_{\rho}$ are: 
\begin{align} \label{OmegaFlagR3}
     \begin{cases}
        \Omega_{\rho}^1 &= \{ (\ell , H) \in \Flag(\R^3) \mid \ell \in \mathcal{C}_{\rho}\} \\
    \Omega_{\rho}^2 &= \{(\ell ,H) \in \Flag(\R^3) \mid \ell \notin \mathcal{C}, H \cap \mathcal{C}_{\rho} \ne \emptyset \} \\
     \Omega_{\rho}^3 &= \{(\ell, H) \in \Flag(\R^3) \mid H \cap \partial \overline{\mathcal{C}}_{\rho}=\emptyset \} 
    \end{cases} 
\end{align}

We now discuss the recipe from the introduction of this appendix in the context of $\SL(3,\R)$-Hitchin representations.

By the Labourie-Loftin parametrization, given a Hitchin representation $\rho: \pi_1S \rightarrow \SL(3,\R)$, there is naturally associated a pair $(\Sigma, q_3)$ consisting of a Riemann surface $\Sigma$ and a holomorphic cubic differential $q_3 \in H^0(\K_{\Sigma}^3)$ on that Riemann surface \cite{Lab06, Lof01}.  
In particular, using the Hitchin section and the non-abelian Hodge correspondence, $\rho$ is the holonomy of the flat connection built from the Higgs bundle $\mathcal{H}(\rho) = (\V, \Phi)$ in the $\SL(3,\R)$-Hitchin section, where $\V = \K_{\Sigma} \oplus \mathcal{O} \oplus \mathcal{K}_{\Sigma}^{-1}$ and $\Phi=\Phi(q_3)$ is described by the following diagram:
\begin{equation*}
\begin{tikzcd}
\mathcal{K}_{\Sigma} \arrow[r, "1"] & \mathcal{O} \arrow[r, "1"]&\mathcal{K}_{\Sigma}^{-1} \arrow[ll,bend right,swap, "q_3"]
\end{tikzcd}.
\end{equation*}

Let $\X $ be the $\SL(3,\R)$-Riemannian symmetric space, which can be identified with the space of Euclidean metrics on $\R^3$. 
In this case, we make the following choices: 
\begin{itemize}
    \item The map $u: \tilde{\Sigma} \rightarrow \X$ is the harmonic map associated to $\rho$ via $\NAH_{\Sigma}$. 
    \item $ \Psi_0 = \Phi_0+\Phi_0^* = \begin{pmatrix} 0 & 1^* & 0 \\
    1 & 0 & 1^* \\
    0 & 1 & 0 \end{pmatrix} $. 
\end{itemize}
We recall Baraglia's insight that the tautological section $s \in \Gamma(\Sigma, \V)$, given by $s(p) = 1 \in \mathcal{O}|_p$, corresponds to the hyperbolic affine sphere $\sigma: \tilde{S} \rightarrow \R^3$ associated to $\rho$ via Labourie-Loftin \cite{Bar10}. 
The solution to Hitchin's equations for this Higgs bundle is diagonal, and obtains the form $h = \diag(\frac{1}{g},1, g)$, for $g \in \Gamma(\K \overline{K} )$ a Riemannian metric on $\Sigma$, where $g$ is the Blashke metric of the affine sphere $\sigma$. The translation between $u, \sigma$ and the Higgs bundle $(\V,\Phi)$ extremely convenient, and is analogous to the case of $G=\Gtwosplit$. The harmonic metric $h$ corresponds to $u$ and the bundle-valued endomorphism 
$[\K \stackrel{1^*}{\leftarrow} \mathcal{O} \stackrel{1}{\rightarrow} \K^{-1} ]$ interprets as the differential $\mathrm{I} = d\sigma$, whereby $\Psi_0= \Sym(\mathrm{I})$ is the symmetrization of $\mathrm{I}$ with respect to the Euclidean metric $h$ on the real locus $\V^{\R} \subset \V$. 

Here is the parallelism property in this case. Once more, it is an exercise in identifications. 

\begin{proposition}[$\Psi_0$-Parallelism]\label{ParallelismPair}
Let $X$ be a local vector field on $\Sigma$. Let $\nabla^g$ be the Levi-Civita connection of $g$ and $\nabla^\X$ the connection on $\T\X$. Conflating $u^*\nabla^\X$ with $\nabla^\X$, we have 
\[ (\nabla^\X\circ \Psi_0)(X)= (\Psi_0\circ\nabla^g)(X ).\]
\end{proposition}

Following the proposed recipe, we have the following result. Let $\rho \in \Hit(S, \SL(3,\R))$ be Hitchin. Write $\mathcal{F}=\Flag(\R^3)$ and define $\overline{M} \rightarrow \tilde{\Sigma}$ to be a fiber sub-bundle of 
$\tilde{\Sigma} \times \Flag(\R^3)$ with fiber 
\[ \overline{M}|_x = \{ F \in \Flag(\R^3) \mid F \in \mathcal{B}_{\mathcal{F}}(\Psi_0|_x) \}.\]
The $(G,X)$-manifolds of interest shall be of the form $M = \pi_1S \backslash \overline{M}$.

\begin{theorem}[$\Flag(\R^3)$-structures for $\Hit(S,\SL(3,\R)$ {\cite{RT25}}]
Let $\rho:\pi_1S \rightarrow \SL(3,\R)$ be Hitchin. The manifold $\overline{M}$ has fiber $\mathbb{S}^1 \sqcup \mathbb{S}^1 \sqcup \mathbb{S}^1$ and the tautological developing map $\dev: \overline{M} \rightarrow \Flag(\R^3)$ by $(x,F)\mapsto F$ is a diffeomorphism onto $\Omega_{\rho}$. 
\end{theorem}

The component $\Omega^2_{\rho}$ in \eqref{OmegaFlagR3} is sometimes called the de Sitter component of $\Omega_{\rho}$, for which this theorem is proven in  \cite[{Proposition 5.13}]{RT25}. For the two other components, the theorem is a consequence of the fact that the projectivization of the affine sphere and the dual affine sphere are exactly the preserved open convex set in $\mathbb{RP}^2$ and its dual, as remarked in \cite[{Section 5.4}]{RT25}. 

\medskip

The geometry and topology of the base of pencil of $\Psi_0$ at $x\in \widetilde{\Sigma}$ are described in \cite[Example 6.9]{Dav25}. Note that in this reference the pencil is described in a different basis, however both pencils are tangent to a totally geodesic copy of $\mathbb{H}^2$ associated with an irreducible representation of $\PSL(2,\R)$, so the pencils are conjugate under $G$ and their corresponding bases are diffeomorphic.

Let $\ell_0\in \mathbb{RP}^2$ and $H_0\in \left(\mathbb{RP}^2\right)^*$ be affine sphere and the linear part of the tangent hyperplane to the affine sphere.
The base of $\Psi_0$ is the union of three circles:
\begin{itemize}
    \item $\lbrace (\ell, H) \in \Flag(\R^3)\mid \ell= \ell_0 \rbrace$,
    \item $\lbrace (\ell, H)\in \Flag(\R^3)\mid  \ell\subset H_0, \,\ell_0\subset H\rbrace$,
    \item $\lbrace (\ell, H)\in \Flag(\R^3)\mid  H=H_0\rbrace$.
\end{itemize}

The quotients of the domain of discontinuity are understood topologically. 
In fact, each quotient $M_i = \rho(\pi_1S)\backslash \Omega_{\rho}^i$ is homeomorphic to $\mathbb{P}(\mathrm{T}^1S)$, the projective unit tangent bundle of $S$. 

\begin{remark}
Hitchin representations in $\SL(3,\R)$ admit geometric structures modeled on $\RP^2$ by \cite{Gol90, CG93}. This structure is directly related to the one we just described, as a developing map for it is just the affine sphere. However, this structure is obtained from a domain of discontinuity that is quite exceptional: it is not associated with a balanced Tits-Bruhat ideal, in the sense of \cite{KLP18}, and as such actually one cannot use a similar constructions for Hitchin representations into $\SL(2n+1,\R)$. However, the base in $\mathbb{RP}^2$ of the the pencil defined by $\Psi_0$ admits two component, one that is a point, and one that is a projective line. If we remove the projective line, we obtain again the description of the convex projective structure by the affine sphere.
\end{remark}

\subsection{\texorpdfstring{$G=\SO_0(2,n+1)$}{G=SO(2,3)} and \texorpdfstring{$X = \Pho(\R^{2,n+1})$}{Pho(R(2,n)}}

In this subsection, we interpret the $(\SO_0(2,n+1), \Pho(\R^{2,n+1})$ structures of \cite{CTT19} for maximal representations $\rho:\pi_1S\rightarrow \SO_0(2,n+1)$ using bases of pencils. \medskip 

To start, let $\rho: \pi_1S \rightarrow \SO_0(2,n+1)$ be a $P_1$-Anosov representation and $\xi^1: \partial_{\infty}\pi_1S\rightarrow \Ein^{1,n}$ the associated limit map. For a point $\ell \in \Ein^{1,n}$, the \emph{Einstein universe} of isotropic lines in $\R^{2,n+1}$, consider the thickening 
\[ K_{\ell} = \{ \omega \in \Pho(\R^{2,n+1}) \mid \ell \subset \omega\} ,\]
which gives rise to a co-compact domain of discontinuity $\Omega_{\rho}$ \cite{GW12}, given by  
\begin{align}\label{OmegaPho2n}
 \Omega_{\rho} := \Pho(\R^{2,n+1}) \setminus \bigcup_{x\in \partial_{\infty}\pi_1S} K_{\xi^1(x)}. 
\end{align}

We now turn to the differential-geometric side.
By \cite{CTT19}, for every maximal representation $\rho:\pi_1S \rightarrow \SO_0(2,n+1)$, there is a unique $\rho$-equivariant (immersed) maximal spacelike surface $\sigma: \tilde{\Sigma}\rightarrow \Ha^{2,n}$ in pseudo-hyperbolic space $\Ha^{2,n},$ the set of negative lines in $\R^{2,n+1}$. In particular, 
this result picks out a distinguished Riemann surface $\Sigma$ on $S$ given the representation $\rho$. 
The $\SO_0(2,n+1)$-Riemannian symmetric space $\X$ identifies with the Grassmannian of spacelike 2-planes $\Gr_{(2,0)}(\R^{2,n+1})$. We recall that the Gauss map of $\sigma$ is the $\rho$-equivariant minimal surface \[u: \tilde{\Sigma} \rightarrow \Gr_{(2,0)}(\R^{2,n+1}), \;\; x \longmapsto d\sigma(\mathrm{T}_x\tilde{\Sigma}) \] 
associated to $\rho$ via $\NAH_{\Sigma, \SO_0(2,n+1)}$. 
In fact, $\sigma$ gives more data. 
For each point $x \in \tilde{\Sigma}$, the map $\sigma$ induces a splitting 
$\R^{2,n+1} = \mathscr{L} \oplus T \oplus N$, where $\mathscr{L} = \sigma(x)$, $T = u(x)$, and $N = (\mathscr{L}_x \oplus T_x)^\bot$, namely, $N$ is the normal space to $T$ in the tangent space $\T_{\sigma(x)}\Ha^{2,n}$. In this case, we set $\mathrm{I} \in \Omega^1(\tilde{\Sigma}, \Hom(\mathscr{L},T))$ to be the differential $d\sigma$ reinterpreted as bundle-valued one-form. We can then consider 
\[ \Psi_0:=\mathrm{I}^* \in \Omega^1(\tilde{\Sigma}, \Hom(T,\mathscr{L})),\] the adjoint of $\mathrm{I}$, and moreover $\Psi_0$ takes also the form $\Psi_0 \in \Omega^1(\tilde{\Sigma}, u^*\T\X))$ under the identification $\T_{U} \Gr_{(2,0)}(\R^{2,n+1}) \simeq \Hom(U,U^\bot)$. 

The relevant parallelism here is as follows: if $g= \sigma^*g_{\Ha^{2,n}}$ is the induced conformal metric on $\Sigma$ associated to $\sigma$, then Proposition \ref{ParallelismPair} holds for the pair $(\Psi_0, g)$. 

There is a beautiful connection between $u, \sigma$, and the associated fibered geometric structures on $\rho(\pi_1S)\backslash \Omega_{\rho}$ in \cite{CTT19}. The following proposition offers a reinterpretation of their result in terms of bases of pencils. Note that the following result differs from the formulation in \cite[Appendix A]{Dav25}, where the complex structure on the Hermitian symmetric space $\X$ was used. 

\begin{proposition}\label{Prop:Pho2nFiber}
For any $x \in \tilde{\Sigma}$, the base of pencil $\mathcal{B}(\Psi_0|_x)$ is naturally identified with $\Pho(\sigma(x)^\bot)$. 
\end{proposition}

\begin{proof}
We recall that if $U = u(x)$, then a tangent vector $\phi \in \T_{u(x)}\Gr_{(2,0)}(\R^{2,n+1})$ pointing towards $\Pho(\R^{2,n+1})$ takes the form of a map $\phi: U \rightarrow U^\bot$ that is an anti-isometry onto its image, and $\phi$ points towards $\omega \in \Pho(\R^{2,n+1})$ given by $\omega = \graph(\phi)$.

Fix an orthonormal basis $(u_1,u_2)$ of $u(x)$. Then $\phi$ is determined by $\phi(u_i) =z_i \in Q_-(U^\bot)$. Now, set $\alpha =\Psi_0(X)$ for some $X \in \mathrm{T}_x\tilde{\Sigma}$. By Proposition \ref{Prop:ModelSpaceMetric}, the pairing 
$\langle \phi, \alpha\rangle_\X$ is given by 
\[\langle \phi, \alpha\rangle = \langle \phi(u_1),\alpha(u_1)\rangle_{q_{2,n+1}} + \langle \phi(u_2),\alpha(u_2)\rangle_{_{q_{2,n+1}}}.\]
In fact, the pencil is given by $\Psi_0|_x \cong \Hom(U,\mathscr{L})$. Hence, there is a basis $(\alpha_1,\alpha_2)$ for the pencil $\Psi_0|_x$ such that $\alpha_i(u_j)=\delta_{ij} \sigma(x)$. It follows that $\phi \in \mathcal{B}(\Psi_0|_x)$ if and only if $\image(\phi) \subset \mathscr{L}^\bot$. In other words, $\phi \in \mathcal{B}(\Psi_0|_x)$ if and only if $\phi$ points towards $\omega \in \Pho(\sigma(x)^\bot)$.
\end{proof}

We now re-state the geometric structure results of \cite{CTT19} in terms of bases of pencils. 
To this end, we again define $\overline{M} \rightarrow \tilde{S}$ with fibers $\overline{M}|_x =\mathcal{B}(\Psi_0|_x) \cong \Pho(\R^{2,n})$ in this case. 

\begin{theorem}[Fibered Photon Structures {\cite[Theorem 5.3]{CTT19}}]
Let $n \geq 2$ be an integer, $\rho: \pi_1S \rightarrow \SO_0(2,n+1)$ be a maximal representation and $\sigma: \tilde{\Sigma} \rightarrow \Ha^{2,n}$ be the unique associated (immersed) maximal spacelike surface.  
The tautological developing map $\dev: \overline{M} \rightarrow \Pho(\R^{2,n+1})$ by $\dev(x, \omega)= \omega$ is a diffeomorphism onto the domain \eqref{OmegaPho2n}, and $\dev(\overline{M}_x)=\Pho(\sigma(x)^\bot)$.  
\end{theorem}

In slightly different language, the theorem asserts that $\rho(\pi_1S)\backslash \Omega_{\rho} \cong M$, where $M = \pi_1S\backslash \overline{M}$. Note that one can also directly build $M$ using Higgs bundles as in \cite{CTT19}.

\subsection{\texorpdfstring{$G=\SO_0(2,3)$}{G=SO(2,3)} and \texorpdfstring{$X = \Ein^{1,2}$}{X=Ein(1,2)}}

In this section, we reinterpret the $(\SO_0(2,n+1), \Ein^{1,2})$-structures of \cite{CTT19} for $\SO_0(2,3)$-Hitchin representations using bases of pencils. \medskip

Hitchin representations are $P_2$-Anosov and thus have a limit map $\xi^2: \partial_{\infty}\pi_1S\rightarrow \Pho(\R^{2,3})$. The notion of thickening is dual to that of the previous section. For $\omega \in \Pho(\R^{2,3})$, define 
\[ K_{\omega} = \{ \ell \in \Ein^{1,2} \mid \ell \subset \omega \}.\]
We then consider the co-compact domain of discontinuity $\Omega \subset \Ein^{1,2}$ defined by \cite{GW12}:
\begin{align}\label{OmegaEin12}
\Omega_{\rho}: = \Ein^{1,2} \setminus \bigcup_{x \in \partial_{\infty}\pi_1S} K_{\xi^2(x)}.
\end{align}

Since $\SO_0(2,3)$-Hitchin representations are also maximal, all of the differential geometry setup of the previous subsection applies. In particular, $\rho: \pi_1S \rightarrow \SO_0(2,3)$ Hitchin admits a unique $\rho$-equivariant maximal spacelike surface $\sigma:\tilde{\Sigma} \rightarrow \Ha^{2,2}$ whose Gauss map $u: \tilde{\Sigma} \rightarrow \Gr_{(2,0)}(\R^{2,3})=\X$ is a minimal surface in $\X$. 

Here, $\Psi_0$ is different than in the previous section. The maximal surfaces $\sigma$ associated to $\SO_0(2,3)$-Hitchin representations are special in the sense that their second fundamental forms are non-vanishing, with two-dimensional image at every point. In this case, we consider the second fundamental form of $\sigma$, by $\sff \in \Omega^1(\tilde{\Sigma}, \Hom(T,N))$. We then define 
\[ \Psi_0 := \sff \in \Omega^1(\tilde{\Sigma}, u^*\T\X).\]

By Labourie \cite{Lab17}, we can associate naturally to $\rho$ a pair $(\Sigma, q_4)$, where $q_4 \in H^0(\K_{\Sigma}^4)$ is a holomorphic quartic differential on $\Sigma$. Then, using the Hitchin section and the non-abelian Hodge correspondence, we have the Higgs bundle $\mathcal{H}(\rho) = \NAH_{\Sigma, \SO_0(2,3)}^{-1}(\rho) = (\V, \Phi)$, where
$\V = \bigoplus_{i=2}^{-2} \K_{\Sigma}^i$. Here, $\Phi = \Phi(q_4)$ obtains the form:
\begin{equation*}
\begin{tikzcd}
\mathcal{K}_{\Sigma}^2 \arrow[r, "1"] & \mathcal{K}_{\Sigma} \arrow[r, "1"]&\mathcal{O} \arrow[r, "1"]&\mathcal{K}_{\Sigma}^{-1} \arrow[r, "1"]\arrow[lll,bend right,swap, "q_4"]&\mathcal{K}_{\Sigma}^{-2}\arrow[lll,bend right,swap, "q_4"]
\end{tikzcd}.
\end{equation*}

The solution to Hitchin's equations yields a harmonic metric $h$ on $\mathcal{H}(\rho)$ that is diagonal and obtains the form 
$h = \diag(h_2,h_1,1, h_1^{-1}, h_2^{-1})$. Here, observe that $g: =h_2^{-1}h_1 \in \Gamma(\K\overline{\K})$ is a Riemannian metric on $\Sigma$, which happens to be conformal. More geometrically, $g = g_{\mathcal{N}}/g_{\mathcal{T}}$, where $\mathcal{N} \cong \K^{-2}, \mathcal{T}\cong \K^{-1}$ are the holomorphic normal and tangent line bundles of $\sigma$ and $g_{\mathcal{N}},g_{\mathcal{T}}$ their metrics. 
The parallelism in this case is that Proposition \ref{ParallelismPair} holds for the pair $(\Psi_0, g)$. 

We now identify the relevant base of pencil geometrically and topologically.  
\begin{proposition}\label{Prop:Ein12Fiber}
For any $x \in \tilde{\Sigma}$, the base of the pencil $\mathcal{B}_{\Ein^{1,2}}(\Psi_0|_x)$ identifies naturally with $\Ein(\mathscr{L}_x\oplus T_x) \cong \mathbb{S}^1$. 
\end{proposition}

\begin{proof}
First, recall that $\phi \in \T_{x}\X$ points towards $\ell \in \Ein^{1,2}$ if and only if $\phi$ is a rank one map whose graph over its cokernel is $\ell$. 

Now, take $\phi \in \T_{u(x)}\X$ pointing towards $\Ein^{1,2}$. Then $\phi$ obtains the form $\phi(u_1) =z_1$, $\phi(u_2)=0$ for some basis $(u_1,u_2)$ of $u(x)$. For any element $\alpha \in \Psi_0|_x$, we have $c\langle \phi,\alpha\rangle_{\X} = \langle \phi(u_1), \alpha(u_1)\rangle_{q_{2,n+1}}$, up to some positive scalar $c \in \R$. 

In this case, $\sff|_x: T_x \times T_x\rightarrow N_x$ has two-dimensional image. 
That is, fixing any tangent vector $X \in T_x$, $\image(\sff(X,\cdot)) = N$. 
Similar to the proof of Proposition \ref{Prop:Pho2nFiber}, one finds $\phi \in \mathcal{B}_{\Ein^{1,2}}(\Psi_0|_x)$ if and only if $\image(\phi) \subset (N_x^\bot \subset u(x)^\bot)= T_x$. In particular, $\phi \in \mathcal{B}_{\Ein^{1,2}}(\Psi_0|_x)$ if and only if $\ell = \graph^*(\phi) \in \Ein(\mathscr{L}_x \oplus T_x)$. 
\end{proof}

Once more, set $\overline{M} \rightarrow \tilde{S}$ to be the smooth $\mathbb{S}^1$-fiber sub-bundle of $\tilde{S}\times \Ein^{1,2}$ with fiber $\overline{M}|_x= \mathcal{B}_{\Ein^{1,2}}(\Psi_0|_x)$. We reach a reinterpretation of $\Ein^{1,2}$-structures from \cite{CTT19}. 

\begin{theorem}[Fibered Einstein Structures {\cite[Theorem 5.6]{CTT19}}]
Let $\rho: \pi_1S\rightarrow \SO_0(2,3)$ be Hitchin. Then the tautological developing map $\dev: \overline{M} \rightarrow \Ein^{1,2}$ by $\dev(x, \ell) = \ell$ is a diffeomorphism onto the domain \eqref{OmegaEin12}. 
\end{theorem}

In particular, the quotient $M = \pi_1S\backslash \overline{M}$ identifies with $\mathrm{T}^1S$ and also with the quotient $\rho(\pi_1S)\backslash \Omega_{\rho}$. 

\begin{remark}
The fact that $\sigma$ was immersed for $\rho$-Hitchin is not essential; one needs only the condition on $\sff$. Indeed, Filip builds fibered $\Ein^{1,2}$-structures for non-Hitchin $\SO_0(2,3)$-representations in \cite{Fil23} by an equivalent developing map to that of \cite{CTT19}.
\end{remark}

\subsection{Summary}\label{Sec:Summary}

Here, we give a Lie-theoretic table that summarizes the differential-geometric construction of fibered geometric structures for the five cases of interest of rank two geometries in terms of bases of pencils.

\begin{table}[ht]
\small
 \renewcommand{\arraystretch}{1.2}
    \begin{tabular}{|c|c|c|c|c|}
    \hline
        $G$ & $X$ & Auxiliary map $\sigma$ & Pencil $\Psi_0$ & Projected metric $g$ \\ \hline
        $\SL(3,\R)$ & $\Flag(\R^3)$ & Hyperbolic affine sphere in $\R^3$ & $\Sym(\mathrm{I})$ & Blashke metric \\ \hline 
        $\SO_0(2,n+1)$ & $\Pho(\R^{2,n+1})$ & $\beta$-maximal spacelike surface in $\Ha^{2,n}$ & $\mathrm{I}^*$ & Induced metric $\sigma^*g_{\Ha^{2,n}}$ \\ \hline
        $\SO_0(2,3)$ & $\Ein^{1,2}$ & $\alpha$-maximal spacelike surface in $\Ha^{2,2}$ & $\sff$ & $g_N/g_{T}$ \\ \hline 
        $\Gtwosplit$ & $\Ein^{2,3}$ &  $J$-holomorphic $\beta$-curve in $\quadric$ & $\rm{I}+\tff^{0,1}$ & Induced metric $-\sigma^*g_{\quadric}$ \\ \hline 
        $\Gtwosplit$ & $\Pho^\times$ & $J$-holomorphic $\alpha$-curve in $\quadric$ & $\sff^*$ & $g_N/g_{T}$  \\ \hline 
    \end{tabular}
    \caption{The auxiliary map $\sigma$, whose appropriate Gauss map is $u:\tilde{\Sigma} \rightarrow \X$. The pencil $\Psi_0$ along $u$ is described in term of differential-geometric data $\sigma$. The conformal Riemannian metric $g$ encodes the parallelism $\nabla^\X\circ \Psi_0 = \Psi_0 \circ \nabla^g$. }
    \label{tab:placeholder}
\end{table}

Here, we provide some brief explanation of the table. We set $(\beta, \alpha)$ to be the pair of short-long roots for $\g \in \{ \mathfrak{so}(2,n+1), \g_2'\}$. For $G=\SO_0(2,n+1)$, $\alpha$ and $ \beta$-maps are maximal spacelike surfaces in $\Ha^{2,n}$ with extra adjectives: 
\begin{itemize}[noitemsep]
    \item $\beta$-maps $\tilde{\Sigma} \rightarrow \Ha^{2,n}$ are \emph{immersed},
    \item $\alpha$-maps $\tilde{\Sigma} \rightarrow \Ha^{2,2}$ are those with \emph{isomorphic $\sff$}.
\end{itemize} 
If $\mathcal{T}, \mathcal{N}$ are the holomorphic tangent and normal line bundles of $\sigma$, then the $\alpha$-condition entails $\mathcal{N}\otimes \mathcal{T}^{-1} \cong \mathcal{K}^{-1}$ and the $\beta$-condition entails $\mathcal{T} \cong \mathcal{K}^{-1}$. 

Similarly, for $G=\Gtwosplit$, the $\alpha, \beta$-curves $\sigma:\tilde{\Sigma} \rightarrow \quadric$ are equivariant alternating $J$-holomorphic curves with extra adjectives: 
\begin{itemize}[noitemsep]
    \item $\beta$-curves $\tilde{\Sigma} \rightarrow \quadric$ are \emph{immersed},
    \item $\alpha$-curves $\tilde{\Sigma} \rightarrow \quadric$ are those with \emph{isomorphic $\sff$}. 
\end{itemize}
For $G \in \{\SO_0(2,n), \Gtwosplit\}$, the subgroup $\hat{H}$ turns out to be the same for both the representations of $\SL(2,\R)$ we consider; namely, it is a maximal torus $T$ in the maximal compact $K$. The map $\hat{\sigma}:\tilde{S}\to G/\hat{H}=G/T$ in these cases is a cyclic surface as in \cite{Bar10, Lab17, CT24} a notion inspired by the $\tau$-primitive harmonic maps to $G/T$ for $G$ a compact Lie group, studied in \cite{BPW95}. The space $G/T$ has special structure and in particular, the tangent bundle $\mathrm{T}(G/T)$ has a line-sub-bundle $[\mathcal{L}_{\gamma}$] for each simple root $\gamma$ of $\g$. 
The pencils $\Psi_0$ can uniformly be viewed as modifications to $du$ that are encoded Lie-theoretically as projections of $d\hat{\sigma}$ to $[\mathcal{L}_{\eta}]$, with $\eta$ determined by the relevant Anosov flavor. In particular, $\eta$ is either $\alpha$ or $\beta$ in each of the four cases for $G \in \{\SO_0(2,n+1), \Gtwosplit\}$. A similar discussion applies to the case of $G=\SL(3,\R)$ and the affine sphere, except in this case, one should project $d\hat{\sigma}$ to the sum $[\mathcal{L}_{\alpha}]\oplus [\mathcal{L}_{\beta}]$, which corresponds to using the $P_{\Delta}$-Anosov condition. 

\section{The five \texorpdfstring{$\sllie_2(\R)$}{sl(2,R)}-subalgebras in \texorpdfstring{$\g_2'$}{g2'}}\label{Appendix:SL2}

 There are exactly five $\mathfrak{s}_2\R$-subalgebras in $\g_2'$ up to the adjoint action, classified in \cite{Dok88} under a much broader classification. In particular, the conjugacy class of the $\sllie_2$-triple $\{E,F,H\}$ is determined by the conjugacy class of the nilpotent $E$. See \cite{CM17} or \cite{BCGBO24} for further details. In this appendix, we explicitly describe these five subalgebras in $\g_2'$ in terms of their action on $\imoct$. Note the subtle difference between the fourth and fifth subalgebras in Table \ref{Table:sl2ing2}. 

 \begin{table}[ht]
 \renewcommand{\arraystretch}{1.2}
    \small 
 	\begin{center}
		\begin{tabular}{| c | c | c  | c | } 
			\hline 
                 Irreducible Subspaces & Splitting Preserved & \; Nilpotent & \; Anosov condition  \\
                \hline 
               7 & $\imoct$ & $e_{-\alpha} +e_{-\beta}$ & $\Delta$-Anosov \\
               \hline 
     3+2+2 & 
        $\langle x_{3}, x_{2} \rangle \oplus \langle x_1, x_0, x_{-1} \rangle \oplus \langle x_{-2}, x_{-3} \rangle$ & $e_{-\beta}$ & $\beta$-Anosov \\
        \hline 
        2+2+1+1+1 &  $\langle x_{3}\rangle  \oplus \langle x_{2}, x_1 \rangle \oplus \langle x_0 \rangle \oplus \langle  x_{-1}, x_{-2} \rangle \oplus \langle x_3 \rangle $ & $e_{-\alpha}$ & $\alpha$-Anosov \\ \hline 
    
        3+3+1 & $\begin{cases} \mathbb{E}_{+1}(\mathcal{C}_{x_0}) \oplus \mathbb{E}_{-1}(\mathcal{C}_{x_0}) \oplus \langle x_0 \rangle \\ 
        \R^{2,1} \oplus \R^{0,1} \oplus \R^{1,2} \end{cases} $
        & $e_{-\beta}+e_{\delta}$ 
        & $\alpha$-Anosov \\ \hline
        3+3+1 & $\R^{1,2} \oplus \R^{1,0} \oplus \R^{1,2}$ & $e_{-\alpha}+e_{\delta}$ & $\alpha$-Anosov \\ \hline
			 \end{tabular} 
        \vspace{1ex}
        \caption{\emph{The five $\sllie_2(\R)$-subalgebras $\mathfrak{s}$ in $\g_2'$ up to the adjoint action and their geometric properties, including the  Anosov condition for Fuchsian representations $\rho_0: \pi_1S \rightarrow (\mathrm{P})\SL(2,\R) \rightarrow \Gtwosplit$ factoring through the corresponding analytic subgroup.
        Unlike in $\sllie(n,\R)$, the dimensions of the irreducible subspaces do not classify the subalgebra $\mathfrak{s}$ up to conjugacy. The fourth subalgebra simultaneously preserves two different splittings. }}
        \label{Table:sl2ing2}
	\end{center} 
\end{table} 

In many cases below, we use the same data from Subsection \ref{Subsec:LieTheory}, namely the Cartan subalgebra $\mathfrak{a} < \g_2'$ from \eqref{g2Cartan} with corresponding root system $\Sigma=\Sigma(\g_2',\mathfrak{a})$ and simple roots $\Delta = \{\alpha, \beta\}$. 

\subsection{Principal \texorpdfstring{$\mathfrak{sl}_2\R$}{sl(2,R)}} 

The principal subalgebra $\mathfrak{s}_{\Delta}$ is the unique $\sllie_2\R$-subalgebra of $\g_2'$ up to conjugation that acts irreducibly on $\imoct$. 
The principal embedding $\iota_{\Delta}: \PSL(2,\R) \hookrightarrow \Gtwosplit$ is explicitly written out and examined in \cite[Section 5.1]{Eva25}. 
In our Higgs bundle notation, a principal nilpotent $E \in \mathfrak{s}$ is given by $\alpha =1, \, \beta = 1, \, \delta =0$.
Representations of the form $\iota_{\Delta} \circ \rho_0: \pi_1S \rightarrow \Gtwosplit$ are $P_{\Delta}$-Anosov and also Fuchsian-Hitchin, where $\rho_0:\pi_1S\rightarrow \PSL(2,\R)$ is Fuchsian. 

\subsection{Short Root  \texorpdfstring{$\mathfrak{sl}_2\R$}{sl(2,R)}}
\label{Subsec:SL2ShortRoot}

Let us denote by $\mathfrak{s}_{\beta}$ the $\sllie_2\R$-subalgebra of the short root $\beta$, namely $\mathfrak{s}_{\beta}= \spann_{\R} \{ E_{-\beta}, E_{\beta}, T_{\beta} \}$. Recall that $\mathfrak{s}_{\beta} $ is a Lie subalgebra of the Levi subalgebra $ \mathfrak{l}_{\alpha}$ of the parabolic subalgebra $\mathfrak{p}_{\alpha}$ associated to the long root $\alpha$. The representation $\iota_1: \mathfrak{s}_{\beta} \rightarrow \gl(\imoct)$ decomposes into irreducible subspaces of dimension $3+2+2$, which, in an appropriate $\R$-cross product basis $X= (x_k)_{k=3}^{-3}$ are as follows:
\[ \imoct =  \langle x_3, x_2 \rangle \oplus \langle x_1, x_0, x_{-1} \rangle \oplus  \langle x_{-2}, x_{-3} \rangle. \]
This splitting is also preserved by the Levi subalgebra $\mathfrak{l}_{\alpha} \cong \gl_2\R$. Denote the corresponding subgroup by $\GL(2,\R)\cong L_{\alpha} < \Gtwosplit$. The analytic subgroup in $\Gtwosplit$ associated to $\mathfrak{s}_{\beta}$ is isomorphic to $\SL(2,\R)$ and realized as the mutual stabilizer of the volume forms on the 2-planes $\langle x_3,x_2\rangle$ and $\langle x_{-2},x_{-3}\rangle$. 

This representation $\mathfrak{s}_{\beta}\hookrightarrow \g_2'$ factors as $\mathfrak{sl}_2(\R)\hookrightarrow \mathfrak{so}(2,2)\hookrightarrow \g_2'$. The  relevant $\SO(2,2)$-subgroup in $\Gtwosplit$ is $\SO(2,2) = \Stab_{\Gtwosplit}(\Ha')$, where $\Ha' <\Oct'$ denotes the \emph{split quaternions}. The $\SO(2,2)$-subgroup acts with dimension 3+4 blocks in $\imoct$, and the invariant splitting is of the form $\R^{1,2} \oplus \R^{2,2}$. In terms of the basis $(x_k)_{k=3}^{-3}$, the invariant subspaces of $\SO(2,2)$ are $\langle x_1,x_0,x_{-1}\rangle$ and $\langle x_3,x_2,x_{-2},x_{-3}\rangle$. 
 The analytic subgroup $\SL(2,\R)_{\beta}$ of $\SO(2,2)$ associated to $\mathfrak{s}_{\beta}$ acts irreducibly on the 3-block, but the 4-block $\R^{2,2}$ decomposes into irreducible 2+2 blocks, which are dual isotropic planes $\langle x_3, x_2\rangle$ and $\langle x_{-2}, x_{-3}\rangle$.
 
Write $\iota_{\beta}: \SL(2,\R) \hookrightarrow \Gtwosplit$ for the corresponding embedding of $\SL(2,\R)$. Representations of the form $\iota_{\beta} \circ \rho_0$ are $P_{\beta}$-Anosov, for $\rho_0:\pi_1S \rightarrow \SL(2,\R)$ Fuchsian. In our Higgs bundle notation, the nilpotent element $E \in \mathfrak{s}_{\beta}$ corresponds to $\alpha =0, \,\beta = 1, \,\delta =0$. 
 
\subsection{Long Root \texorpdfstring{$\mathfrak{sl}_2\R$}{sl(2,R)}}\label{Subsec:SL2LongRoot}

Denote by $\mathfrak{s}_{\alpha}$ the $\sllie_2\R$-subalgebra $\mathfrak{s}_{\alpha} = \spann_{\R}\{ E_{-\alpha}, E_{\alpha}, T_{\alpha} \}$ of the long root $\alpha$. Similar to the previous case, we have $\mathfrak{s}_{\alpha} < \mathfrak{l}_{\beta}$, where $\mathfrak{l}_{\beta}$ is the Levi subalgebra $\mathfrak{l}_{\beta}$ of the parabolic subalgebra $\mathfrak{p}_{\beta}$. The representation $\mathfrak{s}_{\alpha} \rightarrow \gl(\imoct)$ decomposes into irreducible blocks of dimension 2+2+1+1+1, with irreducible splitting from the $\beta$-height grading: 
\[ \imoct =  \langle x_3 \rangle \oplus \langle x_2, x_1 \rangle \oplus \langle x_0 \rangle \oplus \langle x_{-1}, x_{-2} \rangle \oplus \langle x_{-3} \rangle .\] 
The Levi subgroup $L_{\beta}$ is realized as $L_{\beta} = \Stab_{\Gtwosplit}(\langle x_3\rangle) \cap \Stab_{\Gtwosplit}(\langle x_{-3}\rangle) \cong \GL(2,\R)$, acting faithfully on the 2-plane $\langle x_3, x_2\rangle$, seen in the model of Proposition \ref{Lem:NullStiefel} relative to $n :=(x_2,x_1,x_{-3}) \in \mathcal{N}_\R$. The analytic subgroup associated to $\mathfrak{s}_{\alpha}$ is given by $\SL(2,\R)_{\alpha}=\Stab_{L_{\beta}}(x_{-3})=\Stab_{L_{\beta}}(x_{-3})\cap \Stab_{L_{\beta}}(x_3)$. 
The restiction map $\SL(2,\R)_{\alpha} \rightarrow \SL(\langle x_2,x_1\rangle)$ is a Lie group isomorphism. 

The representation $\mathfrak{s}_{\alpha}$ on $\imoct$ factors as $\mathfrak{s}_{\alpha} \hookrightarrow \sllie_3\R \hookrightarrow \g_2'$ through a reducible representation $\sllie_2\R\hookrightarrow \sllie_3\R$. In other words, the subgroup $\mathrm{SL}(2,\R)_{\alpha}$ is realized through the composition $\SL(2,\R) \stackrel{red}{\hookrightarrow} \mathrm{SL}(3,\R) \hookrightarrow \Gtwosplit$. Recall from \eqref{eq:QuadricStabilizers} that for any unit timelike vector $x_0 \in Q_-(\imoct)$, one has $\Stab_{\Gtwosplit}(x_0) \cong \mathrm{SL}_3(\R) $. In particular, this 2+2+1+1+1 splitting is also a refinement of the 3+3+1 spitting of the fourth subalgebra of Table \ref{Table:sl2ing2}.

In our Higgs bundle notation, the nilpotent element $E \in \mathfrak{s}_{\alpha}$ is described by $\alpha =1, \, \beta = 0, \, \delta =0$. Fuchsian representations in $\SL(2,\R)$ included into $\Gtwosplit$ via $\SL(2,\R)_{\alpha}$ are $P_{\alpha}$-Anosov.

\subsection{(Principal \texorpdfstring{$\mathfrak{sl}_2\R$}{sl(2,R)} in \texorpdfstring{$\sllie_3\R \hookrightarrow \g_2'$}{sl(3,R)--> g2'}) = (Stabilizer of timelike vector in \texorpdfstring{$\mathfrak{so}(2,2) \hookrightarrow \g_2'$}{so(2,2)->g2'})}  \label{Appendix:SL2Number4}

Our next $\mathfrak{sl}_2\R$-subalgebra, denoted $\mathfrak{s}_4$, embeds in $\g_2'$ as  
$\sllie_2\R \stackrel{princ}{\hookrightarrow} \sllie_3\R \hookrightarrow \g_2' $, where the $\sllie_3\R$ is the $\sllie_3\R$-subalgebra generated by the root vectors of the long roots. 
This representation $\mathfrak{s}_4 \hookrightarrow \gl(\imoct)$ decomposes into irreducible subspaces of dimension $1+3+3$. However, this subalgebra happens to preserve \emph{two different} such splittings. 

Fixing a background $\R$-cross product basis $(x_k)_{k=3}^{-3}$, one such splitting preserved is  
\begin{align}\label{s4FirstSplitting}
    \imoct =  \langle x_0\rangle \oplus \mathbb{E}_{+1}(\mathcal{C}_{x_0}) \oplus  \mathbb{E}_{-1}(\mathcal{C}_{x_0}) = \langle x_0 \rangle \oplus \, \langle x_{3}, x_{-1}, x_{-2} \rangle \oplus \langle x_{-3}, x_{1}, x_2 \rangle. 
\end{align}
The splitting \eqref{s4FirstSplitting} preserved by $\mathfrak{s}_4$ consists of a timelike line and two dual isotropic 3-planes, which are the $(+1)$ and $(-1)$-eigenspaces of the cross product endomorphism $\mathcal{C}_{x_0}$ of $x_0$. 

To describe the other splitting, consider the subgroup $H_4 := \Stab_{\Gtwosplit}(\mathbf{l}) \cap \Stab_{\Gtwosplit}( \langle \,\mathbf{j}, \mathbf{li}, \mathbf{lk}\,\rangle)$, which has irreducible subspaces 
\[ \langle \,\mathbf{j}, \mathbf{li}, \mathbf{lk}\,\rangle \oplus \langle \,\mathbf{j}, \mathbf{k}, \mathbf{lj} \,\rangle \oplus \langle \mathbf{l} \rangle \cong \R^{1,2} \oplus \R^{2,1} \oplus \R^{0,1}. \]
One checks $H_4 \cong \SO_0(1,2)\cong \PSL(2,\R)$ via the null Stiefel triplet model from Proposition \ref{Prop:FirstStiefel}. We claim that, in fact, the subalgebra $\Lie(H_4)$ is conjugate to $\mathfrak{s}_4$. 
Note that $H_4 < \Stab_{\Gtwosplit}(\mathbf{l})\cong \SL(3,\R)$. Moreover, one finds that $H$ acts irreducibly $\mathbb{E}_{+1}(\mathcal{C}_{\mathbf{l}})$. Indeed, to see this, note that any element $y \in \mathbb{E}_{+1}(\mathcal{C}_{\mathbf{l}})$ can be written $y = x+\mathbf{l}x$ for $x \in \langle \mathbf{j},\mathbf{li},\mathbf{lk}\rangle$. Hence, $H_4$ acts irreducibly on $\mathbb{E}_{+1}(\mathcal{C}_{\mathbf{l}})$ as a consequence of the irreducible action on the 3-plane $\langle \mathbf{j},\mathbf{li},\mathbf{lk}\rangle $. This means $\Lie(H_4)$ is then a principal $\mathfrak{sl}_2\R$-subalgebra in $\mathfrak{sl}_3(\R)$, which then coincides with $\mathfrak{s}_4$ up to conjugation by $\SL(3,\R)$ due to the uniqueness of principal $\sllie_2\R$-subalgebras in $\sllie_3\R$. 

Let us denote $\iota_4: \mathrm{PSL}_2\R \hookrightarrow \mathrm{SL}_3(\R) \hookrightarrow \Gtwosplit$ as the associated inclusion of $H_4$. Representations $\iota_4 \circ \rho_0$ are $P_{\alpha}$-Anosov, where $\rho_0:\pi_1S \rightarrow \PSL(2,\R)$ is Fuchsian. In our Higgs bundle notation, the nilpotent $E$ is described by $\alpha =0, \, \beta = 1, \, \delta = 1$. See Remark \ref{Remk:BetaSpecialCases} for a uniformizing Higgs bundle. 

\subsection{Stabilizer of spacelike vector in \texorpdfstring{$\mathfrak{so}(2,2) \hookrightarrow \g_2'$}{so(2,2)->g2'}}\label{Appendix:SL2Number5}

The final $\mathfrak{sl}_2\R$-subalgebra $\mathfrak{s}_5$ includes as $\mathfrak{s}_5 \stackrel{\cong}{\rightarrow} \mathfrak{so}(1,2) \stackrel{red}{\hookrightarrow}\mathfrak{so}(2,2) \hookrightarrow \g_2'$. The associated Lie subgroup $H_5 < \Gtwosplit$ can be described explicitly by 
\[ H_5 = \Stab_{\Gtwosplit}( \langle \mathbf{i}, \mathbf{l}, \mathbf{li}\rangle) \cap \Stab_{\Gtwosplit}(\mathbf{j}).\]
We recall here that $\Stab_{\Gtwosplit}( \langle \mathbf{i}, \mathbf{l}, \mathbf{li}\rangle)\cong \SO(2,2)$ since $\mathbb{H}'\cong \spann \{1,\mathbf{i},\mathbf{l},\mathbf{li}\}$ is a copy of the split quaternions in $\Oct'$ and $\Stab_{\Gtwosplit}(\mathbf{j})\cong \SU(1,2)$ by \eqref{eq:QuadricStabilizers}. We claim that $H_5 \cong \SO_0(2,1)\cong \PSL(2,\R)$. To see this, use the model point $p_0 =(\mathbf{i},\mathbf{j},\mathbf{l}) \in V_{(+,+,-)}(\imoct)$ from Proposition \ref{Prop:FirstStiefel} and consider the action of $H_5$ on $p_0$. 

The splitting preserved by $H_5$ is as follows:
\[ \imoct = \langle \mathbf{i}, \mathbf{l}, \mathbf{li}\rangle \oplus \langle \mathbf{k, \mathbf{lj}, \mathbf{lk}} \rangle \oplus \langle \mathbf{j} \rangle \cong \R^{1,2} \oplus \R^{1,2} \oplus \R^{0,1}.\]
On the other hand, $H_5$ is an $\SO_0(1,2)$-subgroup of $\SU(1,2)$. Recall that $\mathcal{C}_{\mathbf{j}}$ defines a complex structure on $[\mathbf{j}^\bot \subset \imoct] \cong \R^{3,3}$ due to \eqref{DCP}. From this perspective, we see that 
\[ V_1^\C = V_1 \oplus \mathbf{j} V_1 = V_1 \oplus V_2, \] 
and similarly $V_2^\C = V_1 \oplus V_2$. Thus, each of $V_1$ and $V_2$ is a totally real subspace of $(V_1 \oplus V_3, \mathcal{C}_\mathbf{j}) \cong \C^{1,2}$. 

If $\iota_5: \PSL(2,\R)\hookrightarrow \Gtwosplit$ is the inclusion of $H_5$, then representations $\iota_5 \circ \rho_0$ are $P_{\alpha}$-Anosov, for $\rho_0: \pi_1S\rightarrow \PSL(2,\R)$ Fuchsian. In our Higgs bundle notation, a nilpotent $E_5$ representing $\mathfrak{s}_5$ is given by $\alpha =1$, $\beta =0$, $\delta =1$. See Remark \ref{Remk:SpecialAlphaBundles} for a uniformizing Higgs bundle, and note the subtle differences in the uniformizing Higgs bundles of $H_4$ and $H_5$-Fuchsians. 

\section{Regularity of pencils}\label{Sec:RegularityPencil}

In this brief subsection, we discuss some useful auxiliary results on regularity of pencils in the $\Gtwosplit$-symmetric space. 

The following proposition clarifies when the simple roots vanish on the Cartan projection of a general element $\Psi \in \g_2'$. 
We label the simple roots $\Delta =\{\alpha, \beta\}$, with $\beta$ the short root, so that $\Gtwosplit/P_{\beta} \cong \Ein^{2,3}$ and $\Gtwosplit/P_{\alpha} \cong \Pho^\times$. In fact, we shall take the model data $(\mathfrak{a}, \Delta)$ as in Subsection \ref{Subsec:LieTheory}. 
Recall that an element $X \in \mathfrak{a}$ is called $\sigma$-\emph{regular}, for $\sigma \in \mathfrak{a}^*$, if and only if $\sigma(X) \neq 0$. 

\begin{proposition}\label{Prop:RegularityinG2'}
    The Cartan projection $\mu(\Psi)\in \mathfrak{a}$ of a non-zero semi-simple $\Psi \in \mathfrak{g}_2'$ of the Lie algebra of $\Gtwosplit$ with characteristic polynomial $\chi(\Psi)=X^7-AX^5+BX^3-CX$, is determined up to a multiplicative scalar by the value of the invariant:
\[I(\Psi)=\frac{54 C}{A^3}\in [0,1].\]
Moreover, 
\begin{enumerate}[noitemsep]
    \item $\Psi$ is $\beta$-regular if and only if $I(\Psi)\neq 0$.
    \item $\Psi$ is $\alpha$-regular if and only if $I(\Psi)\neq 1$.
\end{enumerate}

\end{proposition}

\begin{proof}

The quantity $I$ depends only of the conjugacy class of $\Psi$ as it is defined in terms of the coefficients of the characteristic polynomial, and is homogeneous, i.e. $I(\Psi)=I(\lambda\Psi)$ for all $\lambda\neq 0$.\medskip

The conjugacy class of the semi-simple element $\Psi$ is uniquely determined by the non-zero coefficients of the invariant polynomials, which are $A,B,C$ for elements of the Lie algebra of $\Gtwosplit$, but for elements of $\g_2'$, $B=A^2/4$, so $I$ determines the projective class of $\mu(\Psi)$.

\medskip

The projective class of $\mu(\Psi)$ lies in a segment because $\Gtwosplit$ has rank $2$. Let us  express the model Cartan subalgebra $\mathfrak{a}$ as in \eqref{g2Cartan}. The extreme points of this segment are the diagonal matrices corresponding to the coroots 
$T_{2\alpha+3\beta}=\text{diag}(1,1,0,0,0,-1,-1)$ and $T_{\alpha+2\beta}=\text{diag}(2,1,1,0,-1,-1,-2)$, the elements on the boundary of the Weyl chamber, which span exactly the projective classes that are not $\alpha$-regular, respectively $\beta$-regular. The value of $I$ for these are respectively $1$ and $0$.
\end{proof}

\begin{figure}[ht]
    \centering
\begin{tikzpicture}
    \fill[blue!50] (0,0) -- (0:3.3cm) -- (30:3.3cm) -- cycle;
    \foreach\ang in {0,60,120,...,300}{
     \draw[->,thick] (0,0) -- (\ang:1.731cm);
    }
    \foreach\ang in {30,90,...,330}{
     \draw[->,thick] (0,0) -- (\ang:3cm);
    }
    \node[right] at (90:3cm) {$t_{\alpha}$};
    \node[below right] at (300:1.731cm) {$t_{\beta}$};
    \node[below] at (0:1.731cm) {$t_{2\beta+\alpha}$};
    \node[above left] at (30:3cm) {$t_{3\beta+2\alpha}$};

  \end{tikzpicture}
\caption{Illustration of the coroots in $\mathfrak{a}$, and the model Weyl chamber $\mathfrak{a}^+$.}
    \label{fig:G2regularity}
\end{figure}

Going forwards, we fix the basis $(e_k)_{k=3}^{-3}$ from \eqref{BaragliaBasis}. We view $\g_2' < \g_2^\C$ as the subalgebra fixing the standard real subspace $\imoct \subset \imoct^\C$. Concretely, this subalgebra is the fixed point set of the involution $\lambda:\g_2^\C \rightarrow \g_2^\C$ given by $\lambda(A) = Q_0\overline{A}Q_0$, where the matrix $Q_0$ is:  
\[ Q_0 = \begin{pmatrix}
    &&&&&&1\\
    &&&&&1&\\
    &&&&1&&\\
    &&&1&&&\\
    &&1&&&&\\
    &1&&&&&\\
    1&&&&&&\\
\end{pmatrix} \]
Now, we consider a matrix $\Psi \in \g_2'$ of the following form, for some complex numbers $a, b \in \C$: 

\begin{equation}\label{PsiMatrix} \Psi=\begin{pmatrix}
    0& b^* &0&0&0&0&0\\
    b& 0 &a^*&0&0&0&0\\
    0& a&0&\sqrt{2}ib^*&0&0&0\\
    0& 0 &-\sqrt{2}ib&0&\sqrt{2}ib^*&0&0\\
    0& 0 &0&-\sqrt{2}ib&0&a^*&0\\
    0& 0 &0&0&a&0&b^*\\
    0& 0 &0&0&0&b&0
       \end{pmatrix}.
    \end{equation}
Note that $\Psi \in \g_2^\C$ by Subsection \ref{Subsec:LieTheory}, and $\Psi \in \g_2'$ by direct inspection. 
A straightforward calculation shows the value of this invariant $I$ for elements $\Psi$ having the form above is exactly:
\begin{align}\label{I_Quantity}
   I(\Psi)=27\frac{|ab^2|^2+|b|^6}{\left(|a|^2+3|b|^2\right)^3}.
\end{align}

We determine when such an element $\Psi$ is regular for the long root $\alpha$. 

\begin{proposition}\label{Prop:Alpha2RegularPencil}
The element $\Psi$ in \eqref{PsiMatrix} is $\alpha$-regular if and only if $a\neq 0$.
\end{proposition}

\begin{proof}
We have $I(\Psi)=1$ if and only if the following quantity vanishes:  
\[ \left(|a^2|+3|b^2|\right)^3-27\left(|a^2b^4|+|b^6|\right)=|a^6|+9|a^4b^2|. \] 
This is non-zero if and only if $a\neq 0$.
\end{proof}

\setstretch{1.0}
\bibliographystyle{alpha}
\bibliography{bib}

\end{document}